\documentclass[master,       
               twoside,        
               BCOR10mm,       
               ngerman,english, 
               ]{GAUBM2}

\usepackage{setspace}  
\usepackage{babel}     
\usepackage{calc}      
\usepackage[T1]{fontenc}       
\usepackage[latin1]{inputenc}  

\usepackage{amsmath,amssymb,mathrsfs,amsthm} 

\usepackage{lmodern} 
\usepackage[comma,numbers,sort&compress]{natbib}
\usepackage{booktabs}                      
\usepackage{longtable}                     
\usepackage{array}                         

\usepackage{textcomp,gensymb}

\usepackage{units}   

\usepackage{enumerate}

\usepackage[all]{nowidow}

\usepackage[pdftex]{graphicx}

\usepackage{algorithm2e}

\usepackage{floatrow}

\usepackage{float}

\usepackage[labelfont=bf,format=plain, font=footnotesize]{caption}

\usepackage{wrapfig}

\usepackage{subfig}
\usepackage{sidecap}

 \usepackage[pdfstartview=FitH,      
             breaklinks=true,        
             bookmarksopen=true,     
             bookmarksnumbered=true  
             ]{hyperref}
             
\usepackage{framed}
\usepackage{pstricks}
\usepackage{tabularx}

\newcommand{\mN}{\mathbb{N}}
\newcommand{\mZ}{\mathbb{Z}}

\newcommand{\mP}{\mathbb{P}}
\newcommand{\mR}{\mathbb{R}}
\newcommand{\mC}{\mathbb{C}}
\newcommand{\mK}{\mathbb{K}}

\newcommand{\mX}{\mathbb{X}}
\newcommand{\mY}{\mathbb{Y}}
\newcommand{\mH}{\mathbb{H}}
\newcommand{\mV}{\mathbb{V}}

\newcommand{\mL}{\mathbb{L}}

\newcommand{\mW}{\mathbb{W}}

\newcommand{\mG}{\mathbb{G}}

\newcommand{\sC}{\mathscr{C}}
\newcommand{\sL}{\mathscr{L}}
\newcommand{\sS}{\mathscr{S}}

\newcommand{\cF}{\mathcal{F}}
\newcommand{\cG}{\mathcal{G}}
\newcommand{\cH}{\mathcal{H}}

\newcommand{\cM}{\mathcal{M}}
\newcommand{\cN}{\mathcal{N}}

\newcommand{\cP}{\mathcal{P}\,}

\newcommand{\cR}{\mathcal{R}}
\newcommand{\cS}{\mathcal{S}}
\newcommand{\cT}{\mathcal{T}}

\newcommand{\fD}{\ensuremath{\mathfrak D}}

\newcommand{\bN}{\ensuremath{\boldsymbol{N}}}
\newcommand{\bR}{\ensuremath{\boldsymbol{R}}}
\newcommand{\bI}{\ensuremath{\boldsymbol{I}}}

\newcommand{\bpsi}{\ensuremath{\boldsymbol{\psi}}}

\newcommand{\bE}{\ensuremath{\boldsymbol{E}}}
\newcommand{\bB}{\ensuremath{\boldsymbol{B}}}
\newcommand{\bj}{\ensuremath{\boldsymbol{j}}}
\newcommand{\bx}{\ensuremath{\boldsymbol{x}}}

\newcommand{\by}{\ensuremath{\boldsymbol{y}}}
\newcommand{\ba}{\ensuremath{\boldsymbol{a}}}
\newcommand{\bxi}{\ensuremath{\boldsymbol{\xi}}}
\newcommand{\bn}{\ensuremath{\boldsymbol{n}}}

\newcommand{\bbeta}{\ensuremath{\boldsymbol{\beta}}}
\newcommand{\bdelta}{\ensuremath{\boldsymbol{\delta}}}

\newcommand{\Text}[1]{\text{\textnormal{#1}}}
\newcommand{\Textbf}[1]{\textnormal{\textbf{#1}}}

\newcommand{\supp}{\Text{supp}}
\newcommand{\rect}{\Text{rect}}
\newcommand{\tria}{\Text{tria}}

\newcommand{\Frechet}{Fr\'{e}chet }

\newcommand{\I}{\Text{i}}
\newcommand{\E}{\Text{e}}

\newcommand{\MTEXT}[1]{\;\;\;\;\;\text{#1}\;\;\;\;\;}
\newcommand{\ip}[2]{\langle #1, #2 \rangle} 
\newcommand{\D}{\Text{d}}
\newcommand{\norm}[1]{\| #1 \|} 
\newcommand{\Lp}[2]{L^{#1}(#2)}

\newcommand{\Hs}[2]{H^{#1}(#2)}
\newcommand{\HsR}[1]{\Hs{#1}{\mR^m}}
\newcommand{\Wp}[3]{W^{#1,#2}(#3)}
\newcommand{\WpR}[2]{\Wp{#1}{#2}{\mR^m}}
\newcommand{\WtwoR}[1]{\WpR{#1}{2}}

\newcommand{\seqn}[2]{(#1_#2)_{#2 \in \mN}}
\newcommand{\Ck}[2]{\sC^{#1}(#2)}
\newcommand{\Cinf}[1]{\Ck{\infty}{#1}}
\newcommand{\Ckc}[2]{\sC^{#1}_{\Text c}(#2)}
\newcommand{\Cinfc}[1]{\Ckc{\infty}{#1}}

\newcommand{\closure}[2][3]{%
  {}\mkern#1mu\overline{\mkern-#1mu#2}}

\newcommand{\N}{N}
\newcommand{\PN}{\underline{\N}}
\newcommand{\Pdelta}{\underline{\delta}}
\newcommand{\Pbeta}{\underline{\beta}}

\newcommand{\Or}{\mathcal{O}}
\newcommand{\bDd}[1]{\mathcal{\tilde D}_{#1}}

\newcommand{\bDF}[1]{\mathcal{\tilde D}^{(\text F)}_{#1}}
\newcommand{\bDff}[1]{\mathcal{\tilde D}^{(\infty)}_{#1}}
\newcommand{\bDFlat}[1]{\mathcal{\tilde D}^{(\text F)}_{#1,\overline{2}}}

\newcommand{\bMF}{\mathcal{\tilde M}^{(\text F)}}

\newcommand{\nF}{w^{(\Text{F})}}

\newcommand{\bmF}{\tilde m^{(\Text{F})}}
\newcommand{\CR}{\cR_{\Text c}}
\newcommand{\CRdis}{\cR_{\Text c, \Text{dis}}}
\newcommand{\CF}{\cF_{\Text c}}
\newcommand{\NF}{N_{\Text{F}}}
\newcommand{\cbd}{c_{\beta/\delta}}

\newcommand{\cc}[1]{\closure{#1}}
\newcommand{\adj}[1]{#1^\ast}
\newcommand{\ObjDom}{\Omega_{\Text{obj}}}
\newcommand{\SuppDom}{\Omega_{\Text{supp}}}
\newcommand{\ProjDom}{\Omega_{\Text{proj}}}
\newcommand{\LatDom}{\Omega_{\Text{lat}}}

\newcommand{\Sdash}[1]{\mathscr{S}'(\mR^{#1})}
\newcommand{\Sdashc}[1]{\mathscr{S}_{\Text{c}}'(\mR^{#1})}

\newcommand{\A}{\mathscr{A}}
\newcommand{\FFT}{\Textbf{FFT}}
\newcommand{\DFT}{\Textbf{DFT}}

\newcommand{\id}{\Text{id}}

\DeclareMathOperator*{\diam}{diam}
\DeclareMathOperator*{\argmin}{argmin}

\DeclareMathOperator*{\Int}{int}
\DeclareMathOperator*{\KL}{KL}

\newcommand{\UHP}{\mH_+}

\newtheorem{df1}{Definition}[section]
\newtheorem{dfth1}[df1]{Definition and Theorem}
\newtheorem{th1}[df1]{Theorem}
\newtheorem{alg1}{Algorithm}[chapter]
\newtheorem{lem1}[df1]{Lemma}
\newtheorem{cor1}[df1]{Corollary}
\newtheorem{ex1}[df1]{Example}
\newtheorem{prob}{Problem}[chapter]
\newtheorem{res}{Result}[chapter]

\newtheorem{rem1}{Remark}[section]

\newenvironment{pf}[1][Proof:]{\begin{trivlist}\item[\hskip \labelsep {\bfseries #1}]}{\qed\end{trivlist} } 

\newcommand{\thmref}[1]{Theorem~\ref{#1}}
\newcommand{\lemref}[1]{Lemma~\ref{#1}}
\newcommand{\defref}[1]{Definition~\ref{#1}}
\newcommand{\cref}[1]{Corollary~\ref{#1}}
\newcommand{\sref}[1]{$\S$~\ref{#1}}
\newcommand{\figref}[1]{Figure~\ref{#1}}
\newcommand{\probref}[1]{Problem~\ref{#1}}
\newcommand{\chapref}[1]{Chapter~\ref{#1}}
\newcommand{\exref}[1]{Example~\ref{#1}}

\newcommand{\resref}[1]{Result~\ref{#1}}

\newcommand{\aref}[1]{Appendix~\ref{#1}}
\newcommand{\algref}[1]{Algorithm~\ref{#1}}
\newcommand{\tabref}[1]{Table~\ref{#1}}
\newcommand{\noemph}[1]{{#1}}

 \numberwithin{equation}{section}
 \numberwithin{figure}{chapter}
 \numberwithin{table}{chapter}

\begin{document}
\ThesisAuthor{Simon}{Maretzke}
\PlaceOfBirth{Gifhorn}
\ThesisTitle{Regularized Newton Methods for simultaneous Radon Inversion and Phase Retrieval in Phase Contrast Tomography}
	    {Regularisierte Newton-Verfahren zur simultanen Radoninversion und Phasenrekonstruktion in der Phasenkontrast-Tomographie}
\FirstReferee{Professor\ Dr.\ Tim\ Salditt}
\Institute{Institute for Numerical and Applied Mathematics / \\ Institute for X-Ray Physics}
\SecondReferee{Professor\ Dr.\ Thorsten\ Hohage}
\ThesisBegin{12}{07}{2014}
\ThesisEnd{11}{01}{2015}
\frontmatter
\maketitle
\cleardoublepage
\begin{abstract}
Promoted by the advent of coherent synchrotron light sources, 
phase contrast tomography allows to resolve three-dimensional variations of an unknown sample's complex refractive index
from scattering intensities recorded at different incident angles of an X-ray beam. By  diffractive free-space propagation of the transmitted wave field, this method is sensitive not only to absorption but also to refractive phase shifts induced by the specimen, permitting three-dimensional nanoscale imaging of quasi-transparent samples such as biological cells. However, the reconstruction of the specimen structure from the observed data constitutes an algorithmically challenging nonlinear ill-posed inverse problem, mainly due to the characteristic loss of phase information in the detection of the wave field.

In this work, regularized Newton methods are developed for the solution of this tomographic phase retrieval problem, based on a detailed analysis of its mathematical structure. We consider both the  near-field- or Fresnel regime characterized by a moderate propagation length between sample and detector and the far-field limit of large detector distances, where propagation is governed by the Fourier transform. In the former setting, excellent numerical reconstructions are obtained via the chosen Newton-type approach, supplemented by novel theoretical results stating that measurements from a single detector distance are sufficient to uniquely recover both refraction and absorption of a sample. The proposed algorithm simultaneously performs tomographic- and phase reconstruction, which is found to stabilize the latter by exploiting correlations between the diffraction patterns recorded under different incident angles.

\end{abstract}

\cleardoublepage
\onehalfspacing
\tableofcontents



\begin{nomenclature}
\vspace{-1em}
\subsection*{Physical parameters}
  \begin{tabularx}{\textwidth}{lXl} 
       \textbf{Notation} & \textbf{Description} & \textbf{Definition} \\
    \toprule
       $L$ & Sample thickness & \figref{fig:SetupIdeal} $\;\;\;\;\;$ \\
       $d$ & Sample-detector distance & \figref{fig:SetupIdeal} \\
       $k$ & Wavenumber of the X-rays & \eqref{eq:Helmholtz} \\
       $\Delta x$ & Pixel or voxel size & \sref{SS:ContrastFormationNF}, \sref{SS:DiscrGeneral}  \\
       $\NF$ & Fresnel number & \eqref{eq:FresnelNumber}  \\
       $n$ & Refractive index of the sample, $n= 1- \delta + \I \beta$ & \eqref{eq:Wave} \\
       $N $ & Refractive decrement, $N= 1-n$ & \eqref{eq:ndeltabeta} \\
       $\delta$ & Refractive part of $n$  & \eqref{eq:ndeltabeta} \\
       $\beta$ & Absorptive part of $n$  & \eqref{eq:ndeltabeta} \\
       $\tilde \Psi$ & Paraxial wave field (envelope) & \eqref{eq:DefParaxWave}  \\
       $P$ & Probe- or illumination function & \eqref{eq:ProjSol}  \\
       $O_{(0)}$ & Object transmission function (normalized) & \eqref{eq:OTF}, \eqref{eq:ForwardOpNF}  \\
       $I$ & Detected intensities  & \eqref{eq:ContrastDetectorImage}  \\
        $\theta$ & Tomographic incident angle  & \figref{fig:SetupIdeal}, \sref{SS:TomoRadon}  \\
  \end{tabularx}
  
  \vspace{-.5em}
\subsection*{Spaces and Domains}
  \begin{tabularx}{\textwidth}{lXl} 
       $\supp(u)\;\;\;\:$ & Support of a function or distribution $u$& \sref{SS:Lp}, \sref{SS:TempDistri} $\;\;\;\;\;$\\
       $\Lp{p}{\Omega}$ & Lebesgue $L^p$-space on a set $\Omega \subset \mR^m$ & \sref{SS:Lp} \\
       $\sS(\mR^m)$ & Schwartz space in $\mR^m$ & \sref{SS:Schwartz} \\
       $\sS'(\mR^m)$ & Tempered distributions, dual of $\sS(\mR^m)$ & \sref{SS:TempDistri} \\
       $\Sdashc{m}$ & Compactly supported distributions & \sref{SS:TempDistri} \\
       $H^s(\mR^m)$ & Sobolev space of order $s \geq 0$ & \eqref{eq:DefHs} \\
       $\mX$ & Object space, domain of a forward operator & \sref{SS:IRNM-Motivation}  \\
       $\mY$ & Image- or data space of a forward operator & \sref{SS:IRNM-Motivation}  \\
       $\ObjDom$ & Cylindrical object domain, $\supp(N) \subset \ObjDom$ & \eqref{eq:ObjDomain} \\
       $Z^m$ & Image space of $\CF,\CR$, $Z^m = [0;2\pi) \times \mR^{m-1}$ & Def. \ref{def:Radon} \\
  \end{tabularx}

\subsection*{Operators}
  \begin{tabularx}{\textwidth}{lXl} 
       \textbf{Notation} & \textbf{Description} & \textbf{Definition} \\
    \toprule
	$\Re(f)$  & Real part (pointwise for vectors/functions $f$) &  \\
	  $\Im(f)$ & Imaginary part (pointwise) &  \\
       $\norm{\cdot}_{2} $ & Euclidean  norm in $\mC^m$ &  \\
       $ \norm{\cdot}_{\infty}$ & Maximum norm in $\mC^m$ &  \\
       $\ip{\cdot}{\cdot}_{\mH}$ & Inner product of a Hilbert space $\mH$ & \sref{S:OpsAdjoints} \\
       $\norm{\cdot}_{\mV}$ & Norm of a vector space $\mV$ & \sref{S:OpsAdjoints} \\
       $\norm{\cdot}$ & Operator norm or discrete object norm  & \eqref{eq:defOpNorm}, \eqref{eq:ObjScaleNorm} \\
       $F$ & Forward operator & \eqref{eq:ForwardOpNF}, \eqref{eq:ForwardOpFF} \\
       $F'[N]$ & \Frechet derivative & \eqref{eq:FrechetForwOp}, Def. \ref{def:FrDer} $\, $ \\
       $\bDF{d}$ & Fresnel propagator (for the envelope field $\tilde \Psi$) & \eqref{eq:EnvelopeFresnelProp}  \\
       $\bMF_d$ & Fresnel multiplication in Fourier space & \eqref{eq:EnvelopeFresnelProp}  \\
       $\cF$ & Fourier transform & Def. \ref{def:FT}  \\
       $\cR$ & 2D Radon transform & Def. \ref{def:Radon}  \\
       $\CF,\CR$ & Cylindrical Fourier- and Radon transforms & \eqref{eq:CylTransforms}  \\
  \end{tabularx}
  
\subsection*{Sub- and Superscripts}
  \begin{tabularx}{\textwidth}{lXl}
       $F_d, I_d$ & Near-field operator/intensity (detector dist. $d$)  & \eqref{eq:ForwardOpNF} \\
       $F_\infty, I_\infty$ & Far-field operator/intensity (detector dist. $ \infty$) & \eqref{eq:ForwardOpFF} \\
       $N^\dagger, I^\dagger$ & Exact object and corresponding intensity data & \sref{S:GenIll-posed} \\
       $I^{\Textbf{err}}$ & Noisy intensity data & \sref{S:GenIll-posed} \\
       $F_{\Text{dis}}, \mX_{\Text{dis}}\;\;\,$ & Discretized operator or space &  \sref{SS:DiscrGeneral} \\
       $\cF_{\underline m}(f)$ & Transform w.r.t. variables $x_1,\ldots x_m$  of $f$ & Def. \ref{def:FT} \\
       $\cF_{\overline m}(f)$ & Transform w.r.t. variables $x_m, x_{m+1}, \ldots $ of $f$ & Def. \ref{def:FT} \\
        $T^\ast, f^\ast$ & Adjoint of $T$ or Schwarz-reflected function $f$ & Def. \ref{def:OpsAdjoints}, \eqref{eq:DefSchwarzRefl} \\
	$\mV_{\mR}$ & Real analogue of a complex Banach space $\mV$ & \sref{S:OpsAdjoints} \\
  \end{tabularx}

\end{nomenclature}

\mainmatter   


\begin{chapter}{Introduction}\label{C:Intro}

 Over the past decades, classical X-ray tomography, better known as \noemph{computed tomography} (CT), has become a workhorse of noninvasive medical diagnosis. 
 CT scanners measure the attenuation experienced by X-rays traversing a whole patient or single organs under different incident angles \cite{Cormack1963CT,Hounsfield1973CT}.
 Based on the mathematical theory of the \emph{Radon transform} and its inversion \cite{Radon1917} laid down almost a hundred years ago, three-dimensional images are reconstructed from these measurements.
 
 Notably, the small wavelengths of hard X-rays in the order of one nanometer or less would also allow for tomographic imaging of \noemph{micro-} or even \noemph{nanoscale} structures according to Rayleigh's criterion for optical resolution \cite{Rayleigh1879Resolution}. 
 Indeed, so-called micro-CT scans have been successfully applied to \noemph{in vivo} tomography of small animals, resolving features down to a size of less than 100 micrometers \cite{Holdsworth2002InVivoMicroCT,Badea2008InVivoMicroCT}. However, when it comes to imaging of even smaller objects, such as single biological cells, an in principal desirable property of X-rays becomes problematic: if it was not for the partial \emph{transparency} of biological tissues for X-rays, no residual intensities would ever be measured behind a patient from which CT images could be reconstructed. For instance, the opacity of a human body for \emph{visible} light evidently rules out transmission radiography with this type of radiation. However, the desirable transparency of an entire human torso necessarily implies that the \noemph{absorption} of X-rays by a single bacterium of one micrometer thickness is far too small to yield sufficient contrast in the recorded intensities behind the specimen.
 
 One might thus conclude that nanoscale light-element samples are simply too transparent to lend themselves to X-ray tomography. Yet, note that \noemph{non-absorbing} does not mean \noemph{non-interacting}: everyday experiences of the \emph{refraction} of visible light by glass or water teach us that even perfectly transparent materials leave some traces in transmitted radiation. The imprint of sample information within the traversing X-rays is described by its spatially varying \emph{refractive index}, typically written in the form $n = 1 - \delta + \I \beta$ in the considered hard X-ray regime. Indeed, it turns out that the induced refraction, governed by the quantity $\delta \sim 10^{-6}$, is usually two to three orders of magnitude larger than the absorptive part $\beta$ for typical biomolecules \cite{BartelsDiss,Henke1993TypicalDeltaBeta,Mayo2002quantitative}. Accordingly, refraction-sensitive imaging methods may achieve considerably improved contrast compared to purely absorption-based radiography by classical CT, permitting the desired resolution of nanoscale quasi-transparent structures. As the governing parameter $\delta $ manifests itself in the form phase shifts within the transmitted X-ray wave field, this approach is denoted as \emph{phase contrast imaging}. Its applicability to two-dimensional imaging of micro- and nanoscale specimen, measuring the \emph{projection} of $\delta $ along the X-ray's incident direction, has been experimentally demonstrated for both synthetic- and biological structures down to the size of bacterial cells \cite{Wilkins1996,Cloetens1996PCI,Pogany1997noninterferometric,Mayo2002quantitative,Miao2003EColiFarfield}. Performing tomographic reconstruction via the aforementioned Radon inversion from an ensemble of projection images obtained for different X-ray incident angles, \emph{phase contrast tomography} permits 3D imaging of such samples by resolving the refractive decrement $\delta$ \cite{Cloetens1999,Mayo2003PCT}.
 
 However, there are two major difficulties associated with this fascinating imaging technique: for once, since it is based on interference effects, the approach typically requires highly \emph{coherent} X-rays as provided by third generation synchrotron light sources or free-electron lasers \cite{Nugent2010coherent} but not by state-of-the-art lab sources. On the other hand, the induced refractive phase shifts on which the imaging method is based cannot be inferred directly from the detected intensities behind the sample - as opposed to the absorption in CT. The fact that available X-ray detectors may only measure wave intensities but not the phase of the incident radiation indeed gives rise to a \emph{phase retrieval problem}, i.e.\ the lost phase information has to be implicitly recovered in the image reconstruction from the observed data. Most prominently, this physical limitation of the measurement process implies that the phase shifts encoding the refractive index of the traversed sample are completely invisible if the detector is placed directly behind the specimen, see \cite[sec. 4.4]{PaganinXRay}. Only due to the \emph{diffraction} experienced by the transmitted X-rays as they propagate to a distant detector are the imprinted phase perturbations in the wave field partially encoded into observable intensities \cite{Nugent1996Propagationbased,Paganin1998}. Alternatively, the required phase-sensitivity of the measurements may be achieved by interferometric techniques \cite{Bonse1965interferometric,Momose1995interferometric,Wilkins1996}.
 
 In this work, we study three-dimensional imaging by propagation-based phase contrast tomography. Discontinuity of the tomographic Radon inversion, i.e.\ its noise-amplifying property, but also the characteristic loss of phase information, which allows for possibly ambiguous or unstable reconstructions, render this an ill-posed problem in the sense of \citet{Hadamard1902}. We distinguish between the \emph{near-field} case of moderate propagation lengths and imaging from intensity data recorded in the \emph{far-field} limit of large distances between the specimen and the detector. The latter setting is better known as the established technique of \emph{coherent diffractive imaging} \cite{Miao1998FourierNonPeriodic,MiaoNature1999CDIFirstExp}, for which a vast amount of theoretical studies on phase retrieval ambiguities have been published ever since the pioneering works of \citet{Akutowicz1956I,Akutowicz1957II} and \citet{Walther1963}. Most importantly, non-unique reconstructions of two-dimensional images, which are not related by simple geometrical transformations, are found to be ``pathologically rare'' \cite{Barakat1984,Fienup19782DPhaseRetrFeasible}. Reconstruction algorithms based on convex optimization have been designed to account for the remaining ambiguities \cite{Fienup1982HIO,Fienup1986algorithm,Luke2005RAAR}.
 
 In near-field imaging, uniqueness theory is on a less advanced stage: it has been shown that projections of the complex refractive index $n = 1 - \delta + \I \beta$ may be recovered uniquely from intensities recorded at \emph{two} different detector distances \cite{Jonas2004TwoMeasUniquePhaseRetr}. On the other hand, it is commonly argued \cite{Nugent2007TwoPlanesPhaseVortex,Burvall2011TwoPlanes} that a single measurement is not sufficient for this although numerical results for near-field phase contrast tomography suggest that unique reconstructions might be possible in this setting \cite{Ruhlandt2014}. Uniqueness of near-field phase retrieval is therefore analyzed in this work.
 
 The focus, however, lies on designing tailored reconstruction algorithms for phase contrast tomography. Ill-posedness and nonlinearity of the problem are accounted for by constructing \emph{regularized Newton methods} \cite{Bakushinskii1992IRNM}, the potential of which for (non-tomographic) phase retrieval has been demonstrated in \cite{BartelsDiss,Hohage2013}. Retaining the nonlinearity, this approach promises a larger regime of applicability and increased accuracy compared to direct near-field reconstruction techniques based on the contrast transfer function  \cite{Cloetens1999,Cloetens1999Diss,BartelsDiss} or transport-of-intensity-equations \cite{Reed1983TIEPhaseRetr,Nugent1996Propagationbased,Wilkins1996,Paganin1998}. The latter are valid only in the limits of weakly scattering specimen or small propagation distances, respectively. Moreover, these techniques often require measurements from up to four detector distances \cite{Krenkel2014BCAandCTF} and typically assume some coupling between $\delta$ and $ \beta$, e.g.\ proportionality or vanishing absorption $\beta = 0$. The principal motivation for this work is to overcome these limitations, ideally reconstructing both $\delta$ and $ \beta$ from intensities measured at a single detector distance. An essential feature of the Newton-type reconstruction methods presented here is that they perform the Radon inversion and phase retrieval \emph{simultaneously}. Thereby, tomographic consistency conditions are imposed already in the phase reconstruction which has been shown to promote stability and accuracy in other algorithms \cite{Chapman2006,Barty2008ceramicfoam,Ruhlandt2014}. The idea is simply to incorporate the a priori knowledge that all diffraction patterns result from one and the same specimen - merely illuminated under different angles. 
 
 In \chapref{C:PhysProb}, a  physical model for phase contrast tomography is developed as considered in this work, describing the encoding of the specimen structure in measurable intensities. The derived mathematical formulation is further analyzed in \chapref{C:Analysis}, investigating ill-posedness of the inverse reconstruction problem with special focus on uniqueness of the involved phase retrieval. \chapref{C:NumMeth} is dedicated to the construction of regularized Newton methods for far-field- and near-field tomography, yielding our principal \algref{alg:PCT}. Its performance is investigated in \chapref{C:NumRes}, discussing numerical reconstruction results for both simulated and experimental data. \aref{C:MathPre} introduces the basic mathematical concepts and tools applied in this work. 
 
\end{chapter}


\begin{chapter}{Physical Model}\label{C:PhysProb}

  
  This chapter introduces the physical theory of image formation in propagation-based phase contrast tomography, closely following the presentation in \cite{PaganinXRay}. Starting with a brief discussion of an exemplary experimental setup and its idealized analogue considered in this work, a mathematical description of the problem is derived by reviewing the theory of monochromatic electromagnetic waves in inhomogeneous media and introducing different approximations. The overall objective is the formulation of a forward operator, mapping the spatial structure of an unknown specimen onto the diffraction patterns which are detected under different tomographic incident angles of the X-rays. An introduction of the mathematical tools used in this work can be found in \aref{C:MathPre}.

\begin{section}{Experimental Realization and Idealized Model} \label{S:RealIdeal}
  
  \figref{fig:SetupExp} shows a sketch of the GINIX setup \cite{Kalbfleisch2011GINIX} (G\"ottingen Instrument for Nano-Imaging with X-rays) as an exemplary experimental realization of phase contrast tomography: an undulator forces accelerated free electrons from a synchrotron storage ring onto wiggling trajectories by alternating dipole magnets, resulting in the emission of highly brilliant and coherent X-rays. A cascade of optical elements controls width, spectrum and intensity of the incident beam. The latter is focused by a pair of Kirkpatrick-Baez mirrors whose elliptical profiles define the focal point. Optionally, X-ray waveguides may be placed in the focal plane in order to improve the coherence of the illumination, see \cite[sec. 4.3.3]{BartelsDiss}. The emanating cone beam from behind the focus illuminates an unknown specimen. By interaction with the incident radiation, the sample structure is encoded in the scattering intensities measured at a distant charged-coupled device (CCD) detector. An evacuated flight tube in the beam line reduces undesirable absorption and scattering from air or residual particles as the transmitted X-rays propagate to the detector. Rotation of the specimen holder about a lateral axis yields diffraction patterns for different incident angles of the X-ray beam, encoding the three-dimensional structure of the object in question. For details concerning the experimental setup and the involved optical components, see for instance \cite{Nielsen2011ModernXray,BartelsDiss,PaganinXRay}.
  
  The physical model of phase contrast tomography studied in this work is based on the \emph{idealized} setup sketched in \figref{fig:SetupIdeal}. In particular, it is assumed that the incident X-rays are perfectly monochromatic, i.e.\ are of a fixed wavenumber $k$, and that their propagation and interaction with the sample are governed by classical electrodynamics in a medium of spatially varying refractive index $n$. In the following sections, these idealizations are supplemented with further approximations in order to obtain a mathematical description for the encoding of structural object information in the detected intensity data. Although the physical setting in \figref{fig:SetupIdeal} is three-dimensional, we consider the more general case of $m \in \mN$ lateral dimensions, denoted by $\bx$, plus the axial $z$-direction. 

    \begin{figure}[htb!]
	  \centering
	  \subfloat[Experimental realization (GINIX setup): an undulator emits highly brilliant and coherent X-rays which are collimated onto a focal point by a pair of elliptically-shaped Kirkpatrick-Baez mirrors, illuminating the specimen at 1 or 2 by a cone beam. This gives rise to diffraction patterns observed at a distant detector. (Source: \cite{MartinMaster,Kalbfleisch2011GINIX}, modified)]{\includegraphics[width=.9\textwidth]{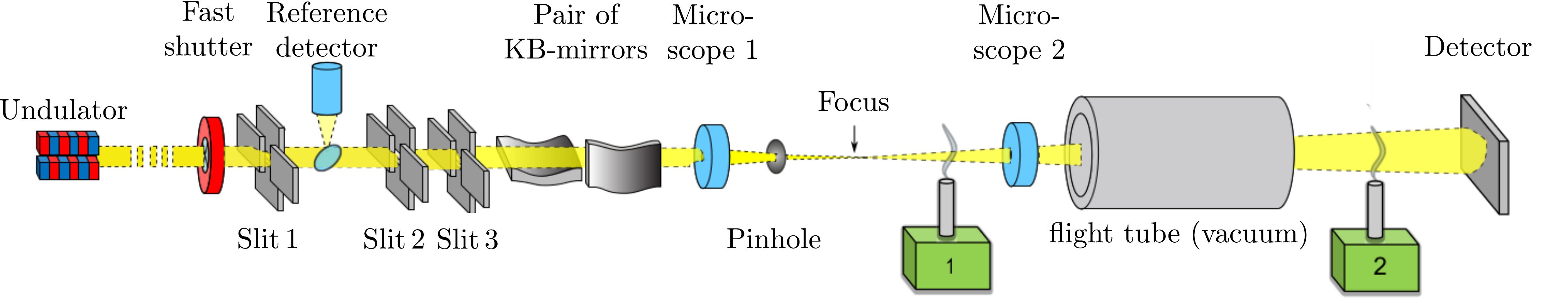}\label{fig:SetupExp}}
	  \\
	  \subfloat[Idealized model: Incident monochromatic electromagnetic illuminate a compact specimen characterized by a spatially varying refractive index $n = 1-\delta + \I \beta$. The interaction results in perturbed wave fronts in the \emph{exit-surface} $E_0$ which manifests themselves in the intensity profiles recorded in the \emph{detector plane} $E_d$. (Source: \cite{AikeMaster}, modified)]{\includegraphics[width=.9\textwidth]{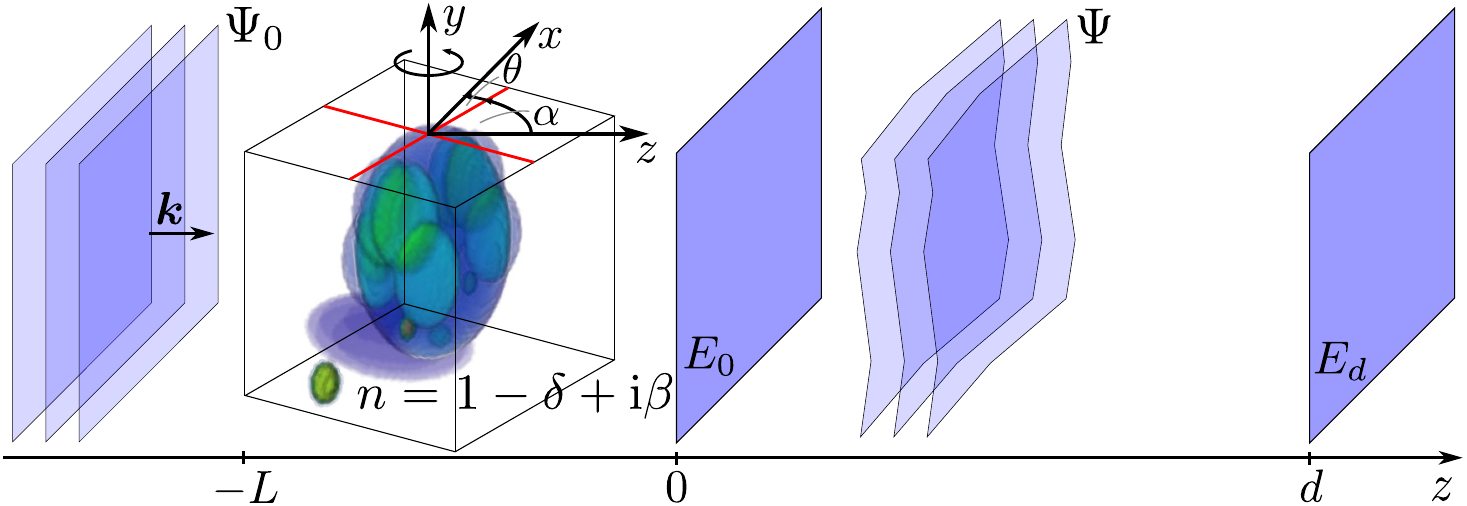}\label{fig:SetupIdeal}}
	  \caption{Exemplary experimental setup for propagation-based phase contrast tomography and its idealized analogue considered in this work.}
	  \label{fig:Setup}
  \end{figure}
  
\end{section}

\begin{section}{Helmholtz Equation and Paraxial Approximation} \label{S:WaveEqs}

	In order to derive a model for X-ray propagation and interaction with matter, we consider Maxwell's equations for isotropic, non-conducting, linear materials in the absence of net electric charges or currents:
	\begin{subequations} \label{eq:Max}
	\begin{align}
	 \nabla \cdot (\varepsilon\boldsymbol{E}) &= 0 \label{eq:Max1} \\
	  \nabla \cdot \boldsymbol{B}  &= 0 \label{eq:Max2} \\
	  	  \nabla \times \boldsymbol{E} + \partial_t \bB  &= 0 \label{eq:Max3}\\
	  	  \nabla \times  \frac 1 {\varepsilon\mu}  \boldsymbol{B} - \partial_{t} \boldsymbol{E} &= 0 \label{eq:Max4} 
	\end{align}
	\end{subequations}
	Under the assumption that the magnetic permeability $\mu $ and dielectric permittivity $\varepsilon$ are slowly varying on lengthscales of the electromagnetic fields, so that derivatives of the material properties can be neglected against those of $\bE$ and $\bB$, \eqref{eq:Max} yields \emph{wave equations} of the form \cite[pp. 66-69]{PaganinXRay} 
	\begin{equation}
	    \left( \varepsilon \mu \partial_t^2  - \nabla^2 \right) \bE  = 0 \MTEXT{and} \left( \varepsilon \mu \partial_t^2  - \nabla^2 \right) \bB  = 0 \label{eq:WaveEH}
	\end{equation}
	It can be shown that all six components may be parametrized by a single scalar but \emph{complex-valued} function $ \Phi$, the squared modulus of which gives the intensity of the electromagnetic field \cite[sec. 8.4]{BornOptics}. According to \eqref{eq:WaveEH}, $\Phi$ is governed by the wave equation
	\begin{equation}
	    \left( \frac{n^2}{c_0^2} \partial_t^2  - \nabla^2 \right) \Phi  =  0 \label{eq:Wave},
	\end{equation}
	where the \emph{refractive index} $n := \sqrt{\varepsilon \mu /(\varepsilon_0 \mu_0)} =  c_0 \sqrt{\varepsilon \mu }  $ has been introduced and $c_0$ denotes the speed of light in vacuum.
	
	Since we consider monochromatic X-rays, the time-dependence in \eqref{eq:Wave} can be eliminated by an ansatz of the form $\Phi(t,\cdot) = \E^{\I  \omega t}  \Psi$. Inserting this ansatz into \eqref{eq:Wave} and defining $k := \frac \omega {c_0}$, we find that $\Psi$ is described by the \emph{Helmholtz equation}
	\begin{equation}
	  \left(  \nabla^2  + n^2 k^2   \right) \Psi  = 0. \label{eq:Helmholtz}
	\end{equation}
	
	In well-controlled experimental settings like those outlined in \sref{S:RealIdeal}, wave dynamics is often strongly anisotropic being characterized by a predominant propagation along the optical $z$-axis up to small angular deviations.
	We exploit this by making the \emph{paraxial approximation}, corresponding to an ansatz of the form
	\begin{equation} \Psi(\bx, z) = \E^{\I k z} \tilde \Psi(\bx, z), \label{eq:DefParaxWave} \end{equation}
	where the envelope $\tilde \Psi$ is assumed to be slowly varying on axial lengthscales $1/k$. Inserting the ansatz into \eqref{eq:Helmholtz}, the contribution $ \partial_z^2 \tilde \Psi$ may thus be neglected against higher orders in $k$, leading to the \emph{paraxial Helmholtz equation}
	\begin{equation}
 	  \left(  2\I k  \partial_{z} +  k^2 (n^2 -1) + \nabla_{\bx}^2   \right) \tilde \Psi = 0 \label{eq:ParaxHelm}  
 	\end{equation}
 	where $\nabla_{\bx}^2$ denotes the Laplacian in the lateral coordinates $\bx$.
 	
 	From hereon, we restrict to paraxial waves and thus always consider the envelope $\tilde \Psi$ as the governing field. We will refer to it as the \emph{wave field} bearing in mind that it has to be supplemented with the plane wave factor $\E^{\I k z}$ (and the time-harmonic one $\E^{\I  \omega t}$) in order to obtain the physical ``waves''. Note that $|\tilde \Psi|^2$ still provides a measure for wave intensities as the supplementary factors are of constant modulus 1.
	
\end{section}

\begin{section}{Free-Space Propagation} \label{S:FreeProp}

	In the sequel, the theory of $\S$\ref{S:WaveEqs} is applied to derive an expression for the propagation of paraxial monochromatic waves \emph{in vacuum}. With respect to the idealized setup in \figref{fig:SetupIdeal}, this will yield a relation between the propagated wave field in the detector plane $E_d$ and the \emph{contact image} at the exit-surface $E_0$.

	\begin{subsection}{Fresnel Propagator} \label{SS:ParaxialFresnel}
	 
 	 Vacuum is characterized by a constant refractive index $n=1$. A model for free-space propagation of paraxial waves is obtained by taking the lateral Fourier transform $\cF_{\underline{m}}$ of \eqref{eq:ParaxHelm}, corresponding to the \emph{angular spectrum approach} discussed in \cite{PaganinXRay}. Writing 
 	 \begin{equation}
 	 \hat \Psi_{\bxi} (z) := \cF_{\underline{m}} ( \tilde \Psi ) (  \bxi, z )
 	 \end{equation}
 	 and noting that $\cF_{\underline{m}} ( \nabla^2_{\bx} \tilde \Psi ) ( \bxi, z )   = - \bxi^2 \hat \Psi_{\bxi}(z)$ according to \eqref{eq:FTDeriv}, this yields
 	\begin{equation}
 	   \left(  2 \I k \partial_z  - \bxi^2    \right)  \hat \Psi_{\bxi} = 0  \MTEXT{and thus} \hat \Psi_{\bxi}(z) = \underbrace{\exp\left( - \frac{ \I     \bxi^2 z}{2 k}  \right)}_{=: \bmF_z(\bxi)} \hat \Psi_{\bxi}(0)  \label{eq:AngSpec}
 	 \end{equation}
 	 for all Fourier modes $\hat \Psi_{\bxi}$. By inverse Fourier transform $\{ \hat \Psi_{\bxi} \}_{\bxi \in \mR^m} \stackrel{\cF_{\underline{m}}^{-1}} \mapsto \tilde \Psi$, \eqref{eq:AngSpec} yields an expression for the propagation of general paraxial wave fields
 	\begin{equation}
	 \tilde \Psi_d = \bDF{d}( \tilde \Psi_0 ) := \cF^{-1}  \bMF_d \cF ( \tilde \Psi_0 ) \MTEXT{with} \bMF_d: \hat f \mapsto \bmF_d \cdot \hat f, \label{eq:EnvelopeFresnelProp} 
	\end{equation}
	where  we have set $\tilde \Psi_z := \tilde \Psi(\cdot, z)$ for notational convenience. We denote $\bDF{d}$ as the \emph{Fresnel propagator} although this term is typically used for the corresponding propagator $\E^{\I k d} \bDF{d}$ of the total wave field given by \eqref{eq:DefParaxWave}.
	
	According to \eqref{eq:EnvelopeFresnelProp}, the propagated wave field at $z=d$ in \figref{fig:SetupIdeal} is related to the exit-wave at $z=0$ 
	by a simple multiplication in Fourier space. 
	This remains valid if the paraxial approximation is relaxed and the complete Helmholtz equation \eqref{eq:Helmholtz} is retained in the derivation. The Fresnel propagator can then be derived by expanding the unitary propagation factor $\exp \left( \I d  (k^2-\bxi ^2)^{\frac 1 2} \right)$, which is obtained in this case, to quadratic order in $\bxi / k$. See \cite[pp. 15 f.]{PaganinXRay} for details. Accordingly, the paraxial approximation provides an accurate description of free-space propagation  whenever the Fourier spectrum $\cF( \Psi)$ of the \emph{total} wave field is strongly peaked around the dominant wavevector $k \boldsymbol e_z$.
	
	
	Moreover, note that
	\begin{equation}
	 \bDF{d}: \Lp 2 {\mR^{m }} \to \Lp 2 {\mR^{m }}
	\end{equation}
        defines a \emph{unitary operator} by \eqref{eq:EnvelopeFresnelProp} as a composition of the unitary maps $\cF, \cF^{-1}$ and a multiplication with a function of modulus 1 (see \exref{ex:Adjoints}-\Text{(b)} and \sref{S:FT}). Physically, this property corresponds to \emph{energy conservation} of the wave field as it is propagated between lateral planes over a distance $d$.
 	
	\end{subsection}

	\begin{subsection}{Convolution Formulation} \label{SS:ConvForm}
	
	According to \sref{SS:ParaxialFresnel}, propagation of paraxial waves in vacuum corresponds to a multiplication with the factor $\bmF_d $ in Fourier space. In a distributional sense, it holds that \cite[pp. 11 f.]{BartelsDiss}
	\begin{equation*}
	 \cF^{-1} \left( \bmF_d   \right) ( \bx )  = \underbrace{\E^{- \I m \pi / 4 }}_{=:\nu^{\frac m 2}} \left(\frac{k }{ d} \right)^{\frac m 2} \underbrace{\exp\left(   \frac{ \I k  \bx^2}{ 2d }  \right)}_{=:\nF_d(\bx)}.
	\end{equation*}
	This expression can be used to rewrite the Fresnel propagator defined in \eqref{eq:EnvelopeFresnelProp}, using the convolution theorem \eqref{eq:ConvThm} and the definition of the Fourier transform:
	
	\begin{align}
	 \bDF{d}(\tilde \Psi_0)(\bx ) &= \cF^{-1} \left( \bmF_d \cdot \cF (\tilde \Psi_0)  \right) (\bx ) = (2 \pi)^{-\frac m 2 }\left( \cF^{-1} \left(  \bmF_d \right)   \ast  \tilde \Psi_0   \right) (\bx)  \nonumber \\
	 &=  \left(\frac{ \nu k }{ 2\pi d} \right)^{\frac m 2} \exp \left(  \frac{\I k \bx^2 }{2 d} \right) \int_{\mR^{m}} \left[ \tilde \Psi_0(\by) \cdot  \exp \left(  \frac{\I k \by^2 }{2 d} \right) \right]  \;  \exp \left(- \frac{\I k \bx \cdot \by }{ d} \right) \; \D y \nonumber \\
	 &= \left(\frac{ \nu k }{ d} \right)^{\frac m 2} \nF_d(\bx) \cdot \cF  \left(     \nF_d \cdot \tilde \Psi_0  \right)  \left(  \frac{ k \bx }{ d } \right). \label{eq:FresnelConv}
	\end{align}
	 According to \eqref{eq:FresnelConv}, propagation of a wave field is - up to a multiplication with the unitary factor $\nF_d$ and rescaling - essentially achieved by a lateral Fourier transform.
	 
	 \end{subsection}

	\begin{subsection}{Fresnel Number and Far-Field Limit} \label{SS:FresnelFarfield}
	In the following, it is assumed that the wave field $\tilde \Psi_0$ to be propagated is non-negligible only within a region of lengthscale $b$, e.g.\ due to confinement by a suitable aperture, such that the dimensionless \emph{Fresnel number}
	\begin{equation}
	 	N_{\text F} := \frac{b^2}{\lambda d}  =   \frac{k b^2}{2 \pi d} \label{eq:FresnelNumber}
	\end{equation}
	is $\ll 1$. Then for all relevant contributions to the Fourier integral in \eqref{eq:FresnelConv}, the unitary factor $\nF_d$ is close to unity and thus may be suppressed. Accordingly, free-space propagation in this limit reduces to a Fourier transform of the original wave field - up to a prefactor of constant modulus:
	\begin{equation}
	 \bDF{d}(\tilde \Psi_0)(\bx ) \stackrel{N_{\text F} \ll 1 } \approx \left(\frac{ \nu k }{ d} \right)^{\frac m 2} \exp\left( \frac{  \I k \bx^2}{ 2d } \right)  \cdot \cF \left(\tilde \Psi_0 \right) \left(  \frac{ k \bx }{ d } \right) =:  \bDff{d}(\tilde \Psi_0)(\bx ) \label{eq:Fraunhofer}
	\end{equation}
	The approximation of \eqref{eq:Fraunhofer} is denoted as the \emph{far-field} or \emph{Fraunhofer diffraction} formula \cite[p. 12]{BartelsDiss}. Notably, its validity is not restricted to the regime of paraxial waves considered here.
	In the general case of Fresnel numbers $N_{\text F} \gtrsim 1 $, the complete Fresnel propagator, given by \eqref{eq:EnvelopeFresnelProp}  or \eqref{eq:FresnelConv}, has to be retained. Throughout this work, this setting is referred to as the \emph{near-field-}, \emph{holographic-} or \emph{Fresnel regime}.
	
	Taking $b$ as an arbitrary reference lengthscale and setting $\bxi' := \frac b {2\pi} \bxi$, the Fresnel number yields a dimensionless form of the propagating factor $\bmF_d$:
	\begin{equation}
	 \bmF_d(\bxi) =  \exp\left( -\frac{\I   d \bxi^2  }{2 k }   \right) = \exp\left( -\I  \frac{2\pi^2 d }{kb^2} \cdot  \frac{b^2\bxi^2}{4\pi^2}  \right) = \exp\left( - \frac{ \I  \pi \bxi'^2 }{\NF}  \right) \label{eq:FresnelNumNoDim}
	\end{equation}
	From \eqref{eq:FresnelNumNoDim}, it can be seen that wave structures of lengthscale $b$ or larger will essentially be preserved under propagation for $\NF \gg 1$ as the exponential factor is close to unity for the corresponding Fourier frequencies. For $\NF \lesssim 1$, on the other hand, features of size smaller or equal to $b$ get significantly distorted as they are propagated over a distance $d$. The Fresnel number thus constitutes a dimensionless measure for propagation effects onto a wave field.

	\end{subsection}
%

\end{section}

\begin{section}{Interaction with Matter} \label{S:Interact}

      By the description of vacuum wave propagation established in \sref{S:FreeProp}, we obtain a map relating the contact image in \figref{fig:Setup} to the resulting wave field at the detector. In the following section, this description is supplemented with a model for the interaction between the incoming radiation and the sample, governing the encoding of the latter's spatial structure in measurable data.

      \begin{subsection}{Refractive Index in the X-Ray Regime} \label{SS:XrayRefrIdx}
      
      Classically, electrons in matter perceive the presence of a traversing electromagnetic wave in the form of an alternating electric field. Bound in the electrostatic potential of the cores, these can therefore be viewed as forced harmonic oscillators. The oscillating electrons in the material - representing accelerated charged particles - in turn emit electromagnetic radiation themselves which is superimposed upon the incident electromagnetic wave. In the picture of a forced harmonic oscillator, it can be seen that both the amplitude and the phase shift of the radiated field depend on the quotient $\omega / \omega_0$ of the driving frequency $\omega$ versus the resonance frequency $\omega_0$ of the bound electron. This frequency dependence of the material response manifests itself globally in the refractive index $n$. See \cite[sec 2.3]{BornOptics} for details.
      
      For hard X-rays with photon energies $\sim \unit[10]{keV}$ propagating in light-element materials, the excitation is usually well above resonance, i.e.\ $\omega / \omega_0 \gg 1$ \cite{Cloetens1999Diss}. For this reason, the forced oscillations carry a phase shift of $\pi$ with respect to the excitation, which leads to an overall phase difference of $\frac{3\pi} 2 $ in the irradiated response wave. Superimposing this contribution upon the incident field gives rise to a phase speed of the total wave that is greater than in vacuum, corresponding to a refractive index
      $n$ with a real part that is slightly smaller than unity \cite[sec. 3.1]{Nielsen2011ModernXray}. We account for this by writing 
      \begin{equation}
       n = 1 - \delta + \I \beta =  1 - \N \label{eq:ndeltabeta}
      \end{equation}
      where $\delta$ and $\beta$ are real and non-negative. The imaginary part $\beta$ parametrizes absorption by the material in addition to the refractive effects.
      For the model protein $\text H_{50} C_{30} N_9 O_{10} S$, the refractive index at $\unit[13.8]{keV}$ is given by \cite{Henke1993TypicalDeltaBeta}
      \begin{equation*}
      \delta \approx 1.6\cdot 10^{-6} \MTEXT{and} \beta \approx 1.9\cdot 10^{-9}.
      \end{equation*}
       These orders of magnitude are typical of light element samples such as biological tissues. The small $\beta$ corresponds to negligible absorption on micrometer-scales. Therefore, solely absorption-based imaging methods - such as classical computed tomography (CT, see for instance \cite{KakTomography}) - are unsuited for X-ray nanoscopy. On the contrary, the refractive decrement $\delta$, giving rise to phase shifts rather than to attenuation of the transmitted radiation, is usually neglected in these techniques. However, it is typically about 1000 times larger and may thus yield reasonable contrast for nanoscale structures.
      
      Note that the adopted description in terms of the refractive index $n$ may only account for \emph{coherent} scattering processes, i.e.\ such which do not change the frequency $\omega$ of the incident radiation. This neglects in particular incoherent quantum effects like inelastic \emph{Compton scattering}. However, quantum theory indeed shows that these are of lesser significance for X-ray phase contrast imaging \cite{Slowik2013QuantumPhaseContrast}. Based on the \emph{classical} model outlined above, the refractive decrement $\delta$ introduced in \eqref{eq:ndeltabeta} can be related to the electron density $\rho_{\text e}$ in the scattering medium if the excitation frequency $\omega$ is well above resonance:
      \begin{equation}
       \delta \approx \frac{2 \pi \rho_{\text e} r_0}{k^2} \label{eq:deltaVSrho}
      \end{equation}
      Here, $r_0 = \unit[2.82\cdot 10^{-15}]{m}$ denotes the Thompson scattering length and $k$ is the wavenumber of the incident wave \cite[p. 63]{Nielsen2011ModernXray}. According to \eqref{eq:deltaVSrho}, refraction-sensitive imaging methods such as phase contrast tomography are \emph{quantitative} in that they measure the physical observable $\rho_{\text e}$. Indeed, it can be shown \cite{Nielsen2011ModernXray} that 
      \begin{equation} \delta \propto Z \omega^{-2} \MTEXT{whereas} \beta \propto Z^4 \omega^{-4}, \end{equation} 
      $Z$ denoting the atomic number of a single-element medium. For light-element samples ($Z$ small) and hard X-rays ($\omega$ large), these relations emphasize the necessity of imaging methods which are sensitive to the refractive phase shifts induced by $\delta$.
      
      \end{subsection}

      \begin{subsection}{Image Formation and the Projection Approximation} \label{SS:ProjApprox}

      We now return to the ideal setup sketched in \figref{fig:SetupIdeal} and consider the scattering interaction of the incident wave with the object in the domain $-L \leq z \leq 0$. As argued in \sref{SS:XrayRefrIdx}, the refractive index in the hard X-ray regime is always close to unity. Thus, we may approximate
            \begin{equation}
       n^2 - 1 =  - 2(\delta - \I \beta) + (\delta - \I \beta)^2  \approx -2(\delta - \I \beta)  = 2 (n-1) = -2 \N. \label{eq:n-approx}
      \end{equation}
       It is assumed that the incident wave satisfies the paraxial approximation, i.e.\ is described by \eqref{eq:ParaxHelm}, and that the thickness $L$ of the object is sufficiently small such that diffraction inside the material can be neglected. This is the \emph{projection approximation} and amounts to neglecting the lateral Laplacian in \eqref{eq:ParaxHelm} and thus to a description by geometrical optics. See \cite[sec. 2.2]{PaganinXRay} for details. Together with \eqref{eq:n-approx}, this yields
      	  \begin{equation}
	  \left(   \partial_z + \I   k \N    \right) \tilde \Psi =  \frac{\I } {2k} \nabla_{\bx}^2 \tilde \Psi \approx 0 \label{eq:ProjODE}.  
	\end{equation}
	Solving this ordinary differential equation and setting $\tilde \Psi_z:= \tilde \Psi(\cdot ,z)$, we obtain
	      	  \begin{equation}
        \tilde \Psi_0 = \underbrace{\tilde \Psi_{-L}}_{=: P }   \cdot  \underbrace{\exp \left( -\I k \PN  \right)}_{=: O } \MTEXT{with} \PN(\bx) := \int_{-L}^0 \N (\bx, z ) \; \D z \label{eq:ProjSol}
	\end{equation}
	
	According to \eqref{eq:ProjSol}, the interaction between radiation and matter is completely described by a multiplication of the \emph{probe-} or \emph{illumination function} $P$, representing the incident wave field, with the \emph{object transmission function} (OTF) $O$ \cite{Thibault2009ProbeObjectFct}. The latter is given by an exponential of line integrals in $z$-direction, i.e.\ by \emph{projections} of the sample along the optical axis. Physically, the wave field thus behaves like parallel non-interacting rays, which merely accumulate phase shifts and attenuation as they pass through matter. This is illustrated in \figref{fig:ProjApprox}. The axial integration can be interpreted as the formation of a \emph{shadow-} or \emph{contact image} of the specimen described by the spatially varying density $\N$, which is imprinted upon the transmitted wave field via \eqref{eq:ProjSol}.
	\begin{figure}[hbt!]
	 \centering
	 \includegraphics[width=.6\textwidth]{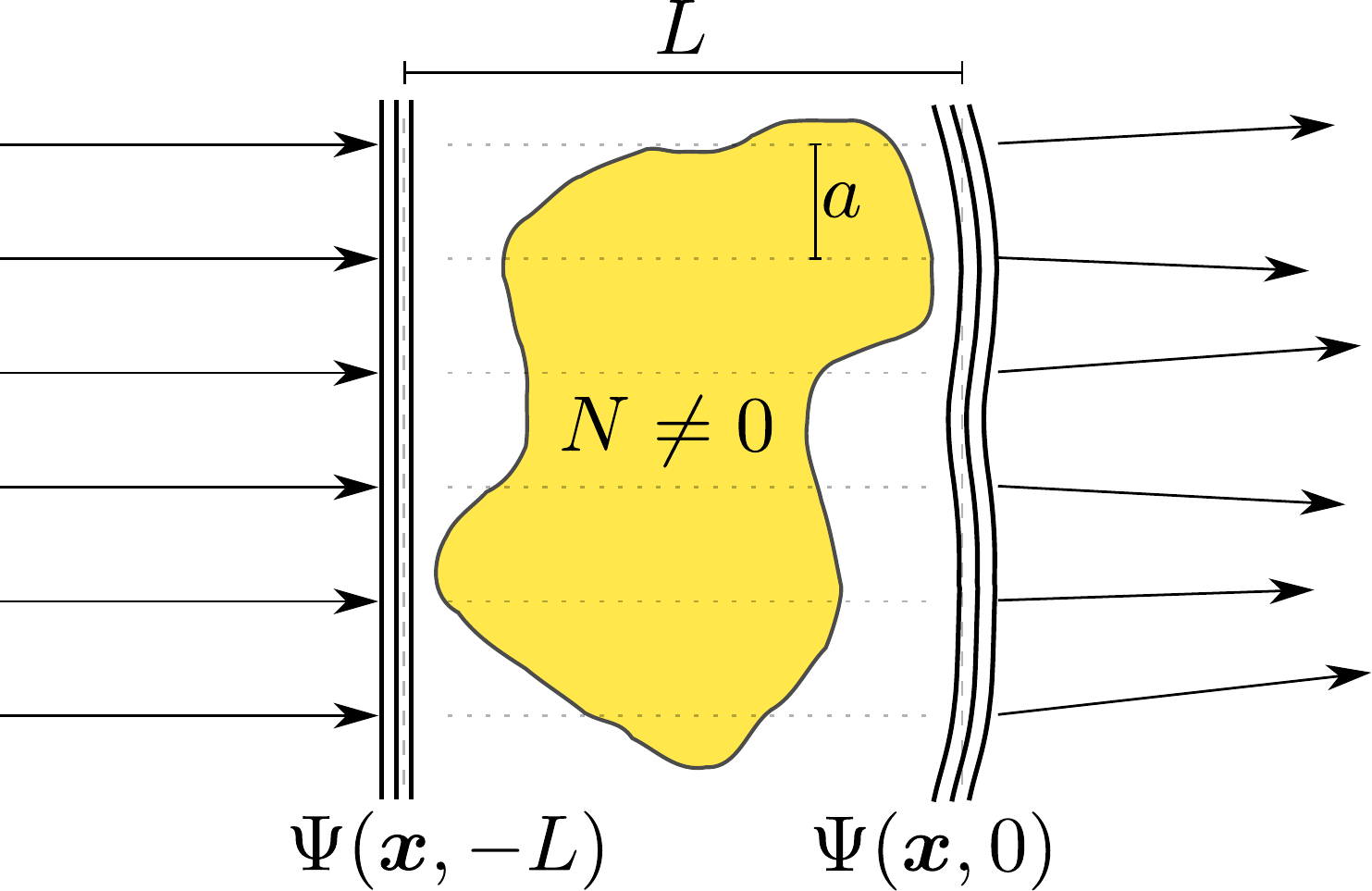}
	 \caption{Illustration of X-ray scattering in the projection approximation: objects of small thickness $ L \ll 2a^2k$, parametrized by a refractive index $n = 1-N$ satisfying $k L N \lesssim 1$, interact with incident monochromatic radiation of wavenumber $k$ as if the latter was composed of parallel non-interacting rays. The interaction reduces to an accumulation of phase and attenuation along the ray trajectories which leads to a perturbed wave field $\Psi$ at the exit-surface, the \emph{contact image} formed according to \eqref{eq:ProjSol} (Source: \cite{MartinMaster}, modified, inspired by \cite{BartelsDiss,PaganinXRay}).\label{fig:ProjApprox}}
	\end{figure}
	
	As can be seen from its Fourier space representation $\cF \nabla_{\bx}^2 = -\bxi^2$ (compare \eqref{eq:FTDeriv}), the neglected lateral Laplacian acts most significantly on small lengthscales. Hence, modeling the scattering within the framework of the projection approximation necessarily results in a lower limit for the attainable lateral resolution. Indeed, assuming plane wave illumination $P = 1$ and that $\PN   = \Or(L\N)$ varies only on lateral lengthscales $ \gtrsim a$, we obtain for $\tilde \Psi_0$ given by \eqref{eq:ProjSol}
	\begin{equation}
	\frac{ \nabla_{\bx}^2 \tilde \Psi_0}{ \tilde \Psi_0} =  - \I \underbrace{k \nabla_{\bx}^2 \PN}_{\lesssim \Or(kLa^2N)} - \underbrace{k^2 (\nabla_{\bx} \PN)^2}_{\lesssim \Or(k\PN \cdot (k L a^2N ))} \lesssim \Or( kLa^2N ).  \label{eq:ProjApproxEst}
	\end{equation}
	Here, it is further assumed that $k \PN \lesssim 1$, corresponding to arguments of at most order 1 in the exponential in \eqref{eq:ProjSol}\footnote{This excludes regimes of strong absorption and excessive phase wrapping (see \sref{SS:PhaseWrap}) which are typically unsuited for phase contrast imaging anyway.}.
	According to \eqref{eq:ProjApproxEst}, the neglected right hand side of \eqref{eq:ProjODE} is $\lesssim \frac 1 2 La^2 N \tilde \Psi$, whereas the retained second summand on the left hand side is $\sim k N \tilde \Psi$. Comparing these, it is found that the projection approximation is consistent for sufficiently weakly interacting objects of thickness
	\begin{equation}
	  L \ll 2a^2k.\label{eq:ProjApproxBound}
	\end{equation}
	On lateral lengthscales $a \lesssim ( L/(2k) )^{\frac 1 2}$, application of \eqref{eq:ProjSol} will result in blurry reconstructions of the projected refractive index $\PN$ due to the neglected diffusive coupling of neighboring ``rays'' induced by $\nabla_{\bx}^2$. In order to resolve these scales, the object may be decomposed into a sequence of thinner slices, applying the projection approximation within each of these but retaining diffractive effects by Fresnel propagation between the slices \cite{Cloetens1999Diss}.
	
	\end{subsection}

      \begin{subsection}{Phase Wrapping} \label{SS:PhaseWrap}
      
	From \eqref{eq:ProjSol} it can be seen that the real part $\delta$ of the specimen's refractive index manifests itself in the form of phase shifts of the transmitted X-rays compared to a propagation in vacuum. This effect is due to an increased phase speed in matter, i.e.\ faster propagation of the wave fronts. If the accumulated phase shifts partly exceed $2 \pi$, they are no longer represented uniquely in the exit wave field since wave front displacements by a multiple of the wavelength $\lambda = \frac{2\pi }k$ cannot be detected. This \emph{phase wrapping} problem is illustrated by \figref{fig:PhaseWrapSketch}.
      \begin{figure}[hbt!]
	 \centering
	 \includegraphics[width=.5\textwidth]{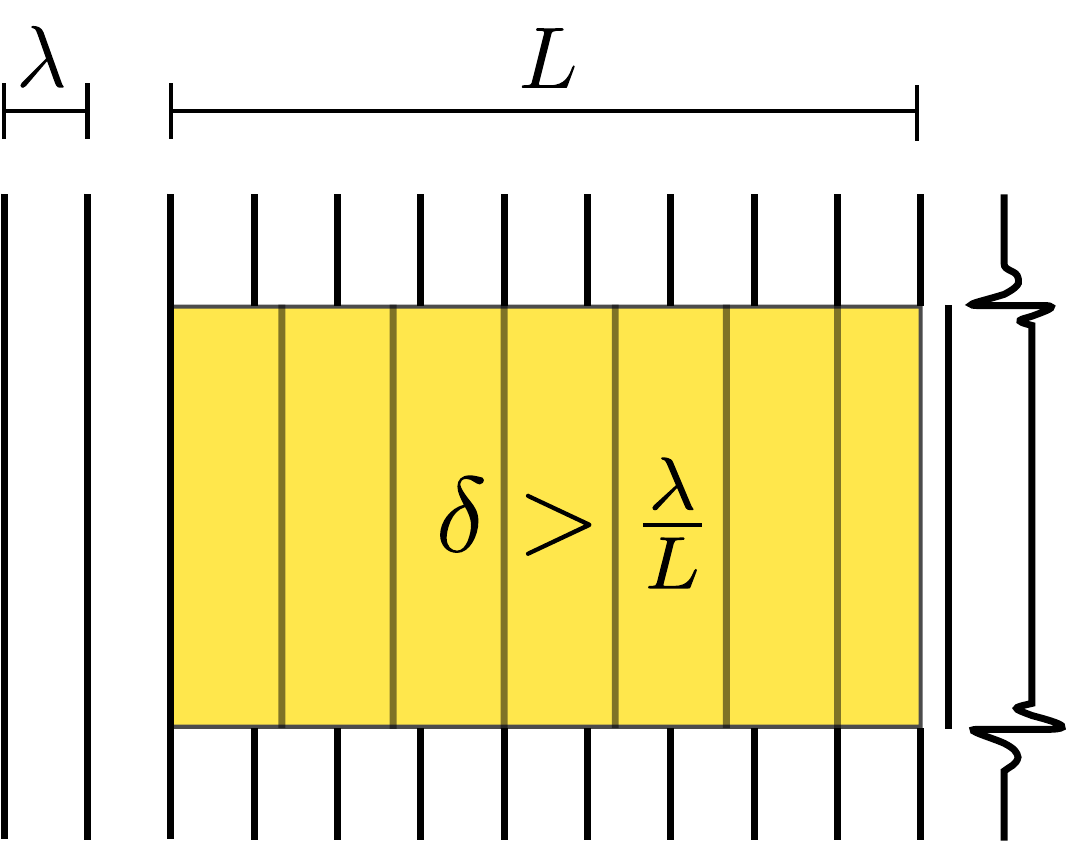}
	 \caption{Illustration of the phase-wrapping problem: an object (box-shaped for simplicity) with a refractive index $n = 1-\delta + \I \beta$ s.t.\ $\delta > 0$ is illuminated by monochromatic radiation of wavelength $\lambda$ which is represented by plane wave fronts. The wave fronts propagate faster within the material, resulting in a phase shift at the exit-surface. For $\delta > \frac \lambda L$, the accumulated phase may exceed $2\pi$ and thus cannot be identified uniquely in the exit wave: in the depicted example, the total phase shift is $\frac{14 \pi} 3$, whereas the relative phase discontinuity at the exit-surface is only $\frac{2 \pi} 3$ \label{fig:PhaseWrapSketch}}
	\end{figure}
	
	Mathematically, phase wrapping is reflected by the $2 \pi$-periodicity of the exponential in \eqref{eq:ProjSol} in the imaginary part of its argument. The effect has to be accounted for in image reconstruction whenever we have for some $\bx$
	\begin{equation}
	  \int_{-L}^0 \delta (\bx, z ) \; \D z > 2\pi k = \lambda. \label{eq:PhaseWrapCond1}
	\end{equation}
	A sufficient condition for phase wrapping to be absent is thus given by
	\begin{equation}
	  \norm{\delta}_{L^{\infty}} = \max_{\bx, z} |\delta(\bx , z )| < \frac \lambda L = \frac{ 2\pi }{kL}. \label{eq:PhaseWrapCond2}
	\end{equation}
	Accordingly, the projections $\Pdelta(\bx) := \int_{-L}^0 \delta (\bx, z ) \; \D z$ can be recovered uniquely from the exit wave field \eqref{eq:ProjSol} for sufficiently thin and weakly refracting objects. For stronger, moderately phase-wrapping objects, reasonable reconstructions may be achieved by heuristic phase unwrapping algorithms \cite{Ying2006PhaseUnwrapping}: assuming natural discontinuities of $\Pdelta$ to be small, lateral phase jumps of magnitude $\geq \pi $ are eliminated by adding integer multiples of $2\pi$ to the \emph{a priori} reconstructed guess.
	
	\end{subsection}

\begin{subsection}{Special Objects} \label{SS:SpecialObjects}

	Beyond the description of image formation for a general refractive index $N$ given by \eqref{eq:ProjSol}, it is useful to introduce a few special cases providing simplified parametrizations of certain specimen: 
	\begin{itemize}
	 \item \emph{Pure phase objects:} As argued in \sref{SS:XrayRefrIdx}, thin, light-element samples often give rise to negligible absorption. This can be accounted for by the approximation
	 \begin{equation}
	  N = \delta - \I \beta \approx \delta. \label{eq:PurePhaseObj}
	 \end{equation}
	 \item \emph{Pure absorption objects:} The opposite limiting case considered in classical computed tomography where refractive effects are negligible
	 \begin{equation}
	  N = \delta - \I \beta \approx -\I \beta. \label{eq:PureAbsObj}
	 \end{equation}
	 \item \emph{Single-material objects:} A specimen composed of a single material, merely varying in density, is characterized by a certain ratio between phase shifts and absorption. Thus there is a constant $\cbd$ such that 
	 \begin{equation}
	  N = \delta - \I \beta = (1 - \I \cbd) \delta. \label{eq:SingleMatObj}
	 \end{equation}
	\item \emph{Weak objects:} Beyond the above cases, the scattering object may be sufficiently weak for the projection $\PN$ to satisfy
	 \begin{equation}
	  \norm{\PN}_{L^{\infty}} = \max_{\bx} |\PN(\bx )| \ll  \frac 1 k .  \label{eq:WeakObjectCond}
	\end{equation}
	This limit legitimizes a linearization of the object transmission function:
	\begin{equation}
	O = \exp(-\I k \PN ) \approx 1 -\I k \PN.   \label{eq:WeakObjectLin}
	\end{equation}

	\end{itemize}

	\end{subsection}

 \end{section}

\begin{section}{Phase Contrast} \label{S:PhaseContrast}

\begin{subsection}{X-ray Detectors and the Phase Problem} \label{SS:Detection}

	In experimental realizations of phase contrast tomography as depicted in \figref{fig:SetupExp}, the scattered radiations are observed for instance by a CCD detector, possibly coupled to a scintillator in order to convert incident X-ray photons to visible light. See \cite[sec. 3.6.2]{PaganinXRay} for an overview on detecting devices. In a classical picture, the electromagnetic field of coherent X-rays of a wavelength $\lambda \sim \unit[10^{-10}]{m}$ oscillate with a frequency $\sim\unit[10^{18}]{Hz}$. Hence, temporal dynamics of the wave field are too fast to be measured by any existing technology \cite[p. 44]{PaganinXRay}. This, however, implies that spatial phase variations in the scattered wave field, representing time lags in the oscillations in the order of $10^{-18}$ seconds or less, cannot be measured. Instead, X-ray detectors are sensitive only to time-averaged wave \emph{intensities}. This defect is known as the \emph{phase problem} of optics, playing a significant role not only in X-ray tomography but also in crystallography \cite{Millane1990,Cruickshank1961PhaseProbCrystal} and electron microscopy \cite{Misell1973PhaseProbEMicroscopy}, for example.
	
	Mathematically, the phase problem implies that not the complex wave field $\Psi$ itself but only its squared modulus $I_d:= |\tilde \Psi_d|^2$, giving the intensity in the plane $z=d$, is accessible by measurements. For the contact image at $z= 0$ in the scattering experiment in \figref{fig:SetupIdeal}, this yields
	\begin{align}
	 I_0  &= |\tilde \Psi_0|^2 = |P |^2 \cdot |O |^2 = |P |^2 \cdot   \exp \left( 2k \Im \left(\PN \right) \right) \label{eq:ContrastContactImage}
	\end{align}
	according to \eqref{eq:ProjSol}. \eqref{eq:ContrastContactImage} implies that the entire \emph{refractive} information, represented by $\Pdelta = \Re(\PN)$, would be irretrievably lost if the intensities were detected in the exit-surface $E_0$. On the contrary, the attenuation $ \Pbeta  = - \Im (  \PN)$ is fully retained in the data. This is the regime of \emph{absorption contrast} characterized by the limit $\NF \to \infty$ of vanishing diffraction, in which classical CT scanners operate \cite[p. 5]{Natterer}.
	
	\emph{Phase contrast} imaging aims at resolving the projected refractive decrement $\Pdelta$ in order to overcome the limitation of X-ray radiography to macroscale, sufficiently absorbing objects. In this work, the required phase-sensitivity is achieved by the diffractive encoding of phase information via propagation of the contact image to a distant detector plane $E_d$ with the near-field- and far-field propagators introduced in \sref{SS:ParaxialFresnel} and \sref{SS:FresnelFarfield}. In this setting, the detected intensities are given by
	\begin{align}
	 I_d &= \left| \bDd{d}\left(\tilde \Psi_0 \right) \right|^2 = \left| \bDd{d}\left(P \cdot O  \right)  \right|^2 =  \left| \bDd{d}\left(P \cdot  \exp \left( - \I k  \PN  \right)  \right)  \right|^2 \label{eq:ContrastDetectorImage} 
	\end{align}
	
	\end{subsection}

\begin{subsection}{Contrast Formation in the Near-Field} \label{SS:ContrastFormationNF}
	
	In order to gain an insight into contrast formation in the near-field regime, characterized by $\bDd{d} = \bDF{d}$ in \eqref{eq:ContrastDetectorImage}, we assume plane wave illumination $P = 1$ and a \emph{weak object}, i.e.\ apply the linearization of the OTF given by \eqref{eq:WeakObjectLin}. Then, by \eqref{eq:ContrastDetectorImage} the observed intensities on the detector are
	\begin{align}
	 I_d \approx \left|1 - \I k \bDF{d}(\PN)\right|^2  = 1 - \I k \left( \bDF{d}(\PN) - \cc{\left[\bDF{d}(\PN)\right]}\right) + k^2 \left|  \bDF{d}(\PN)\right|^2 \label{eq:WeakObjRealCTF}
	\end{align}
	The quadratic contribution in $\PN$ is negligible within the weak object approximation. Using the definition of the Fresnel propagator \eqref{eq:EnvelopeFresnelProp} and $\cc{\left[\bDF{d}(N)\right]} = \bDF{-d}(\cc{\PN})$, the \emph{contrast transfer function} (CTF) is obtained by Fourier transforming \eqref{eq:WeakObjRealCTF}:
	\begin{equation}
	 \cF( I_d- 1 )(\bxi)  =-2 k \left( \sin\left(\frac{  d \bxi^2  }{2 k } \right) \cF(\Pdelta)(\bxi) + \cos\left(\frac{  d \bxi^2  }{2 k } \right) \cF(\Pbeta)(\bxi)  \right) \label{eq:CTF}
	\end{equation}
	For a detailed derivation, see \cite{Guigay1977CTF}.
	  \begin{figure}[hbt!]
	 \centering
	 \includegraphics[width=.6\textwidth]{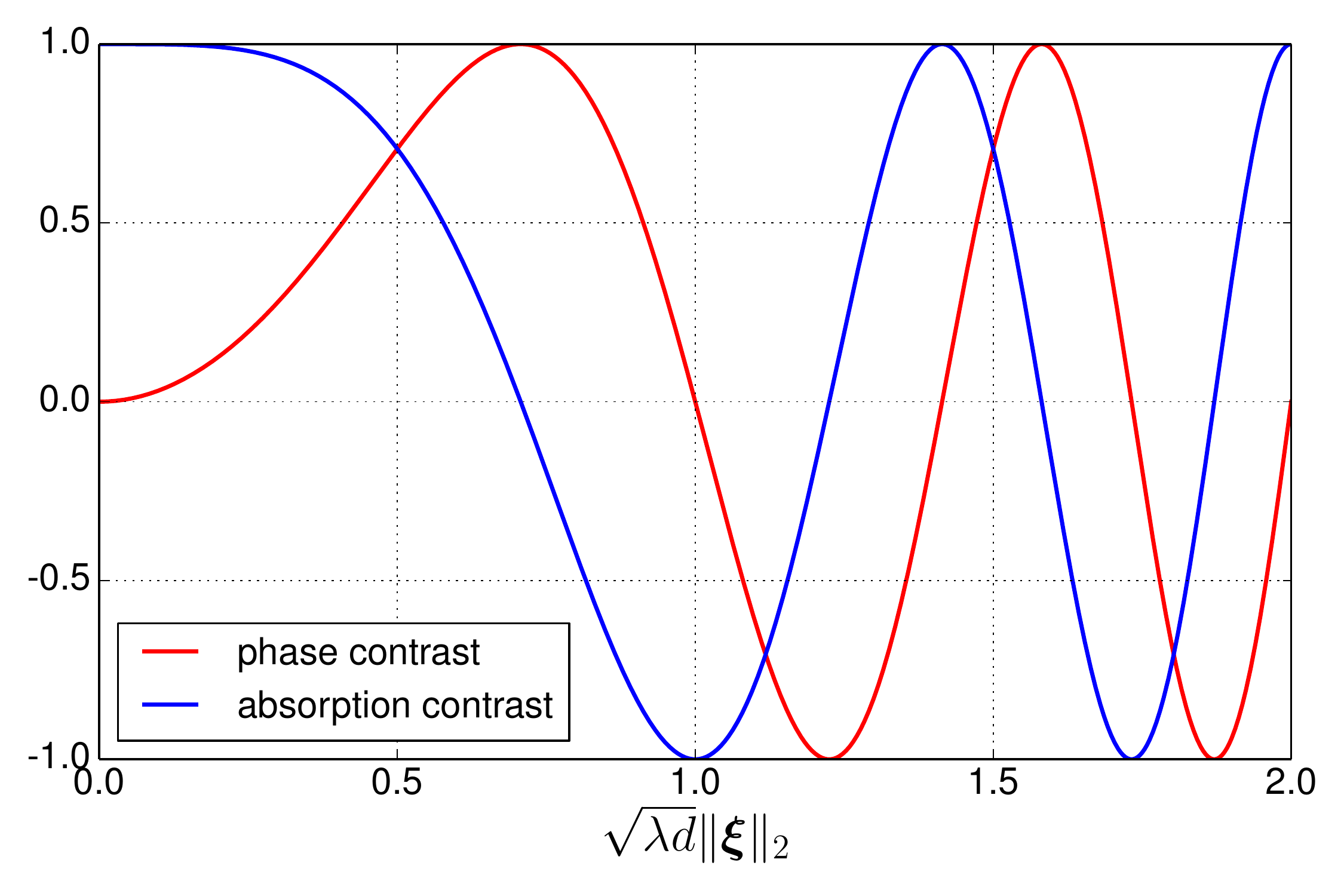}
	 \caption{Contrast transfer function (CTF): For weak objects, near-field phase contrast imaging is governed by \eqref{eq:CTF}, resulting in a manifestation of the projected refraction $\Pdelta$ (phase contrast) and absorption $\Pbeta$ (absorption contrast) in wave intensities that is oscillating in the Fourier frequencies $\bxi$. The zeros of these curves correspond to underrepresented Fourier modes of $\Pdelta$ or $\Pbeta$ in the intensity data for given propagation distance $d$ and wavelength $\lambda = \frac{2\pi}k$ (replicates \cite[Fig. 1]{Pogany1997noninterferometric})  \label{fig:CTF}}
	\end{figure}
	
	From \eqref{eq:CTF} it can be seen that, by propagation over the detector distance $d > 0$, the refractive decrement $\Pdelta$ manifests in a measurable intensity pattern, denoted as a \emph{hologram}, along with the absorptive part $\Pbeta$. Moreover, the derived Fourier space representation reveals that the achieved contrast in both phase and absorption are oscillatory in the Fourier frequencies $\bxi$ due to the sine- and cosine prefactors in \eqref{eq:CTF}. This oscillation is visualized in \figref{fig:CTF}. The setup-dependent zeros of these factors at
	\begin{subequations} \label{eq:CTFZeros}
	 \begin{align}
	\frac{ d \bxi^2 }{2 k } &\in \pi \mZ \MTEXT{(phase contrast)} \label{eq:CTFZerosPhase} \\
	\frac \pi 2 +  \frac{ d   \bxi^2 }{2 k } &\in \pi \mZ \MTEXT{(absorption contrast)} \label{eq:CTFZerosAbs}
	\end{align}
	\end{subequations}
	correspond to spatial frequencies $\bxi$ of $\Pdelta$ and $\Pbeta$ which are underrepresented in the intensity data and which thus cannot be reconstructed with reasonable accuracy. In particular, note that phase contrast is weak at the lower end of the spectrum, especially for small propagation distances. This implies that a certain minimum propagation distance is required in practice to achieve reasonable contrast. Otherwise, only sharp edges in the projected refraction $\Pdelta$, corresponding to high Fourier frequencies, will be visible in the holograms whereas bulk regions are merely represented at all.
	  \begin{figure}[hbt!]
	 \centering
	 \includegraphics[width=\textwidth]{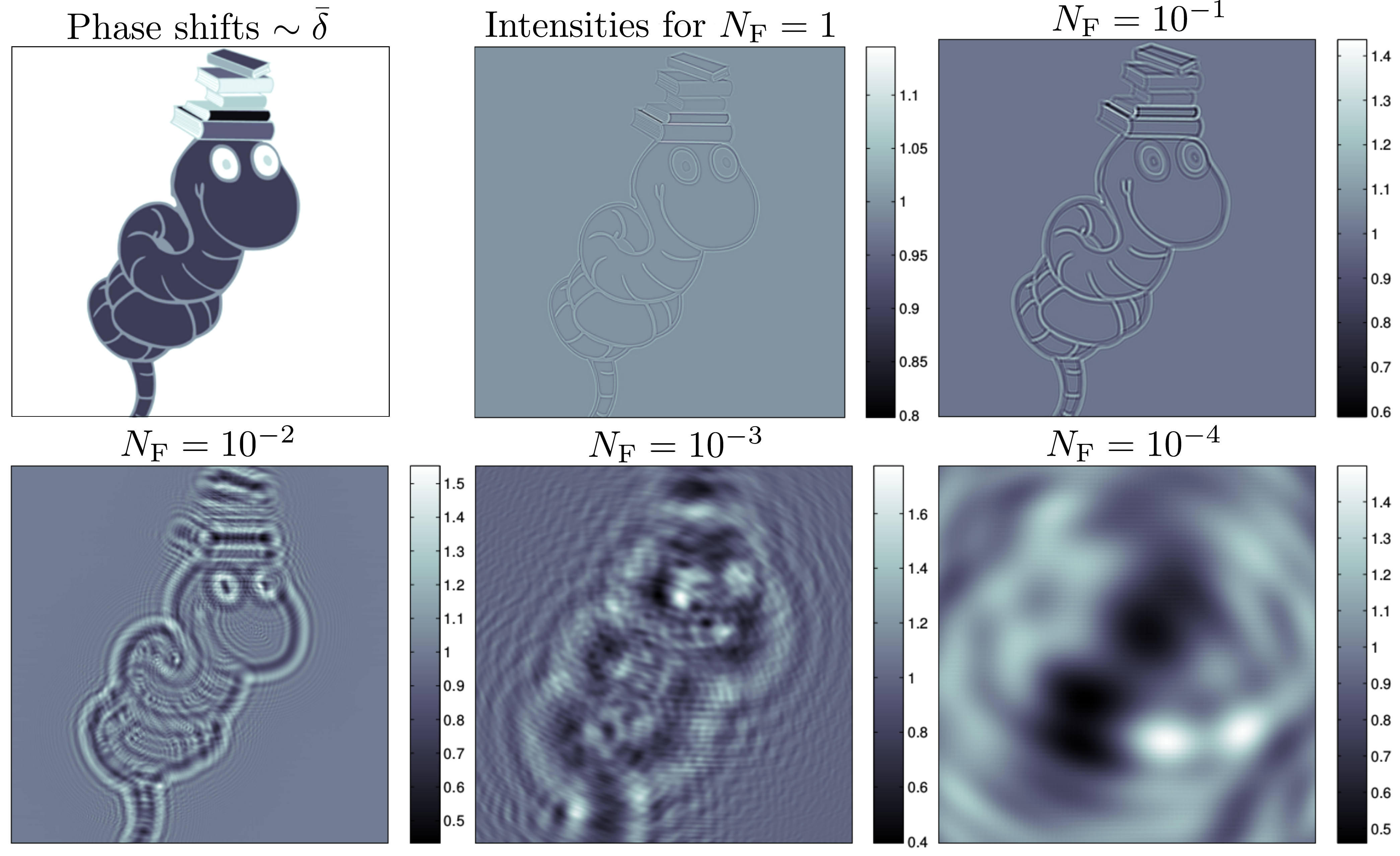}
	 \caption{Simulated holograms of a test object for different Fresnel numbers $\NF =  {k b^2}/{(2 \pi d)}$ where the lengthscale $b$   is taken as the aspect length $\Delta x$ of a single pixel. The projected refractive decrement $\Pdelta$ is assumed to be given by the upper-left image of size $256 \times 256$. For small propagation distances $d \propto \NF^{-1}$, only sharp edges in the contact image manifest themselves in the intensity patterns. As $\NF$ decreases, more and more fringes appear by diffraction of the propagated wave field, encoding phase information. For small Fresnel numbers $\NF \to 0$, corresponding to the far-field limit $d \to \infty$, the structure of the original contact image can no longer be identified in the propagated holograms  (Source: \cite{MartinMaster}, modified)  \label{fig:ContrastTransferDemo}}
	\end{figure}
	
	This effect is illustrated by \figref{fig:ContrastTransferDemo}, showing simulated intensity data for a given phase image at different propagation distances. Here, the latter are expressed in terms of the Fresnel number $\NF \propto d^{-1}$ (compare \sref{SS:FresnelFarfield} and \eqref{eq:FresnelNumber}), where the lengthscale $b$ is taken as the aspect size  $\Delta x$ of a single pixel. From a simple shadow image of the object's edges at $\NF = 1$, i.e.\ quasi-geometrical optics representation, the propagated data becomes wavy due to diffractive effects as $\NF$ decreases, showing more and more fringes encoding the phase information of the contact image. We denote this as the regime of \emph{holographic imaging} being the one of principal interest in this work, represented by the cases $\NF = 10^{-2}$ and $\NF = 10^{-3}$ in \figref{fig:ContrastTransferDemo}. For even smaller Fresnel numbers (here: $\NF = 10^{-4}$), features of the contact image can no longer be identified in the holograms, which corresponds to the transition to the far-field case discussed in \sref{SS:ContrastFormationFF}.
	
	Nevertheless, \eqref{eq:CTF} may be solved for $\Pdelta$ or $\Pbeta$, provided that the other field is known as is the case in particular for pure phase- or single material objects (see \sref{SS:SpecialObjects}). By suitable regularization around the zeros of the CTF, this yields a reconstruction method for the projected refractive index as outlined for instance in \cite{Cloetens1999}. Its regime of applicability may be enlarged from the weak object case assumed here, to weakly absorbing samples which give rise to slowly varying phase shifts \cite{Turner2004FormulaWeakAbsSlowlyVarPhase,Guigay1977CTF}. As the reconstruction of $\Pdelta$ (and $\Pbeta$) from intensity data implicitly requires the recovery of the missing phase information of the wave field, these CTF techniques can be viewed as examples of \emph{phase retrieval} methods.
	
\end{subsection}

\begin{subsection}{Contrast Formation in the Far-Field} \label{SS:ContrastFormationFF}

	The limit of large propagation distances $d$ is characterized by the far-field propagator $\bDd{d} = \bDff{d}$ defined in \eqref{eq:Fraunhofer}. According to \eqref{eq:ContrastDetectorImage}, the measurable intensity data in this case is given by
	  \begin{equation}
	 I_d(\bx) = \left| \bDff{d}\left(P  \cdot O \right)(\bx)   \right|^2 = \left(\frac{k }{  d} \right)^{m}\left|   \cF \left(P \cdot O\right) \left(  \frac{ k \bx }{ d } \right)  \right|^2.\label{eq:ContrastFF}
	\end{equation}
	Thus, up to a suitable rescaling of the magnitudes and the lateral coordinate, far-field imaging, also called \emph{coherent diffractive imaging}, measures the squared modulus of the contact image's Fourier transform \cite{Miao1998FourierNonPeriodic,MiaoNature1999CDIFirstExp}. Evidently, such data is sensitive not only to the modulus of $O$, determined by the absorption $\Pbeta$, but also to its phase, encoding the projected refraction $\Pdelta$. However, other than in the near-field case, there is no simple approach to see which part of the information can be reconstructed from the intensity patterns. This question is therefore postponed to \sref{SS:PhaseRetrFF}.
	
	 \begin{figure}[hbt!]
	 \centering
	 \includegraphics[width=.5\textwidth]{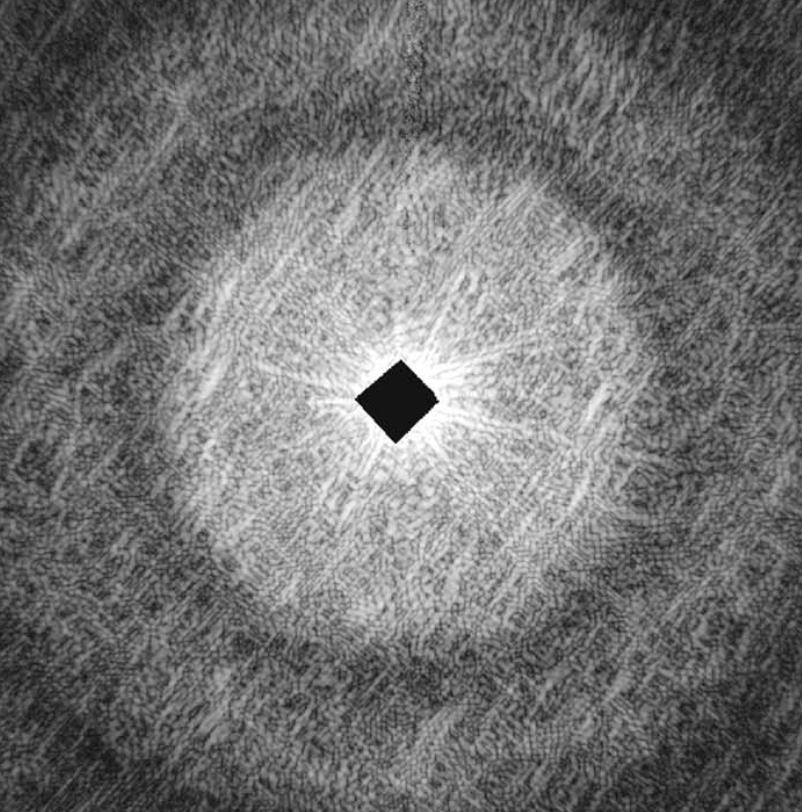}
	 \caption{Example of an experimentally observed far-field phase contrast image \cite{Chapman2006}. The intensities in the center of the diffraction pattern, corresponding to low Fourier frequencies dominated by the unscattered probe beam, have to be blocked by a suitable beam stop in order not to damage the detector. In practice, these are therefore inaccessible to measurements. \label{fig:FarfieldBeamStop}}
	\end{figure}
	
	The scaling operations in \eqref{eq:ContrastFF} can be suppressed in the mathematical model by introducing the far-field intensity $I_\infty\left(  \bxi \right):= \left(\frac{k }{  d} \right)^{m} I_d(\bx),\,\bxi := \frac{ k \bx }{ d } $, so that
	\begin{equation}
	 I_\infty =  \left|   \cF \left(P \cdot O\right)   \right|^2 = \left| \cF\left(P\right)  +  \cF\left(P\cdot \left[\exp \left( - \I k  \PN  \right) - 1 \right] \right)   \right|^2 \label{eq:FarfieldIntensity}
	\end{equation}
	The first summand on the right hand side of \eqref{eq:FarfieldIntensity} yields the contribution of the unscattered probe beam, which usually varies on much larger lengthscales than the projected specimen $\PN$. Consequently, $\cF\left(P\right)$ will be negligible except for a neighborhood of the origin representing low Fourier frequencies. In particular, in the ideal case of plane wave illumination $P=1$, this contribution reduces to a single Dirac-delta-peak at 0. In practice, the resulting intensities around the center of a far-field diffraction pattern often exceed the dynamic range of the CCD detector and may even damage the latter if not blocked by a suitable \emph{beam stop} \cite{Thibault2010XrayDiffMicrBeamstop}. Consequently, these contributions will not be represented accurately in experimental far-field data as depicted in \figref{fig:FarfieldBeamStop}.
	
	The second summand in \eqref{eq:FarfieldIntensity}, on the other hand, has essentially the same support and characteristic lengthscales as $\PN$ and will thus typically yield structures in Fourier space that extend to much higher frequencies. These empirical facts allow to neglect the - locally dominant - probe term in \eqref{eq:FarfieldIntensity}, i.e.\ we have for $\bxi \in \mR^{m }$ with $\norm{\bxi}_2$ sufficiently large
	\begin{equation}
	 I_\infty(\bxi) \approx  \left|   \cF\left(P\cdot \left[\exp \left( - \I k  \PN  \right) - 1 \right] \right)(\bxi)   \right|^2.   \label{eq:FarfieldNoProbe}
	\end{equation}
	
	\end{subsection}
	
\end{section}

\begin{section}{Tomography} \label{S:TomoTheory}

	\begin{subsection}{Parametrization by the Radon Transform} \label{SS:TomoRadon}

	So far, the developed physical model of phase contrast imaging does not take into account that the setup in \figref{fig:Setup} allows to rotate the specimen in the $x$-$z$-plane. This changes the incident angle $\alpha$ of the radiation, which propagates along the $z$-axis of a coordinate system fixed in space, with respect to the object's coordinate system as illustrated in \figref{fig:SetupIdeal}. For reasons of convention, we write the corresponding rotated version of the refractive index in terms of the angle $\theta = \frac \pi 2 - \alpha$:
	\begin{align}
	N_{\theta}(\boldsymbol x,z) =  N_{\theta}(x,\by,z) &:= N(\cos(-\alpha) x - \sin(-\alpha) z,\by,\sin(-\alpha) x + \cos(-\alpha) z) \nonumber \\
	&\;= N(\sin(\theta) x + \cos(\theta) z,\by,-\cos(\theta) x +  \sin(\theta) z). \label{eq:RotN}
	\end{align}
	Here, $\bx = (x, \by)$ and $\by \in \mR^{m-1}$ denotes the lateral dimensions which are not affected by the rotation.
	Inserting \eqref{eq:RotN}
	into \eqref{eq:ProjSol} yields a version of the object transmission function which accounts for the variable incident angle of the illumination:
	\begin{align}
	O(\theta, \bx) &:= \exp \left( -\I k \int_{-L}^0 N_{\theta}(\bx, z ) \; \D z \right)\nonumber  \\
	&\;= \exp \left( -\I k \int_{-L}^0 N(\sin(\theta) x + \cos(\theta) z,\by,-\cos(\theta) x +  \sin(\theta) z) \; \D z \right). \label{eq:OTF}
	\end{align}
	Comparing \eqref{eq:OTF} to \eqref{eq:Radon} and \eqref{eq:CylRadon}, it can be seen that the angle-dependent line integrals exactly match the \emph{cylindrical Radon transform} defined in \sref{S:Radon}, i.e.
	\begin{equation}
	O(N) = \exp \left( -\I k \CR(N) \right). \label{eq:OTFbyRadon}
	\end{equation}
	The derived expression \eqref{eq:OTFbyRadon} describes tomographic image formation by means of the Radon transform, valid within the framework of the paraxial- and the projection approximations. Note that the object transmission function is introduced as an operator acting on $N$, as is convenient for the analysis of the following chapters.
	
	\end{subsection}

	\begin{subsection}{The Forward Operators of Phase Contrast Tomography} \label{SS:PhaseTomoOperator}
	
	With the theoretical background provided in \sref{S:FreeProp}, \sref{S:PhaseContrast} and \sref{SS:PhaseTomoOperator}, we are finally in a position to combine the different stages of image formation
		to a complete model for phase contrast tomography:
	\begin{itemize}
	 \item[$\boldsymbol 1$] Scattering interaction within the projection approximation
	 \item[$\boldsymbol 2$] Diffraction of propagating paraxial waves
	 \item[$\boldsymbol 3$] Detection of the propagated scattered wave field
	\end{itemize}
	
	First, consider the near-field case where $\boldsymbol 2$ is governed by the Fresnel propagator. Combining \eqref{eq:OTFbyRadon}, \eqref{eq:ProjSol} and \eqref{eq:EnvelopeFresnelProp}, we find that the detected intensity under the incident angle $\theta \in [0;2 \pi)$ at $\bx \in \mR^m$ in the setup in \figref{fig:Setup} is given by
	\begin{align}
	 I_d(\theta, \bx ) &= \left| \bDF{d} \left[ P \cdot O(N)(\theta, \cdot) \right] (\bx) \right|^2 \nonumber  \\
	 &= \left| \bDF{d} \left[ P \cdot \exp \left( -\I k \CR(N)(\theta,\cdot) \right) \right] (\bx) \right|^2. \label{eq:PhaseTomoIntensities}
	\end{align}
	Equation \eqref{eq:PhaseTomoIntensities} defines the \emph{forward operator} of near-field phase contrast tomography, which maps the field $N$ parametrizing the spatial structure of the specimen onto measurable data:
	\begin{align}
	 F_d: N \mapsto I_d &=  \left| \bDF{d}(P) + \bDFlat{d} \left( P \cdot \left[ \exp \left( -\I k \CR(N) \right) -1 \right] \right)  \right|^2 \nonumber \\
	  &=  \left| \bDF{d}(P) + \bDFlat{d} \left( P \cdot O_0 (N) \right)  \right|^2 \label{eq:ForwardOpNF}
	\end{align}
	Here, the \emph{normalized object transmission functions} $O_0 := O -1$ is inserted and the subscript $\overline 2$ indicates application of the propagator in the lateral coordinate $\bx$. 
	
	Analogously, the expression \eqref{eq:FarfieldNoProbe} for the far-field intensities, motivated in \sref{SS:ContrastFormationFF}, may be supplemented with the tomographic object transmission function to obtain a forward operator for far-field phase contrast tomography:
	\begin{align}
	 F_\infty: N \mapsto I_\infty =  \left| \cF_{\overline 2}\left( P \cdot \left[ \exp \left( -\I k \CR(N) \right) -1 \right] \right)  \right|^2  =   \left| \cF_{\overline 2}\left( P \cdot O_0 (N) \right)  \right|^2\label{eq:ForwardOpFF}
	\end{align}
	
	In addition to the general operators defined above, it is instructive to consider $F_d$ and $F_\infty$ in the weak object limit $N \ll 1$ (see \sref{SS:SpecialObjects}) represented by
	\begin{equation} O_0 (N) \approx \I k \CR(N). \end{equation}
	For plane wave illumination $P = 1$, this yields by \eqref{eq:EnvelopeFresnelProp} and the Fourier slice theorem  \eqref{eq:FourierSliceCyl}
	\begin{subequations}\label{eq:WeakFwOps}
	\begin{align} 
	F_\infty(N) &\approx  k^2 \left| \CF (N)   \right|^2 \label{eq:WeakFwOpFF} \\
	 F_d(N) &\approx 1 + 2 k \Im \left( \cF_{\overline 2}^{-1} \bMF_{d, \overline 2} \CF (N)  \right) +  \underbrace{k^2 \left| \cF_{\overline 2}^{-1} \bMF_{d, \overline 2} \CF (N)  \right|^2}_{= \Or(\norm{N}^2)}   \label{eq:WeakFwOpNF}
	\end{align}
	\end{subequations}
	Accordingly, the detected far-field intensities in this setting essentially correspond to squared modulus of a Fourier transform of the object $N$ in cylindrical coordinates. The near-field operator $ F_d$, on the other hand, reduces to its linearization up to quadratic errors in $N$ and likewise gives rise to a cylindrical Fourier transform.

	Rather than predicting the tomographic holograms to be observed for a known specimen by application of $F_d$ or $F_\infty$, however, the physically relevant operation is to recover the sample from recorded intensity data by inverting these maps. The objective of this work is to solve this \emph{inverse problem}:
	\begin{prob}[Inverse Problem of propagation-based Phase Contrast Tomography]\label{prob:1}
		From intensity data $I^\dagger_\ast$ given by \eqref{eq:ForwardOpNF} or \eqref{eq:ForwardOpFF}, reconstruct the specimen's refractive index $n^\dagger = 1- N^\dagger$.
	\end{prob}
	The subsequent chapters are dedicated to the analysis of \probref{prob:1} and to the development of stable algorithms for a numerical solution.

	\end{subsection}
	
\end{section}

\end{chapter}


\begin{chapter}{Analysis of the  Problem}\label{C:Analysis}

  In \chapref{C:PhysProb}, a physical model of phase contrast tomography has been introduced
  and formulated as a (nonlinear) inverse problem. In the following sections,
  the derived operator equation is analyzed, showing \Frechet differentiability as well as ill-posedness
  of reconstruction Problem \ref{prob:1} and
  establishing sufficient conditions for uniqueness of its solution. The analysis is based on the mathematical definitions and theorems introduced in \aref{C:MathPre}.

\begin{section}{Well-Definedness and \Frechet Differentiability} \label{S:WellDefFrechet}

  Before turning to the \emph{inverse} reconstruction problem, we show well-behavedness of the forward operators defined in \sref{SS:PhaseTomoOperator} in a mathematical sense by proving their \Frechet differentiability on suitable domains. Exploiting the nested definitions of the operators in \eqref{eq:ForwardOpNF} and  \eqref{eq:ForwardOpFF}, this can be achieved step-wise by virtue of the differentiation rules in \thmref{thm:FrechetProps}.
  
 \begin{subsection}{Object Transmission Function} \label{SS:FrechetOTF}
  
  As indicated in the setup sketched in \figref{fig:SetupIdeal}, we assume that the object lies within the cylindrical domain
  \begin{equation}
  \ObjDom:= \left\{ (x,\by,z) \in \mR^{m+1}: \; \norm{(x,z)}_2^2 \leq \frac L 2, \; \by \in \left[-\frac{L_y} 2;\frac{L_y} 2\right]^{m-1}  \right\}, \label{eq:ObjDomain}
  \end{equation}
  and that the field $N = 1-n \in \Lp \infty {\ObjDom}$ parametrizing the specimen is bounded and supported in $\ObjDom$. By \defref{def:Radon} and \cref{cor:CylRadonFourProps}, this implies $\CR(N) \in \Lp \infty {\ProjDom}$ where the projection domain $\ProjDom$ is given by
  \begin{equation}
  \ProjDom := [0; 2\pi) \times \LatDom \MTEXT{where} \LatDom  := \left[-\frac{L } 2;\frac{L } 2\right] \times \left[-\frac{L_y} 2;\frac{L_y} 2\right]^{m-1}
  \end{equation}
  On these domains, the normalized OTF $O_0$ is \Frechet differentiable:
  \vspace{1em}
  \begin{lem1}[\Frechet Derivative of the Object Transmission Function] \label{lem:FrechetOTF} 
   \begin{equation*} O_0 : \Lp \infty {\ObjDom} \to  \Lp \infty {\ProjDom}; \; N \mapsto \exp \left( -\I k \CR(N) \right) -1\end{equation*}
   is \Frechet differentiable. For $N, h_N \in \Lp \infty {{\ObjDom}}$, the derivative is given by
      \begin{equation} O_0'[N]h_N = -\I k  \exp \left( -\I k \CR(N) \right) \cdot \CR(h_N).  \label{eq:OTFDerivative} \end{equation}
      Moreover $O_0'[N]$ has a unique extension to an operator $\Lp 2 {\ObjDom} \to  \Lp 2 {\ProjDom}$.
  \end{lem1}
  \begin{pf}
   $O_0$ is of the form $O_0 = G \circ H$ where
   \begin{equation*}G: \Lp \infty {\ProjDom} \to  \Lp \infty {\ProjDom}; \; g \mapsto \exp(g) -1 \end{equation*}
   and $H: \Lp \infty {{\ObjDom} } \to  \Lp \infty {\ProjDom}; \; f \mapsto - \I k \CR f$ is a bounded linear operator. Hence, by \exref{ex:Frechet}-(c) and \thmref{thm:FrechetProps}-(d,e), $G$ and $H$ are differentiable with
    \begin{equation*} G[g]h_g = \exp(g)\cdot h_g \MTEXT{and} H[f]h_f = H(h_f) = - \I k \CR h_f \end{equation*}
    for all $g,h_g \in \Lp \infty {\ProjDom} , f,h_f \in \Lp \infty {{\ObjDom}}$. According to the chain rule \eqref{eq:ChainRule}, this implies that $O_0$ is \Frechet differentiable with derivative given by \eqref{eq:OTFDerivative}.
    
    By \cref{cor:CylRadonFourProps}, $\CR$ is also bounded as an operator $\Lp 2 {\ObjDom} \to  \Lp 2 {\ProjDom}$. For $N \in \Lp \infty {{\ObjDom}}$, the prefactor on the right hand side of \eqref{eq:OTFDerivative} is in $\Lp \infty {\ProjDom}$ so that
    \begin{equation*} \Lp 2 {\ProjDom} \to \Lp 2 {\ProjDom}; \; g \mapsto -\I k  \exp \left( -\I k \CR(N) \right) \cdot g \end{equation*}
    is likewise continuous. Hence, there exists a bounded extension $O_0'[N]: \Lp 2 {\ObjDom} \to  \Lp 2 {\ProjDom}$. This extension is unique by \thmref{thm:DenseIncl}-(c).
  \end{pf}
  
\end{subsection}

 \begin{subsection}{Superposition of the Probe Field and Propagation} \label{SS:FrechetProbeProp}
  
  We assume that the probe wave field $P$ in \eqref{eq:ForwardOpNF} and  \eqref{eq:ForwardOpFF} is given by the superposition of a plane wave component, parametrized by a constant $c_P$, plus a bounded $L^1$-perturbation, i.e.
  \begin{align}
 P = P_0 + c_P  \in \mP &:= \Lp 1 {\mR^m} \cap \Lp \infty {\mR^m} \oplus 1 \cdot \mC \nonumber \\
 &\; = \{ f + c \cdot 1 : f \in \Lp 1 {\mR^m} \cap \Lp \infty {\mR^m}, \, c \in \mC \}. \label{eq:ProbeDomain},
  \end{align}
  Here, $1: \mR^m \to \mR$ denotes the one-function in $\mR^m$. In particular, the constructed setting includes the important special cases of both incident ideal plane waves and more realistic illumination by a \emph{Gaussian beam}. The latter is characterized by a lateral intensity profile that is everywhere Gaussian, see \cite[sec. 3.1]{Teich1991Photonics} for details.
  
  As a next step in the analysis of the forward operators, we consider generalized propagators mapping the normalized OTF for a single angle onto the corresponding scattered component of the complex wave field at the detector:
  \begin{align}
   D_d  :   \vartheta \mapsto  \bDF{d} \left( P \cdot \vartheta    \right) \MTEXT{and} D_\infty: \vartheta \mapsto \cF \left( P \cdot \vartheta   \right) \label{eq:O-P-Propagator}
  \end{align}
 The operators $D_d$ and $D_\infty$ represent the near- and far-field case, respectively.
   \vspace{1em}
  \begin{lem1}[Boundedness of the Generalized Propagators] \label{lem:BoundedGenProps} 
   For $P \in \mP$, the maps given by \eqref{eq:O-P-Propagator} are well-defined bounded linear operators
   \begin{equation*} 
    D_d, D_\infty:  \Lp 2 {\LatDom} \to   \Lp 2 {\mR^m} \cap \Lp \infty {\mR^m}.
   \end{equation*}
  \end{lem1}
  \begin{pf}
  The multiplication $M: \vartheta \mapsto P \cdot \vartheta $ is linear and preserves the support domain $\LatDom$. Moreover, with $P = P_0 + c_P \in \Lp \infty {\mR^m} \oplus 1$, we obtain for all $o \in \Lp 2 {\LatDom}$
   \begin{equation*}
    \norm{ M(\vartheta) }_{ \Lp 2 {\LatDom}} = \norm{ c_P \vartheta + P_0 \cdot \vartheta }_{ \Lp 2 {\LatDom}} \leq  (|c_P| + \norm{ P_0 }_{ \Lp \infty {\mR^m}} )  \norm{ \vartheta }_{ \Lp 2 {\LatDom}},
   \end{equation*}
 i.e.\ $M : \Lp 2 {\LatDom} \to \Lp 2 {\LatDom}$ is bounded. Now, according to \thmref{thm:LpBoundedEmbed}, there is a continuous embedding
   \begin{equation*}
    \Lp 2 {\LatDom} \to \Lp 1 {\mR^m} \cap \Lp 2 {\mR^m}
   \end{equation*}
   as $ \LatDom $ is bounded. By \thmref{thm:FTBounded} and \cref{cor:FourierL2}, the Fourier transform $\cF$ is bounded as an operator $\Lp 1 {\mR^m} \cap \Lp 2 {\mR^m} \to \Lp \infty {\mR^m} \cap \Lp 2 {\mR^m}$. Hence,
   \begin{equation*}
    D_\infty = \cF \circ M: \Lp 2 {\LatDom} \to   \Lp 2 {\mR^m} \cap \Lp \infty {\mR^m}
   \end{equation*}
    is well-defined and continuous. According to \eqref{eq:FresnelConv}, we have $\bDF{d} = \cM_2 \circ \cF \circ \cM_1$  where the $\cM_j$ are multiplications with functions of constant modulus. Hence, the same result holds true for
     \begin{equation*}\pushQED{\qed} 
	D_d =   \bDF{d} \circ M: \Lp 2 {\LatDom} \to   \Lp 2 {\mR^m} \cap \Lp \infty {\mR^m}. \qedhere \popQED
      \end{equation*}
 \renewcommand{\qedsymbol}{} 
  \end{pf}
  \vspace{-2em}
  
  As in the case of $\bDFlat{d}$ and $\cF_{\overline 2}$ in \eqref{eq:ForwardOpNF} and \eqref{eq:ForwardOpFF}, respectively, we write
  \begin{equation}
   D_{d,\overline 2}(\psi) := \bDFlat{d}(P\cdot \psi) \MTEXT{and} D_{\infty, \overline 2}(\psi) := \cF_{\overline 2}(P\cdot \psi) \label{eq:GenPropsTomo1}
  \end{equation}
  for the operators which propagate the complete tomographic, i.e.\ $\theta$-dependent, OTF by application of $D_d, D_\infty$ in the lateral coordinates. 
By \lemref{lem:BoundedGenProps}, this defines bounded linear operators
\begin{equation} 
    D_{d,\overline 2},  D_{\infty, \overline 2}:  \Lp 2 {\ProjDom} \to   \Lp 2 {Z^{m+1}} \cap \Lp \infty {Z^{m+1}} \label{eq:GenPropsTomo2},
   \end{equation}
      where $Z^{m+1} = [0;2 \pi) \times \mR^m$ as in \sref{S:Radon}.
  
  In the near-field case described by \eqref{eq:ForwardOpNF}, the scattered wave field is superimposed with the unscattered part of the probe beam $\bDF{d}(P)$. Note that
  \begin{equation*}
   \bDF{d}: \mP \subset \Lp 1 {\mR^m}  \oplus 1 \cdot \mC \to \Lp \infty {\mR^m}  
  \end{equation*}
  is bounded since $\bDF{d}:\Lp 1 {\mR^m} \to \Lp \infty {\mR^m}$ is bounded, as argued in the proof of \lemref{lem:BoundedGenProps}, and the constant part is simply reproduced under propagation. Thus,
  \begin{equation}
   S: \Lp \infty {Z^{m+1}} \to \Lp \infty {Z^{m+1}}; \; \psi \mapsto \bDF{d}(P) + \psi \label{eq:ProbeSuperposeOp}
  \end{equation}
  is well-defined and \Frechet differentiable where the derivative $S'[\psi]h_\psi = h_\psi$ equals the identity and may therefore be trivially extended to $\Lp p {Z^{m+1}}$ for any $p \in [1;\infty)$.
  
  \end{subsection}

 \begin{subsection}{Total Forward Operators} \label{SS:TotalForwOp}
  
  The intermediate results from \sref{SS:FrechetOTF} and \sref{SS:FrechetProbeProp} enable us to finally prove \Frechet differentiability for the forward operators of near- and far-field phase contrast tomography. Note that all function spaces have to be treated as \emph{real} Banach- or Hilbert spaces, indicated by the subscript $\mR$ introduced in \sref{S:OpsAdjoints}, in order to obtain differentiability of the squared modulus operation, see \exref{ex:Frechet}-(b). The differentiability result reads as follows:
  \vspace{1em}
  \begin{th1}[\Frechet Differentiability of the Forward Operators] \label{thm:ForwardOpsFrechet} 
   For $P \in \mP$, the forward operators given by \eqref{eq:ForwardOpNF} and  \eqref{eq:ForwardOpFF} are well-defined and \Frechet differentiable on the domains
   \begin{align*}
        F_d&: \Lp \infty {\ObjDom}_\mR \to \Lp \infty {Z^{m+1}}_\mR \\
     F_\infty&: \Lp \infty {\ObjDom}_\mR \to \Lp 1 {Z^{m+1}}_\mR \cap \Lp \infty {Z^{m+1}}_\mR
   \end{align*}
   with $Z^{m+1} = [0;2 \pi) \times \mR^m$ as in \sref{S:Radon}. For $N, h_N \in \Lp \infty {\ObjDom}_{\mR}$, the derivatives are
   \begin{subequations} \label{eq:FrechetForwOp}
    \begin{align}
        F_d'[N]h_N &= 2k^2 \Re \{ \; \cc{  \left[\bDFlat{d} \left( P \cdot  \exp \left( -\I k \CR(N) \right)  \right) \right] } \nonumber \\
   &\;\;\;\;\;\;\;\;\;\;\;\;\;\: \cdot \bDFlat{d} \left( P \cdot \exp \left( -\I k \CR(N) \right) \cdot \CR(h_N) \right) \;  \} \label{eq:FrechetForwOpNF} \\
     F_\infty'[N]h_N &= 2k^2 \Re \{ \; \cc{ \left[ \cF_{\overline 2} \left( P \cdot \left[\exp \left( -\I k \CR(N) \right) -1\right] \right) \right] }\nonumber \\
   &\;\;\;\;\;\;\;\;\;\;\;\;\;\: \cdot \cF_{\overline 2} \left( P \cdot \exp \left( -\I k \CR(N) \right) \cdot \CR(h_N) \right) \;  \} \label{eq:FrechetForwOpFF}
   \end{align}
   \end{subequations}
   Moreover, for any $N \in \Lp \infty {\ObjDom}_{\mR}$, there exist unique bounded linear extensions
    \begin{equation}
        F_d'[N], F_\infty'[N]: \Lp 2  {\ObjDom}_\mR \to \Lp 2  {Z^{m+1}}_\mR \cap \Lp \infty  {Z^{m+1}}_\mR.  \label{eq:FrechetForwOpDoms}
   \end{equation}
  \end{th1}
  \begin{pf}
  The forward operators can be decomposed as
  \begin{align*}
   F_d  = A \circ S \circ D_{d,\overline 2} \circ O_0 \MTEXT{and}  F_\infty  = A  \circ  D_{\infty, \overline 2} \circ O_0,
  \end{align*}
  where $A: f \mapsto |f|^2$ denotes the pointwise squared modulus operator. According to \exref{ex:Frechet}-\Text{(b)}, $A$ is \Frechet differentiable both on $\Lp \infty {Z^{m+1}}_\mR$ and as a map $\Lp 2 {Z^{m+1}}_\mR \to \Lp 1 {Z^{m+1}}_\mR$. Moreover, the boundedness of the propagators $ D_{d,\overline 2},  D_{\infty, \overline 2}$ defined by \eqref{eq:GenPropsTomo1} and \eqref{eq:GenPropsTomo2} remains true when restricted to 
  \begin{equation*} \Lp \infty {\ProjDom}_{\mR}   \to   \Lp 2 {Z^{m+1}}_{\mR} \cap \Lp \infty {Z^{m+1}}_{\mR} \end{equation*}
  by the continuous embeddings in \thmref{thm:LpBoundedEmbed}.
  In combination with \lemref{lem:FrechetOTF}, this implies that
  \begin{align*}
   F_\infty: \Lp \infty {\ObjDom}_\mR &\stackrel{O_0}\to \Lp \infty {\ProjDom}_\mR \stackrel{ D_{\infty, \overline 2}} \to   \Lp 2 {Z^{m+1}}_\mR \cap \Lp \infty {Z^{m+1}}_\mR \\
   &\stackrel{A} \to   \Lp 1 {Z^{m+1}}_\mR \cap \Lp \infty {Z^{m+1}}_\mR
  \end{align*}
  is well-defined and \Frechet differentiable according to \thmref{thm:FrechetProps}. The derivative is obtained by the chain rule \eqref{eq:ChainRule}, yielding  for all $N, h_N \in \Lp \infty {\ObjDom}_\mR$
  
\begin{align}
   F'_\infty[N]h_N  &= A'[ D_{\infty, \overline 2} \circ O_0(N)] \circ  D_{\infty, \overline 2}  \circ O_0'[N] (h_N) \nonumber \\
   &= 2k^2 \Re \left\{ \cc{\left[ D_{\infty, \overline 2} \circ O_0(N)\right]} \cdot  D_{\infty, \overline 2}(O_0'[N] h_N) \right\}  \label{eq:FrechetFormalFF}\\
   &= 2k^2 \Re \{ \; \cc{ \left[ \cF_{\overline 2} \left( P \cdot \left[\exp \left( -\I k \CR(N) \right) -1\right] \right) \right] } \nonumber \\
   &\;\;\;\;\;\;\;\;\;\;\;\;\;\: \cdot \cF_{\overline 2} \left( P \cdot \exp \left( -\I k \CR(N) \right) \cdot \CR(h_N) \right) \;  \}. \nonumber 
  \end{align}
  \vspace{-3em}
  
  The principal difference in the case of the near-field operator $F_d$ lies in the superposition of the unscattered probe contribution, induced by the composition with the operator $S$ defined in \eqref{eq:ProbeSuperposeOp}. Since $S$ is differentiable, so is
    \vspace{-.5em}
   \begin{align*}
   F_d: \Lp \infty {\ObjDom}_\mR &\stackrel{O_0}\to \Lp \infty {\ProjDom}_\mR \stackrel{ D_{d,\overline 2}} \to   \Lp 2 {Z^{m+1}}_\mR \cap \Lp \infty {Z^{m+1}}_\mR \stackrel{S} \to  \Lp \infty {Z^{m+1}}_\mR \\ &\stackrel{A} \to    \Lp \infty {Z^{m+1}}_\mR
  \end{align*}
  by \thmref{thm:FrechetProps}-(c,e) and \lemref{lem:BoundedGenProps}. Applying these, we obtain
      \vspace{-1em}
  \begin{align}
   F'_d[N]h_N  &= A'[ S \circ D_{d,\overline 2} \circ O_0(N)] \circ \overbrace{S'[ D_{d,\overline 2} \circ O_0(N)]}^{ = \text{ identity in }\Lp \infty {Z^{m+1}}_\mR  } \circ  D_{d,\overline 2}  \circ O_0'[N] (h_N) \nonumber  \\
   &= 2k^2 \Re \left\{ \cc{\left[S \circ  D_{d,\overline 2} \circ O_0(N)\right]} \cdot  D_{d,\overline 2}(O_0'[N] h_N) \right\} \label{eq:FrechetFormalNF} \\
   &= 2k^2 \Re \{ \; \cc{ \left[ \bDFlat{d} \left( P \cdot  \exp \left( -\I k \CR(N) \right)  \right) \right] } \nonumber \\
   &\;\;\;\;\;\;\;\;\;\;\;\;\;\: \cdot \bDFlat{d} \left( P \cdot \exp \left( -\I k \CR(N) \right) \cdot \CR(h_N) \right) \;  \} \nonumber 
  \end{align}
  for all $N, h_N \in \Lp \infty {\ObjDom}_\mR$, i.e.\ the expression given in \eqref{eq:FrechetForwOpNF}.
  
  As for the extensions of the \Frechet derivatives, note that for all $N \in  \Lp \infty {\ObjDom}_\mR$
  \begin{equation*}
  D_{\ast,2} \circ O_0'[N]: \Lp 2 {\ObjDom}_\mR \to \Lp 2 {Z^{m+1}}_\mR \cap \Lp \infty {Z^{m+1}}_\mR 
  \end{equation*}
  is a well-defined bounded linear operator according to \lemref{lem:FrechetOTF} and \lemref{lem:BoundedGenProps}. Moreover, the left hand factors in \eqref{eq:FrechetFormalFF} and \eqref{eq:FrechetFormalNF} are in $\Lp \infty {Z^{m+1}}_\mR$ for $N \in \Lp \infty {\ObjDom}_\mR$ so that the multiplication with these defines bounded linear maps
    \begin{equation*} 
    \Lp 2 {Z^{m+1}}_\mR \cap \Lp \infty {Z^{m+1}}_\mR \to\Lp 2 {Z^{m+1}}_\mR \cap \Lp \infty {Z^{m+1}}_\mR 
    \end{equation*}
  By continuity of $\Re: \Lp p {Z^{m+1}}_\mR \to \Lp p {Z^{m+1}}_\mR$ (see \exref{ex:Adjoints}-(c)), these observations imply the existence of the bounded extensions given in \eqref{eq:FrechetForwOpDoms}. By \thmref{thm:DenseIncl}, $\Lp \infty {\ObjDom}_\mR$ is dense in $\Lp 2 {\ObjDom}_\mR$ so that the latter are unique.
  \end{pf}
  \vspace{1em}
  
  Note that the necessity for the larger image space $\Lp \infty {Z^{m+1}}_\mR$ in \thmref{thm:ForwardOpsFrechet} for the near-field forward operator $F_d$ arises from the contributions of the unscattered probe beam. For the chosen class of probe functions $P \in \mP$, the plane wave component in the incident illumination is infinitely extended in the lateral dimensions, i.e.\ does not vanish at infinity. Consequently,
  \vspace{-.5em}
  \begin{equation} F_d(\Lp \infty {\ObjDom}_\mR) \not \subset  \Lp p {Z^{m+1}}_\mR \MTEXT{for any} p < \infty.  \vspace{-.5em}\end{equation}
  A stronger result with respect to the image space of $F_d$ might be obtained by subtracting the constant empty-beam intensities $|\bDF{d} ( P ) |^2$ from the scattering contributions which encode the desired object information $N$. However, this step is omitted here to retain notational simplicity.
  
  The $L^2$-extensions of the \Frechet derivatives allows to study these on \emph{Hilbert spaces}. In particular, this permits to define \emph{adjoints} being of significance for the reconstruction methods introduced in \chapref{C:NumMeth}. As differentiability implies continuity, we further obtain the following corollary:
    \vspace{.5em}
  \begin{cor1}[Continuity of the Forward Operators] \label{cor:ForwardContinuity}
   The forward operators $F_d$ and $F_\infty$ on the spaces studied in \thmref{thm:ForwardOpsFrechet}  are \emph{continuous}.
  \end{cor1}

\end{subsection}

\end{section}

\begin{section}{Well-Posedness and Ill-Posedness} \label{S:GenIll-posed}

  As shown in \sref{S:TomoTheory}, phase contrast tomography gives rise to an inverse problem
  \vspace{-.5em}
  \begin{equation}
   F(N^\dagger ) = I^\dagger \label{eq:InverseProblem} \vspace{-.5em}
  \end{equation}
  where $F \in \{F_d, F_\infty\}$ is the forward operator mapping the unknown refractive index $1- N^\dagger$ of the specimen onto the corresponding scattering intensity data $I^\dagger$. A general characterization of the solvability of such problems is due to Hadamard \cite{Hadamard1902}:
  
  \vspace{.5em}
  \begin{df1}[Well-Posedness and Ill-Posedness \text{\cite{Hadamard1902}}] \label{def:IllWell}
  A problem of the form \eqref{eq:InverseProblem} is called \emph{well-posed} if for all admissible data
  \begin{itemize}
   \item[\Text{(a)}] A solution  exists
      \item[\Text{(b)}] The solution is unique
        \item[\Text{(c)}] The solution depends continuously on the data
  \end{itemize}
  Otherwise, it is called \emph{ill-posed}.
  \end{df1}
  \vspace{1em}
  It is evident that a problem for which conditions $\Text{(a)}$ and/or $\Text{(b)}$ are violated does not allow for a reasonable reconstruction of - in our case - the refractive index. On the other hand, if these are satisfied, then the inverse map $F^{-1}$ exists. Hence, $N^\dagger $ can be reconstructed uniquely from the \emph{exact} data $I^\dagger$. Condition $\Text{(c)}$ in \defref{def:IllWell} corresponds to continuity of $F^{-1}$.
  Its significance arises from the physical fact that any realistic experiments are subject to noise and other measurement errors. In the considered setting, this implies that only a systematically perturbed and/or noisy version of the intensities
  \begin{equation}
   I^{\Textbf{err}} = I^\dagger + \Textbf{err} \label{eq:NoisyData}
  \end{equation}
   can be observed. Now, the claimed continuity ensures that small observations errors $\Textbf{err}$ result in small deviations of the corresponding reconstruction $N^{\Textbf{err}} := F^{-1} ( I^{\Textbf{err}} )$ from the exact solution $N^\dagger$ in the corresponding topologies. Conversely, if part (c) of \defref{def:IllWell} is violated, then $\Textbf{err} \to 0$ need not imply $N^{\Textbf{err}} \to N^{\dagger}$, i.e.\ measurement errors may be amplified by arbitrary factors in the reconstruction.
  
  The \emph{forward problem} of phase contrast tomography, i.e.\ the assignment of scattering intensities $I$ for given specimen data $N$, is implemented by the maps $F_d, F_\infty$. From the well-definedness of these (see \thmref{thm:ForwardOpsFrechet}), it is immediately clear that unique intensity patterns $I \in \Lp \infty {Z^{m+1}}_{\mR}$ exists for all objects $N \in \Lp \infty {\ObjDom}_{\mR}$. Moreover, $I$ depends continuously on $N$ according to \cref{cor:ForwardContinuity}. In the language of \defref{def:IllWell}, this yields the following result:
  \vspace{1em}
  \begin{res}[Well-Posedness of the Forward Problem] \label{res:Well-posedForward}
   The forward problem of phase contrast tomography, given by the evaluation of \eqref{eq:ForwardOpNF} or \eqref{eq:ForwardOpFF}, is well-posed.
  \end{res}
  \vspace{1em}
  
  Now we turn to the \emph{inverse} \probref{prob:1}, i.e.\ the reconstruction of $N$ from observed non-ideal intensities $I$.
  By considering \eqref{eq:ForwardOpNF} and \eqref{eq:ForwardOpFF}, it can be seen that this requires the subsequent solution of essentially four subproblems:
  \vspace{1em}
   \begin{itemize}
   \item[$\boldsymbol 1 $] Reconstruct the exit-wave fields $P\cdot O(N)$ from $I = \left| \bDd{d}\left( P \cdot O(N) \right) \right|^2$
         \item[$\boldsymbol 2 $] Recover the object transmission function $O(N)$ from $P\cdot O(N)$
      \item[$\boldsymbol 3 $] Solve the object transmission function $O(N)$ for the projections $\CR(N)$
        \item[$\boldsymbol 4 $] Invert the cylindrical Radon transform $\CR$ to obtain $N = \delta - \I \beta$
  \end{itemize}
    \vspace{1em}
    
  Step $\boldsymbol 1 $ amounts to solving a non-linear \emph{phase retrieval problem} of the form discussed in \sref{S:PhaseContrast}. The ill-posedness of this subproblem is studied in \sref{S:PhaseRetrieval}.
  
  Part $\boldsymbol 2 $ is conceptually the simplest one, as it boils down to a division by the probe function $P$. This operation is unique and continuous, i.e.\ well-posed, if
  \begin{equation}
   |P(\bx)| > \varepsilon \MTEXT{for all} \bx \in \mR^{m} \label{eq:ProbeBoundedBelow}
  \end{equation}
for some lower bound $\varepsilon > 0$. This holds for instance for plane wave illumination, i.e.\ for $P$ constant. On the contrary, the reconstruction of $O(N)$ is discontinuous if $P$ approaches zero and even non-unique if $P$ vanishes on some open set $U \subset \mR^{m}$.

  By \eqref{eq:OTFbyRadon}, step $\boldsymbol 3 $ is equivalent to the inversion of a pointwise exponential. As discussed in \sref{SS:PhaseWrap}, this operation is well-posed with respect to the imaginary phase $-\I k \CR(\delta)$ if and only if phase-wrapping is absent. Likewise, the reconstruction of the real-part of the exponent, corresponding to the projected attenuation $- k \CR(\beta)$, gives rise to large error amplifications only in the case of very strong absorption. Both excessive absorption and phase-wrapping may be precluded experimentally by suitably choosing the wavelength of the incident X-rays for a given sample.
  
  Both step $\boldsymbol 2 $ and $\boldsymbol 3 $ thus give rise trivial - but not necessarily insignificant - sources of ill-posedness.
  For this reason, a further study of these is omitted.
  On the other hand, the final reconstruction step $\boldsymbol 4 $  given by a inverse cylindrical Radon transform  is subject to further analysis in \sref{S:RadonIll-posed}.

\end{section}

\begin{section}{Radon Inversion} \label{S:RadonIll-posed}

The objective of this section is a characterization of the inverse cylindrical Radon transform $\CR$, which has to be evaluated as part of \probref{prob:1}, in terms of \defref{def:IllWell}. The principal tool for this analysis is provided by the Fourier Slice \thmref{thm:FourierSlice}. Owing to the slice-wise definition of $\CR$, it is sufficient to consider the two-dimensional Radon transform $\cR$. Throughout this section, we regard $\cR$ as a map on compactly supported $L^2$-functions, i.e.\ as an operator
\begin{equation}
 \cR: \Lp 2 \Omega \to \Lp 2 {Z^2} \label{eq:RadonDoms}
\end{equation}
for some measurable and bounded $\Omega \subset \mR^{2}$.

\begin{subsection}{Existence and Consistency} \label{SS:RadonSurjective}

As a first step, we study which functions $g \in \Lp 2 {Z^2}$ have a preimage under the Radon transform, corresponding to the \emph{existence} of solutions for the inverse problem $\cR(f) = g$. From a physical perspective it is evident that the projections $\cR(f)(\theta, \cdot)$ of a single two-dimensional object $f \in \Lp 2 \Omega$ for different incident angles $\theta \in [0; 2\pi)$ cannot be entirely independent, but may be highly correlated if the difference in $\theta$ is small. In Fourier space, these correlations can be understood by virtue of \thmref{thm:FourierSlice}: the projections are such that their lateral Fourier transforms may be consistently arranged to a polar Fourier transform of a function $ f \in \Lp 2 \Omega$. In particular, the Fourier ``slices'' of any $g = \cR(f)$ must be consistent at 0 where they all intersect in the polar grid. In real space, this give rise to the \emph{Helgason-Ludwig consistency conditions}, which are even sufficient in suitable Schwartz-spaces \cite{Ludwig1966Radon}:
  \vspace{1em}
  \begin{th1}[Helgason-Ludwig Consistency Conditions \cite{Helgason1965Radon,Ludwig1966Radon}] \label{thm:HelgasonLudwig} 
   Let $\Omega \subset \mR^{2}$ be measurable and bounded. For a function $g \in \Lp 2 {Z^2}$ to have preimage $f \in \Lp 2 \Omega$ such that $g = \cR(f)$, it is necessary that
    \begin{itemize} 
         \item[$\Text{(a)}$] $g$ is compactly supported
     \item[$\Text{(b)}$] $g(\theta + \pi, -x) = g(\theta, x)$ for all $(\theta, x) \in Z^2$
          \item[$\Text{(c)}$] For all $k \in \mN_0$,  $\int_{\mR } g(\theta, x) x^k \; \D x$ is almost everywhere equal to a polynomial of degree $\leq k$ in $(\cos\theta, \sin\theta) $.
    \end{itemize}
  \end{th1}
  \vspace{1em}
  Compactness of the support is evident from the definition of the Radon transform. Condition $\Text{(b)}$ simply states that the projections will be reproduced up to a reflection under a rotation of the incident angle by exactly $180^\circ$. Hence, the sinogram of a function $f$ is uniquely defined already by all incident angles $\theta \in [0;\pi)$.
  
  Condition $\Text{(c)}$ in \thmref{thm:HelgasonLudwig} is a bit more involved. For $k = 0$, it  states that $\int_{\mR } g(\theta, x)  \; \D x$ must not depend on $\theta$. This holds for $g = \cR(f)$ since
  \begin{equation*}
   \int_{\mR }  \cR f (\theta, x)  \; \D x = \int_{\mR^2 } f(x \bn_{\theta} + y \bn_{\theta}^\perp) \; \D y \D x = \int_{\mR^2 } f(x, y) \; \D y \D x.
  \end{equation*}
  The constraint for $k=1$ implies that the center of mass of the projections must move on a sinusoidal curve as $\theta$ is varied. Analogously, the conditions for $k \geq 2$ define admissible variations in $\theta$ of higher-order moments.

  Part (b) can be incorporated experimentally by restricting the measurement to half of the sinogram. Likewise, the consistency conditions (c) for the $k$-th lateral moments may be exploited to partially recover the projection data for incident angles that cannot be measured \cite{KudoSaito1991SinogramRecovery}. On the other hand, \thmref{thm:HelgasonLudwig}-(c) corresponds to a significant and complicated restriction of the image space $\cR( \Lp 2 \Omega ) \subset \Lp 2 {Z^2}$. Consequently, noisy or systematically perturbed data $g^{\Textbf{err}} = \cR(f) + \Textbf{err}$ will almost surely violate the Helgason-Ludwig consistency conditions for any realistic error model, i.e.\ will not admit for an exact solution $\cR^{-1} ( g^{\Textbf{err}} )$. In other words, $\cR$ is \emph{not surjective} for any simple restriction of its image space which implies that Radon inversion violates the well-posedness-condition (a) in \defref{def:IllWell}.
  
  This problem has to be accounted in numerical implementations, for example by choosing iterative methods such as the \emph{Algebraic Reconstruction Technique} \cite{Kaczmarz1937ART,Gordon1970ART}. In this approach, the current iterates $f_j$ are subsequently projected to match the tomographic data $g^{\Textbf{err}}(\theta_k, \cdot)$ for a new incident angle $\theta_k$. By suitably regularizing this method, data inconsistencies in the sense of \thmref{thm:HelgasonLudwig} may be suppressed.
\end{subsection}

\begin{subsection}{Uniqueness} \label{SS:RadonInjective}

The next step is to investigate whether the Radon transform $\cR$ is \emph{injective},
i.e.\ whether a function may be reconstructed uniquely from its sinogram $\cR(f)$. Here, the answer is positive and follows from the fact that a measurement of the Radon transform is equivalent to a sampling in Fourier space by \thmref{thm:FourierSlice}:
  \vspace{1em}
  \begin{th1}[Injectivity of the Radon Transform \text{\cite[p. 11]{Natterer}}] \label{thm:RadonInjective}
   The Radon transform $\cR: \Lp 2 \Omega \to \Lp 2 {Z^2}$ is injective and any $f \in \Lp 2 \Omega$ is uniquely determined by the projections $\cR f_{| V \times \mR} $ restricted to an arbitrary open set of angles $V \subset [0 ;  2 \pi )$.
  \end{th1}
  \begin{pf}
   Let $ f \in \Lp 2 \Omega$ such that $\cR f_{| V \times \mR} = 0$ for some open set $V \subset [0 ;  2 \pi )$. Then \thmref{thm:FourierSlice} implies that
   \begin{equation*}
   \cF_{\Text p }f _{| V \times \mR} = \cF_2 ( \cR f_{| V \times \mR} ) = 0,
   \end{equation*}
   i.e.\ $\cF(f)$ vanishes on the wedge $W = \{ x \bn_{\theta} : x \in \mR, \theta \in V \}$ with non-empty interior $\Int(W)$ by \eqref{eq:PolarFT}. As $ \Omega$ is bounded, $f$ can be interpreted as a distribution of compact support. Thus, $\cF f$ has a unique extension to an entire function in $\mC^2$ by \thmref{thm:PaleyW}. Via Taylor-expansion in $\Int(W)$, this implies $\cF(f) = 0$ and hence $f= 0$ according to \cref{cor:FourierL2}. By linearity of $\cR$, this proves the claim. 
  \end{pf}
  \vspace{1em}
  In addition to mere injectivity, \thmref{thm:RadonInjective} states that the measurement of projections for an arbitrarily small - but continuous - interval of incident angles is sufficient for the unique reconstruction of compactly supported $L^2$-functions. Note, however, that the analytic continuation in Fourier space applied in the proof is highly sensitive to errors and thus can hardly be used in practical reconstructions. Nonetheless, the result bears some significance for the completion of ``missing wedges'', i.e.\ small angular sections from which no projections may be measured due to experimental constraints, e.g.\ due to obstructing instruments. The good news in \thmref{thm:RadonInjective} is that this does not preclude unique reconstructions in principal.

\end{subsection}

\begin{subsection}{Ill-Posedness} \label{SS:RadonContiInv}
 
 According to \thmref{thm:RadonInjective}, the inverse Radon transform
 \begin{equation}
   \cR^{-1}: \Lp 2 {Z^2} \supset \cR ( \Lp 2 \Omega ) \to \Lp 2 \Omega \label{eq:inverseRadonDoms}
 \end{equation}
 exists. This operator can even be described by explicit inversion formulae, see \cite[sec. II.2]{Natterer} for an overview.
 A particularly straightforward option is to inverse Fourier transform the polar Fourier data obtained from $\cR(f)$ in the spirit of \thmref{thm:FourierSlice}. Setting $f \in \Lp 2 \Omega$ arbitrary, $g = \cR(f) \in \Lp 2 {Z^2}$ and
 \begin{equation} \Xi_{ 1 / 2  } (g)(\theta, \xi) := \xi^{\frac 1 2 }\cF_2(g)(\theta, \xi) \stackrel{\eqref{eq:FourierSlice}}=  \xi^{\frac 1 2 }\cF_{\Text p }(f)(\theta, \xi)  \label{eq:xi12Op} \MTEXT{for all} (\theta, \xi ) \in Z^2, \end{equation}
 this yields by a transformation from Cartesian to polar coordinates
 \begin{align}
  \norm{f}^2_{\Lp 2 \Omega} &\stackrel{\text{Cor. \ref{cor:FourierL2}}}= \norm{\cF f}^2_{\Lp 2 {\mR^2}} = \int_{\mR^2} \left| \cF f(\bxi) \right|^2 \; \D \xi  \stackrel{\bxi = \sigma \bn_\theta}= \int_{0}^\infty \int_{0}^{2\pi} \sigma | \underbrace{\cF f(\sigma \bn_\theta)}_{= \sqrt{2\pi}\cF_{\Text p } f (\theta, \sigma) } |^2 \; \D \theta \D \sigma \nonumber \\
  &\;\;\stackrel{\eqref{eq:xi12Op}}=  (4\pi)^{- 1} \int_{\mR} \int_{0}^{2\pi}  |  \Xi_{ 1 / 2  } (g) (\theta, \sigma) |^2 \; \D \theta \D \sigma    = (4\pi)^{- 1} \norm{ \Xi_{ 1 / 2  } (g) }_{\Lp 2 {Z^2}}^2. \nonumber 
 \end{align}
 
 Accordingly, the assignment $\Xi_{ 1 / 2  } (g) \mapsto f$ is isometric in $L^2$-norm up to a constant factor. On the other hand, the operator $\Xi_{ 1 / 2  }$ defined by \eqref{eq:xi12Op} is \emph{unbounded} in $\Lp 2 {Z^2}$, weighting every Fourier component in the lateral coordinate with the square root of the corresponding frequency $\xi^{\frac 1 2 }$. This implies that $\cR^{-1}: g \mapsto f$ amplifies noise of frequency $\xi$ by this slowly but unboundedly growing factor. In this sense $\cR^{-1}$ is \emph{discontinuous} on the domains in \eqref{eq:inverseRadonDoms} so that the problem of Radon inversion violates part (c) of \defref{def:IllWell}. The characteristic noise amplification is associated with a weak smoothing of the corresponding forward operation $f \mapsto \cR(f)$, suppressing high frequency components of a signal $f$ in its sinogram. This can alternatively be seen from the singular value decomposition of $\cR$ as derived for instance in \cite[sec.~IV.3]{Natterer}.
 
 The following result summarizes the observations of this section:
  \vspace{1em}
  \begin{res}[Ill-Posedness of the Radon Inversion] \label{res:Ill-posedRadon}
   Radon Inversion of compactly supported signals is unique - even if the projections are not measured for all incident angles. However, the problem may not have a solution for inconsistent projections and is mildly ill-posed as data noise of frequency $\xi$ is amplified by factors $\sim \xi^{\frac 1 2 }$. 
  \end{res}
  \vspace{1em}
 
\end{subsection}

\end{section}

\begin{section}{Phase Retrieval} \label{S:PhaseRetrieval}

As discussed in \sref{S:GenIll-posed}, our principal tomographic imaging \probref{prob:1} involves a phase retrieval problem which arises from the loss of phase information in the detection of the propagated wave field, see \sref{SS:Detection}. In addition to the ill-posedness of the other subproblems studied in the preceding sections, this characteristic defect in the measurement process may be expected to prevent unique reconstructions. This section is therefore dedicated to the theory of phase retrieval in the considered settings of near- and far-field phase contrast imaging, discussing sufficient conditions that ensure unique recovery of the missing phase information.

\begin{subsection}{Abstract Formulation} \label{SS:PhaseRetrAbstract}

 The inherent phase retrieval step in \probref{prob:1} amounts to reconstructing the contact image $P \cdot O_0(N)(\theta, \cdot)$ for any incident angle $\theta$ from the corresponding scattering intensities.
 In the near-field case governed by \eqref{eq:ForwardOpNF}, the latter are given by
 \begin{equation*}
 I_d(\theta, \cdot)  = \left| \bDF{d}(P_d) + \bDFlat{d} \left( P \cdot O_0 (N)(\theta, \cdot) \right)  \right|^2.
\end{equation*}
  Introducing dimensionless coordinates $\bxi := \left( k / d \right)^{\frac 1 2}\bx $ for all $\bx \in \mR^{m}$ and
 \begin{equation}
 \psi (\bxi ) :=  \nu P  \cdot O_0(N)(\theta, \bx ), \;\;\;\; P_d(\bxi) := \bDF{d}(P) (\bx), \;\;\;\; I(\bxi) := I_d(\theta, \bx), 
  \end{equation}
  this expression can be rewritten using the convolution formulation of the Fresnel propagator \eqref{eq:FresnelConv}. Setting  $\nF(\bxi )  = \exp( \I \bxi^2 /2 ) $, this yields for all $\bxi \in \mR^{m}$
  \begin{equation}
 I (\bxi) =  \left|  \exp \left( - \frac{\I \bxi^2} 2 \right)  P_d(\bxi) +  \cF \left(\nF \cdot \psi  \right)(\bxi)  \right|^2. \label{eq:PhaseRetrAbstractNF}
\end{equation}
As the rescaling operations and coordinate transforms are invertible, we obtain the following abstract formulation of the near-field phase retrieval problem:
\vspace{1em}
	\begin{prob}[Near-Field Phase Retrieval]\label{prob:NF}
		For $\A \subset \Sdash{m}$ and known probe function $P_d$, reconstruct $\psi \in \A$ from  intensity data $I$ given by \eqref{eq:PhaseRetrAbstractNF}.
	\end{prob}
\vspace{1em}
\noindent Note that we allow for \emph{distributional} solutions $\psi \in \Sdash{m}$ although \eqref{eq:PhaseRetrAbstractNF} might not even be well-defined for such. The choice of the subset $\A \subset \Sdash{m}$ restricts the admissible solutions, corresponding to a certain \emph{a priori knowledge} of $\psi$.

Now consider the \emph{far-field} setting. According to \eqref{eq:ForwardOpFF}, the measured intensities under an incident angle $\theta$ in this case can be written in the abstract form
\begin{equation}
 I  =  \left|  \cF \left(  \psi \right)   \right|^2 \label{eq:PhaseRetrAbstractFF}
\end{equation}
by setting $\psi := P \cdot O_0(N)(\theta, \cdot ) $ and  $I (\bxi) := I_\infty( \theta, \bxi )$. Hence, the problem of far-field phase retrieval can be stated as follows:
 \vspace{1em}
	\begin{prob}[Far-Field Phase Retrieval]\label{prob:FF}
	For $\A \subset \Sdash{m}$, reconstruct $\psi \in \A$ from  intensity data $I$ given by \eqref{eq:PhaseRetrAbstractFF}.
\end{prob}
\vspace{1em}
\noindent Phase retrieval in far-field imaging thus amounts to the reconstruction of a function from the squared modulus of its Fourier transform. Variants of this abstract problem arise also in many other imaging contexts, as for instance in crystallography \cite{Millane1990,Cruickshank1961PhaseProbCrystal} and electron microscopy \cite{Misell1973PhaseProbEMicroscopy}. \probref{prob:FF} is therefore addressed in \sref{SS:PhaseRetrFF}.

It can be seen from \eqref{eq:PhaseRetrAbstractNF} that the near-field phase retrieval \probref{prob:NF} is closely related: here, the product of $\nF \cdot \psi$ is Fourier transformed and superimposed as a \emph{perturbation} upon the probe contributions in the first summand of \eqref{eq:PhaseRetrAbstractNF}. The intensity data is given by the squared modulus of this superposition of known background probe field and unknown perturbation. The principal difference to the far-field case thus lies in the presence of a known reference signal. Implications for the uniqueness of \probref{prob:NF} are explored in \sref{SS:PhaseRetrNF}.

By construction, the \emph{contact image} $\psi$ to be reconstructed in \eqref{eq:PhaseRetrAbstractNF} and \eqref{eq:PhaseRetrAbstractFF} is essentially given by $O_0(N)(\theta, \cdot)$ and thus has compact support whenever $N$ is compactly supported, i.e.\ for specimen of bounded spatial extent. We will widely restrict to this case assuming $\psi \in  \A \subset \Sdashc{m}$ in  \probref{prob:NF} and \probref{prob:FF}.
From \thmref{thm:PaleyW}, it then follows that $\cF \left( \psi \right)$ and $\cF ( \nF \cdot \psi )$ are \emph{entire functions} (see \sref{S:SchwartzLp}) which simplifies analysis considerably. In particular, \eqref{eq:PhaseRetrAbstractNF} and \eqref{eq:PhaseRetrAbstractFF} are well-defined in this case as the Fourier transforms are indeed $\sC^\infty$-functions.

The following section therefore introduces some notions from the theory of entire functions of a single variable, providing a classical approach to the uniqueness theory of phase retrieval applicable to the abstract problems motivated here.

 \end{subsection}

 \begin{subsection}{Preliminaries: Entire Functions of one Complex Variable} \label{SS:EntireFunctions}



The given overview on entire functions is adapted from the manuscript \cite{maretzke2014uniqueness}, which has been submitted to \emph{Inverse Problems}. The introduction in the latter is based on a more detailed treatment in \cite{BoasEntire,ConwayOneComplex1,FreilingEntire}.

In general, \emph{entire functions} are maps 
\begin{equation*} f: \mC^m \to \mC^m \end{equation*}
which are everywhere analytic, i.e.\ characterized by a globally convergent Taylor series. For simplicity, we restrict to the univariate case $m=1$, i.e.\ assume that $f$ is an entire function in $\mC$.

An important characterization is given by the growth behavior. Therefore, we set
 \begin{align}
                M_f(r) &:= \max_{\xi \in \mC: |\xi| = r} |f(\xi)| \MTEXT{and} m_f(r) := \min_{\xi \in \mC: |\xi| = r} |f(\xi)|. \label{eq:DefmM}
        \end{align}
Asymptotic bounds for $M_f$ give rise to the definition of its \emph{order} $\lambda_f$ and \emph{type} $\tau_f$ of the entire function $f:\mC \to \mC$:
\begin{subequations} \label{eq:DefOrderType}
\begin{align}
 \lambda_f &:= \begin{cases} 0 &\text{for $f$ constant} \\
\limsup_{r \to \infty} \frac{ \log \log M_f(r) }{ \log r } &\text{else}  \end{cases} \label{eq:DefOrder} \\
\tau_f &:= \begin{cases}  0  &  \text{if $\lambda_f =0$} \\
\limsup_{r \to \infty} r^{-\lambda_f} \log M_f(r) & \text{else} \end{cases} \label{eq:DefType}
\end{align}
\end{subequations}
If $f$ is of finite order and type $0 < \lambda_f, \tau_f < \infty$, we write \begin{equation} f(\xi) = \Or(\exp(\tau_f|\xi|^{\lambda_f})). \end{equation} Moreover, we say that order 1 entire functions are of \emph{exponential order}. According to \thmref{thm:PaleyW}, the Fourier transform of any compactly supported tempered distribution is an entire function of at most exponential order and finite type.

Alternatively, an entire function $f:\mC \to \mC$ may be characterized by its zeros in the complex plane. If $f$ is not identically zero, its roots counted by their multiplicity form an at most countably infinite sequence
\begin{equation} Z_f := \{a_j\}_{j \in J} \subset \mC \setminus \{0\}, J \subset \mN \end{equation}
with no accumulation point in $\mC$. Note that we exclude a possible zero in the origin from this definition and that $Z_f$ is assumed to be monotonically increasing in modulus. Classifying the asymptotic behavior of $Z_f$, we define the \emph{convergence exponent} $\rho_f \in [0;\infty]$ and \emph{rank} $p_f \in \mN_0 \cup \{\infty\}$ of $f$ by
\begin{subequations} \label{eq:DefConvExpRank}
\begin{align}
 \rho_f &:= \inf \{ \rho \geq 0 : \sum_{j \in J} |a_j|^{-\rho} < \infty \}  \label{eq:DefConvExp} \\
 p_f &:= \min \{ p \in \mN_0 : \sum_{j \in J} |a_j|^{-(p+1)} < \infty \}.  \label{eq:DefRank}
\end{align}
\end{subequations}
Note that the definitions do not depend on $f$ directly so that the notions can be generalized to arbitrary ordered sequences in $\mC$.

It turns out that order and rank are closely related. Combined with Weierstrass' factorization theorem for holomorphic functions, which is discussed for instance in \cite{ConwayOneComplex1}, this observation leads to Hadamard's factorization theorem:
	\vspace{1em}\begin{th1}[Hadamard's factorization theorem \cite{BoasEntire,ConwayOneComplex1}]\label{thm:Hadamard}
		Let $f: \mC \to \mC$ be an entire function of finite order $\lambda_f$ and not identically zero. Then $f$ has rank $p_f \leq \lambda_f$ and it admits a factorization
                \begin{equation} f(\xi) = \xi^m \, \exp(q_f(\xi)) \prod_{j \in J} E_{p_f} \left( \frac{\xi}{a_j} \right)  \MTEXT{for all} \xi \in \mC, \label{eq:Hadamard} \end{equation}
                where $m \in \mN_0$ is the order of the zero at $\xi = 0$, $q_f$ a polynomial of degree $\leq \lambda_f$ and
                \begin{equation*} E_n(z) = (1-z)\exp \left( \sum_{j = 1}^n \frac{z^j} j \right). \end{equation*}
		The product in \eqref{eq:Hadamard} converges uniformly on any compact subset $K \in \mC$. 
	\end{th1}\vspace{1em}
A converse variant of \thmref{thm:Hadamard} is also true: for any sequence of zeros of finite convergence exponent, \emph{canonical products} of the form \eqref{eq:Hadamard} define entire functions of finite order:
	\vspace{1em}\begin{th1}[Borel \cite{FreilingEntire}]\label{thm:OrderCanProducts}
		Let $\{a_j\}_{j \in J} \subset \mC \setminus \{0\}, J \subset \mN$ be a possibly finite sequence of monotonically increasing modulus, finite rank $p$ and convergence exponent $\rho$.
		Then a product of the form \eqref{eq:Hadamard} defines an entire function $f$ for any polynomial $q_f$ and $m \in \mN_0$. Moreover, the order of $f$ is
		\begin{equation*}
		 \lambda_f = \max \{ \deg(q_f), \rho \}.
		\end{equation*}
	\end{th1}\vspace{1em}
For convenience, we define the \emph{Schwarz reflection} $\adj{f}$ of a entire function $f$ by
\begin{equation}
 \adj{f} (\xi ) := \cc{f (\cc{\xi})} \MTEXT{for all} \xi \in \mC. \label{eq:DefSchwarzRefl}
\end{equation}
Note that $\adj{f}$ is entire, has the same order, type, convergence exponent and rank as $f$ and $f^{\ast\ast} = f$. With this notation, \thmref{thm:Hadamard} allows to quantify the amount of information gained by measuring the squared modulus $|f|^2$ of an entire function $f: \mC \to \mC$ on a segment of the real line.

\vspace{1em}\begin{lem1}[Phase Retrieval for Entire Functions \cite{Akutowicz1956I,Akutowicz1957II,Walther1963}]\label{lem:InfoModulus2}
		Let $f$, $\tilde f$ be entire functions of finite order $\lambda_{f} \geq \lambda_{\tilde f}$ such that for some $U \subset \mR$ open 
		\begin{equation*}
		 |f|_{|U}^2 = |\tilde f|_{|U}^2.
		\end{equation*}
		Then there exist entire functions $f_1,f_2$ of order $\leq \lambda_{f}$ such that
		\begin{equation}
		 f = f_1 \cdot f_2  \qquad \text{and} \qquad \tilde f = f_1 \cdot \adj{f_2}. \label{eq:FactCharacterize}
		\end{equation}
		Moreover, if $\{a_j\}_{j \in J} \subset \mC \setminus \{0\}$ are the non-zero roots of $f$, there exists a polynomial $Q$ of degree $\leq \lambda_f$ with imaginary coefficients and $K\subset J$ such that for all $\xi \in \mC$
		\begin{align} 
		\tilde f(\xi) =  \left(  \exp(  Q(\xi) ) \prod_{j \in K}\frac{  E_{p_f}\left( \xi /\cc{a_j} \right) }{ E_{p_f} \left( \xi / a_j \right) } \right) f(\xi) \label{eq:HadamardChar} .
		\end{align}
		Conversely, if $f_1$ and $f_2$ are entire functions of order $\lambda$, then $f$ and $\tilde f$ defined by \eqref{eq:FactCharacterize} are entire functions of order $\leq \lambda$ satisfying $|f|_{|\mR}^2 = |\tilde f|_{|\mR}^2$.
	\end{lem1}\vspace{1em}
	\begin{pf}
	 It is sufficient to consider the case $f \neq 0$. By \thmref{thm:Hadamard}, $f$ admits a factorization of the form \eqref{eq:Hadamard}. Noting that $|f^2|_{|U} = (f  \cdot  f ^\ast)_{|U}$ and that $f^\ast$ is entire, we find that $|f|^2$ has an extension to an entire function $F = f  \cdot  f ^\ast$ of order $\leq \lambda_f$. As such, $F$ is uniquely determined by its values on $U$ (e.g.\ by Taylor expansion) and thus coincides with the respective entire extension $\tilde F$ of $|\tilde f|_{|U}^2$. By the factorization of $f$, we obtain for all $\xi \in \mC$
		\begin{equation*}
		 \tilde F(\xi) = F(\xi) = f(\xi) \cc{f(\cc \xi)} = \xi^{2m} \, \exp(2\Re(q_f)(\xi)) \prod_{j \in J} E_{p_f} \left( \frac{\xi}{a_j} \right)  \cdot E_{p_f} \left( \frac{\xi}{\cc{a_j}} \right). 
		\end{equation*}
	In particular, we find that $F$ uniquely determines all zeros of $f$ modulo complex conjugation, as well as the real parts of the coefficients of $q_f$ and $m \in \mN_0$. Consequently, $\tilde f$ may differ from $f$ at most by a subset $K \subset J$ of ``flipped'' zeros and a multiplicative factor $\exp(Q)$, where $Q$ is a polynomial of degree $\leq \lambda_{f}$ with purely imaginary coefficients. Thus, the Hadamard factorization of $\tilde f$ is given by
 	\begin{align}  \tilde f(\xi) &= \xi^m \, \exp(q_f(\xi) + Q(\xi) ) \prod_{j \in J \setminus K} E_{p_f} \left( \frac{\xi}{a_j} \right) \cdot \prod_{j \in K} E_{p_f} \left( \frac{\xi}{\cc{a_j}} \right)  \nonumber \\
 	&=\left(  \exp(  Q(\xi) ) \prod_{j \in K}\frac{  E_{p_f}\left( \xi /\cc{a_j} \right) }{ E_{p_f} \left( \xi / a_j \right) } \right) f(\xi) \label{eq:HadamardCharPf}\\
 	&= \underbrace{ \left( \xi^m \exp\left(q_f(\xi) + \frac{Q(\xi)} 2  \right) \prod_{j \in J \setminus K} E_{p_f} \left( \frac{\xi}{a_j} \right) \right)}_{=: f_1(\xi)} \cdot \underbrace{\left(   \exp \left( \frac{Q(\xi)} 2  \right) \prod_{j \in K} E_{p_f} \left( \frac{\xi}{\cc{a_j}} \right)  \right)}_{=: \adj{f_2}(\xi)}. \nonumber \end{align}
 	Since the convergence exponents of the subsequences $\{a_j\}_{j \in J \setminus K}$ and  $\{a_j\}_{j \in K}$ are at most as large as that of the total one, $f_1$ and $f_2$ are entire functions of order $\leq \lambda_f$ according to \thmref{thm:OrderCanProducts}. Noting that $Q^\ast = -Q$, we further obtain for all $\xi \in \mC$
 	\begin{align*}
 	   f_1(\xi) f_2(\xi) &= \xi^m \exp\left( q_f(\xi) + \frac{(Q + \adj Q)( \xi)} 2 \right) \left( \prod_{j \in J \setminus K} E_{p_f} \left( \frac{\xi}{a_j} \right) \right)\cdot \left(   \prod_{j \in K} E_{p_f} \left( \frac{\xi}{\cc{\cc{a_j}}} \right)  \right)  \\ 
 	  &= \xi^m \exp(q_f(\xi))  \prod_{j \in J } E_{p_f} \left( \frac{\xi}{a_j} \right) = f(\xi).
 	\end{align*}
 	This proves the first claim. The second claim is shown by the second equality in \eqref{eq:HadamardCharPf}. For the converse statement, we simply note that for all $x \in \mR$
 	   \begin{equation*}\pushQED{\qed} 
	  |f(x)|^2 =  (f \cdot f^\ast)(x) = (f_1 \cdot f_2 \cdot \adj{f_1} \cdot \adj{f_2} ) (x)   = (\tilde f \cdot \tilde f^\ast)(x) =| \tilde f (x) |^2. \qedhere \popQED
	  \end{equation*}
	   \renewcommand{\qedsymbol}{}
	\end{pf}\vspace{-1em}
  By \thmref{thm:PaleyW}, \lemref{lem:InfoModulus2} has immediate consequences for the phase retrieval problems considered in  \sref{SS:PhaseRetrAbstract}, as will be discussed in \sref{SS:PhaseRetrFF} and \sref{SS:PhaseRetrNF}.

According to the definition in \eqref{eq:DefOrderType}, adding a function $g$ with $\lambda_g < \lambda_f$ to $f$ may neither change its order nor its type. Likewise, it is clear that multiplication with $g$ cannot \emph{increase} any of these parameters. The following lemma shows that they may neither \emph{decrease} if $g$ is of at most exponential order and not identically zero:
	\vspace{1em}\begin{lem1}[Decay bounds for low order entire functions \cite{BoasEntire}]\label{thm:MaxDecayExpOrder}
		Let $f$ be an entire function of order $0 \leq \lambda_f \leq 1$ that is not identically zero and let $\varepsilon > 0$. Then
		\begin{equation*}
		 \limsup_{r \to \infty} m_f(r) M_f(r)^{1+\varepsilon} > 0
		\end{equation*}
		In particular, if $f \leq \Or(\exp(\tau_f|\xi|^{\lambda_f}))$, then $\limsup_{r \to \infty} m_f(r) \E^{(\tau_f+\varepsilon) r} = \infty$.
	\end{lem1}\vspace{1em}
The essential message of \lemref{thm:MaxDecayExpOrder} is that non-vanishing factors of at most exponential order may never weaken super-exponential growth.

 \end{subsection}

 \begin{subsection}{Far-Field Phase Retrieval} \label{SS:PhaseRetrFF}

 In the following, we analyze \probref{prob:FF}, i.e.\ the reconstruction of a function or tempered distribution $\psi$ from the squared modulus of its Fourier transform. The focus is on uniqueness of the solution for compactly supported contact images $\psi \in \A \subset  \Sdashc{m}$, which is investigated both for $m =1$ and in the higher dimensional case.
 
 \subsubsection{Phase Retrieval without Constraints}
 
 According to \cref{cor:FourierL2}, the Fourier transform maps $\Lp 2 {\mR^m}$ bijectively onto itself. Hence, it is clear that the solutions of \probref{prob:FF} for $ \A = \Lp 2 {\mR^m}$, i.e.\ for general square-integrable functions of not necessarily compact support, are highly non-unique \cite{Klibanov1995}:
 \vspace{1em}
 \begin{th1}[Non-Uniqueness for general $L^2$-functions] \label{thm:Nonunique-L2}
  For $\A = \Lp 2 {\mR^m}$, let $\psi$ be a solution of \probref{prob:FF}. Let $u: \mR^n \mapsto \mC$ be measurable with $|u| = 1$ and
  \begin{equation*}
   \psi_u := \cF^{-1}\left( u \cdot \cF(\psi)  \right)
  \end{equation*}
  Then $\psi_u $ solves \probref{prob:FF}, i.e.\ $\psi_u \in \A $ and $| \cF( \psi_u) |^2 =  | \cF( \psi  ) |^2 $.

  \end{th1}
\begin{pf}
 By construction, we have $|\cF(\psi_u)|^2 = |u|^2 |\cF(\psi)|^2 = |\cF(\psi)|^2$.
 According to \cref{cor:FourierL2}, this also implies $\psi_u \in \A$ since $\cF(\psi_u)$ is measurable and
   \begin{equation*}\pushQED{\qed} 
  \norm{\psi_u}_{\Lp 2 {\mR^m}}^2 = \int_{\mR^m} |\cF (\psi_u)|^2 \; \D x = \int_{\mR^m} |\cF (\psi)|^2 \; \D x = \norm{\psi}_{\Lp 2 {\mR^m}}^2 \stackrel{\psi \in \Lp 2 {\mR^m}} < \infty. \qedhere \popQED
 \end{equation*}
 \renewcommand{\qedsymbol}{}
\end{pf}
  \vspace{-1em}
  The result shows that \probref{prob:FF} has an uncountably large number of different solutions in $\Lp 2 {\mR^m}$ whenever the measured intensities are not identically zero.
  
   \subsubsection{Trivial Ambiguities} \label{SSS:TrivialAmbiguities}

    According to \thmref{thm:Nonunique-L2}, far-field phase retrieval of $L^2$-functions without further constraints is \emph{not feasible}. In the following, we therefore restrict to compactly supported signals, assuming $\psi \in \A \subset \Sdashc{m}$ as motivated in \sref{SS:PhaseRetrAbstract}.
    
  By \thmref{thm:PaleyW}, $\A$ is then mapped onto entire functions of at most exponential order (cf. \eqref{eq:DefOrder}). If $\psi \in \A$, we thus obtain $\psi_u \notin \A$ for the alternate solutions in \thmref{thm:Nonunique-L2} whenever $u \cdot \cF(\psi)$ is not entire or of super-exponential order. Nevertheless, there are certain choices of $u$ which retain these properties, i.e.\ for which the assignment $\psi \mapsto \psi_u$ preserves compactness of the support. In particular, this is true for $u\in \{s_{\beta_0}, t_{\bbeta_1}, r \}$ with $\beta_0 \in \mR$, $\bbeta_1 \in \mR^m$ defined by
\begin{subequations}
  \begin{align}
   s_{\beta_0}(\bxi) &:= \exp( \I \beta_0 ) \\
   t_{\bbeta_1}(\bxi) &:= \exp( \I \bbeta_1 \cdot \bxi ) \\
   r(\bxi) &:=   \cc{\cF( \psi )(\bxi)} / \cF( \psi )(\bxi)
  \end{align}
\end{subequations}
for all $\bxi \in \mR^m$. According to \thmref{thm:PropFT}, the  corresponding real-space transformations are given by:
  \vspace{1em}
\begin{itemize}
 \item[(a)] \emph{Scaling} by a unitary constant: $\psi \mapsto \psi_{s_{\beta_0}} = \exp( \I \beta_0 ) \psi$
     \item[(b)] \emph{Translation} by a constant shift $\bbeta_1$: $\psi \mapsto \psi_{t_{\bbeta_1}} = \psi(\cdot - \bbeta_1) $
  \item[(c)]  \emph{Reflection} in the origin and complex conjugation: $\psi \mapsto \psi_r = \cc{\psi(-\cdot)} $
\end{itemize}
  \vspace{1em}
Hence, when merely a compact support of the solution of \probref{prob:FF} is assumed, uniqueness may only hold up to these \emph{``trivial ambiguities''}. Moreover, note that many properties of $\psi$ are preserved under the transformations (a), (b) and (c). This implies that these ambiguities may not be easily overcome by imposing additional constraints such as positivity or a certain regularity in a Sobolev space sense (compare \sref{S:Sobolev}):
\vspace{1em}
\begin{rem1}\label{rem:TrivAmbiguities}
 If $\psi$ is of regularity $H^s(\mR^m)$ for $s \geq 0$, then $\psi_u \in H^s(\mR^m)$ for $u\in \{s_{\beta_0}, t_{\bbeta_1}, r \}$, i.e.\ the transformations \Text{(a), (b)} and \Text{(c)} preserve regularity. Moreover, if $\psi$ is real-valued or positive, so are $\psi_{t_{\bbeta_1}}$ and $\psi_r$.
\end{rem1}
\vspace{1em}
On the other hand, the translational symmetry (b) can be ruled out if the exact support of the solution $\psi$ (or more generally its convex hull) is known and imposed as a constraint. Then any shifted version of $\psi$ would violate this restriction. If the support is non-pointsymmetric in addition, then the \emph{twin-images} $\psi_{\Text{twin}, \bx_0}$ defined by reflections on some point $\bx_0 \in \mR^m$
\begin{equation}
 \psi_{\Text{twin}, \bx_0}(\bx) = \cc{\psi(\bx_0-\bx)} 
\end{equation}
would neither be compatible with the support constraint. Hence, the only remaining trivial ambiguity would be the scaling symmetry (a).

In practice, a priori information on the exact support of a specimen is rarely accessible. In order to overcome the translational symmetry in this case, iteratively updated support estimates have to be incorporated in reconstruction algorithms. An example is given by the \emph{Shrinkwrap Algorithm} \cite{Fienup1986algorithm} which has been successfully applied to reconstruct experimental far-field data \cite{Marchesini2003PhysRevB}.

\subsubsection{Holographic Constraints}

A different technique to rule out the above ambiguities is by \emph{perturbational-} or \emph{holographic-} approaches. If the objective function $\psi$ can be written as a sum of a known and - in a suitable sense - dominant part $b$ plus a perturbation $h$, the fixed support location of $b$ breaks the translational symmetry. This ansatz even has the potential to break the phase retrieval ambiguities altogether as illustrated by the following example:
\vspace{1em}
\begin{ex1}[Speckle Holography \text{\cite{Bates1973SpeckleHolo,Millane1990}}] \ \\ \label{ex:Speckle}
	     Let $\psi = (2 \pi)^{n/2} \delta_0 + h \in \Sdashc{m}$ be such that the distance of $\supp(h)$ to the Dirac delta $\delta_0$ in the origin is greater than the diameter of $\supp(h)$. Then $\psi$ can be reconstructed uniquely up to twin-image symmetry
	\end{ex1}
\begin{pf} Using $\cF(\delta_0) = (2 \pi)^{-n/2}$, we obtain
 	     \begin{equation*} |\cF(\psi)|^2 = 1 + \cF(h) + \cc{\cF(h)} + \cF(h) \cdot \cc{\cF(h)}  \end{equation*}
	     By the convolution theorem \eqref{eq:ConvThm}, this implies with $h_-(\bx) := \cc{h(-\bx)}$
	     \begin{equation} \cF^{-1}(|\cF(h)|^2 -1) = h + h_- + h \ast h_-.  \label{eq:Speckle} \end{equation}
	    Due to the assumed distance between $\supp(h)$ and the origin, the supports of $h, h_-$ and $ h \ast h_-$ do not overlap. Hence, all of these may be recovered from the data $|\cF(\psi)|^2$ using \eqref{eq:Speckle}. However, the contributions from $h$ and $h_-$ cannot be distinguished from one another, which gives rise to the remaining twin-image ambiguity.
\end{pf}
\vspace{1em}
\exref{ex:DiracDelta} demonstrates that - in spite of the above ambiguities - establishing uniqueness of the Fourier-data phase retrieval \probref{prob:FF} is not a hopeless endeavor. In fact even the remaining twin-image symmetry in the considered holographic setting may be broken if the Dirac delta in \exref{ex:DiracDelta} is replaced by a non-pointsymmetric reference function.

\subsubsection{Uniqueness Theory in 1D}

In the following, we focus on compactly supported signals in a single dimension $\psi \in \A \subset \Sdashc{}$. By the theory of \sref{SS:EntireFunctions}, this allows a complete characterization of the ambiguities in phase retrieval beyond the study of the ``trivial'' ones in the preceding paragraph. This characterization is due to \citet{Akutowicz1956I,Akutowicz1957II} and \citet{Walther1963}, following directly from \lemref{lem:InfoModulus2} and \thmref{thm:PaleyW}:

\vspace{1em}\begin{th1}[Phase Retrieval from 1D Fourier Data \cite{Akutowicz1956I,Akutowicz1957II,Walther1963}]\label{thm:PhaseAmbiguities1D}
		Let $ \psi \neq 0$ solve \probref{prob:FF} for $\A =  \Sdashc{}$. Let  $\{a_j\}_{j \in J} \subset \mC \setminus \{ 0 \}$, $J \subset \mN$ denote the complex zeros of $\cF(\psi)$ counted by their multiplicity. Then $\tilde \psi \in \A$ solves \probref{prob:FF} if and only if
		\begin{align}
 \cF(\tilde \psi)(\xi )  = \left( \exp(-\I (\beta_0 + \beta_1 \xi ) )  \prod_{j \in K} \frac{ 1-  \xi / \cc{a_j} }{ 1-  \xi / a_j }  \right) \cF(\psi)(\xi)  \label{eq:PhaseSolutions}
\end{align}
for some $J \subset K$, $\beta_0,\beta_1 \in \mR$.
	\end{th1}\vspace{1em}
	
In particular, \thmref{thm:PhaseAmbiguities1D} includes the trivial ambiguities associated with scaling, translational and twin-image symmetries, represented by the exponential prefactor and the case $K = J$, respectively. Beyond this, however, many more alternate solutions may be constructed by taking the product in \eqref{eq:PhaseSolutions} only over an arbitrary subset $K \subset J$. Since the multiplication of $\cF(h_1)$ with the factor
\begin{align}
 \prod_{j \in K} \frac{ 1-  \xi / \cc{a_j} }{ 1-  \xi / a_j }     \label{eq:ZeroFlippingFactor}
\end{align}
replaces the zeros  $\{a_j\}_{j \in K}$ by their complex conjugates $\{\cc{a}_j\}_{j \in K}$, i.e.\ reflects them on the real axis from one complex half plane into the other, this transformation is called ``zero-flipping''. Note that it is \emph{without effect} for any zeros on the real line. On the other hand, if $\cF( \psi)$ has non-real zeros of which only a part is flipped by the product in \eqref{eq:PhaseSolutions}, then the constructed alternate solution $\tilde \psi$ in \thmref{thm:PhaseAmbiguities1D} will in general not be related to $\psi$ by simple geometrical transformations. Both cases occur naturally as illustrated by the following example:
\vspace{1em}
\begin{ex1} \label{ex:BumpFunctions}
Define the bump functions $b_{\exp}, b_{\rect}: \mR \to \mR$ by
 	\begin{equation}
 	b_{\rect} (x) := \begin{cases}
	                  (2 \pi )^{\frac 1 2} &\text{for }x \in [- 1 ; 1 ] \\
	                  0       &\text{else}
	                 \end{cases}, \;\;\;\;
 b_{\exp} (x) := \begin{cases}
	                  (2 \pi )^{\frac 1 2} \exp(x) &\text{for }x \in [0;1] \\
	                  0       &\text{else}
	                 \end{cases}
	  \label{eq:DefBaseFunc}
	\end{equation}	
Then the solution to \probref{prob:FF} is unique for $I = |\cF(b_{\rect})|^2$, $\A = \Sdashc{}$ up to translation and a global factor of modulus 1.
For $I = |\cF(b_{\exp})|^2$, the phase retrieval problem has an infinite number of real-valued solutions $\tilde b \in \Sdashc{}$ with support in $[0;1]$ that are of the same regularity as $b_{\exp}$.
\vspace{1em}
\end{ex1}
\begin{pf}
 The Fourier transforms of $b_{\rect}, b_{\exp}$ are given by
	\begin{equation}
	 \cF ( b _{\rect} ) (\xi) :=  \frac{2\sin(\xi)}{\xi}  \MTEXT{and}  \cF( b _{\exp} ) (\xi)  =  \frac{\I (\E^{1 - \I \xi} - 1)}{\xi - \I}. \label{eq:FourBaseFunc}
	\end{equation}
Accordingly, the zeros of $\cF ( b _{\rect} )$ are exactly given by $\pi \mZ \subset \mR$. Hence, the zero-flipping factor in \eqref{eq:ZeroFlippingFactor} is always 1, ruling out all non-trivial ambiguities. The product in \thmref{thm:PhaseAmbiguities1D} reduces to the exponential prefactor representing the remaining symmetry transformations in the claim.

On the other hand, the complex roots of $\cF( b _{\exp} )$ are obtained as
	\begin{equation}
	 \E^{1 - \I \xi} = 1 \;\; \Leftrightarrow \;\; 1 - \I \xi \in 2\pi \I \mZ \;\; \Leftrightarrow \;\; \xi \in 2\pi \mZ - \I  \label{eq:ZerosBaseFunc}
	\end{equation}
	We see that all zeros $\seqn{a}{j}$ have an imaginary part $\Im(a_n) = -1$, i.e.\ lie in the lower complex half-plane $\mH_-$. By flipping any finite subset $\{a_j\}_{j\in J}$ of these into $\mH_+$, an alternate solution $\tilde b \in \Sdashc{}$ to \probref{prob:FF}  is obtained, given by
	\begin{equation}
	 \cF(\tilde b)(\xi)  = \left( \prod_{j \in K} \frac{ 1-  \xi / \cc{a_j} }{ 1-  \xi / a_j } \right) \cF( b _{\exp} ) (\xi) \MTEXT{for all} \xi \in \mR. \label{eq:ZeroFlipbexp}
	\end{equation}
	As regularity in a Sobolev space sense (see \sref{S:Sobolev}) is determined by the decay behavior of the Fourier transform and $| \cF(\tilde b)|^2 = |  \cF( b _{\exp})|^2$, we have $\tilde b \in H^s(\mR)$ if and only if $b _{\exp} \in H^s(\mR)$ for all $s \geq 0$. Moreover, the prefactor in \eqref{eq:ZeroFlipbexp} is of algebraic growth for $K$ finite. Since the convex hull of the support is determined by the \emph{exponential} growth behavior of the Fourier transform according to the estimate \eqref{eq:PaleyWBound} in \thmref{thm:PaleyW}, $\tilde b$ is supported in $[0;1]$. Moreover, whenever all zeros are flipped or retained in pairs $(a_j, -\cc{a}_j)$, then the symmetry
	\begin{equation}
	\cc{ \cF(\tilde b)(-\xi) } = \cF(\tilde b)(\xi) \MTEXT{for all} \xi \in \mR 
	\end{equation}
	is preserved. Consequently, $\tilde b$ is real-valued in this case according to \eqref{eq:FTReflect}.
\end{pf}
\vspace{1em}
    
Note that almost all entire functions have infinitely many non-real zeros. In this respect, $b_{\rect}$ represents the special case in \exref{ex:BumpFunctions} whereas the construction of alternate solutions applied to $b_{\exp}$ is applicable to almost all choices $\psi \in \Sdashc{}$. Moreover, it should be emphasized that this partial ``zero-flipping'' indeed yields signals that may differ significantly in shape from the original one. This is demonstrated in \figref{fig:ZeroFlip} for a discrete version of $ b _{\exp}$. Hence, we may state:
  \vspace{1em}
  \begin{res}[Phase Retrieval Ambiguities for 1D Fourier Data] \label{res:1DFTPhaseRetrIllposed}
   Phase retrieval from Fourier data, given by \probref{prob:FF}, is highly non-unique for objects $\psi \in \A = \Sdashc{}$. These ambiguities may \emph{not} be overcome by imposing real-valuedness, support in a particular interval or a certain regularity.
  \end{res}
  \vspace{1em}
  
  On the other hand, the characterization in \thmref{thm:PhaseAmbiguities1D} may be exploited to derive sufficient conditions which allow for a unique reconstruction of the object. The idea is to restrict the set of admissible solutions $\A \ni \psi$ such that all complex zeros of $\cF(\psi)$ may be located unambiguously in either of the complex half planes, for instance by ensuring that $\cF(\psi)$ is nonzero in $\mH_\pm$ or that zeros necessarily occur in complex-conjugate pairs. This rules out the zero-flipping ambiguity so that uniqueness up to trivial transformations is achieved. Some examples of such sufficient criteria are given by \cite{Klibanov1995}:
  \vspace{1em}
	\begin{itemize}
	 \item[$\Text{(a)}$] $\psi \in \Sdashc{}$ is real-valued and symmetric or antisymmetric w.r.t some $a \in \mR$
	 \item[$\Text{(b)}$] $\psi \in \Lp 1 {\mR}$ is non-negative and non-increasing on $[a;b] \supset \supp(\psi)$
	 \item[$\Text{(c)}$] For some $[a;b] \subset \mR$, $\seqn{c}{j} \subset \mC$, $a = x_0 < x_1 < \ldots \leq b$, $h \in L^1([a;b])$ and $n \in \mN_0$, the $n$-th derivative of $\psi \in \Sdashc{}$ is of the form
	 \begin{equation}
	  \psi^{(n)}  = \sum_{j \in \mN} c_j \delta_{x_j} + h   \MTEXT{with} |c_1| \geq \norm{h}_1 + \sum_{j \in \mN \setminus \{1 \} } |c_j|. \label{eq:Crit3}
	 \end{equation}
 	\end{itemize}
	\vspace{1em}
 The constant bump $b_{\rect}$ in \exref{ex:BumpFunctions} is uniquely determined according to $\Text{(a)}$. On the other hand, if follows from $\Text{(b)}$ that none of the zero-flipped versions of $b_{\exp}$ may be non-negative and non-decreasing as is confirmed by \figref{fig:ZeroFlip}. In $\Text{(c)} $, uniqueness is obtained from the knowledge of the dominant singularity of $\psi$  given by the Dirac delta $\delta_{x_0}$ at the lower boundary of the support. A generalized variant of this criterion has been applied to establish uniqueness in X-ray reflectivity measurements \cite{Hohage2008XrayRefl} - another experimental setup which gives rise to \probref{prob:FF}. However, it should be noted that for $n \geq 1$, i.e.\ in the case of prescribed singularities in a derivative and not in the function itself, $\psi^{(n-1)}(x)$ must decay to 0 for $x \to \infty$ owing to the compactness of the support. By \eqref{eq:Crit3}, this implies that $c_1^{-1} \psi^{(n-1)}$ must be monotonically decreasing in $(a; \infty)$. Hence, all criteria $\Text{(a)}$, $\Text{(b)}$, and $\Text{(c)}$ make strong structural assumptions on the signal $\psi$.
 \begin{figure}
     \centering
    \subfloat{\includegraphics[width=0.32\textwidth]{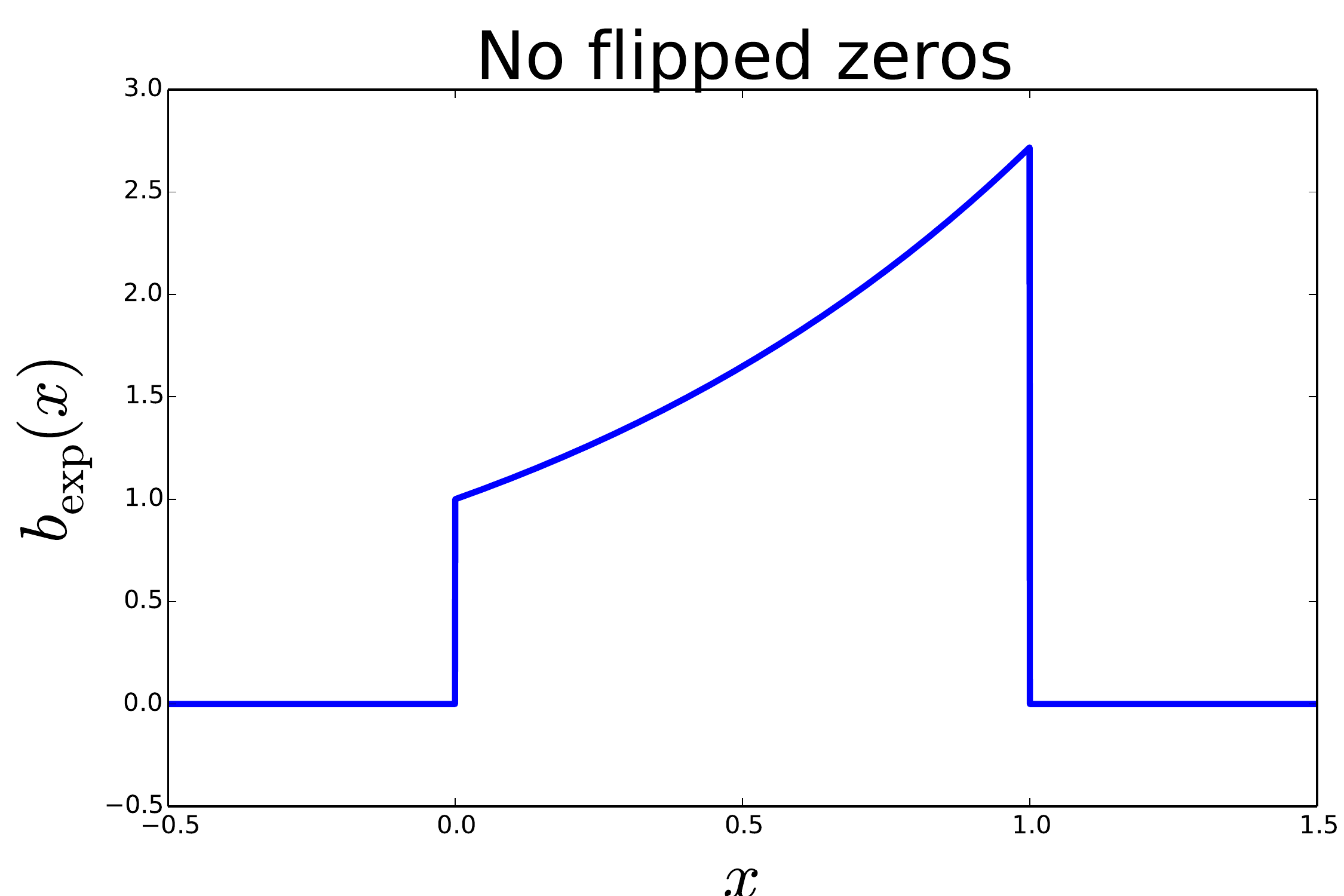}}
     \hfill
    \subfloat{\includegraphics[width=0.32\textwidth]{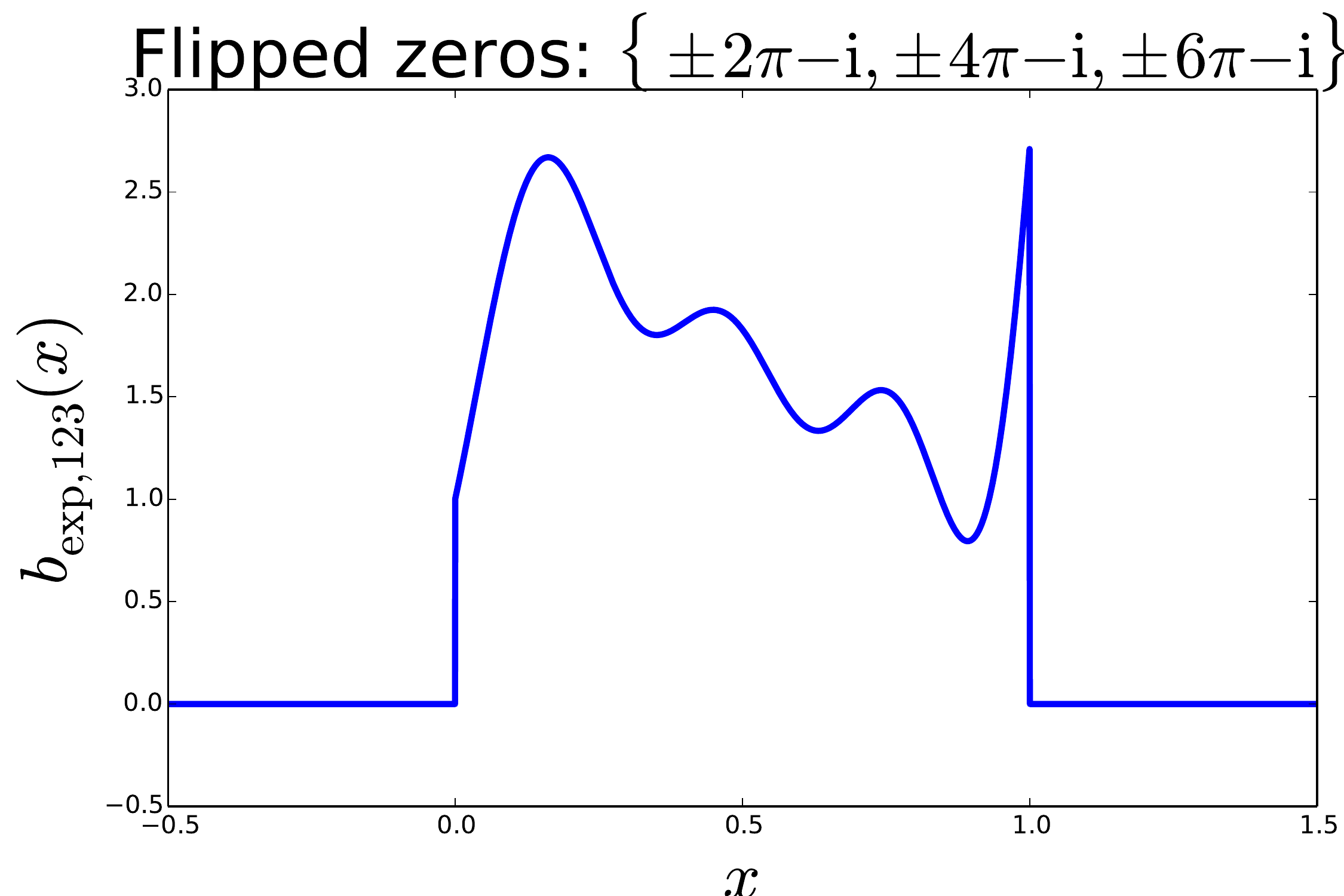}}\vspace{-1em}
    \\
    \subfloat{\includegraphics[width=0.32\textwidth]{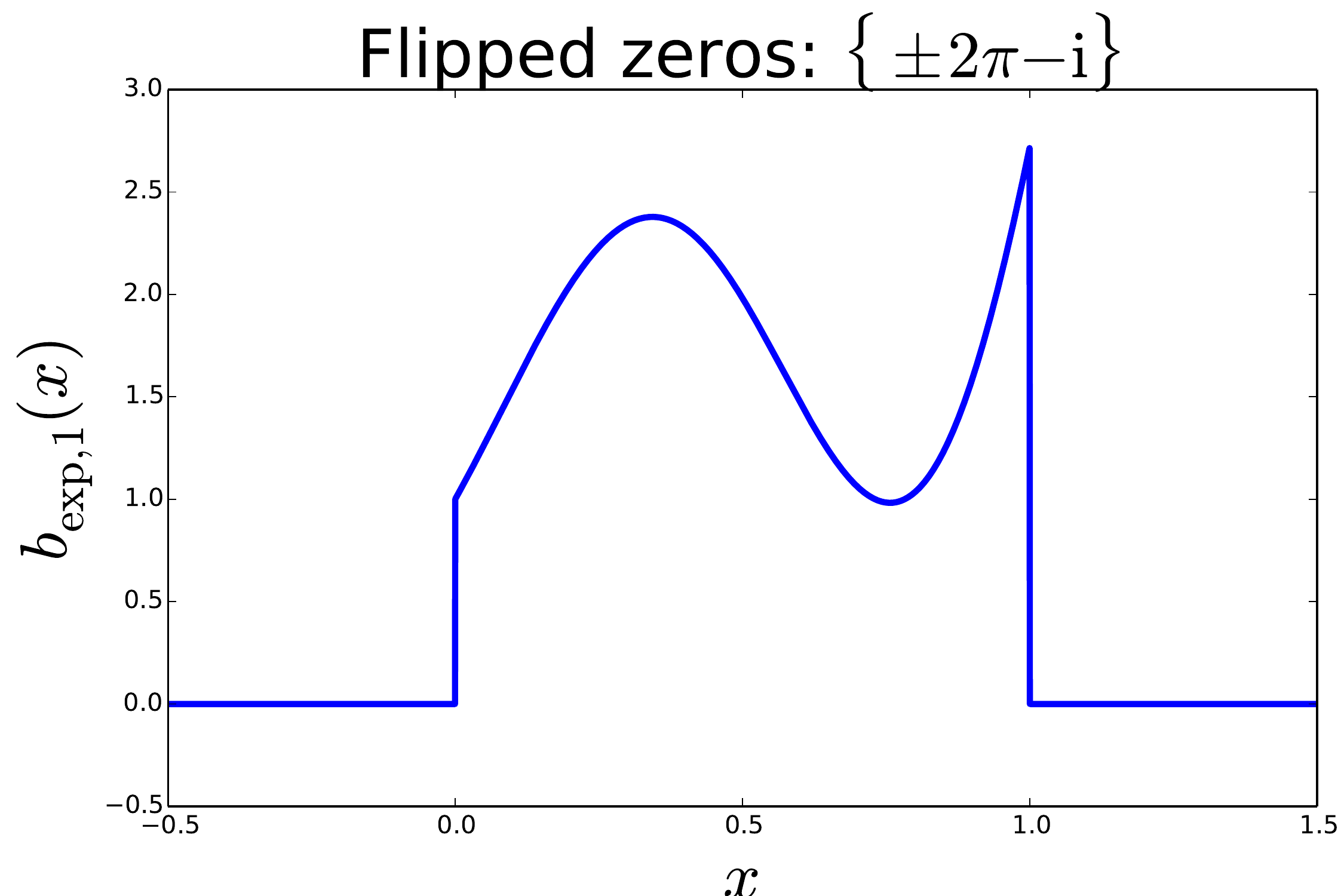}}
    \hfill
    \subfloat{\includegraphics[width=0.32\textwidth]{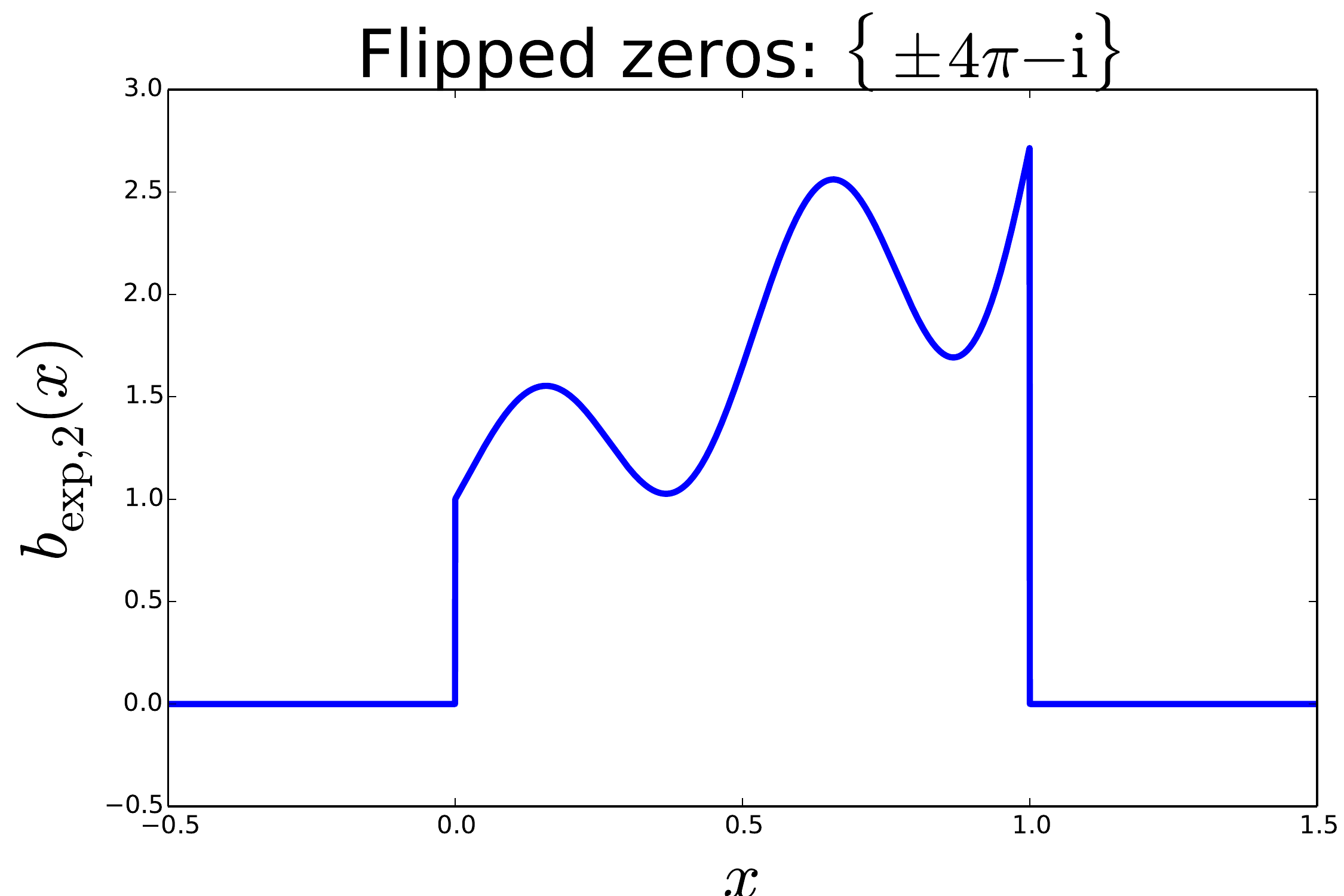}}
    \hfill
    \subfloat{\includegraphics[width=0.32\textwidth]{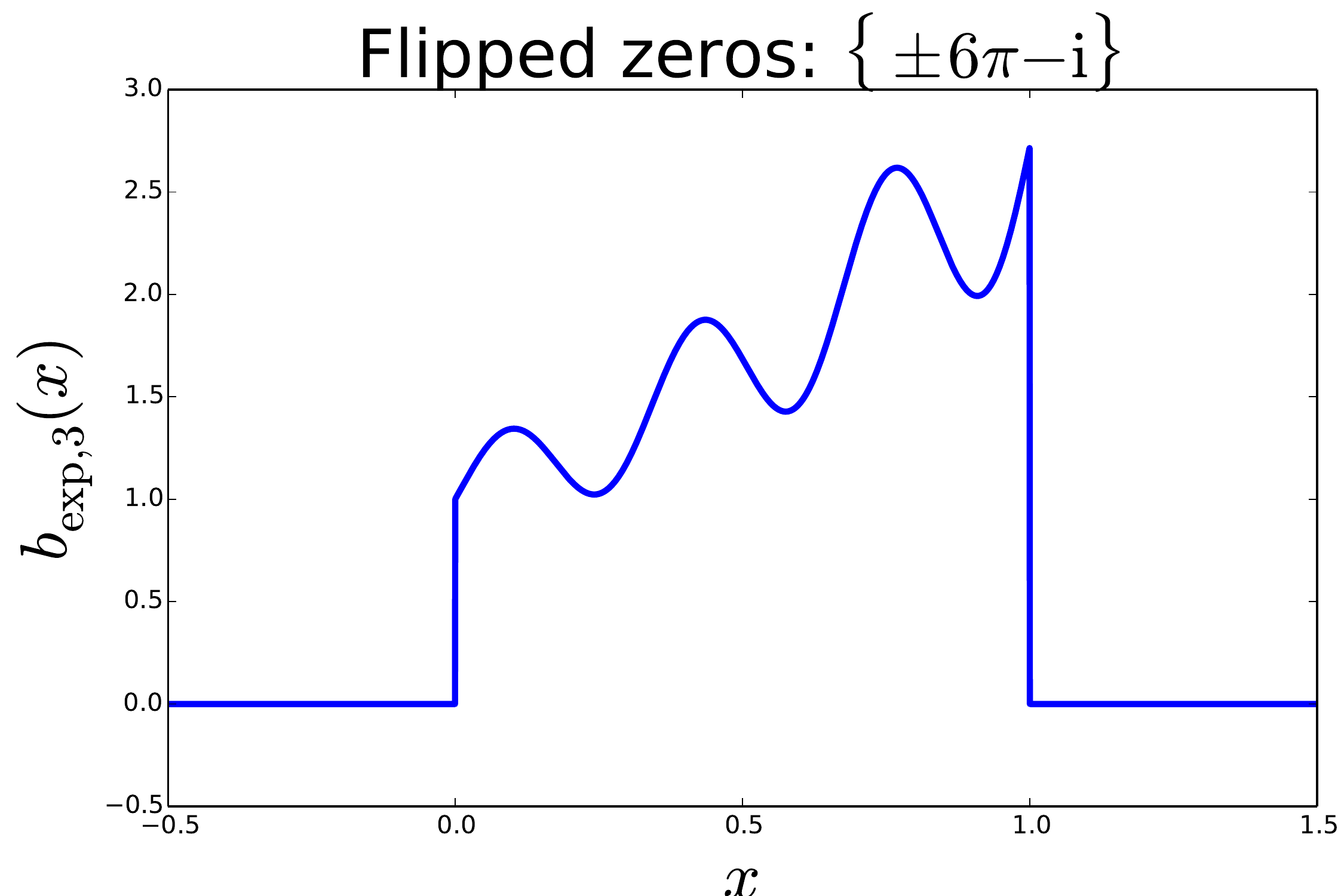}}\vspace{-1em}
    \\
    \subfloat{\includegraphics[width=0.32\textwidth]{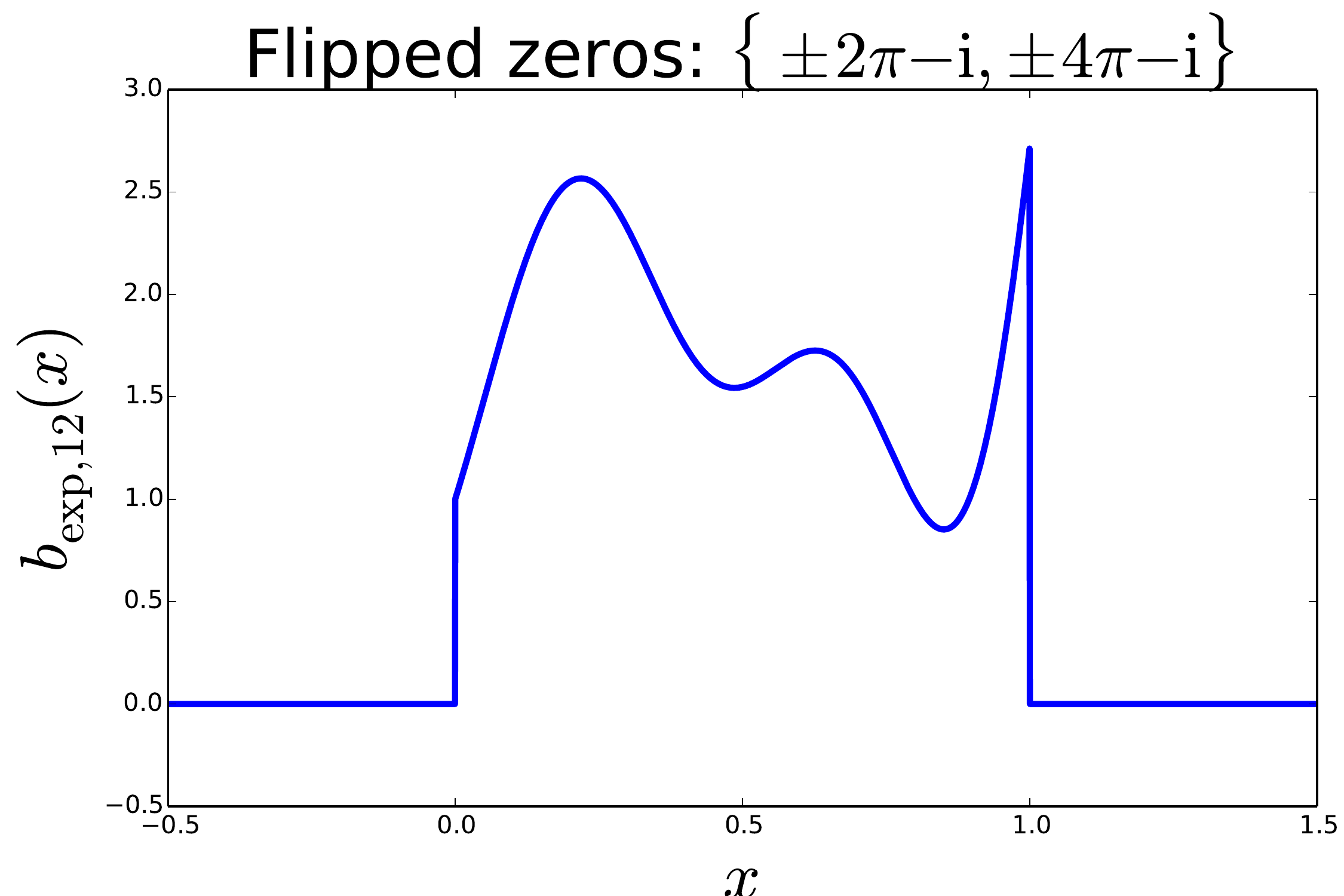}}
    \hfill
    \subfloat{\includegraphics[width=0.32\textwidth]{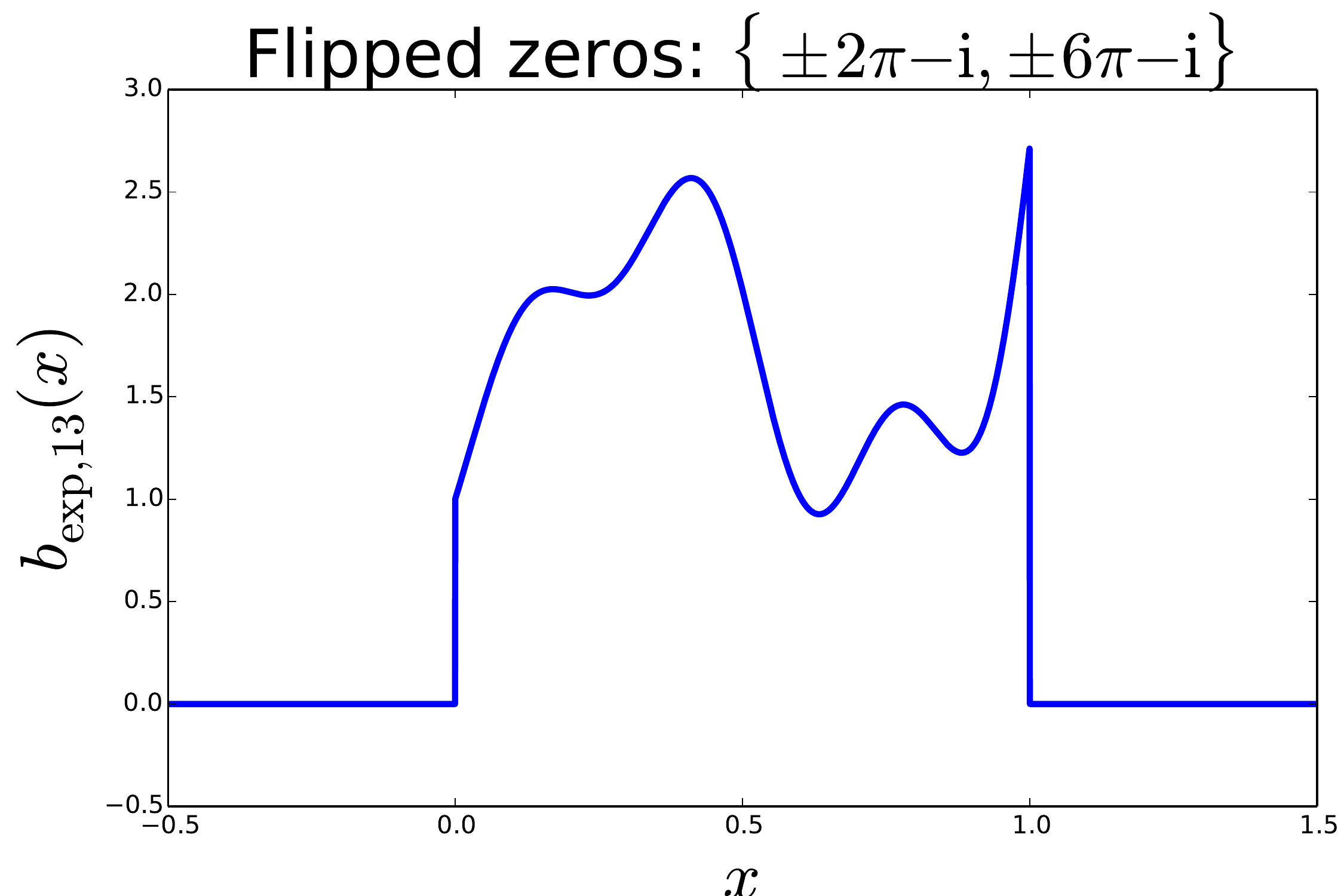}}
    \hfill
    \subfloat{\includegraphics[width=0.32\textwidth]{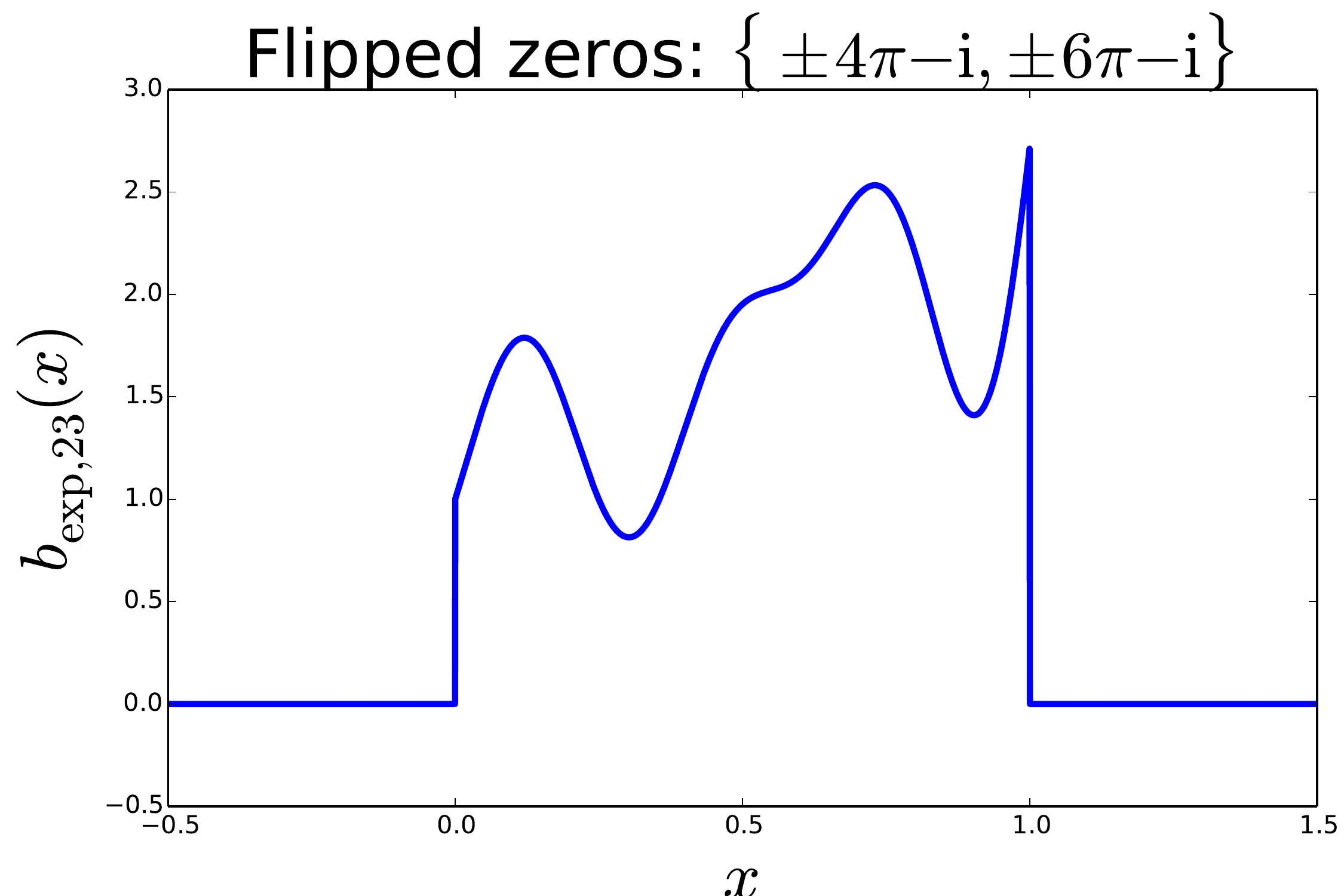}}
    \caption{Simulation of the zero-flipping construction in \exref{ex:BumpFunctions}. The discrete Fourier transform of $b_{\exp}$ (top-left) is multiplied with different combinations of the zero-flipping factors $\frac{ (1-  \xi / \cc{a_j})(1 +  \xi /  a_j) }{ (1-  \xi / a_j)(1 +  \xi /  \cc{a_j}) }$ for $a_j = 2 \pi j - \I, j \in \mN$, yielding alternate solutions to the phase retrieval problem. All of these are real-valued (and positive) and have the same regularity and support as $b_{\exp}$. \label{fig:ZeroFlip}}
    \end{figure}
    \captionsetup[subfigure]{justification=justified,singlelinecheck=false,position=bottom}

 \subsubsection{Uniqueness Result in 1D}
 
 	In the sequel, we derive a new criterion, following the \emph{holographic} approach of establishing uniqueness in phase retrieval by superimposing some known profile upon the unknown object to be reconstructed. For the reference signal, we choose the exponential ramp $b_{\exp}$ defined in \eqref{eq:DefBaseFunc}. We will consider perturbations $h = \psi - b_{\exp} \in W^{1,1}(\mR)$ that are small in Sobolev norm (see \defref{def:SobSpace} in \sref{S:Sobolev}) and completely contained within the support of $b_{\exp}$. The latter enforces $\supp(\psi) = [0;1]$ and thus rules out translational variance. Moreover, the asymmetry of $b_{\exp}$ will break the reflectional twin-image symmetry. The manifestation of this asymmetry in its Fourier transform \eqref{eq:FourBaseFunc}, gives rise to the estimate
	\begin{equation}
	 (1+ |\xi|) |\cF( b _{\exp} ) (\xi)| \geq \frac{1+ |\xi|}{|\xi - \I|} \left| \exp(1+ \Im(\xi) - \I \Re(\xi)) -1 \right| \stackrel{\Im(\xi) \geq 0}{\geq} \exp( \Im(\xi) ) (\E -1) \label{eq:ExpBaseEstimate}
	\end{equation}
	for all $\xi \in \mH_+ = \{ z \in \mC: \Im(z) \geq 0\}$. Note that no analogue holds for the rectangular bump $b_{\rect}$ due to its real zeros. On the other hand, \eqref{eq:ExpBaseEstimate} implies that $\cF( b _{\exp} + h )$ may not have zeros in $\mH_+$ for suitably ``small'' perturbations $h \in \Sdashc{}$. This observation leads to the following uniqueness criterion:
	\vspace{1em}
	\begin{th1}[Uniqueness Criterion for Phase Retrieval from 1D Fourier Data]\label{thm:UniquePhaseFF1D}
 		Let $h \in W^{1,1}(\mR)$ s.t.\ $\supp(h) \subset [0;1]$ and $\norm{h}_{W^{1,1}(\mR)} < \E -1$. Set $\psi:= b _{\exp} + h$. Then $h$ is uniquely determined by $|\cF(\psi)|_{|U}^2$ on an arbitrary open set $U \subset \mR$.
	\end{th1}
	\begin{pf}
	 $\psi$ has compact support, so that $|\cF(\psi)|^2$ is an entire function and thus uniquely determined in $\mC$ by its values in $U \subset \mR$ as argued in the proof of \lemref{lem:InfoModulus2}. Moreover, for all $\xi \in \UHP$ we have the estimate
	\begin{align}
	 (1+ |\xi|) |\cF(h) (\xi)| &\;\;\;\;\;= \;\;\;\;\;\left|\int_{\mR} \exp(-\I \xi x) h \; \D x \right| + \left| \int_{\mR} (-\I \xi) \exp(-\I \xi x) h \; \D x \right| \nonumber \\
	 &\stackrel{h, h' \in \Lp 1 {\mR}}= \left|\int_{\mR} \exp(-\I \xi x) h \; \D x \right| + \left| \int_{\mR} \exp(-\I \xi x) h' \; \D x \right| \nonumber \\
	 &\;\;\;\;\; \leq  \;\;\;\;\; \int_{\mR} \underbrace{\exp(\Im(\xi)x )}_{\leq \exp(\Im(\xi)) \, \forall \, x \in \supp(h) } (|h| + |h'|) \; \D x  \nonumber \\
	 &\;\;\;\;\;\leq \;\;\;\;\; \exp(\Im(\xi)) \int_{\mR}  (|h| + |h'|) \; \D x = \exp(\Im(\xi)) \norm{h}_{W^{1,1}(\mR)} \nonumber \\
	 &\;\;\;\;\;< \;\;\;\;\; \exp( \Im(\xi) ) (\E -1) \label{eq:1DFTPhaseRetrEstimate}.
	\end{align}
	By inequality \eqref{eq:ExpBaseEstimate}, this implies
	\begin{equation*} |\cF(\psi)(\xi)| \geq |\cF(b _{\exp}) (\xi)| - |\cF(h)(\xi)| > 0 \MTEXT{for all} \xi \in \UHP, \end{equation*}
	i.e.\ $\cF(\psi)$ has no zeros in the upper half plane. Hence, $\cF(\psi) $ is uniquely determined up to an exponential factor $\E^{\I(\beta_0 + \beta_1 \xi)}, \beta_i \in \mR$. However, by the Sobolev embedding theorem \thmref{thm:SobEmbed}, $h$ is continuous. In particular, this yields
	\begin{equation*} h(0) = \lim_{x \to 0} h(x) =0 \MTEXT{and} h(1) = \lim_{x \to 0} h(x) =0, \end{equation*}
	so that $\psi(0) = 1$ and $\psi(1) = \E$. This fixes the value of the multiplicative constant $\E^{\I \beta_0}$. The remaining linear exponent exponential factor $\E^{\I \beta_1 \xi }$ in Fourier space corresponds to a translation in real space and is thus uniquely determined by the support constraint $\supp(\psi) \subset [0;1]$ due to the nonzero boundary values $\psi(0)$ and $\psi(1)$.
	\end{pf}
	\vspace{1em}
	
	Note that, apart from a certain degree of regularity, \thmref{thm:UniquePhaseFF1D} does not make any structural assumptions on the perturbation $h$ - as opposed to the criteria $\Text{(a)}$, $\Text{(b)}$, and $\Text{(c)}$ discussed above. All that is needed is smallness of $h$ compared to $b_{\exp}$, which unfortunately has to be measured in the $W^{1,1}$-norm rather than in some $L^p$-sense, being apparently violated by the distortions in \figref{fig:ZeroFlip}. As a benefit, however, $h$ does not have to satisfy any form of monotonicity, symmetry or real-valuedness.
	In fact, we can interpret the result in the sense that we may uniquely reconstruct any \emph{complex-valued} compactly supported function $h \in W^{1,1}(\mR)$ by superimposing a suitably scaled version of the exponential ramp $b _{\exp}$ as a reference signal.

\begin{subsubsection}{Uniqueness in Higher Dimensions}

	Despite the mathematical beauty of the one-dimensional theory outlined in the preceding sections, the physically relevant setting for this work is phase reconstruction of \emph{two}-dimensional contact images $\psi$ from the corresponding far-field intensities. Moreover, it is seen from \eqref{eq:WeakFwOpFF} that phase contrast tomography of weakly scattering objects may even be interpreted as a \emph{three}-dimensional phase retrieval problem. This motivates the study of \probref{prob:FF} in higher dimensions $m \geq 2$.
	
	Accordingly, let $\psi \in \Sdashc{m}$ arbitrary. A straightforward approach to generalize the 1D theory to this case lies Fourier-transforming $\psi$ in all but the first coordinate:
	\begin{align}
	 \psi_{\bxi_y}(x) := \cF_{\overline 2} ( \psi)(\cdot, \bxi_y) \MTEXT{so that} \cF(\psi)(\xi_x, \bxi_y ) = \cF(  \psi_{\bxi_y} ) ( \xi_x) \label{eq:1DReductionFourier}
	\end{align}
	 for all $(\xi_x, \bxi_y) \in \mR \times \mR^{m-1} = \mR^m$, $x \in \mR$. Then it follows from \eqref{eq:DefFTSdash} and the definition of the support of a distribution (see \sref{S:SchwartzLp}) that $\psi_{\bxi_y}$ is compactly supported for all $\bxi_y \in \mC^{m-1}$. Moreover, the family $\{ \psi_{\bxi_y} \}_{\bxi_y \in \mR^{m-1}} \subset \Sdashc{}$ is composed of solutions to the 1D phase retrieval problems 
        \begin{align}
	 |\cF ( \psi_{\bxi_y} )|^2 = |\cF(\psi)(\cdot , \bxi_y )|^2, \label{eq:1DReductionPhaseProb}
	\end{align}
	which are exactly of the form considered above.
	In particular, alternate solutions $\{ \tilde \psi _{\bxi_y} \}_{\bxi_y \in \mR^{m-1}}\subset \Sdashc{}$ can be constructed by the zero-flipping construction outlined in \thmref{thm:PhaseAmbiguities1D}. However, by this construction,
	\begin{equation*}
	 \check \psi: (\xi_x, \bxi_y) \mapsto  \cF( \tilde \psi _{\bxi_y} )( \xi_x)
	\end{equation*}
	will in general \emph{not} define an entire function in $\mC^m$ and will thus not yield an alternate solution $\tilde \psi := \cF^{-1} ( \check \psi )$ to \probref{prob:FF} for $\A \subset \Sdashc{m}$.
	
	The geometrical reason is that the isolated complex zeros of 1D entire functions are replaced by smooth manifolds in higher dimensions, so-called \emph{zero-sheets}, which need to be ``flipped'' as a whole in order preserve their smoothness.
	Notably, the seemingly isolated zeros of $\cF(\psi _{\bxi_y})$ all belong to one and the same zero-sheet of $\cF( \psi )$ for \emph{almost any} $\psi \in \Sdashc{m}$ \cite{Millane1990}. Algebraically, this is due to the (ir-)reducibility of the corresponding entire functions in Fourier space as outlined by \citet{BruckSodin1979PhaseAmbiguity2D}: by \thmref{thm:Hadamard}, Hadamard factorizations in 1D always reduce to an infinite product of the \emph{primary factors} $E_{p_f} (\xi/a_j)$, i.e.\ monomials scaled with an exponential. This gives rise to the non-trivial ambiguities in
	\thmref{thm:PhaseAmbiguities1D}. On the other hand, almost all  polynomials in $\mC^m$ for $m \geq 2$ are irreducible \cite{Hayes1982Reducible}, i.e.\ cannot be decomposed into polynomials of a smaller degree. Consequently, if the infinite products in the Hadamard factorizations of the $h_{\bxi_y}$ are regarded as entire functions of all variables $(\xi_x, \bxi_y)$, these usually no longer factorize in a non-trivial manner \cite[sec. 10.3]{Hurt2001PhaseRetr}. \citet{Barakat1984}, however, showed that the existence of such factorizations is necessary for non-trivial phase retrieval ambiguities also in $m \geq 2$ dimensions. This leads to the following startling conclusion:
	\vspace{1em}
	\begin{res}[Phase Retrieval Ambiguities for $m$-dimensional Fourier Data \text{\cite{Barakat1984,Millane1990}}] \label{res:NDFTPhaseRetrAmbiguities}
	In $m > 1$ dimensions, phase retrieval of compactly supported signals $\psi \in \A \subset \Sdashc{m}$ from Fourier intensities (\probref{prob:FF}) is almost always unique up to trivial ambiguities. In other words, ``multiplicity of solutions is pathologically rare''.
	\end{res}
	\vspace{1em}

	\citet{Fienup19782DPhaseRetrFeasible} was the first to observe this ``dramatic'' reduction of non-uniqueness in numerical phase reconstructions of two-dimensional images. Nevertheless, note that arbitrarily ambiguous solutions may also occur in $\mR^m$. For instance, the $m$-dimensional exponential ramp, defined by
	\begin{equation}
	 b_{\exp, m} (x_1, x_2, \ldots , x_m ) :=  b_{\exp }(x_1) \cdot  b_{\rect}(x_2) \cdot \ldots  \cdot  b_{\rect}(x_n)  \label{eq:DefBaseMD}
	\end{equation}
	for all $ (x_1, x_2, \ldots , x_m ) \in \mR^m$, gives rise to the same variety of non-trivial ambiguities as its 1D analogue in \exref{ex:BumpFunctions}. Admittedly, this is due to its simple product structure. Yet, the example certainly demonstrates the necessity of \emph{deterministic} uniqueness criteria in higher dimensions in order to avoid excessive ill-posedness of \probref{prob:FF} in the vicinity of signals which cannot be reconstructed uniquely.
	
	In the following, we therefore derive a multidimensional analogue of \thmref{thm:UniquePhaseFF1D} using $b_{\exp, m}$ as a reference signal. The principal idea for this generalization lies in the reduction to the 1D case given by \eqref{eq:1DReductionFourier}: by entire analyticity of $\bxi_y \mapsto \psi_{\bxi_y}$, \probref{prob:FF} admits a unique solution $\psi \in \A \subset \Sdashc{m}$ if \eqref{eq:1DReductionPhaseProb} uniquely determines $\psi_{\bxi_y} \in \Sdashc{}$ for all $\bxi_y \in V$ in some open set $V \subset \mR^{m-1}$. In particular, it is sufficient if $h_{\bxi_y}$ is unique in the limit $\norm{\bxi_y}_2 \to \infty$. According to this argument, the \emph{absolute} smallness of the perturbation $h$ assumed in \thmref{thm:UniquePhaseFF1D} can be relaxed to an \emph{asymptotic} smallness of $h_{\bxi_y}$ for suitably large $\bxi_y$, obtained by imposing a certain Sobolev regularity of $h$ (cf. \sref{S:Sobolev}):
 
	\vspace{1em}
	\begin{th1}[Uniqueness Criterion for $m$-dimensional Phase Retrieval]\label{thm:UniquePhaseFFMD}
		For $m \geq 2$, let $h \in H^{\frac 3 2}(\mR^m)$ with $\supp(h) \subset \Omega := [0;1] \times [-1;1]^{m-1}$ and let $\psi := b_{\exp, m} + h$. Then $h$ is uniquely determined by $|\cF(\psi)|^2_{|U}$ on an arbitrary open set $U \subset \mR^m$.
	\end{th1}
	\begin{pf}
	 $\psi$ has compact support, so that $|\cF(\psi)|^2$ is an entire function and thus uniquely determined in $\mC^m$ by its values in $U \subset \mR^m$. For the remainder of the proof given here, we assume $m=2$. The general statement is a technical corollary following from the dimension-reduction in \eqref{eq:1DReductionFourier}, proven in \aref{A:MDFourPhaseRetr}. 
	 
	  For convenience, we define
	 \begin{equation*}
	  S_{> \frac 1 2} := \left\lbrace \xi \in \mR : |\sin(\xi)| > \frac 1 2\right\rbrace = \bigcup_{k \in \mZ} \left(\left(k + \frac 1 6\right) \pi; \left(k + \frac 5 6\right) \pi \right).
	 \end{equation*}
	  Then \eqref{eq:FourBaseFunc} and the estimate \eqref{eq:ExpBaseEstimate} imply for all $\xi_x \in \UHP$, $\xi_y \in  S_{> \frac 1 2}$
	 \begin{align}
	  |\cF(\psi)(\xi_x,\xi_y)| &\geq |\cF( b _{\exp} )(\xi_x) \cF( b _{\rect} )(\xi_y)| - |\cF(h)(\xi_x,\xi_y)| \nonumber \\
	  &> \frac{\exp( \Im(\xi_x) ) (\E -1)}{(1+|\xi_x|) |\xi_y|} - |\cF(h) (\xi_x, \xi_y)|. \label{eq:Unique2DEst0}
	 \end{align}
	 Setting $ h_{\xi_y}(x) := \cF(h(x,\cdot))(\xi_y)$, we have $h_{\xi_y} \in H^s(\mR) \subset W^{1,2}(\mR)$ for any $\xi_y \in \mC$ by \lemref{lem:HsDecay}. Since $h$ is compactly supported, this implies $h_{\xi_y} \in W^{1,1}(\mR)$ according to \thmref{thm:LpBoundedEmbed}. Applying the estimate \eqref{eq:1DFTPhaseRetrEstimate}, we thus obtain for all $\xi_x \in \UHP, \xi_y \in \mR$
	\begin{align}
	\E^{-\Im(\xi_x)}   &(1+ |\xi_x|) \cdot |\cF(h) (\xi_x, \xi_y)| \leq     \norm{ h_{\xi_y}}_{W^{1,1}(\mR)} \nonumber \\
	 =  \;  &\ip{|h_{\xi_y} |}{1_{[0;1]}}_{L^2(\mR)} + \ip{ |\partial_x h_{\xi_y} |}{1_{[0;1]}}_{L^2(\mR)} \leq   \norm{h_{\xi_y}}_{L^2(\mR)} + \norm{\partial_x h_{\xi_y} }_{L^2(\mR)}  \nonumber \\
	  \leq \;  &\frac 1 2   \norm{h_{\xi_y}}_{W^{1,2}(\mR)} \leq C \norm{h_{\xi_y}}_{H^{1}(\mR)}. \label{eq:Unique2DEst1}
	\end{align}
	Here, we have applied Cauchy--Schwarz inequality to the inner product of $h_{\xi_y}$ with the indicator function $1_{[0;1]}$ of its support. The constant $C > 0$ results from the equivalence of the $W^{1,2}(\mR)$- and $H^1(\mR)$-norms stated in \thmref{thm:SobFourier}. According to \lemref{lem:HsDecay}, the map $\xi_y \mapsto \norm{h_{\xi_y}}_{H^{1}(\mR)}$ is continuous. Furthermore, we have the bound
	\begin{align}
	 \int_\mR (1 + |\xi_y|^2 )^{\frac 1 2 } \norm{h_{\xi_y}}_{H^{1}(\mR)}^2 \;\D \xi_y &= \int_{\mR^2} (1 +|\xi_y|^2 )^{\frac 1 2 } (1 + |\xi_x|^2 ) |\cF(h) (\xi_x,\xi_y)|^2 \;\D \xi_x \D \xi_y \nonumber \\
	 &\leq \int_{\mR^2} (1 + |\xi_x|^2 +|\xi_y|^2)^{\frac 3 2 }  |\cF(h) (\xi_x,\xi_y)|^2 \;\D \xi_x \D \xi_y  \nonumber \\
	 &= \norm{h}_{H^{\frac 3 2} (\mR^2)}^2 < \infty. \label{eq:Unique2DEst2}
	\end{align}
	
	Now let $h_1, h_2 \in H^{\frac 3 2}(\mR^2)$ be two solutions to the phase retrieval problem $|\cF(\psi_j)|^2 =  |\cF(\psi)|^2$ for $\psi_j := b_{\exp, m} + h_j$ with supports in $[0;1] \times [-1;1]$ and define
	\begin{equation*} H: \xi_y \mapsto \max \{ \norm{ h_{1, \xi_y}}_{H^{1}(\mR)}, \norm{h_{2, \xi_y}}_{H^{1}(\mR)} \}. \end{equation*}
	Then $\xi_y \mapsto (1 + |\xi_y|^2 )^{\frac{1}{4}} H(\xi_y)$ is continuous and square-integrable according to \eqref{eq:Unique2DEst2}. Continuity implies that the set
	\begin{align*}
	 V &:= \left\lbrace \xi_y \in \mR :  H(\xi_y) < C^{-1}  (1 + |\xi_y|^2 )^{- \frac 1 2} \right\rbrace \cap S_{> \frac 1 2} 
	\end{align*}
	is open. Moreover, it must be non-empty, as otherwise
	\begin{align*}
	 \int_\mR (1 + |\xi_y|^2 )^{\frac 1 2} |H(\xi_y)|^2 \;\D \xi_y &\geq \int_{S_{> \frac 1 2}} (1 + |\xi_y|^2 )^{\frac 1 2} |H(\xi_y)|^2 \;\D \xi_y \\
	 &\geq  \int_{S_{> \frac 1 2}} C^{-2}(1 + |\xi_y|^2 )^{-\frac 1 2 }   \;\D \xi_y = \infty,
	\end{align*}
	in contradiction to \eqref{eq:Unique2DEst2}.
	
	Combining the estimates \eqref{eq:Unique2DEst0} and \eqref{eq:Unique2DEst1}, we find that for all $\xi_x \in \UHP$, $\xi_y \in  V$ and $j \in \{1,2\}$
	\begin{align*}
	\E^{-\Im(\xi_x)}   (1+ |\xi_x|) \cdot |\cF(\psi) (\xi_x, \xi_y)|  &> \frac{  \E -1 }{  |\xi_y|} - \E^{-\Im(\xi_x)}     (1+ |\xi_x|) \cdot |\cF(h_j) (\xi_x, \xi_y)|\\
	\geq   \frac 1 {|\xi_y|} - C \norm{h_{j, \xi_y}}_{H^{1}(\mR)}  &\geq  \frac 1 {|\xi_y|} - C H(\xi_y)   \stackrel{\xi_y \in V}{\geq}  \frac 1 {|\xi_y|} - \frac 1 {(1 + |\xi_y|^2)^{\frac 1 2 }}  \geq 0,
	\end{align*}
	i.e.\ $\cF(\psi_j)(\cdot ,\xi_y)$ has no zeros in the upper complex half plane for $\xi_y \in V$. Hence, by \thmref{thm:PhaseAmbiguities1D}, $\cF(\psi_1)(\cdot ,\xi_y)$ and $\cF(\psi_2)(\cdot ,\xi_y)$ may differ at most by an exponential factor $\exp( \I(\beta_0(\xi_y) + \beta_1(\xi_y)\xi_x) )$. Notably, $h_j  \in H^{\frac 3 2}(\mR^2) $ is continuous according to \thmref{thm:SobEmbed} and supported in $[0;1] \times [-1;1]$. Thus, we necessarily have
	\begin{equation*} 
	  \lim_{x \to x_0 } \cF( \psi_j (x, \cdot)) =  b _{\exp}(x_0) \cF( b _{\rect}  ) \MTEXT{for}  x_0 \in \{0,1\}.
	\end{equation*}
	This may only hold for all $j\in \{1, 2 \}$ if the exponential factor is unity, i.e.\ if
	\begin{equation*}
	 \cF(\psi_1)(\cdot ,\xi_y) = \cF(\psi_2)(\cdot ,\xi_y) \MTEXT{for all} \xi_y \in V.
	\end{equation*}
	Hence, we obtain $\cF(\psi_1)_{|W} = \cF(\psi_2)_{|W}$ on the open set $W := \mR \times V \subset \mR^2$. Since $\cF(\psi_j)$ is entire, this implies $\cF(\psi_1) = \cF(\psi_2)$ everywhere and therefore $h_1 = h_2$, proving uniqueness of the solution to the phase retrieval problem.
	\end{pf}
	\vspace{1em}
	
	We emphasize that the uniqueness stated in \thmref{thm:UniquePhaseFFMD} is absolute and deterministic, holding for \emph{any} complex-valued $h \in H^{\frac 3 2}(\mR^m)$ supported in $\Omega$ and not just for \emph{almost all} signals or modulo trivial ambiguities. Several other criteria bear either of these defects and/or make additional structural assumptions on regularity, symmetry or  real-valuedness of the reconstructed object, see for instance \cite{Bates1984Uniqueness,Klibanov20062DUniquePurePhase}.
	Another intriguing feature is that, as opposed to the 1D analogue in \thmref{thm:UniquePhaseFF1D}, the perturbation $h$ need not be small compared to the reference signal. All that is necessary is a certain degree of regularity, which yields the asymptotic behavior of $h_{\xi_y}$ for large $\xi_y \in \mR$ required in the proof. Note that the functions $\psi = b_{\exp,m} + h$ for which uniqueness  holds indeed form a \emph{dense} affine subspace of $L^2(\Omega)$ as
	\begin{equation*} \Ckc{\infty}{\Omega} \subset H^{\frac 3 2}(\mR^m) \cap L^2(\Omega) \subset L^2(\Omega) \end{equation*}
	is dense by \thmref{thm:DenseIncl}. Accordingly, we may state: 
	
	\vspace{1em}
	\begin{res}[Unique Phase Retrieval from Fourier Data on a Dense Set] \label{res:PhaseRetrFFDenseUnique}
	  	For $m >~1$ and $\Omega = [0;1] \times [-1;1]^{m-1}$, there exists a dense subset $\A \subset L^2(\Omega)$ such that \probref{prob:FF} is uniquely solvable, i.e.\ any square-integrable signal $\psi$ supported in $\Omega$ may be approximated arbitrarily well by functions for which phase retrieval from Fourier data is unique.
	\end{res}
	\vspace{1em}
	
	It is furthermore noteworthy that, although it might seem constructed and artificial, the reference signal in \thmref{thm:UniquePhaseFF1D} may indeed be implemented in experimental setups: within the framework of the projection approximation (see \sref{SS:ProjApprox}), this can be achieved by placing a phase-shifting plate of rectangular cross-section and exponentially varying thickness in the incident beam such that the unknown specimen in \figref{fig:Setup} lies entirely in its ``shadow''. On the other hand, note that uniqueness may not be ensured by inserting a plate of uniform thickness as the resulting reference signal - a constant rectangular ``bump'' in the phase shifts, i.e.\ a multiple of the support's indicator function - does not break the twin-image symmetry. 
	
	Nevertheless, imposing a support constraint by a reference of constant magnitude is viable alternative to the exponential ramps studied here - for \emph{non-pointsymmetric} supports. For illustration, we consider a scaled indicator function of a triangle in $\mR^2$
	\begin{equation}
	 b_{\tria}(x,y) = \begin{cases}  
	                   (2 \pi)^{\frac 1 2 } &\text{for } x \in [0;1], y \in [-x;x] \\
	                   0 &\text{else}
	                  \end{cases}
	  \label{eq:DefBaseTria}
	\end{equation}
	From \eqref{eq:FourBaseFunc}, it then follows for all $x \in [0;1]$, $\xi_y \in \mC$
			\vspace{-.5em}
	\begin{equation}
	 \cF(b_{\tria}(x,\cdot))(\xi_y) = \frac { 2 \sin( \xi_y x  ) } { \xi_y } = - \frac{ \I }{ \xi_y } \left( \exp( \I  \xi_y x ) - \exp( -\I  \xi_y x ) \right).
	\end{equation}
		In particular, we obtain in the limit $\xi_y = - \I a, a \to \infty$
		\vspace{-.5em}
	\begin{equation}
	 \cF(b_{\tria}(x,\cdot))(\xi_y) \sim \begin{cases} \frac 1 a  \exp(a x)  &\text{for } x \in [0;1]  \\
	                   0 &\text{else} \end{cases}
	\end{equation}
	Accordingly, choosing $b_{\tria}$ as a reference signal results in a similar Fourier space representation as the exponential ramp studied in \thmref{thm:UniquePhaseFF1D}. By making suitable assumptions on the perturbation $h$ bounding its asymptotic growth in Fourier space, this might be exploited to derive uniqueness results for triangular support ``bumps'' by similar techniques as applied in the above proof.
	\end{subsubsection}
	
\subsubsection{Ill-Posedness}

	In the course of this section, we have seen that \probref{prob:FF} is ill-posed in the sense of \defref{def:IllWell} as solutions may be severely nonunique even if a compact support is assumed. However, even in settings where uniqueness can be established, far-field phase reconstruction remains \emph{discontinuous} with respect to reasonable error metrics and thus violates \defref{def:IllWell}-(c). This can be seen by considering
			\vspace{-.25em}
	\begin{equation}
	 F: \Lp 2 \Omega \to \Lp 1 {\mR^m};  \; \psi  \mapsto | \cF(\psi)|^2 \label{eq:FwPhaseRetrFF}
	\end{equation}
	for $m \geq 2$ and $\Omega = [0;1] \times [-1;1]^{m-1}$.
 	By \exref{ex:Frechet}-(c) $F$ is continuous.
	According to Result \ref{res:PhaseRetrFFDenseUnique}, there exists a dense set $\A \subset \Lp 2 \Omega$ such that the restriction $F_{|\A}$ is injective. Thus, the inverse $ (F_{|\A})^{-1} : F(\A) \to \A$ exists. However, since ambiguity persists in the closure $\closure{\A} =  \Lp 2 \Omega$, we may in general construct
	\begin{equation*}
	 \psi \in \A, \tilde \psi \in \Lp 2 \Omega \MTEXT{such that} F(\psi) = F(\tilde \psi) \MTEXT{and} \psi \neq \tilde \psi.
	\end{equation*}
	For $(\psi_j)_{j \in \mN} \subset \A$ with $\psi_j \to \tilde \psi$, continuity of $F$ then implies $F(\psi_j) \to F(\tilde \psi) =   F(\psi) $. Yet, we have $\psi_j \not \to \psi$ by construction, which shows that $(F_{|\A})^{-1}$ is not continuous. We may thus conclude that far-field phase retrieval is severely ill-posed.
	
	This remains true for \emph{near-field} phase retrieval which is studied in the sequel. Non-uniqueness for compact objects, however, turns out to be ruled out completely in this setting by the unscattered probe beam providing a natural reference signal.

  \end{subsection}

   \begin{subsection}{Near-Field Phase Retrieval} \label{SS:PhaseRetrNF}

   We now proceed to the analysis of \probref{prob:NF}, i.e.\ to phase retrieval from near-field data. A first insight is provided by the contrast transfer function (CTF) introduced in \sref{SS:ContrastFormationNF}, which represents a linearization of \eqref{eq:PhaseRetrAbstractNF} in the contact image $\psi \propto \Pdelta - \I \Pbeta$ for plane wave illumination $P_d=1$: according to the derived expression \eqref{eq:CTF} and \cref{cor:FourierL2}, the information encoded in the intensities is sufficient to uniquely reconstruct \emph{either} absorption $\Pbeta$ or phase shifts $\Pdelta$ if the other part is known.
   Both components, i.e.\ arbitrary complex-valued contact images $\psi$, can thus be uniquely recovered whenever \emph{two} intensity measurements for different propagation distances are available. This even remains true if the nonlinearity in \eqref{eq:PhaseRetrAbstractNF} is retained \cite{Jonas2004TwoMeasUniquePhaseRetr}. On the other hand, it is commonly argued \cite{Burvall2011TwoPlanes} that a \emph{single} diffraction pattern is not sufficient for unique reconstructions of complex contact images.
   
   \subsubsection{Counter-Example and Mathematical Setting}
   
   The existence of ambiguities in near-field phase retrieval is indeed confirmed by an explicit counter-example \cite{Nugent2007TwoPlanesPhaseVortex}: if an exit wave $\Psi_0 \in \Lp{2}{\mR^2}$ can be written as
   \begin{equation}
    \Psi_0 ( \rho, \varphi ) = A(\rho ) \exp(\I m \varphi )  \label{eq:PhaseVortex}
   \end{equation}
   in polar coordinates $(\rho, \varphi)$, then the sign of $m \in \mZ$, governing the direction of the ``phase vortex'' described by the exponential factor, cannot be retrieved from the propagated intensities $|\bDF d (\Psi_0)|^2$. However, it should be noted that the vortical phase variations extend infinitely in space. Consequently, an exit wave of the form \eqref{eq:PhaseVortex} may never result from a compactly supported contact image $\psi \in \Sdashc{m}$ imprinted upon incident plane waves $P = 1$, for instance, as the superposition $\Psi_0 = P + \psi$ is constant outside the support of $\psi$, i.e.\  in particular non-vortical. Uniqueness in this physically relevant case is thus \emph{not} ruled out by the counter-example.
   In the sequel, we therefore analyze \probref{prob:NF} for $\A := \Sdashc{m}$.
   
   Considering the argument of the squared modulus in \eqref{eq:PhaseRetrAbstractNF} for plane wave illumination $P_d = 1$
   \begin{equation}
    f( \bxi ):= \exp \left( - \frac{\I \bxi^2} 2 \right)  +  \cF \left(\nF \cdot \psi  \right)(\bxi), \label{eq:ArgMod2NF}
   \end{equation}
   we find that $f: \mC^m \to \mC$ defines an entire function for any $\psi \in \A$. Hence, \probref{prob:NF} amounts to the recovery of an entire function $f$ from $|f|^2$ - just like in the far-field case studied in the preceding section. The principal difference lies in the first summand in \eqref{eq:ArgMod2NF} whose characteristic form arises from the Fresnel propagator. This term defines an entire function of \emph{order two} (compare \sref{SS:EntireFunctions}) growing quadratic-exponentially in $\mC^m$ and is independent of the signal $\psi$ to be reconstructed. On the other hand, the second $\psi$-dependent term in \eqref{eq:ArgMod2NF} may grow at most exponentially by \thmref{thm:PaleyW}. Accordingly, near-field phase retrieval is characterized by a natural holographic reference term of a very specific form. In the analysis of the far-field case in \sref{SS:PhaseRetrFF}, where such terms in turn had to be introduced somewhat artificially, it has been found that these may establish uniqueness in suitable settings.
   
   In the one-dimensional case $m=1$, a complete characterization of the phase retrieval ambiguities in the reconstruction of order-2 entire functions is provided by \lemref{lem:InfoModulus2}. Thereby, we obtain any such function $\tilde f: \mC \to \mC$ satisfying $|f|^2 = |\tilde f|^2$ in the above setting. The crucial point is that the constructed $\tilde f$ need not be \emph{consistent} with the specific structure in \eqref{eq:ArgMod2NF}: if $\tilde f$ corresponds to an alternate solution $\tilde \psi$ of \probref{prob:NF}, then we necessarily have
   \begin{equation}
    f - \tilde f = \cF \left(\nF \cdot (\psi - \tilde \psi)  \right) \label{eq:ArgMod2NFDiff}
   \end{equation}
   as the order-2 reference term in \eqref{eq:ArgMod2NF} must not change under the assignment $\psi \mapsto \tilde \psi$. In particular, $f - \tilde f$ must be an entire function of at most order \emph{one}, which restricts the generality of the order-2 function $\tilde f$ significantly. On the other hand, $f$ and $\tilde f$ are related by factorizations $f = f_1 \cdot f_2$ and $\tilde f = f_1 \cdot f_2^\ast$
   according to \lemref{lem:InfoModulus2}.
   
   \subsubsection{Uniqueness Results}
   
   Based on these observations, the theory of \sref{SS:EntireFunctions} enables us to show that the only consistent solution to all of these constraints is $f = \tilde f$. Generalizing the above setting, this yields the following uniqueness theorem for near-field phase retrieval of compactly supported objects, the proof of which is adapted from the manuscript~\cite{maretzke2014uniqueness}:
   
   \vspace{1em}\begin{th1}[Uniqueness Result for Near-Field Phase Retrieval \text{\cite{maretzke2014uniqueness}}] \label{thm:NFUnique}
 For $w \in \Ck{\infty}{\mR^m}$ everywhere nonzero, $\alpha \in \mC \setminus \mR$ and $\check P \in \Sdashc{m} \setminus \{ 0 \}$ set
 \begin{equation}
  F: \Sdashc{m} \to \Ck{\infty}{\mR^m}; \; F(\psi) = | \cF(\check P )\exp(\alpha (\cdot)^2) + \cF(w \cdot \psi) |^2 \label{eq:defF} 
 \end{equation}
 Then $F$ is well-defined and injective. Moreover, any $ \psi\in \Sdashc{m}$ is uniquely determined by data $F(\psi)_{|U}$ restricted to an arbitrary open set $U \subset \mR^n$.
\end{th1}
\vspace{1em}
	The well-definedness in \thmref{thm:NFUnique} follows from the fact that the argument of the squared modulus \eqref{eq:defF} defines an entire function as discussed above. The remainder of the proof is split into two parts: as a first step, uniqueness is shown for the one-dimensional case $m=1$:
\vspace{1em}
	\begin{pf}[Proof of \thmref{thm:NFUnique} for $m =1$:]
	 Let $\psi, \tilde \psi \in \Sdashc{}$ s.t.\ $F(\psi)_{|U} = F(\tilde \psi)_{|U}$. Define
	 \begin{equation*} 
	  f(\xi):= \cF(\check P ) ( \xi ) \exp(\alpha \xi^2) + \cF(w \cdot \psi)(\xi)
	 \end{equation*}
	  for all $\xi \in \mC$ and $\tilde f$ analogously, so that $F(\psi) = |f|_{|\mR}^2$ and $ F( \tilde \psi) = |\tilde f|_{|\mR}^2$. Since $\alpha \neq0 $ and $\check P \in \Sdashc{} \setminus \{ 0 \}$, $f$ and $\tilde f$ are entire functions of order 2 by \thmref{thm:PaleyW} and \lemref{thm:MaxDecayExpOrder}, matching the setting of \lemref{lem:InfoModulus2}. Accordingly, we have
	 \begin{equation*}
	  f = f_1 \cdot f_2 \MTEXT{and} \tilde f =  f_1 \cdot f_2^\ast
	 \end{equation*}
	  for some entire functions $f_1,f_2$ of order $\leq 2$. Moreover,
	  \begin{equation*}
	 \cF( w \cdot (\psi - \tilde \psi) )  = f - \tilde f = f_1 \cdot (f_2  - f_2^\ast) =: g.,
	 \end{equation*}
	 where $\psi - \tilde \psi \in \Sdashc{}$ is of compact support. Thus, $g$ is an entire function of at most exponential order according to \thmref{thm:PaleyW} and therefore of rank $p_g \leq 1$ by \thmref{thm:Hadamard}. 
	 
	 We show that $\psi = \tilde \psi$ by contradiction. Accordingly, assume $\psi - \tilde \psi \neq 0$. Then $f_1$ and $f_2  - f_2^\ast$ are nonzero factors of order $\leq 2$ of $g$. Consequently, the rank of $f_1$ must be  smaller or equal $p_g \leq 1$ because its zeros $\{a_j\}_{j \in I} \subset \mC \setminus \{ 0 \}$ form a subset of those of $g$, $\{a_j\}_{j \in J}$, where $I \subset J \subset \mN$. This implies that the Hadamard factorization of $f_1$ can be written in the form
	 \begin{align} 
	 f_1(\xi) &= \xi^m \, \exp(\mu_0 + \mu_1 \xi + \mu \xi^2) \prod_{j \in I} E_{1} \left( \frac{\xi}{a_j} \right) \nonumber \\
	 &= \exp( \mu \xi^2)  \underbrace{\left( \xi^m \, \exp(\mu_0 + \mu_1 \xi) \prod_{j \in I} E_{1} \left( \frac{\xi}{a_j} \right) \right)}_{ =:  f_0(\xi)} \label{def:f0}
	 \end{align}
	   for some $\mu_0, \mu_1, \mu \in \mC$. By the same argument as with the rank of $f_1$, the convergence exponent $\rho_{f_0}$ of $f_0$, determined by $\{a_j\}_{j \in I}$, can be at most as large as  $\rho_g$. On the other hand, an application of \thmref{thm:OrderCanProducts} to the Hadamard factorization of $g$ yields $\rho_g \leq \lambda_g \leq 1$. By \thmref{thm:OrderCanProducts}, this implies that $f_0$, as defined in \eqref{def:f0}, is of at most exponential order.
	 
	     Substituting \eqref{def:f0} into the factorizations of $f$ and $\tilde f$ setting $\eta := \Re ( \mu ),\, \gamma := -\Im ( \mu )$, we find that for all $\xi \in \mC$
	 \begin{align*}
	  f(\xi) &= \exp( \mu \xi^2) f_0(\xi)  f_2(\xi) = \exp(-\I  \gamma \xi^2) f_0(\xi)  \exp( \eta \xi^2)f_2(\xi)  \\ 
	  \tilde f(\xi)   &= \exp( \mu \xi^2) f_0(\xi)   \cc{f_2(\cc\xi)} = \exp(-\I  \gamma \xi^2) f_0(\xi)  \left(\cc{\exp(\eta \cc\xi^2) f_2(\cc\xi)} \right)
	 \end{align*}
	 These equalities show that the factor $\exp( \eta \xi^2)$ may be absorbed in $f_2$, as it is invariant under Schwarz reflection $\adj{}$. Thus, we may assume $\eta = 0$ without loss of generality. This implies for all $\xi \in \mC$
	 \begin{align*}
	  g(\xi ) &= \exp( - \I \gamma \xi^2) f_0(\xi)  (f_2(\xi) - \adj{f_2}(\xi)) 
	 \end{align*}
	 and  by multiplication with $\adj{f_0}$ and application of $\adj{}$
	 \begin{equation}
	  f_0(\xi ) \adj{g}(\xi ) = -\exp( 2\I \gamma \xi^2) \adj{f_0}(\xi) g(\xi ) \label{eq:fvsfast}.
	 \end{equation}
	 $f_0 \cdot \adj{g}$ and $\adj{f_0} \cdot g$ are both nonzero entire functions of order $\leq 1$, whereas  $\xi \mapsto \exp( 2\I \gamma \xi^2)$ is of order $2$ for any $\gamma \neq 0$. According to \lemref{thm:MaxDecayExpOrder}, this super-exponential growth could not be compensated by the remaining at most exponential order factors on the right hand side of \eqref{eq:fvsfast}, so that the only possibility for \eqref{eq:fvsfast} to hold for all $\xi \in \mC$ is $\gamma = 0$.
%
%
%
	 
	 Recalling the definition of $f$ and $\tilde f$ and setting $a:=\cF(\check P)$, $b:= \cF( w \cdot \psi )$, $\tilde b:= \cF( \tilde  w \cdot \psi )$ and $e(\xi) := \exp(\alpha \xi^2)$ for all $\xi \in \mC$, this and the preceding results imply
 	 \begin{align}
 	  f_0  \cdot  f_2  &= \, f\, =a \cdot e +   b \label{eq:Aa} \\
 	   \adj{f_0} \cdot  f_2 &= \adj{\tilde f}  =  \adj{a} \cdot \adj{e} +   \adj{\tilde b}. \label{eq:Ab}
 	 \end{align}
 	 By multiplication of \eqref{eq:Aa} and \eqref{eq:Ab} with $\adj{f_0}$ and $f_0$, respectively, we obtain
 	 \begin{equation}
 	   \adj{f_0} \cdot ( a \cdot e  +  b ) =  \adj{f_0} \cdot f_0 \cdot f_2 =  f_0  \cdot (  \adj{a}  \cdot \adj{e} +     \adj{\tilde b} ) \label{eq:B}
 	 \end{equation}
	   For $c \in \{-1,1\}$, consider the diagonals in the complex plane
 	 \begin{equation*}
 	  D_c := \{ z \in \mC :  \Re(z) = c \Im(z) \}
 	 \end{equation*}
 	 and let $s$ denote the sign of $\Im( \alpha )$ (recall that $\Im(\alpha ) \neq 0$ is assumed). Then we have 
  	 \begin{align*} 
  	  |e(\xi)| &= \begin{cases} \exp(-|\Im(\alpha)| |\xi|^2) &\text{for }\xi \in D_s \\
                               \exp( |\Im(\alpha)| |\xi|^2) &\text{for }\xi \in D_{-s} \end{cases} \\
  	  |\adj{e}(\xi)| &= \begin{cases} \exp(|\Im(\alpha)| |\xi|^2) &\text{for }\xi \in D_s \\
                               \exp( -|\Im(\alpha)| |\xi|^2) &\text{for }\xi \in D_{-s} \end{cases}.
  	 \end{align*}
 	 Since all of the remaining factors in \eqref{eq:B} are non-vanishing entire functions of at most exponential order, this implies that the right hand side of \eqref{eq:B} is $\Or(\exp(|\Im(\alpha)| |\xi|^2))$ in $D_s$, whereas the left hand side grows at most exponentially along this diagonal. Contradiction!

 	 Accordingly, the initial assumption $\psi \neq \tilde \psi$ must be wrong. By generality of $\psi, \tilde \psi \in \Sdashc{}$, this proves injectivity of the operator $F$ in the case $m=1$. 
	\end{pf}
\vspace{1em}

 In \sref{SS:PhaseRetrFF}, higher dimensional phase retrieval has been related to the 1D uniqueness theory by partially Fourier-transforming the objective function (see \eqref{eq:1DReductionFourier}), giving rise to a family of compactly supported distributions in one dimension. Here, the general statement of \thmref{thm:NFUnique} is obtained by a similar dimension reduction argument, combined with an application of the 1D result that has already been shown in the first step:
 
 \vspace{1em}
\begin{pf}[Proof of \thmref{thm:NFUnique}:]
Injectivity has already been proven in the case $m=1$, so that we may restrict ourselves to  $m \geq 2$. 

Let $\psi, \tilde \psi \in \Sdashc{m}$ such that $F(\psi)_{|U} = F(\tilde \psi)_{|U}$ for some $U \subset \mR^m$ open. Like in the 1D case, $F(\psi)$ and $F(\tilde \psi)$ have extensions to entire analytic functions in $\mC^m$ by \thmref{thm:PaleyW}, so that $F(\psi) = F(\tilde \psi)$ everywhere. Let $\cF_{\overline 2}: \Sdash{m} \to \Sdash{m}$ denote the Fourier transform in all variables but the first. For $\bxi_y \in \mR^{m-1}$, we set
\begin{subequations} \label{eq:pfthm1-1}
\begin{align}
 P_{0,\bxi_y} &:= \cF_{\overline 2} ( \check P) (\cdot, \bxi_y)\exp ( \alpha  \bxi_y^2 ), \\
  \psi_{ \bxi_y} &:= \cF_{\overline 2} ( w \cdot \psi ) (\cdot, \bxi_y), \\
    \tilde \psi_{ \bxi_y} &:= \cF_{\overline 2} ( w \cdot \tilde \psi ) (\cdot, \bxi_y). 
\end{align}
\end{subequations}
Then $P_{0,\bxi_y}, \psi_{ \bxi_y}, \tilde \psi_{ \bxi_y} \in \Sdashc{}$ and there exists an open set $V \subset \mR^{m-1}$ such that $P_{0,\bxi_y} \neq 0$ for all $\bxi_y \in V$. By construction, we have for all $\bxi_y \in \mR^{m-1}, \xi_x \in \mR$, $\bxi = (\xi_x, \bxi_y)$
\begin{equation*}
\cF(P_{0,\bxi_y} ) ( \xi_x ) \exp(\alpha \xi_x^2) + \cF(\psi_{ \bxi_y})(\xi_x) =\cF(\check P ) ( \bxi ) \exp(\alpha \bxi ^2) + \cF(w \cdot \psi)(\bxi)
\end{equation*}
and an analogous equality for $\tilde \psi$ and $\tilde \psi_{ \bxi_y}$. This implies by assumption
\begin{align}
|\cF(P_{0,\bxi_y} ) ( \xi_x ) &\exp(\alpha \xi_x^2)  + \cF(\psi_{ \bxi_y})(\xi_x)|^2 = F(\psi)(\bxi) \nonumber \\ 
&= F(\tilde \psi)(\bxi) = |\cF(P_{0,\bxi_y} ) ( \xi_x ) \exp(\alpha \xi_x^2) + \cF(\tilde \psi_{ \bxi_y})(\xi_x)|^2 \label{eq:pfthm1-2}
\end{align}
for all $\bxi_y \in \mR^{m-1}, \xi_x \in \mR$.

The leftmost and rightmost expressions in \eqref{eq:pfthm1-2} are exactly the images of $ \psi_{ \bxi_y}$ and $\tilde \psi_{ \bxi_y}$ under the operator $F$ in the \emph{one-dimensional} setting $m = 1$, $\check P = \check P_{\bxi_y}$ and $w = 1$. By application of \thmref{thm:NFUnique} for $m=1$, \eqref{eq:pfthm1-2} thus implies 
\begin{equation}
 \psi_{ \bxi_y} =  \tilde \psi_{ \bxi_y} \quad \text{for all}\quad\bxi_y \in V \label{eq:pfthm1-3}
\end{equation}
According to \thmref{thm:PaleyW}, $ \bxi_y \mapsto \psi_{ \bxi_y}  $ and $ \bxi_y \mapsto \tilde \psi_{ \bxi_y} $ are entire analytic functions  so that \eqref{eq:pfthm1-3} holds even for $\bxi_y \in \mR^{m-1}$. By bijectivity of $\cF_{\overline 2}$ and invertibility of $w$ in a multiplicative sense, $\psi$ and $\tilde \psi$ can be recovered uniquely from $\{ \psi_{ \bxi_y} \}_{\bxi_y \in \mR^{m-1}}$ and $\{ \tilde \psi_{ \bxi_y} \}_{\bxi_y \in \mR^{m-1}}$, respectively, by inversion of \eqref{eq:pfthm1-1}.

Since these families coincide by the 1D uniqueness result, we obtain $\psi = \tilde \psi$ which proves injectivity of $F$.
\end{pf}
\vspace{1em}

By comparison of \eqref{eq:defF} to \eqref{eq:PhaseRetrAbstractNF}, we find that \probref{prob:NF} for plane wave illumination $P_d=1$ exactly matches the setting of \thmref{thm:NFUnique} for the parameters
\begin{equation}
 \alpha = -\frac{\I} 2, \;\;\;\; w = \nF \MTEXT{and} \check P = (2\pi)^{\frac m 2 }  \delta_0.
\end{equation}
Here, $\delta_0$ denotes the Dirac delta centered at 0, compare \exref{ex:DiracDelta} and \exref{ex:FTDiracDelta}. In this setting, the injectivity statement in \thmref{thm:NFUnique} thus leads to the startling conclusion that \emph{any} compactly supported complex-valued contact image may be uniquely reconstructed from near-field intensity data recorded at a single distance.

Moreover, the freedom in the choice of $\check P$ and $\alpha$ in the uniqueness result makes it applicable to \probref{prob:NF} for a large number of other probe functions $P$. As an example, we consider illumination by a Gaussian beam \cite[sec. 3.1]{Teich1991Photonics}, characterized by a propagated wave field of the form
\begin{equation}
 P_{\Text{Gauss}} (\bxi) = \exp( \gamma_0 + \alpha_0 \bxi ^2 ) \MTEXT{where} \gamma_0, \alpha_0 \in \mC, \Re(\alpha_0) < 0, \Im(\alpha_0) \leq 0 \label{eq:GaussianBeamProbe}
\end{equation}
Inserting this probe contribution into \eqref{eq:PhaseRetrAbstractNF}, the resulting expression is found to match \eqref{eq:defF} for the parameter choices
\begin{equation}
 \alpha = \alpha_0 -\frac{\I} 2, \;\;\;\; w = \nF \MTEXT{and} \check P = (2\pi)^{\frac m 2 } \exp(\gamma_0)  \delta_0.
\end{equation}
Hence, \thmref{thm:NFUnique} yields uniqueness of \probref{prob:NF} also for this more realistic illumination function. The findings are summarized by the following corollary:
	\vspace{1em}
	\begin{cor1}[Uniqueness of Near-Field Phase Contrast Imaging \text{\cite{maretzke2014uniqueness}}] \label{cor:UniqueNFPlaneWave}
	 Near-Field Phase Retrieval of compactly supported images $\psi$, given by \probref{prob:NF} for $\A = \Sdashc{m}$, is uniquely solvable for $P_d \in \{ 1, P_{\Text{Gauss}} \}$ corresponding to illumination with plane waves or Gaussian beams. Moreover, any $\psi  \in \A$ can be reconstructed from intensity data $I_{|U}$ of the form \eqref{eq:PhaseRetrAbstractNF} restricted to an arbitrary open set $U \subset \mR^m$.
	\end{cor1}
	\vspace{1em}

   \subsubsection{Ill-Posedness}
   
   For an investigation of ill-posedness of near-field phase retrieval, it is once more illustrative to consider the linearization of \eqref{eq:PhaseRetrAbstractNF}, corresponding to the  CTF representation \eqref{eq:CTF} to which we already referred at the beginning of this section \sref{SS:PhaseRetrNF}. The zeros of the sinusoidal prefactors plotted in \figref{fig:CTF} correspond to Fourier frequencies which are not represented in the near-field intensities as discussed in \sref{SS:ContrastFormationNF}. Accordingly, these give rise to arbitrary error amplifications in the inversion, i.e.\ discontinuity and thus \emph{ill-posedness} of phase retrieval, even if only the real- \emph{or} the imaginary part of the contact image $\psi\propto \Pdelta - \I \Pbeta$ is to be reconstructed.
   
   On the other hand, \cref{cor:UniqueNFPlaneWave} implies that a unique phase reconstruction from \emph{exact} data is still possible for arbitrary complex-valued contact images $\psi$, provided that these are compactly supported. In this case, however, the problem is more severely ill-posed since  uniqueness breaks down for non-compact supports according to the ``phase vortex'' counter-example discussed above - as opposed to the invertibility of the CTF \eqref{eq:CTF} with respect to \emph{either} $\Pdelta$ or $\Pbeta$. Consequently, any \emph{stability estimate}, by which the discontinuity of near-field phase retrieval of general complex-valued images might be bounded, would need to incorporate the support size in a suitable sense. Unfortunately though, the non-constructive proof of \thmref{thm:NFUnique} does not give any hint on how this might be achieved in detail.

  \end{subsection}

\end{section}

\begin{section}{Uniqueness of Phase Contrast Tomography} \label{S:IntermConcl}
 
 In this chapter, we have studied regularity and ill-posedness of the inverse reconstruction problem of phase contrast tomography, given by \probref{prob:1}. The analysis of the different subproblems in the preceding sections now enables us to deduce statements for the complete forward operators \eqref{eq:ForwardOpNF} and \eqref{eq:ForwardOpFF}.
 
 In the near-field case, the intermediate results from \sref{S:GenIll-posed}, \sref{S:RadonIll-posed} and \sref{S:PhaseRetrieval} indeed imply that the tomographic reconstruction is uniquely solvable for known probe functions of reasonable shape if phase-wrapping (see \sref{SS:PhaseWrap}) is absent. The latter can be ensured by restricting to
 \begin{equation}
  \fD_F := \{ N \in \Lp \infty {\ObjDom}_\mR : \Re(N) \geq 0 , \; k L \norm{ \Re(N) }_{\Lp \infty {\ObjDom}} < 2 \pi \} \label{eq:DefNonPhaseWrapDom}.
 \end{equation}
 where $L > 0$ denotes the diameter of the object domain $\ObjDom$.

 	\vspace{1em}
	\begin{cor1}[Uniqueness of Near-Field Phase Contrast Tomography \text{\cite{maretzke2014uniqueness}}] \label{cor:UniqueNFPCT}
	 For non-zero propagation distance $d> 0$, wavenumber $k> 0$, object diameter  $L > 0$ and propagated probe field $\bDF d (P) \in \{ 1, P_{\Text{Gauss}} \}$, the forward operator of near-field phase contrast tomography defined by \eqref{eq:ForwardOpNF} is \emph{injective} as a map
	 \begin{equation}
	  F_d : \fD_F \to \Lp \infty {Z^{m+1}}.
	 \end{equation}
	 Moreover, any $N \in \fD_F$ is uniquely determined by the data $F_{d} (N) _{|W}$ on a wedge-shaped set $W:= V\times U$ for $V \subset [0;2\pi)$, $U \subset \mR^m$ open. In particular, \probref{prob:1} is uniquely solvable in this setting.
	\end{cor1}
	\begin{pf}
	 Let $N \in \fD_F$ and  $W:= U\times V$ for $U \subset [0;2\pi)$, $V \subset \mR^m$ be arbitrary. Since $N$ is compactly supported, so is
	 \begin{equation*} \psi_\theta := P \cdot O_0(N)(\theta, \cdot) = P \cdot [ \exp(- \I k \CR(N)(\theta, \cdot) ) - 1] \end{equation*}
        for all $\theta \in [0; 2\pi)$. Recalling that \probref{prob:NF} was derived in \sref{SS:PhaseRetrAbstract} as an abstract formulation of the reconstruction of $\psi_\theta $ from $F_d(N)(\theta, \cdot)$, \cref{cor:UniqueNFPlaneWave} implies that $\psi_\theta$ is uniquely determined by $ F_{d} (N)_{|W}(\theta, \cdot)$ for all $\theta \in V$. By assumption, the probe $P$ is everywhere nonzero so that 
        \begin{equation*}
        \Lp \infty {\mR^m} \to \Lp \infty {\mR^m}; \;   \exp(- \I k \CR(N)(\theta, \cdot) ) \mapsto \psi_\theta
        \end{equation*}
	is injective, i.e.\ likewise uniquely invertible. By construction of $\fD_F \ni N$ and \thmref{thm:RadonContinuity}, we further have for all $\theta \in [0;2 \pi), \bx \in \mR^m$
	\begin{equation*}
	 - \I k  \CR(N)(\theta, \bx) \in \{ z \in \mC : \Im (z) \in (-2\pi ; 0 ] \}
	\end{equation*}
	so that the pointwise exponential  $ \CR(N)  \to \exp(- \I k \CR(N)) $ is invertible, representing the absence of phase-wrapping. Hence, the data $ F_{d} (N) _{|W}$ uniquely determines $\CR(N)_{V\times \mR^m}$, from which $N$ can be uniquely reconstructed according to the injectivity result for the Radon transform in \thmref{thm:RadonInjective}.
	\end{pf}
	\vspace{1em}
	
 The only difference in the far-field case lies in the phase retrieval step as discussed in \sref{S:PhaseRetrieval}. Here, no uniqueness statement of comparable generality as \cref{cor:UniqueNFPlaneWave} could be derived. We therefore have to content ourselves a with less concrete result:
  	\vspace{1em}
 	\begin{cor1}[Uniqueness of Far-Field Phase Contrast Tomography] \label{cor:UniqueFFPCT}
	 Let $k,L$ be as in \cref{cor:UniqueNFPCT} and let $P \in \mP$ with $\supp(P) = \mR^m$ and $\mP$ defined by \eqref{eq:ProbeDomain}. Then \probref{prob:1} has a unique solution $N \in \fD_F$ for far-field intensities $F_\infty(N) = I_\infty$ given by \eqref{eq:ForwardOpFF} whenever the corresponding phase retrieval problems
	 \begin{equation}
	  |\cF(P \cdot O_0(N)(\theta, \cdot))|^2 = I_\infty(\theta, \cdot) \label{eq:FFPCTPhaseRetr}
	 \end{equation}
	are uniquely solvable for all $\theta \in V$ in some $V \subset [0;2\pi)$ open. In this case, $N$ is uniquely determined by intensities $I_{\infty|V\times U}$ on any open set $U \subset \mR^m$. 
	\end{cor1}
	\begin{pf}
	 In the given setting, $P \cdot O_0(N)_{|V\times \mR^m}$ can be uniquely reconstructed from the data because $I_\infty(\theta, \cdot) $ is  uniquely determined by its values on $U \subset \mR^m$ by \thmref{thm:PaleyW}. The remainder of the proof works exactly as in \cref{cor:UniqueNFPCT}.
	\end{pf}
	\vspace{1em}
	According to \sref{SS:PhaseRetrFF}, uniqueness for the individual phase retrieval problems \eqref{eq:FFPCTPhaseRetr} may be established by symmetry or monotonicity assumptions on the solution or by the holographic approach of superimposing a known reference signal as in \thmref{thm:UniquePhaseFFMD}. Alternatively, one may hope for uniqueness based on the observation in \resref{res:NDFTPhaseRetrAmbiguities} that almost all images for $m\geq 2$ can be uniquely reconstructed up to trivial ambiguities. In either case, ambiguity is reduced significantly by the \emph{combination} of Radon inversion and phase retrieval in \cref{cor:UniqueNFPCT}: due to the correlations between the projections for different $\theta$, uniqueness already holds if the phase retrieval problems are uniquely solvable for an arbitrarily small wedge of incident angles. This observations constitutes a major motivation for the numerical reconstruction method introduced in the subsequent \chapref{C:NumMeth}.
	
\end{section}

\end{chapter}


\begin{chapter}{Reconstruction Method}\label{C:NumMeth}


In the preceding \chapref{C:Analysis}, it has been found that phase contrast tomography constitutes an ill-posed inverse problem in a number of different aspects, including possible non-existence or non-uniqueness of exact solutions and in particular discontinuity of the inverse operators. On the other hand, the governing forward operators turn out to be \Frechet differentiable, i.e.\ well-posed and smooth. In this chapter, we introduce regularized Newton-type methods, which allow a numerical solution of the inverse reconstruction problem by exploiting its particular structure.

\begin{section}{Algorithms and the Simultaneous Approach} \label{S:CombinedApproach}
 
 In the analysis of \chapref{C:Analysis}, we have decomposed \probref{prob:1} into different subproblems in order to investigate its ill-posedness. It likewise seems natural to implement the different reconstruction steps independently, i.e.\ phase retrieval, recovery of the sinogram from the object transmission function and Radon inversion (see \sref{S:GenIll-posed}) as subsequent operations in the work flow. 
 
 A major advantage of this separation is efficiency: in the case of weakly interacting samples, for instance, the near-field phase problem may be solved directly by inverting the contrast transfer function (CTF, see \eqref{eq:CTF}). Combined with direct tomographic reconstruction via \emph{filtered backprojection} (see \cite[sec. V.1]{Natterer}), this approach allows for efficient and accurate imaging from experimental data \cite{Cloetens1999,Cloetens1999Diss,BartelsDiss}.
 Another class of direct phase reconstruction methods outlined and applied in \cite{Reed1983TIEPhaseRetr,Nugent1996Propagationbased} is based on a linearization of the \emph{transport-of-intensity equations} (compare e.g.\ \cite[sec. 4.5.2]{PaganinXRay}), which essentially corresponds to a linearization of the sine term in \eqref{eq:CTF} valid in the limit of small propagation distances. Recent enhancements of this approach such as \emph{Bronnikov-Aided-Correction} \cite{Bronnikov2002BronnikovMod,Witte2009BronnikovAC} can be applied to samples which are weakly absorbing but not necessarily weakly refracting. From an experimental point of view, another advantage of the latter methods is their relative insensitivity to polychromaticity which allows X-ray imaging with laboratory sources, as demonstrated e.g.\ in \cite{Wilkins1996,Paganin1998} and more recently by \cite{Krenkel2014BCAandCTF}. For an overview of direct phase reconstruction techniques, see for instance \cite{Burvall2011TwoPlanes}.
 
 Major drawbacks of these direct methods lie in their restriction to the limited range of validity of the underlying linearizations. For instance, propagation distances in X-ray nanoscopy will typically not be small compared to other lengthscales of the setup. Moreover, accurate CTF-reconstructions of a single projection typically require holograms recorded at multiple detector distances \cite{Cloetens1999,Krenkel2014BCAandCTF}. Finally, no equivalents of these phase retrieval methods exist for \emph{far-field} phase retrieval. These facts motivate phase reconstruction by \emph{iterative} methods. The most commonly used essentially go back to the ideas of \citet{Gerchberg1972GS} and \citet{Fienup1982HIO}: the current iterate is projected alternatingly onto \emph{constraint sets} defined by the measured intensities at possibly multiple propagation planes or by available priori knowledge e.g.\ on support or positivity of the solution. Further improvements of this approach such as the \emph{Shrinkwrap Algorithm} \cite{Fienup1986algorithm} provide iterative support adaption or faster convergence as achieved e.g.\ by \emph{Relaxed Averaged Alternating Reflections} \cite{Luke2005RAAR}. Applications of these \emph{convex optimization} techniques to experimental far-field data have been shown to yield good reconstructions of single material objects from a single intensity measurement without further constraints \cite{Marchesini2003PhysRevB,Bartels2012}.
 
 For the present work, the main benefit of iterative methods is their flexibility which allows for \emph{simultaneous} phase retrieval and Radon inversion. In the far-field case, we have seen in \sref{SS:PhaseRetrFF} that ambiguities reduce tremendously from a single dimension to two-dimensional images. Consequently, a similar improvement can be expected to occur in the transition to \emph{three-dimensional} phase retrieval, to which simultaneous tomographic- and phase reconstruction amounts at least in the weak object limit, see \eqref{eq:WeakFwOpFF}. This conjecture is supported by \cref{cor:UniqueFFPCT} stating that not all of the projections need to be uniquely reconstructible from the corresponding far-field intensities in order to ensure uniqueness of the reconstructed 3D sample - as might be expected if phase retrieval was considered as an independent  subproblem.
 
 The mathematical reason for the apparent stabilizing effect of combining the different steps lies in the strong correlations between projections of one and the same 3D object expressed by the Helgason-Ludwig-Consistency-Conditions in \thmref{thm:HelgasonLudwig}. Even in the near-field case, where uniqueness already holds for phase retrieval of the single projections according to \cref{cor:UniqueNFPlaneWave}, the exploitation of these correlations may be beneficial to reduce ill-posedness: holograms recorded at only slightly different incident angles may have a similar effect as measurements at multiple propagation distances. The Radon inversion, on the other hand, may be less susceptible to artifacts if data inconsistencies are accounted for already in phase retrieval.

 For far-field phase contrast tomography, the simultaneous approach has been successfully implemented by interpreting tomographic data in the sense of \eqref{eq:WeakFwOpFF} as Fourier intensities on a cylindrical grid \cite{Chapman2006,Barty2008ceramicfoam}. Interpolated to Cartesian coordinates, the 3D data set is assigned to iterative phase retrieval algorithms of the alternating-projection-type described above. In the near-field case, simultaneous reconstruction has been implemented in the form of the  iterative reprojection phase retrieval  algorithm (IRP) \cite{Ruhlandt2014}: here, the idea is to embed the iterative \emph{Algebraic Reconstruction Technique} for Radon inversion (ART, see \cite{Kaczmarz1937ART,Gordon1970ART} and \cite[sec. V.4]{Natterer}) in Gerchberg-Saxton-type phase retrieval iterations. Thereby, consistency of the projections is  imposed implicitly. This results in significantly improved reconstructions as demonstrated for simulated data \cite{Ruhlandt2014} - especially in the case of general objects for which refraction $\delta$ and absorption $\beta$ have to be reconstructed independently.
 
 In this work, simultaneous phase retrieval and Radon inversion is enforced by a more radical approach. The principal idea is simply to invert the forward operators of phase contrast tomography introduced in \sref{SS:PhaseTomoOperator} \emph{as a whole}, which ensures precise book-keeping of tomographic correlations. Owing to the nonlinearity of the problem, this can only be achieved by iterative methods. Our choice here is given by iteratively regularized Newton methods \cite{Bakushinskii1992IRNM}, which have already been applied to (non-tomographic) far-field phase retrieval problems \cite{Hohage2013}. As opposed to the convex optimization methods discussed above, this approach takes advantage of the \Frechet differentiability proven in \sref{S:WellDefFrechet}, which promises improved convergence. At the same time, the regularization accounts for the various forms of ill-posedness (see \sref{S:GenIll-posed} - \sref{S:PhaseRetrieval}) of the inverse reconstruction problem to be solved.
 
\end{section}

\begin{section}{Regularized Newton-Type Methods} \label{S:IRNM}

\begin{subsection}{Motivation and Setting} \label{SS:IRNM-Motivation}

In the preceding chapters, we have seen that phase contrast tomography  amounts to the solution of an ill-posed nonlinear operator equation of the form
\begin{equation}
 F(f) = g^{\Textbf{err}} \label{eq:GenOpEq}
\end{equation}
for a \Frechet differentiable operator $F: \mX \to \mY$. The right hand side is given by imperfect, noisy observations
\begin{equation} g^{\Textbf{err}} = g^\dagger + \Textbf{err} \MTEXT{with} g^\dagger = F(f^\dagger)  \in \mY  \end{equation}
whereas the ideal data $g^\dagger$ corresponding to the exact solution $f^\dagger \in \mX$ is unknown.

By analogy to nonlinear equations in $\mR$, a straightforward approach for seeking an approximate solution to \eqref{eq:GenOpEq} is by Newton's method, iteratively solving the linearized problems
\begin{equation}
 F(f_k) +  F'[f_k](f_{k+1} - f_k) = g^{\Textbf{err}} \label{eq:GenOpLin}
\end{equation}
in the $k$-th iterate. However, as the nonlinear problem \eqref{eq:GenOpEq} is ill-posed, so are in general the linearizations \eqref{eq:GenOpLin} to be solved in the Newton iterations \cite[p. 285]{Hanke1996Regularization}. For instance, the derivatives obtained in \thmref{thm:ForwardOpsFrechet} still involve the Radon transform $\CR$, which typically does not admit an exact solution for noisy data and whose inverse is unbounded as seen in \sref{SS:RadonSurjective} and \sref{SS:RadonContiInv}. Accordingly, even the single Newton iterates defined by \eqref{eq:GenOpLin} may not have a unique solution for all $k$ and - even if so - will in general not depend continuously on the data. As only noisy data is available, this implies that standard Newton's method is not applicable to the problems considered in this work. 

\end{subsection}

\begin{subsection}{Iteratively Regularized Gauss-Newton method} \label{SS:IRGNM}

A remedy for the ill-posedness of the linearizations \eqref{eq:GenOpLin} is to slightly modify the problem, computing the Newton iterates as the solution to the quadratic minimization problem
\begin{equation}
 f_{k+1} = \argmin_{f \in \mX} \left( \norm{ F(f_k) +  F'[f_k](f - f_k) - g^{\Textbf{err}} }_{\mY}^2 + \alpha_k \norm{ f- f_0 }_{\mX}^2 \right) \label{eq:GenOpQuadMin}
\end{equation}
where $\alpha_k > 0$ is a \emph{regularization parameter} and $f_0 \in \mX$ denotes the initial guess. The iterates \eqref{eq:GenOpQuadMin} define the \emph{iteratively regularized Gauss-Newton method} (IRGNM) proposed by \citet{Bakushinskii1992IRNM}. Essentially, it corresponds to Tikhonov regularization (see for instance \cite[C. 5]{Hanke1996Regularization}) applied to the linearized problem \eqref{eq:GenOpLin}.

The following lemma shows that the practical problems arising from ill-posedness, namely non-existence, non-uniqueness or discontinuity of the inverse, are ruled out by the introduced \emph{regularization term} $\alpha_k \norm{ f- f_0 }_{\mX}^2$ in the minimization problem formulation:
\vspace{1em}
\begin{lem1}[Well-Posedness of the IRGNM \text{\cite[p. 286]{Hanke1996Regularization}}] \label{lem:RegIRGNM}
 Let $\alpha_k > 0 $, $g^{\Textbf{err}} \in \mY$ and $f_0, f_k \in \mX$. Then the quadratic minimization problem \eqref{eq:GenOpQuadMin} has the unique solution
 \begin{equation}
 f_{k+1} = f_k +  T_k^{-1} \left(F'[f_k]^\ast (g^{\Textbf{err}} - F(f_k)) + \alpha_k(f_0 - f_k) \right)   \label{eq:GenOpRegSol}
\end{equation}
where $T_k := F'[f_k]^\ast F'[f_k] + \alpha_k : \mX \to \mX $ is an isomorphism with  $\norm{T_k^{-1}}  \leq \frac 1 {\alpha_k}$. In particular, the computation of the \emph{IRGNM}-iterates is \emph{well-posed}.
\end{lem1}
\begin{pf}
 For given $\alpha_k > 0 $, $g^{\Textbf{err}} \in \mY$ and $f_0, f_k \in \mX$, consider the quadratic functional $\Phi: \mX \to \mR$ defined by the argument on the right hand side of \eqref{eq:GenOpQuadMin}. As $\mX$ and $\mY$ are Hilbert spaces, $\Phi$ is strictly convex and thus has a unique minimizer $f_{k+1} \in \mX$. By \thmref{thm:FrechetProps} and \exref{ex:Frechet}-(a), $\Phi$ is furthermore \Frechet differentiable where the derivative for all $f, h \in \mX$ is given by
 \begin{align}
  \Phi'[f]h &=  2  \ip{F(f_k) +  F'[f_k](f - f_k) - g^{\Textbf{err}}}{F'[f_k] h}_{\mY} +2 \alpha_k  \ip{f- f_0 }{h}_{\mX} \nonumber \\
  &=  2  \ip{F'[f_k]^\ast ( F(f_k) +  F'[f_k](f - f_k) - g^{\Textbf{err}})+  \alpha_k  (f- f_0) }{ h}_{\mX}. \label{eq:RegIRGNM-1}
 \end{align}
 Here, the defining property of the adjoint in \defref{def:OpsAdjoints} has been used. According to \thmref{thm:FrechetProps}-(f), the global minimizer of $\Phi$ is characterized by $\Phi'[f_{k+1}] = 0$ corresponding to a vanishing first argument on the right hand side of \eqref{eq:RegIRGNM-1}, i.e.
 \begin{equation}
  F'[f_k]^\ast ( F(f_k) +  F'[f_k](f_{k+1} - f_k) - g^{\Textbf{err}})+  \alpha_k  (f_{k+1} - f_0) = 0. \label{eq:RegIRGNM-2}
 \end{equation}
 By rearranging \eqref{eq:RegIRGNM-2}, we obtain the solution \eqref{eq:GenOpRegSol}. Note that the operator
 \begin{equation*} T_k := F'[f_k]^\ast F'[f_k] + \alpha_k : \mX \to \mX \end{equation*}
 is an isomorphism with $\norm{T_k^{-1}} < \frac 1 {\alpha_k}$ by the Lax-Milgram Theorem (see \cite[p. 247]{Werner2007FA}), as $T_k$ is bounded and uniformly \emph{positive-definite}. 
 Applying the Cauchy-Schwartz inequality and the definition of the adjoint, the latter is seen from
 \begin{equation}
  \norm{f}_{\mX} \norm{T_k f}_{\mX} \geq \ip{(F'[f_k]^\ast F'[f_k] + \alpha_k)f }{f}_{\mX} = \underbrace{\norm{F'[f_k] f }^2_{\mX}}_{\geq 0} + \alpha_k \norm{ f }^2_{\mX}
 \end{equation}
 for all $f \in \mX$. Accordingly, $T_k^{-1}$ is well-defined and continuous so that the solution of the quadratic minimization problem \eqref{eq:GenOpQuadMin} given by \eqref{eq:GenOpRegSol} is well-posed.
\end{pf}

\end{subsection}

\begin{subsection}{Choice of the Regularization Parameter} \label{SS:RegPar}

By \lemref{lem:RegIRGNM}, the IRGNM is well-defined and may be stably implemented. However, it remains to be investigated whether the iterates $ \{ f_k \}_{k \in \mN_0}  $ provide reasonable approximations to the solution of the nonlinear operator equation \eqref{eq:GenOpEq}. This depends significantly on the choice of the regularization parameters $ \{ \alpha_k\}_{k \in \mN_0} $.

 In the case $\alpha_k = 0$,  minimizers of \eqref{eq:GenOpQuadMin} are exact solutions of \eqref{eq:GenOpLin} whenever such exist. For $\alpha_k > 0$, the  regularization term  $\alpha_k \norm{ f- f_0 }_{\mX}^2$ enforces uniqueness of the iterates and ensures that they depend continuously on the data $g \in \mY$ according to \lemref{lem:RegIRGNM}. Consequently, If $f_{k+1}^{\dagger} \in \mX$ is the subsequent iterate to $f_k$ computed by replacing (only in the $k$-th iteration!) $g^{\Textbf{err}}$ with $g^\dagger$ in \eqref{eq:GenOpQuadMin} and $f_{k+1}$ its analogue from noisy data, then the resulting error can be estimated by 
\begin{equation}
 \norm{ f_{k+1}  - f_{k+1}^{\dagger} } =  \norm{ T_k^{-1} F'[f_k]^\ast \left( g^{\Textbf{err}}- g^{\dagger} \right) }_{\mY} \leq  \frac{ \norm{F'[f_k]^\ast} }{\alpha_k} \norm{\Textbf{err}}_{\mY}.
\end{equation}
Hence, a bounded \emph{data error} induces bounded deviations of the reconstruction where it should be emphasized that the error estimate deteriorates in the limit $\alpha_k \to 0$.

On the other hand, minimizers of \eqref{eq:GenOpQuadMin} are in general no exact least-square solutions to the original problem \eqref{eq:GenOpLin} but deviate by an \emph{approximation error} growing with $\alpha_k$ due to the balancing of the data residual with the regularization term.
In order to minimize the \emph{total} reconstruction error it is thus necessary to balance data- and approximation errors by suitable \emph{parameter choice rules}:
\begin{itemize}
 \item Choose $\alpha_0 $ large enough to preclude excessive step sizes in the initial iterates
  \item Define $ \{ \alpha_k\}_{k \in \mN_0} $ to be monotonically decreasing
  \item Stop the iterations at $k = k_{\Text{stop}}$ before the data error becomes dominant
\end{itemize}
Details depend on the specific operator $F$, exact solution $f^\dagger$ and expected data errors. One strategy for the choice of $k_{\Text{stop}}$, going back to \citet{Morozov1966Discrepancy}, is given by the \emph{discrepancy principle}. This parameter rule is defined by
\begin{equation}
k_{\Text{stop}}  := \min\{ k \in \mN: \norm{ F(f_k) - g^{\Textbf{err}} }_{\mY} \leq \tau \norm{\Textbf{err}} _{\mY} \} \MTEXT{with fixed} \tau \geq 1. \label{eq:Discrepancy}
\end{equation}
 By implementing \eqref{eq:Discrepancy}, the Newton-iterations are thus stopped as soon as the nonlinear residual reaches the order of the error level. This is reasonable because a further reduction of the residual need not yield a better approximation of the \emph{exact} data $g^\dagger = g^{\Textbf{err}} - \Textbf{err}$, whereas smaller regularization parameters $\alpha_k$ would reduce stability to data errors as discussed above. For details and further parameter choice rules, see for instance \cite[C. 4]{Hanke1996Regularization}.

To conclude, we remark that the IRGNM with a suitable stopping rule indeed defines a \emph{regularization method} for nonlinear inverse problems of the form \eqref{eq:GenOpEq}. In particular, the final iterate converges to the exact solution $f^\dagger$ for $\Textbf{err} \to 0$ under reasonable assumptions \cite{Blaschke1997IRGNMConv}. Moreover, explicit convergence rates can be shown, given source conditions for $f^\dagger$ and bounds for the nonlinearity of $F$ \cite{Bakushinskii1992IRNM,Bakushinsky1994IllposedProbs,Blaschke1997IRGNMConv}. However, verifying these assumptions for phase contrast tomography can be expected to turn out cumbersome which is why convergence analysis is omitted in this work.

\end{subsection}
\vspace{-.5em}

\begin{subsection}{Generalized Newton-Type Methods} \label{SS:IRNM-Gen}

By construction, the regularization term in \eqref{eq:GenOpRegSol} limits the deviations of the iterates from the initial guess $x_0$. The choice of the norm in $\mX$ along with $f_0$ thereby allows to impose desirable properties. If $\mX \ni f_0 $ is for instance given by some Sobolev space (compare \sref{S:Sobolev}), then the bounded deviations imply in particular $f_k \in \mX$ for all $k$ ensuring a prescribed regularity of the iterates. However, many desirable constraints like positivity of reconstructed functions may not be imposed by Hilbert space norms. Therefore, it is reasonable to relax the setting of \sref{SS:IRNM-Motivation} to Banach spaces $\mX$ and $\mY$ and consider general \emph{penalty functionals}
\begin{equation} \cH: \mX \to \mR \cup \{ \infty \} \end{equation}
as regularization terms. Likewise, generalized \emph{data fidelity functionals}
\begin{equation} \cS(g \, ; \; \cdot): \mY \to \mR \cup \{ \infty \} \end{equation}
may provide a more meaningful measure for how well the reconstruction explains the observed data. This leads to generalized Newton methods of the form \cite{Hohage2013}
\begin{equation}
 f_{k+1} = \argmin_{f \in \mX} \left( \cS\left(g^{\Textbf{err}} \, ; \; F(f_k) +  F'[f_k](f - f_k) \right) + \alpha_k \cH(f) \right) \label{eq:GenIRNM}
\end{equation}
The minimizers are unique if $\cH$ and $\cS(g^{\Textbf{err}}\, ; \; \cdot)$ are convex and lower semi-continuous and if either of these is strictly convex. Note, however, that the solution of \eqref{eq:GenIRNM} may in general not be expressed in closed form - as achieved in the case of the IRGNM (see \lemref{lem:RegIRGNM}) - but requires general convex optimization techniques.

\end{subsection}


\end{section}

\begin{section}{Application to Phase Contrast Tomography} \label{S:RecMethodPCT}

In this section, we apply the idea of regularized Newton methods, outlined in \sref{S:IRNM} in an abstract form, to derive reconstruction methods for the inverse problem of phase contrast tomography (\probref{prob:1}).

\begin{subsection}{Basic Reconstruction Method}\label{SS:RecMethodBase}
 
 In \thmref{thm:ForwardOpsFrechet}, \Frechet differentiability of the near-field- and far-field forward operators 
 \begin{equation*} F_d, F_\infty: \mX \to \mY ; \; N \mapsto I_{\ast} \end{equation*}
 has been shown on Banach spaces $\mX = \Lp \infty {\ObjDom}_{\mR}$ and $\mY = \Lp \infty {Z^{m+1}}_{\mR}$. We adopt the generalized version of regularized Newton methods introduced in \sref{SS:IRNM-Gen} for the construction of Newton iterations of the form \eqref{eq:GenIRNM}. Thus, what is left is the choice of suitable data fidelity- and penalty functionals $\cS, \cH$ such that
 \begin{itemize}
  \item[$\boldsymbol 1$] Iterates may be computed stably and efficiently by evaluating \eqref{eq:GenIRNM} and...
    \item[$\boldsymbol 2$] provide good approximations of the object to be reconstructed.
 \end{itemize}

 \begin{subsubsection}{Choice of the Data Misfit Functional}
  
  According to \eqref{eq:GenIRNM}, the data fidelity functional $\cS$ compares the data corresponding to the current reconstruction, approximated by the linearized forward operator, with the observations. The overall aim is \emph{not} an exact match with the imperfect and noisy  measurements $I^{\Textbf{err}} = I^\dagger + \Textbf{err}$, but a faithful approximation of the unknown exact data $I^\dagger$. Accordingly, a good choice of $\cS$ needs to take into account the \emph{statistics} of the expected errors $\Textbf{err}$ in order to provide an accurate measure for closeness to the true solution. For stochastic errors, a canonical choice is the \emph{negative log-likelihood}
  \begin{equation}
   \cS(I^{\Textbf{err}}\, ; \; I) := - \log \Textbf{P} ( I^{\Textbf{err}} | I ) + \Text{const} \label{eq:NegLogLikeli} .
  \end{equation}
  where $\Textbf{P} ( I^{\Textbf{err}} | I ) $ denotes the conditional probability of measuring $I^{\Textbf{err}} $ given the exact data is $I$ \cite{Hohage2013}.
  
  In \emph{near-field} phase contrast tomography, the observed intensities are usually so large and uniform over the detector area that the statistical errors, arising from fluctuations in the number of incident photons  plus  instrument noise, can be modeled as additive Gaussian errors $\Textbf{err}$. For this type of noise, \eqref{eq:NegLogLikeli} suggests $L^2$-data fidelity functionals, i.e.\ a suitable choice for solution of \probref{prob:1} in the near-field setting governed by the operator $F_d$ in \thmref{thm:ForwardOpsFrechet} is given by
  \begin{equation}
   \cS(I^{\Textbf{err}}\, ; \; I)  := \norm{I - I^{\Textbf{err}}}_{\Lp 2 {Z^{m+1}}_{\mR}}^2 \MTEXT{for} I,I^{\Textbf{err}} \in \mY  \label{eq:DataMisfitNF} .
  \end{equation}
  
  In the \emph{far-field} case, the detected intensities are typically much smaller and of much greater lateral variation so that it becomes significant that the radiation is actually \emph{quantized} into single incident photons. Detector pixels accordingly count discrete uncorrelated events over some integration time. The resulting probability distribution of the counts is given by \emph{Poisson statistics}
  \begin{equation}
   \Textbf{P} ( I^{\Textbf{err}}_j | I_j ) = \exp(- I_j) \frac{ I_j^{I^{\Textbf{err}}_j} }{ I^{\Textbf{err}}_j !} \label{eq:PoissStats}.
  \end{equation}
  where $ I^{\Textbf{err}}_j$ denotes the number of counts measured at a pixel $j$ and $I_j$ the exact local intensity. For the small intensities observed in far-field imaging, the fluctuations of the counts constitute the dominant statistical errors. Accounting for probabilities of the form \eqref{eq:PoissStats} in \eqref{eq:NegLogLikeli}, the obtained data fidelity functional is the \emph{Kullback-Leibler-Divergence} (cf. \cite[sec. 5.3]{Hanke1996Regularization}, \cite{Hohage2013}), which for $ I,I^{\Textbf{err}} \in \mY $ is given by
  \begin{equation}
    \KL (I^{\Textbf{err}}\, ; \; I) := \int_{Z^{m+1}} \left(I -  I^{\Textbf{err}} - I^{\Textbf{err}} \ln \frac{I}{I^{\Textbf{err}}} \right) \; \D x \D \theta \label{eq:KL}.
  \end{equation}
  In \eqref{eq:KL}, the conventions $\ln(x) = \infty$ for $x < 0$, $\ln(\frac{x}{0}) = \infty$ for $x > 0 $  and $0 \cdot \ln(\frac{0}{x}) = 0$ for $x \geq 0$ are adopted to ensure well-definedness.
  $\mK\mL (I^{\Textbf{err}}\, ; \; \cdot )$ defines a convex, lower semi-continuous functional with a global minimum at $ I = I^{\Textbf{err}}$.
  
  In order to simplify the implementation, we expand the integrand in \eqref{eq:KL} to quadratic order in $I$ about the observations $I^{\Textbf{err}}$. This yields
  \begin{equation}
    \KL (I^{\Textbf{err}}\, ; \; I) = \int_{Z^{m+1}}  \left(  \frac{( I - I^{\Textbf{err}} )^2}{2 I^{\Textbf{err}}} + \Or( ( I - I^{\Textbf{err}} )^3 ) \right) \; \D x \D \theta  \label{eq:KLQuadApprox}
  \end{equation}
  The quadratic term in \eqref{eq:KLQuadApprox} diverges wherever $I^{\Textbf{err}} = 0$ vanishes, corresponding to pixels where zero photon counts have been recorded that are weighted by infinity. We account for this problem by choosing a relaxed version of \eqref{eq:KLQuadApprox} as the data fidelity functional for far-field phase contrast tomography
  \begin{equation}
   \cS (I^{\Textbf{err}}\, ; \; I)  := \norm{\cG_{\mY}^{\frac 1 2 } (I - I^{\Textbf{err}})}_{\Lp 2 {Z^{m+1}}_{\mR}}^2 \MTEXT{with} \cG_{\mY} I := \frac{I}{ 2 \max( I^{\Textbf{err}}, I_{\min} ) } \label{eq:DataMisfitFF} .
  \end{equation}
  Hence, we arrive at a similar expression as \eqref{eq:DataMisfitNF} where the positive-semidefinite \emph{Gramian} $\cG_{\mY}:  \Lp 2 {Z^{m+1}}_{\mR} \to \Lp 2 {Z^{m+1}}_{\mR}$ in the near-field case is simply the identity.
  
 \end{subsubsection}

  \begin{subsubsection}{Choice of the Penalty Functional}
  
     The classical choice for the penalty functional $\cH$ in \eqref{eq:GenIRNM} is $L^2$-regularization, given by $\cH(N) := \norm{N - N_0 }_{\Lp 2 {\ObjDom}_{\mR}}^2$ for the compactly supported objects $N \in \Lp \infty {\ObjDom}_{\mR}$ considered as admissible reconstructions in \probref{prob:1}. Here, we allow for somewhat more general regularization terms of the form
      \begin{equation}
       \cH(N) := \norm{\cG_{\mX}^{\frac 1 2 }  ( N - N_0 ) }_{\Lp 2 {\mR^{m+1}}_{\mR}}^2 \label{eq:PenaltyPCT}
      \end{equation}
      where the Gramian $\cG_{\mX} : \Lp 2 {\mR^{m+1}}_{\mR} \to \Lp 2 {\mR^{m+1}}_{\mR}$ is assumed to be bounded, self-adjoint and uniformly positive-definite, i.e.\ for some $\varepsilon > 0$
      \begin{equation*} \ip{N}{\cG_{\mX} N}_{\Lp 2 {\mR^{m+1}}_{\mR}} \geq \varepsilon \norm{N}_{\Lp 2 {\mR^{m+1}}_{\mR}}^2 \MTEXT{for all} N \in \Lp 2 {\mR^{m+1}}_{\mR}. \end{equation*}
      This implies that the square root $\cG_{\mX}^{\frac 1 2 }$ is well-defined and that $\cG_{\mX}$ is boundedly invertible. $L^2$-regularization simply corresponds to choosing $\cG_{\mX}$ as the identity.
      Concrete regularization terms of the form \eqref{eq:PenaltyPCT} are introduced in \sref{SS:Constraints}.

   \end{subsubsection}

     \begin{subsubsection}{Construction of the Newton-Iterates}

	Having derived suitable penalty- and data fidelity functionals applicable to the inverse problem of phase contrast tomography, we are finally in a position to define the corresponding Newton-steps. Let $N_0 \in \Lp \infty {\ObjDom}_{\mR}$ denote the initial guess for the object to be reconstructed and let $F\in \{F_d, F_\infty \}$ be the near-field or far-field forward operator characterized in \thmref{thm:ForwardOpsFrechet}. Assume that the observed intensities satisfy $I^{\Textbf{err}} - F (N_0) \in \Lp 2 {Z^{m+1}}_{\mR}$. By inserting \eqref{eq:PenaltyPCT} and \eqref{eq:DataMisfitFF} into \eqref{eq:GenIRNM}, a Newton-step for the solution of \probref{prob:1} is obtained as
	\begin{align}
	  N_{k+1} = \argmin_{N \in \tilde \mX} ( \;\;\;\;\; \;\;\; &\norm{\cG_{\mY}^{\frac 1 2 }(F (N_k) +  F' [N_k](N - N_k) -  I^{\Textbf{err}}) }_{\Lp 2 {Z^{m+1}}_{\mR}}^2 \nonumber \\ 
	  +  \; \alpha_k &\norm{\cG_{\mX}^{\frac 1 2 } (N - N_0)}_{\Lp 2 {\mR^{m+1}}_{\mR}}^2  \;\;\;\;\; \;\;\;\;\;\;\;\; \;\;\;\;\;\;\;\; \;\;\;\;\;\;\;\; \;\;\;\;\;\; \;\;\;\;    ). \label{eq:NewtonQuadProbPCT}
	\end{align}
	
	Recall that $\cG_{\mY}$ is simply chosen as the identity in the near-field case. Moreover, note that the set of admissible solutions $\tilde \mX$ is not specified in \eqref{eq:NewtonQuadProbPCT}. For the domain of the nonlinear forward operators $\tilde \mX = \Lp \infty {\ObjDom}_{\mR}$, minimizers might indeed not exist.
	On the other hand, recall that the \Frechet derivatives $F' [N_k]$ have a unique extension to $\Lp 2 {\ObjDom}_{\mR}$ as proven in \thmref{thm:ForwardOpsFrechet}. For $\tilde \mX = \Lp 2 {\ObjDom}_{\mR}$, \eqref{eq:NewtonQuadProbPCT} resembles the Hilbert space setting of the IRGNM in \sref{SS:IRGNM}. Thus, in analogy to \lemref{lem:RegIRGNM}, we obtain that the Newton iterate $N_{k+1}$ can be stably computed by solving a self-adjoint positive-definite linear problem:
	\vspace{1em}
	\begin{th1}[Newton Step for Phase Contrast Tomography] \label{thm:NewtonStepPCT}
	  Let $\tilde \mX = \Lp 2 {\ObjDom}_{\mR}$, $\alpha_k > 0 $, $N_0, N_k \in \Lp \infty {\ObjDom}_{\mR}$ and $  I^{\Textbf{err}} - F(N_k)  \in \Lp 2 {Z^{m+1}}_{\mR}$. Then the quadratic minimization problem \eqref{eq:NewtonQuadProbPCT} has the unique solution
	\begin{align}
	N_{k+1} = N_k +  (\cG_{\mX}^{-1} F'[N_k]^\ast \cG_{\mY} F'[N_k] + \alpha_k)^{-1} ( \,  \cG_{\mX}^{-1} F'[N_k]^\ast \cG_{\mY} (I^{\Textbf{err}} - F(N_k))& \nonumber \\
	+ \; \alpha_k(N_0 - N_k)& \, )  \label{eq:NewtonStepPCT}
	\end{align}
	where $F'[N_k]^\ast$ denotes the adjoint of the extension $F'[N_k]: \Lp 2 {\ObjDom} \to \Lp 2 {Z^{m+1}}$. Moreover, $N_{k+1} \in  \Lp 2 {\ObjDom}_{\mR} $ depends continuously on the data $I^{\Textbf{err}}$.
	\end{th1}
	\begin{pf}
	 Equipping $\mX := \Lp 2 {\ObjDom}_{\mR}$ and $\mY   := \Lp 2 {Z^{m+1}}_{\mR}$ with the inner products
	 \begin{equation*}
	  \ip{f_1}{f_2}_{\mX}  := \ip{f_1}{\cG_{\mX} f_2}_{\Lp 2 {\mR^{m+1}}_{\mR}} \MTEXT{and} \ip{g_1}{g_2}_{\mY}  := \ip{g_1}{\cG_{\mY} g_2}_{\Lp 2 {Z^{m+1}}_{\mR}}
	 \end{equation*}
	 and identifying $f_j = N_j$, $ g^{\Textbf{err}} = I^{\Textbf{err}}$, \eqref{eq:NewtonQuadProbPCT} can be brought to the form considered in \lemref{lem:RegIRGNM}. Hence, a unique minimizer $N_{k+1} \in \mX$ exists and is given by
	 \begin{equation*}
	N_{k+1} = N_k +  (F'[N_k]^{\star} F'[N_k] + \alpha_k)^{-1} \left(F'[N_k]^{\star} (I^{\Textbf{err}} - F(N_k)) + \alpha_k(N_0 - N_k) \right)
	\end{equation*}
	where $F'[N_k]^\star: \mY \to \mX $ denotes the adjoint of $F'[N_k]$ with respect to the $\mX$- and $\mY$-inner products. By construction, we have for all $f \in \mX$, $g \in \mY$
	\begin{align*}
	 \ip{f }{F'[N_k]^\star g }_{\mX} &= \ip{F'[N_k] f }{ g }_{\mY} = \ip{ F'[N_k] f }{ \cG_{\mY} g }_{\Lp 2 {Z^{m+1}}_{\mR}}  = \ip{  f }{ F'[N_k]^\ast \cG_{\mY} g }_{\Lp 2 {\ObjDom}_{\mR}} \\
	 &= \ip{  f }{ \cG_{\mX}^{-1} F'[N_k]^\ast \cG_{\mY} g }_{\mX},
	\end{align*}
	i.e.\ $F'[N_k]^\star  = \cG_{\mX}^{-1} F'[N_k]^\ast \cG_{\mY}$. Inserting this into the derived expression for the iterate $N_{k+1}$, we obtain the Newton step \eqref{eq:NewtonStepPCT} in the claim.
	
	 By \lemref{lem:RegIRGNM}, the map $(I^{\Textbf{err}} - F(N_k)) \mapsto N_{k+1}$ is continuous with respect to the $\mX$- and $\mY$-norms. As the embeddings $\Lp 2 {Z^{m+1}}_{\mR} \hookrightarrow \mY$ and $\mX \hookrightarrow\Lp 2 {\ObjDom}_{\mR}$ are bounded according to the properties of $\cG_{\mX}$ and $\cG_{\mY}$ stated in the preceding paragraphs, the dependence on the data is likewise continuous in $L^2$-norm.
	\end{pf}
	\vspace{1em}
	
	A major drawback of the constructed Newton step \eqref{eq:NewtonStepPCT} is that the computed iterates $N_{k+1} \in\Lp 2 {\ObjDom}_{\mR}$ need not be in the domain $  \Lp \infty {\ObjDom}_{\mR}$ of the nonlinear forward operators $F\in \{F_d, F_\infty \}$ even if $N_{k } \in \Lp \infty {\ObjDom}_{\mR}$. Therefore, the sequence of Newton-steps is not necessarily well-defined. While this technical peculiarity can be expected to obstruct any rigorous analysis significantly, it supposedly does not lead to numerical instabilities in the reconstruction, as $L^2$- and $L^\infty$-spaces coincide in a \emph{discretized}, i.e.\ finite-dimensional setting.
	
	According to \thmref{thm:NewtonStepPCT}, a final ingredient for the solution of phase contrast tomography by regularized Newton methods is given by the \emph{adjoint} of the \Frechet derivatives. These are derived in the following theorem:
	\vspace{1em}
	\begin{th1}[Adjoints of the Forward Operators] \label{thm:AdjFrechetPCT}
	  For $N \in \Lp \infty {\ObjDom}_\mR$, let $F_d'[N], F_\infty'[N]: \Lp 2  {\ObjDom}_\mR \to \Lp 2  {Z^{m+1}}_\mR$ denote the extensions of the \Frechet derivatives in \thmref{thm:ForwardOpsFrechet}. For $g \in \Lp 2  {Z^{m+1}}_\mR$,  their adjoints are given by
	  \begin{subequations} \label{eq:AdjForwOps}
	    \begin{align}
		F_d'[N]^\ast g &= 2k^2 \CR^\ast  \{  \;\; \cc{\left[ P \cdot \exp \left( -\I k \CR(N) \right) \right ]} \nonumber \\
		 &\;\;  \;\;\;\;\;\;\;\;\;\;\;\;\;\; \cdot \bDFlat{-d} \left( \left[\bDFlat{d} \left( P \cdot  \exp \left( -\I k \CR(N) \right)  \right) \right] \cdot \Re(g) \right) \; \} \label{eq:AdjForwOpNF} \\
		F_\infty'[N]^\ast g &= 2k^2 \CR^\ast  \{  \;\; \cc{\left[ P \cdot \exp \left( -\I k \CR(N) \right) \right ]} \nonumber \\
		 &\;\;  \;\;\;\;\;\;\;\;\;\;\;\;\;\; \cdot \cF_{\overline 2}^\ast \left( \left[\cF_{\overline 2} \left( P \cdot \left[ \exp \left( -\I k \CR(N) \right) - 1 \right] \right) \right] \cdot \Re(g) \right) \; \} \label{eq:AdjForwOpFF}
	  \end{align}
	  \end{subequations}

	\end{th1}
	\begin{pf}
	According to \eqref{eq:FrechetForwOpNF}, the Frechet derivatives of $S \in \{ F_d'[N], F_\infty'[N] \}$ are of the form
	\begin{equation*} S = 2k^2 S_3 \circ \cM_{f_2} \circ S_2 \circ \cM_{f_1} \circ S_1 \end{equation*}
	with $S_1 = \CR$, $S_2 \in \{\bDFlat{d}, \cF_{\overline 2} \}$ and $S_3 = \Re$. By the properties of the adjoint given in \thmref{thm:AdjointProps}, this implies
	\begin{equation} S^\ast = 2k^2  S_1^\ast \circ \cM_{\cc{f_1}} \circ S_2^\ast \circ \cM_{\cc{f_2}} \circ S_3.  \label{eq:AbstractAdjPCT} \end{equation}
	Here, we have used that $S_3 = \Re: \Lp 2  {Z^{m+1}}_\mR \to \Lp 2  {Z^{m+1}}_\mR$ is self-adjoint according to \exref{ex:Adjoints}-(c) and that the adjoints of the multiplication operators $\cM_{f_1}, \cM_{f_2}$ simply amount to multiplications with the complex conjugate factors $\cc{f_1}, \cc{f_2}$. Noting furthermore that $\left(\bDFlat{d}\right)^\ast = \bDFlat{-d}$, i.e.\ that the adjoint (and inverse) of the Fresnel propagator corresponds to back-propagation, we obtain the expressions \eqref{eq:AdjForwOpNF} and \eqref{eq:AdjForwOpFF} by substituting the partial operators into \eqref{eq:AbstractAdjPCT}.
	\end{pf}
	\vspace{1em}
	
	With the explicit expressions for the \Frechet derivatives and their adjoints, given in \thmref{thm:AdjFrechetPCT} and \thmref{thm:ForwardOpsFrechet}, the regularized Newton method defined by \eqref{eq:NewtonStepPCT} allows for reconstructions in phase contrast tomography, yielding approximate solutions of \probref{prob:1}. We emphasize once more that this algorithm corresponds to simultaneous phase retrieval and Radon inversion as the solved linearized problems incorporate the complete tomographic data and yield approximations of the unknown object $N^\dagger$ itself - instead of merely projections from which tomographic reconstruction would have to be computed a posteriori. As argued in \sref{S:CombinedApproach}, this can be expected to yield improved reconstruction results compared to sequential implementations of the subproblems.

   \end{subsubsection}

\end{subsection}

    \begin{subsection}{A Priori Constraints} \label{SS:Constraints}
    
    We have seen in \sref{S:PhaseRetrieval} that imposing a priori knowledge on the reconstructed object may facilitate phase retrieval. In the commonly used alternating-projection-type algorithms (compare \sref{S:CombinedApproach}) imposing additional constraints is fairly simple as these just correspond to yet another projection on a further constraint set. In the following, we discuss in which manner the basic regularized Newton method constructed in \sref{SS:RecMethodBase} may be similarly supplemented to incorporate a priori knowledge.
    
    \begin{subsubsection}{Regularity Constraints}

	Imposing regularity of the solution may  suppress noise in the reconstruction and - according to \thmref{thm:UniquePhaseFFMD} - possibly promote unique phase retrieval. A straightforward approach to do so is by choosing the penalty functional in \eqref{eq:GenIRNM} as the squared norm of a suitable Sobolev space $H^s(\mR^{m+1})$ for $s \geq 0$, see \sref{S:Sobolev}. Comparing the definition in \eqref{eq:DefHs} to \eqref{eq:PenaltyPCT}, it is found such constraints are implemented by the Gramian
	\begin{equation}
	 \cG_\mX (N) = \cF^\ast \left( (1 + \norm{\bxi}_2^2)^s \cdot \cF(N) \right). \label{eq:GramSobolev}
	\end{equation}
	Note that $\cG_\mX$ is self-adjoint and strictly positive but $\cG_{\mX}(N) \in \Lp 2 {\mR^{m+1}}$ only holds for the dense subspace $H^{2s}(\mR^{m+1}) \subset \Lp 2 {\mR^{m+1}}$. However, its inverse
	\begin{equation} \cG_{\mX}^{-1} : \Lp 2 {\mR^{m+1}} \to \Lp 2 {\mR^{m+1}}; \; N \mapsto \cF^{-1} \left( (1 + \norm{\bxi}_2^2)^{-s} \cdot \cF(N) \right) \label{eq:GramSobolevInv} \end{equation} 
	is well-defined and bounded. As only the latter needs to be evaluated in the Newton step \eqref{eq:NewtonStepPCT}, regularity constraints of given order $s \geq 0$ may be imposed by equipping the basic Newton method in \thmref{thm:NewtonStepPCT} with the Gramian in \eqref{eq:GramSobolevInv}.
  
   \end{subsubsection}

    \begin{subsubsection}{Support Constraints}

     In some cases, the support of the specimen to be reconstructed may be known more specifically than in the form of the cylindrical domain $\ObjDom$ defined in \eqref{eq:ObjDomain}. A known support in $\SuppDom \subset \ObjDom$ may be accounted for by choosing $\tilde \mX = \Lp 2 {\SuppDom}$ in \thmref{thm:NewtonStepPCT} and considering the modified forward operators
     \begin{equation}
      \tilde F = F \circ \iota : \Lp \infty {\SuppDom}_{\mR} \to \Lp \infty {Z^{m+1}}_{\mR} \label{eq:SuppConstFwOp}
     \end{equation}
     supplemented with the canonical embedding $\iota: \Lp 2 {\SuppDom}_{\mR} \hookrightarrow \Lp 2 {\ObjDom}_{\mR}$. The latter is linear and bounded with $\norm{\iota} = 1$ (both in $L^2$ and $L^\infty$) and corresponds to an extension of functions $N \in \Lp 2 {\SuppDom}_{\mR}$ with 0 in $\ObjDom \setminus \SuppDom$. The resulting \Frechet derivative is
    \begin{equation}
      \tilde F'[N] = F'[\iota(N)] \circ \iota. \label{eq:SuppConstDeriv}
     \end{equation}
     As $\Lp 2 {\SuppDom}_{\mR}$ defines a closed subspace in $\Lp 2 {\ObjDom}_{\mR}$, the adjoint of $\iota$ equals the orthogonal projection $\cP_{ \SuppDom}$ onto $\Lp 2 {\SuppDom}_{\mR}$ as seen in \exref{ex:Adjoints}-(a). Thus,
     \begin{equation}
      \tilde F'[N]^\ast = \cP_{ \SuppDom}\circ F'[N]^\ast. \label{eq:SuppConstAdj}
     \end{equation}
     Accordingly, support constraints can be incorporated into the regularized Newton iterations in \eqref{eq:NewtonStepPCT} by restricting the set of admissible objects to $\Lp 2 {\SuppDom}$ and projecting onto this space after each evaluation of the adjoint $ F'[N]^\ast$. Moreover, for a function $N \in \Lp 2 {\ObjDom}_{\mR}$, these projections simply amount to setting $N = 0$ outside the support $\SuppDom$.

   \end{subsubsection}
   
   \begin{subsubsection}{Non-Absorbing and Single-Material Objects}

   As discussed in \sref{SS:SpecialObjects}, many specimen of interest give rise to negligible absorption, i.e.\ are described by a real-valued refractive index $n= 1 - N = 1-\delta$, or more generally satisfy the single-material approximation of a fixed ratio between absorption $\beta$ and refraction $\delta$. We account for this within the framework of \thmref{thm:NewtonStepPCT}, by introducing modified forward operator $ \tilde F = F \circ \iota_{c|L^\infty_{\mR} (\ObjDom) }$ composed with
   \begin{equation}
    \iota_c : L^2_{\mR} (\ObjDom) \to \Lp 2 {\ObjDom}_{\mR}; \; N_{\mR} \to c  N_{\mR}
   \end{equation}
   and setting the set of admissible objects to all \emph{real-valued} $L^2$-functions  $  L^2_{\mR} (\ObjDom) = \tilde \mX$. Note that $ L^2_{\mR} (\ObjDom)$ is a closed subspace of $\Lp 2 {\ObjDom}_{\mR}$ and $\iota_c $ is the corresponding canonical embedding, scaled with the factor $c \in \mC$ which defines the $\delta$-$\beta$-ratio of the object.
   
   Hence, the mathematical structure of the modified forward operator is widely identical to the case of support constraints discussed in the above paragraph. In particular, we obtain for the modified \Frechet derivative and adjoint
   \begin{equation}
     \tilde F'[N] = F'[\iota_c(N)] \circ \iota_c \MTEXT{and}  \tilde F'[N]^\ast = \Re \circ \cc c F'[\iota_c(N)]^\ast \label{eq:RealConstAdj}.
   \end{equation}
    Here,  it has been used that the orthogonal projection in $\iota_c^\ast = \cP_{L^2_{\mR} (\ObjDom)} \circ \cc{c}$ is simply the point-wise real part as can be seen from the properties $\Re( \Lp 2 {\ObjDom}_{\mR}  ) =  L^2_{\mR} (\ObjDom)$ and $\Re^\ast = \Re = \Re \circ \Re$ shown in \exref{ex:Adjoints}-(c).
    
    By introducing the expressions \eqref{eq:RealConstAdj} into the Newton step \eqref{eq:NewtonStepPCT}, the reconstruction is thus restricted to single-material objects characterized by $N = c  N_{\mR}$ for $N_{\mR}$ real-valued. Furthermore, note that this assumption may be easily combined with a support constraint since the corresponding orthogonal projections commute, meaning that the constraints are perfectly compatible.

   \end{subsubsection}

   \begin{subsubsection}{Positivity Constraints}
   
   As discussed in \sref{SS:XrayRefrIdx}, the refractive index $n = 1 -\delta + \I \beta $ in the hard X-ray regime typically satisfies $\delta, \beta \geq 0$. This motivates a restriction of the space of admissible functions in our regularized Newton method to
   \begin{equation}
    \tilde \mX  = \{ N \in \Lp 2 {\ObjDom}_{\mR}: \Re(N), -\Im(N) \geq 0 \text{ a.e.}\} =:  C_+. \label{eq:PosCone}
   \end{equation}
    The corresponding projection of an object $N = \delta - \I \beta \in \Lp 2 {\ObjDom}_{\mR}$ onto this set, i.e.\ its best approximation in $C_+$ with respect to the $L^2$-norm, is given by 
    \begin{equation}
     \cP_+: \Lp 2 {\ObjDom}_{\mR} \to C_+; \; (\delta - \I \beta) \mapsto \max(\delta, 0) - \I \cdot \max(\beta, 0) \label{eq:PosProj}.
    \end{equation}
    
    Notably, $ \cP_+$ does \emph{not} define a linear projection according to \eqref{eq:PosProj}. Indeed, it can be seen that the map is neither \Frechet differentiable as the pointwise maximum $x \mapsto \max(x, 0)$ corresponds to a truncation which is not even differentiable in $\mR$.    
    The underlying reason for these peculiarities  lies in the different geometry of $C_+$ compared to the constraint sets considered above: as the set $C_+$ is closed only under multiplication with positive scalars, it does not form a linear subspace of $\Lp 2 {\ObjDom}_{\mR}$ but only a convex \emph{cone}. The nonlinear structure of this set gives rise to the non-smooth projection.
    
    Unfortunately, these observations imply that positivity constraints for $\delta$ and $\beta$ may not be incorporated into the derived regularized Newton method for phase contrast tomography, as it requires \Frechet differentiability of the involved operators. In principal, this could however be achieved by adopting the generalized approach of \emph{semismooth Newton methods}, as applied for instance in \cite{Griesse2008Semismooth,Hintermuller2002SSNewton}.

   \end{subsubsection}
  
   \end{subsection}

\end{section}

\begin{section}{Discretization} \label{S:Discretization}

    For a numerical implementation of phase contrast tomography via the regularized Newton method developed in \sref{S:RecMethodPCT}, we need to leave the infinite-dimensional description adopted so far and discretize the problem. The pursued strategy for this is outlined in the following section.

    \begin{subsection}{General Approach} \label{SS:DiscrGeneral}
    
    \subsubsection{Discrete Spaces}
    
    A discrete approximation of the objects $N \in \mX = \Lp \infty {\ObjDom}_{\mR}$ and corresponding intensity data $I = F(N)  \in \mY = \Lp \infty {Z^{m+1}}_{\mR}$ is obtained by sampling these quantities in the $m+1$ dimensions on equidistant cubic \emph{voxels} of edge length $\Delta x$. For convenience, we take the discretized objects $\bN$ not on a cylindrical domain but to be parametrized by a rectangular grid of voxels. This corresponds to the discrete object space
    \begin{equation}
     \bN \in  \mX_{\Text{dis}}  := \mC^{M_x \times M_y \times M_z} \subset \mX, \label{eq:DefXdisc}
    \end{equation}
    where the inclusion is to be understood by identifying the arrays with piecewise constant functions on a suitable cuboid of voxels contained in $\ObjDom$. Here, $M_x, M_y, M_z \in \mN$ denote number of grid points, i.e.\ the resolution in the different dimensions. Note that we focus on the physically relevant case of $m+1 = 3$ spatial dimensions. However, a two-dimensional toy model for phase contrast tomography of objects varying only in the tomographic plane of rotation, i.e.\ in the $x$- and $z$-directions, is readily obtained by setting $M_y =1$ in \eqref{eq:DefXdisc}.
    
    Discretization of the image space $\mY$ arises naturally owing to the fact that the intensities are measured by CCD detectors of finite aspect size, composed of a discrete number of $K_x \cdot K_y \in \mN$ pixels. Diffraction patterns can likewise only be recorded for a finite number of $K_\theta \in \mN$ different incident angles. Accordingly, the real-valued intensity measurements $ \bI$ can be identified with the space
    \begin{equation}
     \bI \in  \mY_{\Text{dis}}  := \mR^{K_\theta \times K_x \times K_y } \subset \mY.
    \end{equation}
    Once more, the two-dimensional toy model corresponds to the choice $K_y = 1$.

    \subsubsection{Discretization of the Operators}
    
    The interpretation of $\bN \in  \mX_{\Text{dis}}, \bI \in  \mY_{\Text{dis}}$ as piecewise constant functions in $\mX$ and $\mY$ induces a discretization of the forward operators $F \in \{F_d, F_\infty \}$ and \Frechet derivatives, as defined in \eqref{eq:ForwardOpNF}, \eqref{eq:ForwardOpFF} and \eqref{eq:FrechetForwOp}, in the following form:
    \begin{itemize}
     \item Interpret pointwise operations ($+$, $\cdot$, $\exp$, $\Re$, etc.) as componentwise on arrays
     \item Replace continuous Fourier transforms $\cF$ by \emph{fast Fourier transforms} (FFTs)
       \item Approximate line integrals in Radon transforms $\CR$ by weighted sums of voxels
    \end{itemize}
   On the discretized spaces $\mX_{\Text{dis}}$ and $\mY_{\Text{dis}}$, the inner products $\Lp 2 {\ObjDom}_{\mR}$ and $\Lp 2 {Z^{m+1}}_{\mR}$ take the form
    \begin{equation}
     \ip{\bN_1}{\bN_2}_{\mX_{\Text{dis}}} := \Re( \bN_1^\ast \cdot \bN_2 ) \MTEXT{and} \ip{\bI_1}{\bI_2}_{\mY_{\Text{dis}}} := \Re( \bI_1^\ast \cdot \bI_2 ). \label{eq:DefIPdisc}
    \end{equation}
    up to multiplicative constants. Accordingly, these are essentially given by Euclidean inner products so that adjoints in the Newton step \eqref{eq:NewtonStepPCT} can be simply be evaluated by applying the conjugate transpose of matrix representations of discretized forward operations (although doing so explicitly is rarely efficient). As discrete Fourier transforms are unitary up to a multiplicative constant with respect to Euclidean scalar products, their adjoints may be implemented via inverse FFTs.

    \subsubsection{Constraints and Error Metrics}
    
    Defining the \emph{support}  of a discretized object $\bN \in \mX_{\Text{dis}}$ as
    \begin{equation}
     \supp(\bN) := \{ \boldsymbol j \in \{1, \ldots, M_x \} \times \{1, \ldots, M_y \} \times \{1, \ldots, M_z \} : \bN_{\boldsymbol j } \neq 0 \},
    \end{equation}
    all $\bN  \in \mX_{\Text{dis}}$ supported in some subset of the voxel grid form a closed subspace $\tilde \mX_{\Text{dis}} \subset \mX_{\Text{dis}}$. The same is true for the set of all real-valued $\bN$. Hence, discrete support- and single-material constraints may be implemented exactly as outlined in \sref{SS:Constraints} via embeddings and orthogonal projections, i.e.\ extensions and truncations.
    
    Introducing specific Hilbert space data fidelity- and penalty terms is possible by suitably choosing the discrete Gramians $\cG_{\mX_{\Text{dis}}}$ and $\cG_{\mY_{\Text{dis}}}$ analogously to the infinite-dimensional setting considered in \sref{SS:RecMethodBase}. In the far-field case governed by Poisson-errors in the data, we choose the discrete version of \eqref{eq:DataMisfitFF} without the factor $\frac 1 2$
    \begin{equation}
     \cG_{\mY_{\Text{dis}}} \bI := \bI \oslash \max( \bI^{\Textbf{err}}, I_{\min} ), \;\;\;\; I_{\min}  > 0 \label{eq:GramYDisc}
    \end{equation}
    where $\oslash$ denotes component-wise division and $\bI^{\Textbf{err}} \in \mY_{\Text{dis}}$ are noisy observations. Regularization by Sobolev norms, motivated in \sref{SS:Constraints}, can be implemented in the discrete setting by the choice
    \begin{equation}
	\cG_{\mX_{\Text{dis}}} \bN = \FFT^\ast( (1 +  \norm{\bxi_{\Text{dis}}}^2)^s \odot \FFT(\bN) ), \;\;\;\; s \geq 0. \label{eq:GramXDisc}
	\end{equation}
     $\bxi_{\Text{dis}} \in \mC^{M_x \times M_y \times M_z}$ is the frequency array corresponding to the $m+1$-dimensional FFTs and $\odot$ denotes component-wise multiplication. $L^2$-regularization or $L^2$-data fidelity terms are implemented by taking $\cG_{\mX_{\Text{dis}}}$ or $\cG_{\mY_{\Text{dis}}}$ as the identity.

    \subsubsection{Resulting Algorithm}
    
    With the discretization of the forward map $F_{\Text{dis}}: \mX_{\Text{dis}} \to \mY_{\Text{dis}}$ outlined above, the evaluation of the Newton iterates \eqref{eq:NewtonStepPCT} reduces to solving the linear problem
    \begin{align}
	\cT_k( \bN_{k+1} - \bN_k) =  \cG_{\mX_{\Text{dis}}}^{-1} F_{\Text{dis}}'[\bN_k]^\ast \cG_{\mY_{\Text{dis}}} (\bI^{\Textbf{err}} - F_{\Text{dis}}(\bN_k)) + \alpha_k(\bN_0 - \bN_k) \label{eq:DiscNewtonStepPCT}
	\end{align}
    for a finite-dimensional self-adjoint positive-definite operator
    \begin{equation}
     \cT_k = \cG_{\mX_{\Text{dis}}}^{-1} F_{\Text{dis}}'[\bN_k]^\ast \cG_{\mY_{\Text{dis}}} F_{\Text{dis}}'[\bN_k] + \alpha_k \id_{\mX_{\Text{dis}}}. \label{eq:DiscNewtonLinProb}
    \end{equation}
    In order to exploit this form, the solution of \eqref{eq:DiscNewtonStepPCT} is computed by the iterative \emph{conjugate gradient method} (CG), see \cite[sec. 6.7]{Saad2003CG} for algorithmic details. Most importantly, only the stable \emph{forward} operations given by $F_{\Text{dis}}$, $F_{\Text{dis}}'[\bN]$ and $F_{\Text{dis}}'[\bN]^\ast$ as well as $ \cG_{\mX_{\Text{dis}}}^{-1}$ and $ \cG_{\mY_{\Text{dis}}}$ have to be evaluated explicitly for this approach. 
    
    The regularization parameter $\alpha_k > 0 $ is taken to be geometrically decreasing, i.e.
    \begin{equation}
     \alpha_k = r_\alpha^k \alpha_0 \MTEXT{for a fixed factor} r_{\alpha} \in (0;1).
    \end{equation}
    As discussed in \sref{SS:RegPar}, nonzero $\alpha_k$ in \eqref{eq:DiscNewtonStepPCT} give rise to approximation errors (even in the case of exact data) since only a perturbed form of the linearized inverse reconstruction problem is solved in each iteration. This is accounted for by implementing the CG-method such that the iterations are stopped as soon as the approximation accuracy corresponding to the current $\alpha_k$ and the error level $\Textbf{err}$ is reached. A detailed description of this approach can be found in \cite{Frommer1999TikhonovCG}.
    
    All in all, the discretized reconstruction method outlined in this section yields the following basic algorithm for our principal goal, the solution of \probref{prob:1}:
	\begin{alg1}[Regularized Newton Method for Phase Contrast Tomography]\label{alg:PCT} \ \\
		 \ \vspace{-.75em}  \\
	\begin{algorithm}[H]
	\SetKwInOut{KwInit}{Initialization}
	\KwData{Intensities $\bI^{\Textbf{err}} \in \mY_{\Text{dis}}$, setup $F  \in \{ F_d, F_\infty \}$, constraints $\tilde \mX_{\Text{dis}} \subset \mX_{\Text{dis}}$, initial guess $\bN_0 \in \tilde \mX_{\Text{dis}}$, Gramians $\cG_{\mX_{\Text{dis}}}, \cG_{\mY_{\Text{dis}}}$, reg. parameters $\alpha_0 > 0$, $r_{\alpha} \in (0;1)$, $k_{\Text{stop}} \in \mN$, stop rule $K: (\bN, \bI^{\Textbf{err}}, \Textbf{err}) \mapsto k$}
	 \ \vspace{-.75em}  \\
	\KwResult{Discrete approximation $\bN_{k_{\Text{stop}}}$ to the solution of \probref{prob:1}}
	    \ \vspace{-.75em}  \\
	    \KwInit{$\:\,\alpha = \alpha_0$\Text{;} \\$ \bN = \bN_0$\Text{;} }
	    \ \vspace{-.75em}  \\
	\For{$k= 0,1,\ldots, k_{\Text{stop}}$}{
	    \   \\
	    $ \;\;\:\, \cT   =\cG_{\mX_{\Text{dis}}}^{-1} F_{\Text{dis}}'[\bN ]^\ast \cG_{\mY_{\Text{dis}}} F_{\Text{dis}}'[\bN ] + \alpha \id_{\mX_{\Text{dis}}} $\Text{;} \\
	    $ \;\;\bN   \stackrel{\Text{CG}}=  \bN + \cT^{-1} \left( \cG_{\mX_{\Text{dis}}}^{-1} F_{\Text{dis}}'[\bN]^\ast \cG_{\mY_{\Text{dis}}} (\bI^{\Textbf{err}} - F_{\Text{dis}}(\bN)) + \alpha (\bN_0 - \bN) \right)$\Text{;} \\
	    \ \vspace{-.75em} \\
	    $\;\;\;\:\,\alpha =    r_\alpha \alpha $\Text{;}  \\
	    $k_{\Text{stop}} = K (\bN, \bI^{\Textbf{err}}, \Textbf{err})$\Text{;}  \\
		    \ \vspace{-.5em}  
	}
	\end{algorithm}
	\vspace{-.5em}
	\end{alg1}
    
   \end{subsection}

              \begin{subsection}{Implementation of the Propagators and Zero-Padding} \label{SS:ZeroPad}
          
          According to the general approach outlined in \sref{SS:DiscrGeneral}, the far-field and near-field propagators in \eqref{eq:ForwardOpFF} and \eqref{eq:ForwardOpNF} can be discretized in the form
          \begin{subequations} \label{eq:DiscProps}
          \begin{align}
           \cF_{\Text{dis}} ( \bpsi) &:= \FFT( \bpsi)  \label{eq:DiscPropFF} \\
           \bDF{d, \Text{dis}} (\bpsi)  &:= \FFT^{-1}\left( \exp\left( - \frac{ \I  \pi \bxi_{\Text{dis}} }{\NF} \right) \odot \FFT( \bpsi) \right) \label{eq:DiscPropNF}
          \end{align}
          \end{subequations}
          for some discrete contact image $\bpsi \in \mC^{J_x \times J_y}$. Here, we have employed the dimensionless  form of the near-field propagating factor from \sref{SS:FresnelFarfield}, governed by the Fresnel number $\NF$ as a single parameter, using the pixel size $\Delta x$ as the lengthscale.
       
    The discrete Fourier transform implicitly assumes a periodic continuation of the input signal. This periodicity may lead to severe artifacts when incorporated into the discretized propagators via \eqref{eq:DiscProps}, as is illustrated in \figref{fig:ZeroPad} for the near-field case: as a wave field is propagated, fringes may leave the computational field of view and will reappear on the opposite boundary. In order to suppress this non-physical effect, contact images $\bpsi \in \mC^{J_x \times J_y}$ are symmetrically \emph{zero-padded} prior to propagation, i.e.\ the propagators in \eqref{eq:DiscProps} are applied to an extended signal $\bpsi^{\Text{(pad)}} \in \mC^{J_{x}^{\Text{(pad)}} \times J_{y}^{\Text{(pad)}}}$ defined by
        \begin{equation}
     \bpsi_{ j_x,j_y }^{\Text{(pad)}} = \begin{cases}
                                  \bpsi_{  j_x - \left\lceil \frac{J_{x}^{\Text{(pad)}} - J_x}2 \right\rceil,j_y - \left\lceil \frac{J_{y}^{\Text{(pad)}} - J_y}2 \right\rceil} &\text{for } 1 +  \left\lceil \frac{J_{\ast}^{\Text{(pad)}} - J_\ast}2 \right\rceil \leq j_\ast \leq J_\ast +  \left\lceil \frac{J_{\ast}^{\Text{(pad)}} - J_\ast}2 \right\rceil \\
                                  0 &\text{else} \\
                                 \end{cases}
    \end{equation}
    Physically, this simulates additional free space around the contact image into which wave features may propagate without encountering periodic boundaries. From a mathematical perspective, zero-padding ensures that the analytical Fourier transforms in the governing forward operators, defined on the \emph{infinite} lateral domain $\mR^m$, are approximated sufficiently accurately by their (periodic) discretizations.
    
    For optimal computational efficiency of the FFTs, the padding sizes $J_{x}^{\Text{(pad)}}, J_{y}^{\Text{(pad)}}$ need to be chosen as a product of small primes, typically a power of two.
    In order to adapt the lateral resolution $K_x \leq J_{x}^{\Text{(pad)}}, K_y \leq J_{y}^{\Text{(pad)}}$ in image space $\mY_{\Text{dis}}$ to the recorded intensity data $\bI^{\Textbf{err}}$, the propagated padded wave fields $\cF_{\Text{dis}} (\bpsi^{\Text{(pad)}})$ or $\bDF{d, \Text{dis}} (\bpsi^{\Text{(pad)}})$  are symmetrically truncated. Retaining resolutions $K_x >J_x, K_y > J_x$  corresponds to \emph{oversampling} in the data, i.e.\ to a larger number of degrees of freedom image space than in object space, which may stabilize phase- and tomographic reconstruction. Zero-padding and truncation operations correspond to mutually adjoint linear operators that have to be incorporated into the discretized forward operators, derivatives and adjoints used in \algref{alg:PCT}.
                \begin{figure}[hbt!]
     \centering
    \subfloat[Without zero-padding]{\includegraphics[width=0.32\textwidth]{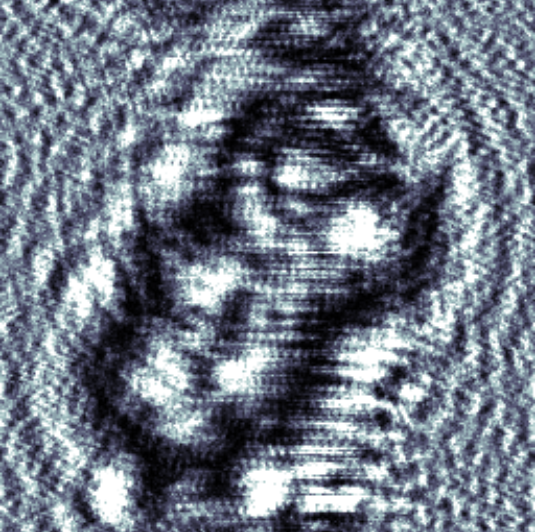}}
    \hfill
    \subfloat[Zero-padded to $512 \times 512$]{\includegraphics[width=0.32\textwidth]{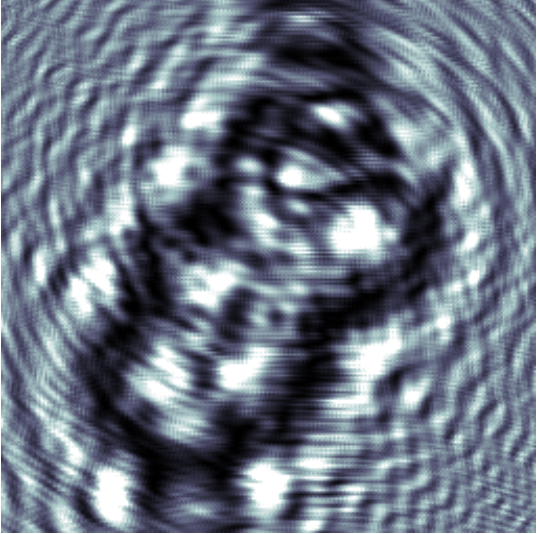}}
    \hfill
    \subfloat[Zero-padded to $1024 \times 1024$]{\includegraphics[width=0.32\textwidth]{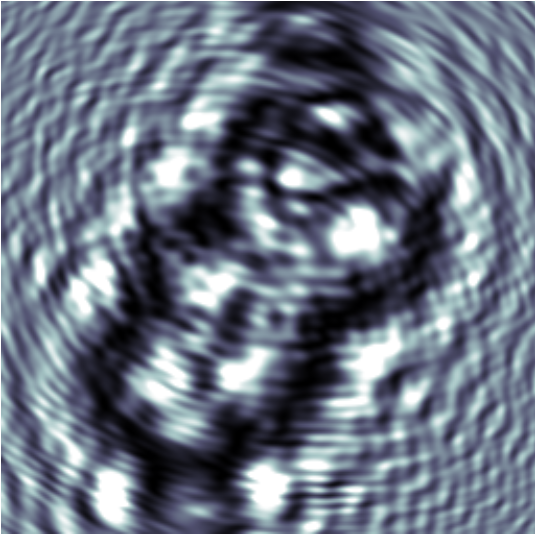}}
    \caption{Simulated near-field holograms for a $256\times 256$-sized section of the test object in \figref{fig:ContrastTransferDemo} at (pixel-length-based) Fresnel number $\NF = 10^{-3}$. From left to right: Fresnel propagation discretized according to \eqref{eq:DiscPropNF} with optional symmetric zero-padding by different factors and truncation of the propagated intensities to the original size. The hologram (a) computed without zero-padding shows severe artifacts caused by outgoing fringes that spuriously reenter the field of view by periodicity of the discrete Fourier transform. This effect is nearly eliminated in (b) owing to the greater computational domain simulated. \label{fig:ZeroPad}}
    \end{figure}

    \end{subsection}

       \begin{subsection}{Complexity and Implementation of the Radon Transform} \label{SS:RadonImplementation}
    
	  The numerical implementation of \algref{alg:PCT} boils down to the computation of a certain number of CG-iterations, each of which essentially requires an evaluation of the discrete \Frechet derivative $F_{\Text{dis}}'[\bN]$  and its adjoint $F_{\Text{dis}}'[\bN]^\ast$. The resulting computational complexity of our regularized Newton-type approach to phase contrast tomography is discussed in the following.
	  
	  For simplicity, we assume that the resolutions in the object- and image spaces $\mX_{\Text{dis}}  = \mC^{M_x \times M_y \times M_z}, \mY_{\Text{dis}}  = \mR^{K_\theta \times K_x \times K_y }$ are of the same order $M \in \mN$ in all dimensions. According to \eqref{eq:FrechetForwOp} and the discretization introduced in \sref{SS:DiscrGeneral} and \sref{SS:ZeroPad}, an evaluation of the derivative $F_{\Text{dis}}'[\bN]$ requires the following arithmetic operations:
    \begin{itemize}
     \item[$\boldsymbol 1$] $\Or(M^3)$: An order one number of componentwise operations 
     \item[$\boldsymbol 2$] $\Or(M^3 \log M)$: An evaluation of the discrete propagators in \eqref{eq:DiscProps}
     \item[$\boldsymbol 3$] $\gtrsim \Or(M^3 \log M)$: A discrete cylindrical Radon transform $\CRdis$
    \end{itemize}
    The complexity of $\Or(M^3 \log M)$ flops for the propagators result from their FFT-based implementation. The adjoints corresponding to $\boldsymbol 1$ and $\boldsymbol 2$ to be computed in the evaluation of $F_{\Text{dis}}'[\bN]^\ast$ are again given by componentwise operations and (inverse) FFTs, i.e.\ require the same number of arithmetic operations.
    If the total number of CG-iterations in the regularized Newton method are independent of the resolution $M$, which is empirically confirmed, the \emph{total} complexity of \algref{alg:PCT} is thus $\Or(M^3 \log M)$ up to the required evaluation of $\CRdis$ and its adjoint.
    
    As mentioned in \sref{SS:DiscrGeneral}, a straightforward idea for the discretization of the Radon transform is to approximate the integrals in \eqref{eq:Radon} by a sum over (bilinearly) interpolated voxel values along the corresponding lines through the grid. See \cite{Beylkin1987DiscreteRadon} for a detailed outline of this approach. This discretization strategy is also pursued in this work, using the standard implementations of the discrete 2D Radon transform $\cR_{\Text{dis}}$ provided by the numerical computing environments \textsc{Matlab} and \textsc{Octave} \cite{Octave}. If $\bN_j \in \mC^{M_x  \times M_z}$ for $j = 1, \ldots, M_y$ denote the slices of a 3D object $\bN \in \mX_{\Text{dis}}$, then its cylindrical Radon transform is obtained via
    \begin{equation}
     \CRdis(\bN)_j = \cR_{\Text{dis}}( \bN_j )  \MTEXT{for all} 1 \leq j \leq M_y,
    \end{equation}
    i.e.\ by slice-wise application of the 2D transform. Due to the independent summations of voxel values along lines in the numerical grid, each of the $\Or(M^3)$ values of $\CRdis(\bN)$ contributes $\Or(M)$ arithmetical operations, giving a total complexity of $\Or(M^4)$ flops for the evaluation of $\CRdis$.
    
    Accordingly, the Radon transforms to be evaluated in each CG-iterations typically constitute the performance-critical part of \algref{alg:PCT}. Its $\Or(M^4)$ complexity for the chosen discretization represents an algorithmic \emph{bottleneck} for the otherwise fast $\Or(M^3 \log M)$ implementation - at least for asymptotically large resolutions $M$. Yet, the required arithmetical operations may in practice be implemented very efficiently by assembling a \emph{sparse matrix} composed of the $\Or(M)$ nonzero integration weights per output component of the \textsc{Matlab}- or \textsc{Octave} Radon transforms $\cR_{\Text{dis}}$:
    \begin{equation} 
      \bR \in \mR^{K_\theta \cdot J_x \times M_x \cdot M_z} \MTEXT{such that} \cR_{\Text{dis}}( \bN_j ) = \bR \boldsymbol \cdot \bN_j
    \end{equation}
    for all $\bN_j \in \mC^{M_x  \times M_z} \cong \mC^{M_x  \cdot M_z}$. The cylindrical Radon transform then corresponds to a sparse matrix-matrix product if the 2D slices are arranged as column vectors:
      \begin{equation} 
      \CRdis(\bN) = \bR \boldsymbol \cdot ( \bN_1 \; \; \bN_2 \; \; \ldots  \; \; \bN_{M_y} ) \MTEXT{for all} \bN \in \mX_{\Text{dis}} \cong \mC^{M_x  \cdot M_z \times M_y}. \vspace{2em}
    \end{equation}
    This sparse matrix representation has the advantage that the adjoint transform $\CRdis^\ast$ may be evaluated simply by applying the transpose of $\bR$. Moreover, sparse matrix-matrix products may be massively parallelized and permit efficient computations also on \emph{graphic cards}, for example, which may be exploited in future.

    However, it is still desirable to reduce the complexity of the Radon transform to the $\Or(M^3 \log M)$ flops required for the remaining algorithm. A promising approach is motivated by the Fourier Slice \thmref{thm:FourierSlice}, stating that analytical 2D Radon transform $\cR$ is equal to a polar Fourier transform \eqref{eq:PolarFT} and an inverse Fourier transform in the lateral coordinate, see \sref{S:Radon}. Accordingly, an alternate discretization of $\cR$ may be obtained via a discrete Fourier transform (DFT)
    \begin{equation}
     \DFT(\bN)(\bxi) = \sum_{\bj \in M_x  \times M_z} \bN_{\bj} \exp( - \I \bj \cdot \bxi )  \MTEXT{for} \bN \in \mC^{M_x  \cdot M_z } \label{eq:DFT}
    \end{equation}
    with frequency vectors $\bxi \in \mG \subset \mR^2$ sampled on a radially equidistant \emph{polar grid}
    \begin{equation} \mG = \{ \sigma (\cos\theta, \sin\theta): \sigma \in \mG_\sigma \subset \mR, \; \theta \in \mG_\theta \subset [0;\pi) \text{ discrete} \} \label{eq:PolarGrid} \end{equation}
    and an inverse FFT in the radial direction. The cylindrical transform once more corresponds to slice-wise application. Note that in the weak object limit, where the cylindrical Fourier transform $\CF$ explicitly enters in the governing operators, see \eqref{eq:WeakFwOps}, this approach reduces the number of required FFTs compared to the standard discretization of $\cR_{\Text{dis}}$. However, it may only be competitive if the polar Fourier transform is implemented by a \emph{fast algorithm} as the naive 2D DFT in \eqref{eq:DFT} is already of complexity $\Or(M^4)$. Unfortunately, no such $\Or(M^2\log M )$  \emph{polar} FFT algorithm is known for the exact evaluation of $\DFT(\bN)$ on the polar grid \eqref{eq:PolarGrid}.
    
    Existing implementations of the polar Fourier transform therefore have to interpolate from oversampled FFTs evaluated on related grids onto the desired polar sampling \cite{Averbuch2006,Fenn2007PolarFFT}. This results in ``fast'' algorithms in the sense of $\Or(M^2\log M )$ complexities, i.e.\ $\Or(M^3\log M )$ in the cylindrical 3D case, yet with a very large prefactor: typically, an evaluation with a low interpolation accuracy already takes $ \sim 100$ times longer than the computation of a Cartesian FFT of comparable size, see for instance \cite{Fenn2007PolarFFT}. Indeed, it has turned out in the preparation of this work that the algorithm proposed in \cite{Averbuch2006} cannot compete with the sparse matrix implementation of the Radon transform up to large resolutions $M \sim 1000$ - even though the programming effort put into the polar FFT was significantly larger. Therefore, this approach is not pursued any further.
    
    A possibly more efficient implementation might however be achieved via the related \emph{pseudo-polar} FFT \cite{Averbuch2001PPFFT}, which may be evaluated in $140 M^2 \log M$ flops, yielding the \emph{exact} discrete Fourier transform on a grid of concentric squares (instead of circles). Yet, note that the particular grid geometry in Fourier space corresponds to a rather unusual sampling of the resulting discretized Radon transform, which is obtained by inverse FFT along the radial direction of the pseudo-polar grid: the sampling in the incident angles $\theta$ is \emph{non-equispaced} whereas the lateral pixel spacing varies with $\theta$. In order to overcome these peculiarities of the pseudo-polar geometry, interpolation would become necessary once more. Accordingly, the potential benefits in computational efficiency would go along with a significant loss in geometrical flexibility of the method if  \algref{alg:PCT} was based on pseudo-polar FFTs. In other words, we arrive at the somewhat undesirable conclusion that the numerical constraints would dictate the experimental setup to a considerable degree.
    
    All in all, it thus seems that the standard, voxel-summation-based discretization of the Radon transform provides the best compromise between accuracy, efficiency and flexibility - despite its \emph{asymptotically} inferior $\Or(M^4)$ complexity. For this reason, we implement phase contrast tomography via the Newton-type \algref{alg:PCT} using the efficient sparse matrix formulation outlined in this section.

   \end{subsection}

\end{section}

\end{chapter}


\begin{chapter}{Numerical Results}\label{C:NumRes}
\vspace{-1em}

%
    The regularized Newton methods developed in the preceding chapter allow for a numerical solution of \probref{prob:1}, corresponding to our principal objective of sample reconstructions in phase contrast tomography.
    In the following, the performance of the derived \algref{alg:PCT} is validated by discussing numerical results obtained by an implementation in \textsc{Matlab}/\textsc{Octave} \cite{Octave}.
        \vspace{-.25em}
    \begin{section}{Far-Field Tomography from Simulated Data} \label{S:NumResFF}
    \vspace{-.25em}
    
    We start our presentation of numerical results with the case of far-field phase contrast tomography, governed by the forward operator given in \eqref{eq:ForwardOpFF}. As the discussion is focused on \emph{qualitative} aspects, we widely restrict to simulations within a two-dimensional toy model in order to simplify visual inspection of the results. All of the described qualitative effects, however,  manifest analogously for far-field reconstructions within the physically relevant 3D geometry.
     \vspace{-.25em}
        \begin{subsection}{Simulation Setup} \label{SS:NumResFF-GenSetup}
        \vspace{-.25em}
     
     The considered 2D toy model  incorporates only the $x$- and $z$-dimensions within the tomographic plane of rotation (compare \figref{fig:SetupIdeal}). Accordingly, the cylindrical Radon transforms in the discrete forward operators, derivatives and adjoints to be evaluated in \algref{alg:PCT} reduce to standard two-dimensional ones applied to planar ``objects'' given by 2D images. The resulting ``contact images'' recorded under different incident angles $\theta$ are simply \emph{one}-dimensional profiles as sketched in \figref{fig:RadonVisualization}. The ensemble of these profiles for different $\theta$, propagated by the 1D Fourier transform, represents the intensity data given by 2D far-field ``sinograms''.
     
     As a two-dimensional test phantom we choose a $256 \times 256$ pixel version of the abstract cell sketch in \figref{fig:CellSketchTestObj}, introduced in \cite{Giewekemeyer2011CellSketch}. By scaling this real-valued object $\bN_{\mR} \in  \mR^{256 \times 256}$ with complex constants $c \in \mC$, we may construct single-material objects
     \begin{equation}
      \bN^\dagger = c\bN_{\mR} \in \mC^{256 \times 256} = \mX_{\Text{dis}}.
     \end{equation}
     of arbitrary absorption-refraction-ratio $\beta/\delta$ and magnitude, see \sref{SS:SpecialObjects}. To leading order, the far-field intensities are independent of the phase of $c$ and thus of $\beta/\delta$ according to the weak object limit of the forward operator given in \eqref{eq:WeakFwOpFF}. For simplicity, we thus restrict to real and positive $c > 0$, i.e.\ to \emph{pure phase objects}.
     
     Note that the \emph{magnitude} of the scaling constant $c$ controls the nonlinearity of the object transmission function $O ( \bN ) = \exp( - \I k \CR (\bN ) )$. We measure this with a scaled maximum norm
      \begin{equation}
      \norm{\bN} = kL \norm{\bN}_{\infty} \label{eq:ObjScaleNorm}
     \end{equation}
     where $k$ is the wavenumber and $L$ the aspect length of the computational domain, i.e.\ the thickness of the object.
     According to \eqref{eq:PhaseWrapCond2}, $\norm{\bN} = 2 \pi$ then defines the transition to strong objects for which phase-wrapping  may occur, whereas $\norm{\bN} \lesssim 0.1$ ensures that the weak object approximation applies, see \sref{SS:PhaseWrap} and \sref{SS:SpecialObjects}.
     
     Apart from the significance as a scale for the strength of an object incorporated in \eqref{eq:ObjScaleNorm}, the parameters $k$ and $L$ have no further qualitative impact on far-field tomography if the lateral coordinate $\bxi_x$ in image space is scaled as in \sref{SS:ContrastFormationFF}. A specification of these is therefore omitted in the numerical study. For simplicity, we further restrict to ideal \emph{plane wave illumination}, setting $P = 1$ in the expressions \eqref{eq:ForwardOpFF}, \eqref{eq:FrechetForwOpFF}, \eqref{eq:AdjForwOpFF} for the forward operator, derivative and adjoint.
	\begin{figure}[hbt!]
	 \centering
	 \includegraphics[width=0.5\textwidth]{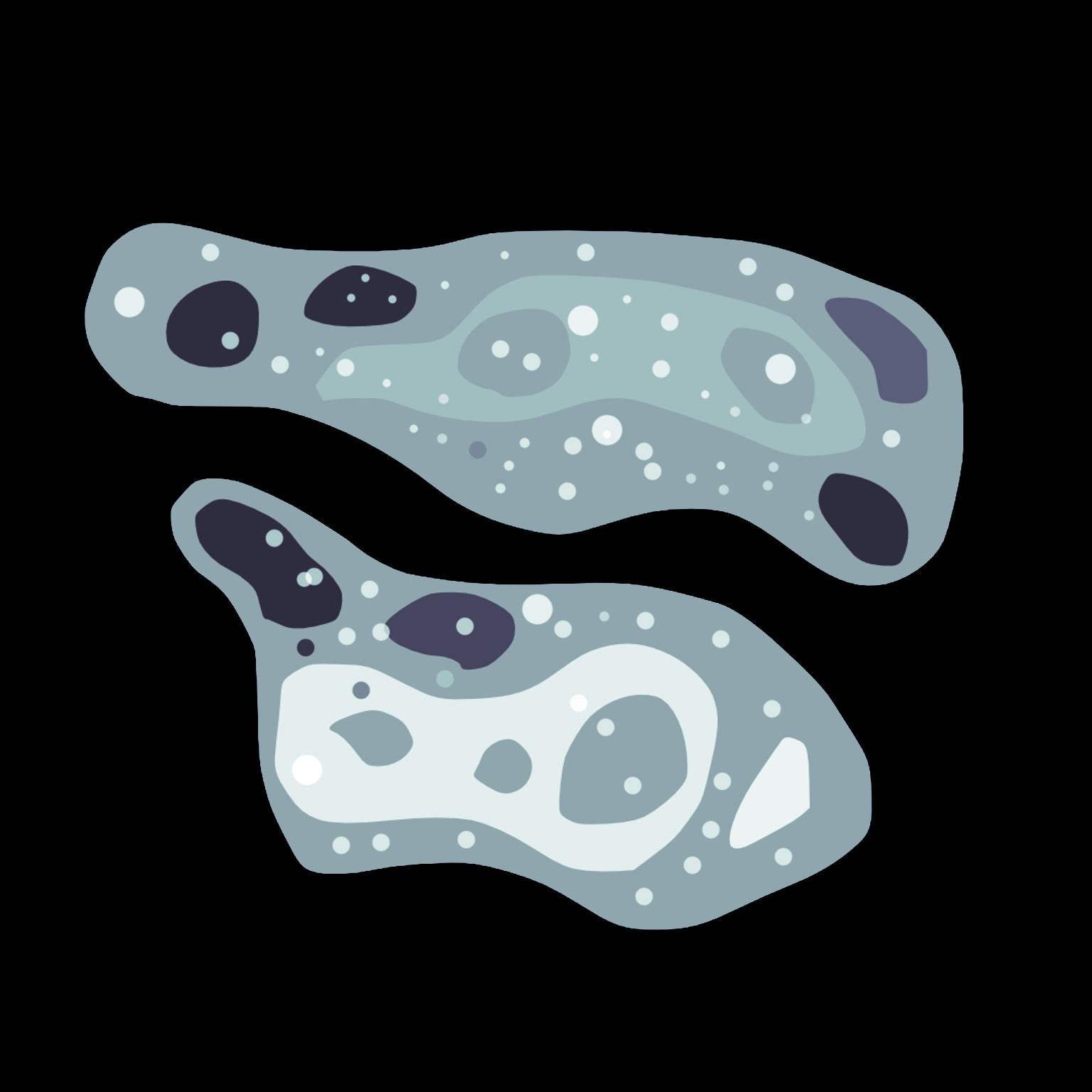}
	 \caption{Sketch of an abstract cell used as a phantom for the 2D far-field test cases of phase contrast tomography. Scaling of the real-valued image by complex constants allows the simulation of arbitrarily strong or weak single-material objects. The tomographic far-field data corresponding to the phantom is given by 1D projections under different incident angles $\theta$, propagated by Fourier transforming in the single lateral coordinate of the sinogram. (Source: \cite{Giewekemeyer2011CellSketch})  \label{fig:CellSketchTestObj}}
	\end{figure}
     
     In order to preclude undersampling issues, we simulate $K_\theta = 512$ incident angles $\theta \in [0; \pi)$ sampled at a lateral resolution of $K_x = 512$ attained by zero-padding of the 1D contact images (see \sref{SS:ZeroPad}), i.e.\ we reconstruct from intensities
     \vspace{-.25em}
     \begin{equation}
      \bI^{\Textbf{err}}  \in \mR^{512 \times 512} = \mY_{\Text{dis}}. \vspace{-.25em}
     \end{equation}
     The data $\bI^{\Textbf{err}} $ corresponding to an exact object $\bN^\dagger \in \mX_{\Text{dis}}$ is generated via
          \vspace{-.25em}
     \begin{equation}
     \bI^{\Textbf{err}} =   \Text{Poi} (  \bI^\dagger) \MTEXT{with} \bI^\dagger := I_0 \cdot F_{\Text{dis}}(\bN^\dagger) \vspace{-.25em} \label{eq:DataGenFF},
     \end{equation}
     assigning the image under the discrete forward operator as the parameter (expectation value) of a Poisson distribution. By the choice of the scaling factor $I_0 > 0$, a certain relative $L^2$-noise level $\varepsilon = \norm{\bI^{\Textbf{err}} - \bI^{\dagger}}_2 / \norm{\Textbf{\bI}^\dagger}_2$ may be prescribed\footnote{Note that the impact of a given noise level $\varepsilon$ onto the reconstruction depends on the resolution.}.
     
     The initial regularization parameter $\alpha_0$ is fixed heuristically such that the CG-method for the first Newton step terminates after $\sim 5$ iterations. This is an indicator for the condition number of the linear problem \eqref{eq:DiscNewtonLinProb} to be solved being neither too small nor too large, which would correspond to over- or underestimated regularization. In order to obtain a balance between the data fidelity and regularization term that is independent of the object's magnitude and intensity factor $I_0$, we further have to apply the scaling $\alpha_0 \sim \norm{\bI^{\Textbf{err}}}_{\mY_{\Text{dis}}} / \norm{\bN^\dagger}^2_{\mX_{\Text{dis}}}$. In the considered setting, a suitable choice is found as
     \begin{equation}
		\alpha_0 = \frac 1 {10} \cdot \frac{ \norm{\bI^{\Textbf{err}}}_{\mY_{\Text{dis}}}^2 }{  \norm{\bN^\dagger}^2_{\mX_{\Text{dis}}} }   \label{eq:NumResFF-alpha0}.
     \end{equation}
     For every subsequent Newton iteration in \algref{alg:PCT}, $\alpha_k$ is reduced by a factor of $r_\alpha = \frac 2 3$. As we are mainly interested in a qualitative validation of the reconstruction method, we further restrict to $L^2$-penalty terms by taking the Gramian $\cG_{\mX_{\Text{dis}}}$ as the identity. The chosen Kullback-Leibler-type data 
     fidelity term is parametrized by setting $\cG_{\mY_{\Text{dis}}}$ according to \eqref{eq:GramYDisc} where we take the truncation at $I_{\min} = I_0$.
     \begin{table}[htb!]
  \centering
  \begin{tabular}{cccccccccc} 
    \toprule
       $\bI^{\Textbf{err}}$  &  $\mX_{\Text{dis}}$ &  $\mY_{\Text{dis}}$ & $\cG_{\mX_{\Text{dis}}}$  & $\cG_{\mY_{\Text{dis}}}$  & $\alpha_0$ & $r_\alpha$  & $\bN_0$ & $P$ & Constraints \\
    \midrule
       \scriptsize{\eqref{eq:DataGenFF}} & $\mC^{256 \times 256}$ & $\mR^{512 \times 512}$  & $\id_{\mX_{\Text{dis}}}$  &  \scriptsize{\eqref{eq:GramYDisc}} &  $ \frac{ \norm{\bI^{\Textbf{err}}}_{\mY_{\Text{dis}}}^2 }{ 10 \norm{\bN^\dagger}^2_{\mX_{\Text{dis}}} } $  & $\frac 2 3$  & $\neq 0 $ & $ 1$ & $\substack{  \text{pure phase obj.} \\ \text{support} } $ \\ 
    \bottomrule
  \end{tabular}
  \caption{Chosen setup parameters for the considered numerical test cases of 2D far-field phase contrast tomography as assigned to \algref{alg:PCT}. The exact test object $\bN^\dagger \in \mX_{\Text{dis}}$ is taken as a scaled and rebinned version of the abstract cell sketch in \figref{fig:CellSketchTestObj}. The probe choice $P=1$ corresponds to the assumption of ideal plane wave illumination. }
  \label{tab:NumResFFSetup}
\end{table}

     In the considered far-field case, the canonical choice of an initial guess $\bN_0 = 0$ for an unknown object leads to immediate stagnation of the Newton method. This is due to the fact that the \Frechet derivative \eqref{eq:FrechetForwOpFF} obtained in \thmref{thm:ForwardOpsFrechet} vanishes identically at $N = 0$. This suggests that the choice of the initial guess is significant for the convergence of the Newton method in general, which is examined in the following numerical examples along with the impact of support constraints on the reconstruction. The general setup for the considered far-field test cases assigned to \algref{alg:PCT} is summarized in \tabref{tab:NumResFFSetup}.
     
     \end{subsection}

        \begin{subsection}{Ab Initio Reconstructions in 2D} \label{SS:NumResFF-NoReference}
     
	  As a first test case, we attempt an \emph{ab initio} reconstruction of the cell phantom in \figref{fig:CellSketchTestObj} as a pure phase object of magnitude $\norm{\bN^\dagger} = \pi$, corresponding to a neither weak nor phase-wrapping object. The idea is not to incorporate any strong a priori knowledge into the reconstruction. Consequently, no support constraint is assumed apart from the restriction to the $256 \times 256$-sized square object domain. 
	  As the initial guess $\bN_0$, we choose a Gaussian of the same peak magnitude as the object to be reconstructed. Without incorporating specific information on the unknown sample, this prescription of an initial centered peak breaks the ``trivial'' symmetry of the involved Fourier phase retrieval problem with respect to translations of the object (compare \sref{SSS:TrivialAmbiguities}-Trivial Ambiguities).

	  
	   \figref{fig:NumResFF-NoSupp} shows the results for a Poisson-noise level of $\varepsilon = 10^{-3}$ in the simulated intensity data after 10 Newton iterations. As seen from the errors plotted in \figref{fig:NumResFF-Conv} (blue curves), the iterates $\bN_k$ no longer improve at this point although the data residual continuous to decrease.
	  Apart from the rough shape, the reconstructed object (\figref{fig:NumResFF-NoSupp-c}) is found to match only poorly with the exact one in (a), the relative $L^2$-error being $\norm{\bN_{10}  - \bN^\dagger}_2 / \norm{\bN^\dagger}_2 \approx 51 \, \% $: structures in the interior of the cell are entirely unidentifiable due to the dominant artifacts in the reconstruction. This observation is contrasted by an accurate fit of the observed intensities data, as seen by comparing figure parts (d) and (f), with a final $L^2$-residual  of $\approx 0.5 \, \%$. 
	  
	  Note that the translational invariance of the far-field (Fourier-)intensities is broken merely by the choice of the initial guess $\bN_0$: if it was not for the initial bias by the prescribed Gaussian peak, the algorithm could reconstruct the object at \emph{any} possible location in the computational domain wherever the shape fits entirely - and would thus stagnate at some blurry intermediate state that is symmetric with respect to all these shifted realizations of the cell phantom. 
	  However, the initial guess $\bN_0$ does not constitute any strict constraint but only acts weakly via the dependence of the Newton iterate $\bN_{k+1}$ on the preceding ones $\bN_k, \bN_0$. Accordingly, as the initial data $F_{\Text{dis}}(\bN_0)$ is far from the simulated intensities $\bI^{\Textbf{err}}$ to be fitted (compare \figref{fig:NumResFF-NoSupp-d}-e), the weakly suppressed translational invariance might still manifest considerably in the course of the Newton iterations.
	  The discrepancy between the good data fit and the artifacts in the reconstructed object may thus be attributed to latent ``trivial'' ambiguities due to translational symmetry.
	  	   \begin{figure}
	  \centering
	  \subfloat[Exact object $\bN^\dagger$ \label{fig:NumResFF-NoSupp-a}]{\includegraphics[height=0.32\textwidth]{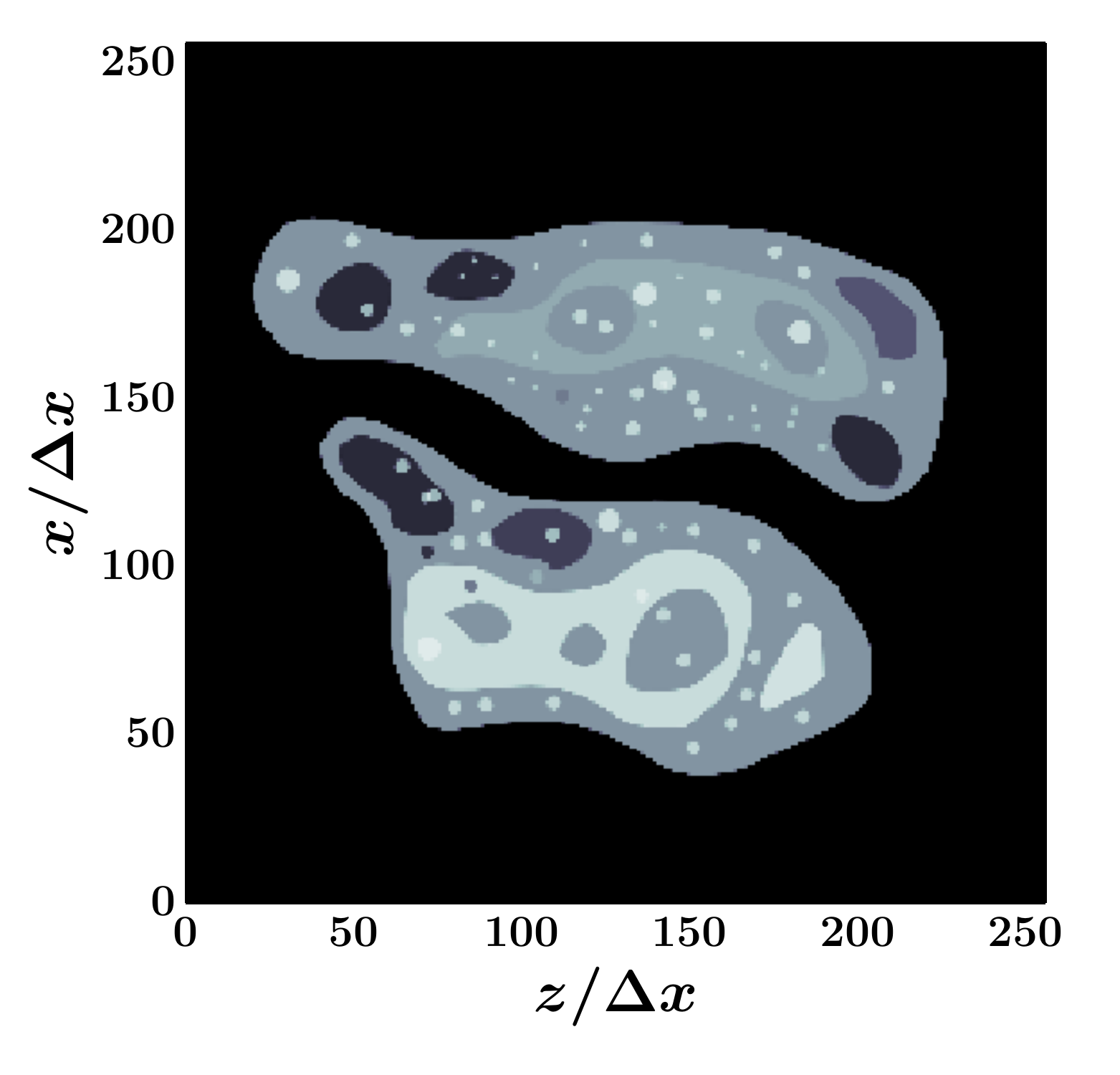}}
	  \hfill
	  \subfloat[Initial guess $\bN_0$\label{fig:NumResFF-NoSupp-b}]{\includegraphics[height=0.32\textwidth]{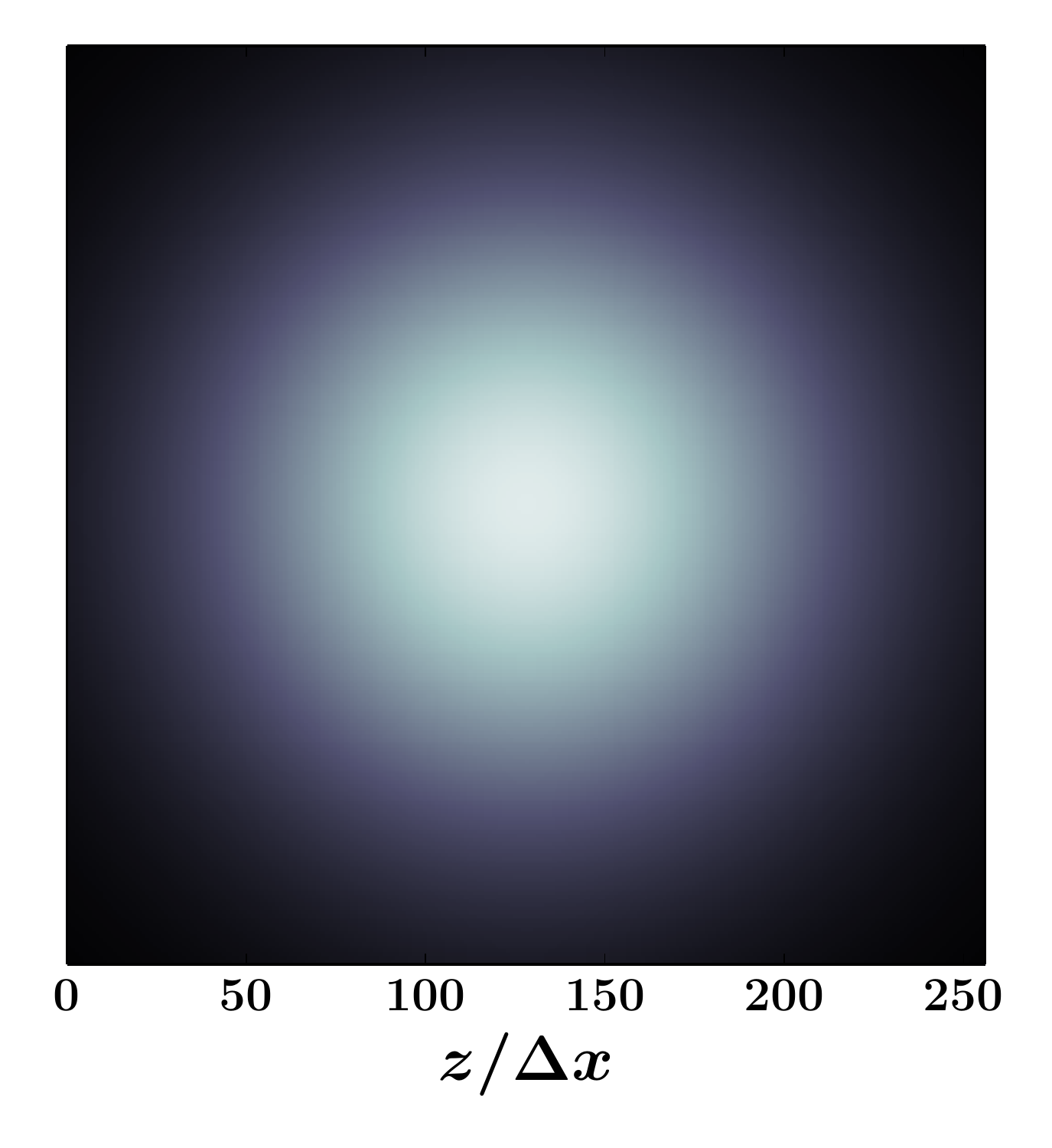}}
	  \hfill
	  \subfloat[Reconstructed object $\bN_{10}$ \label{fig:NumResFF-NoSupp-c}]{\includegraphics[height=0.32\textwidth]{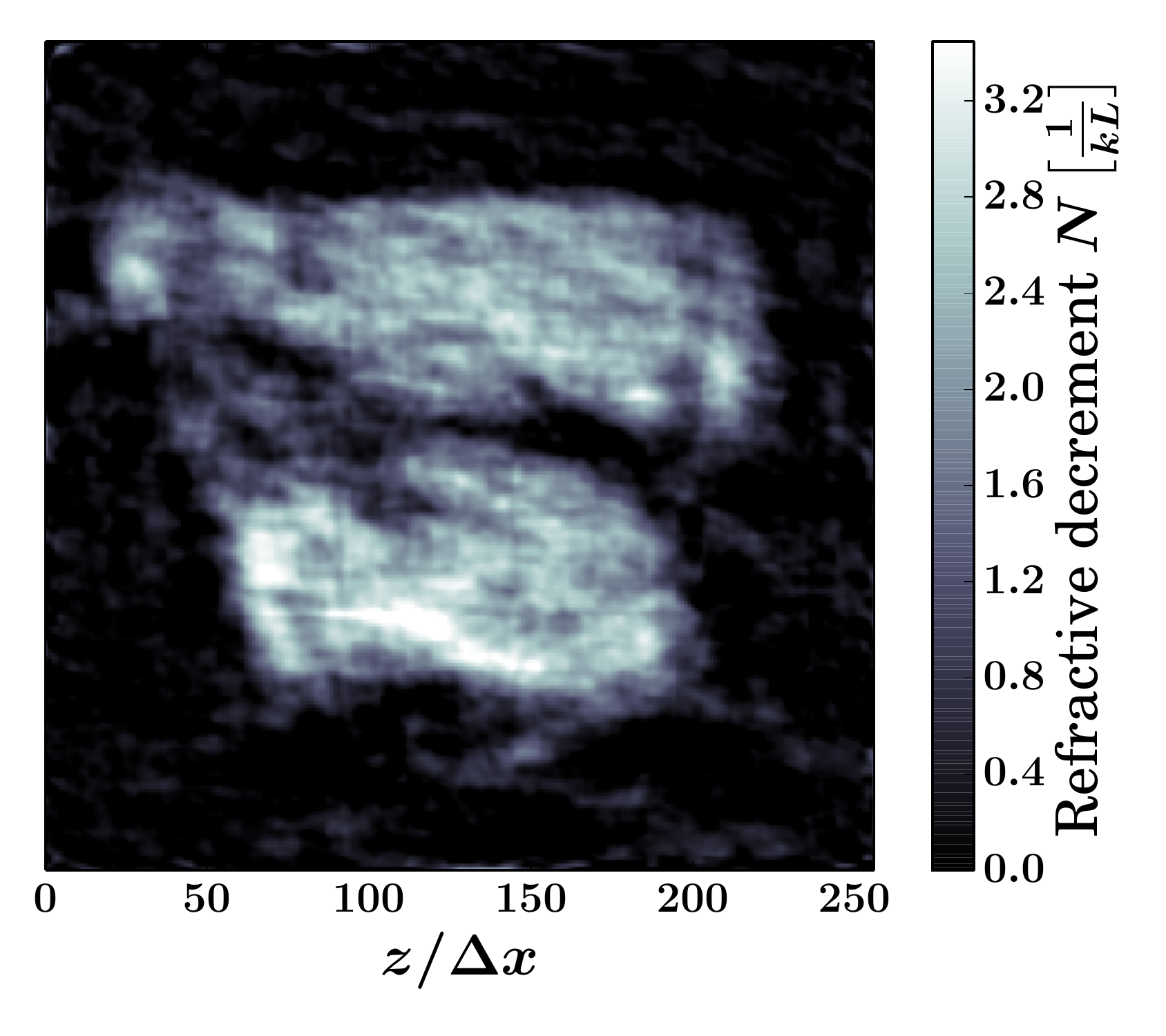}}
	  \\
	  \subfloat[Simulated intensities $\bI^{\Textbf{err}}$ \label{fig:NumResFF-NoSupp-d}]{\includegraphics[height=0.32\textwidth]{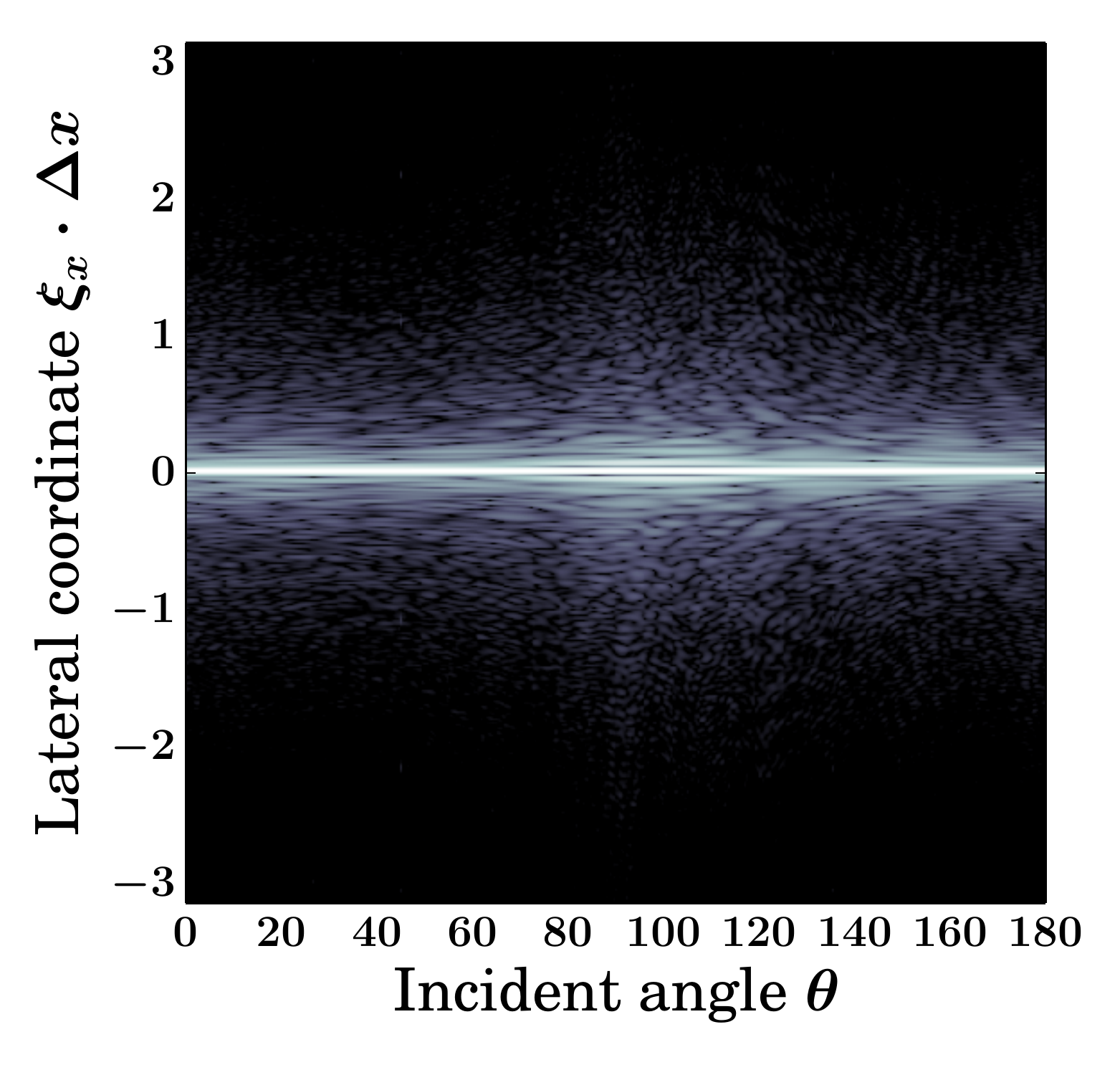}}
	  \hfill
	  \subfloat[Initial data $F_{\Text{dis}}(\bN_0)$ \label{fig:NumResFF-NoSupp-e}]{\includegraphics[height=0.32\textwidth]{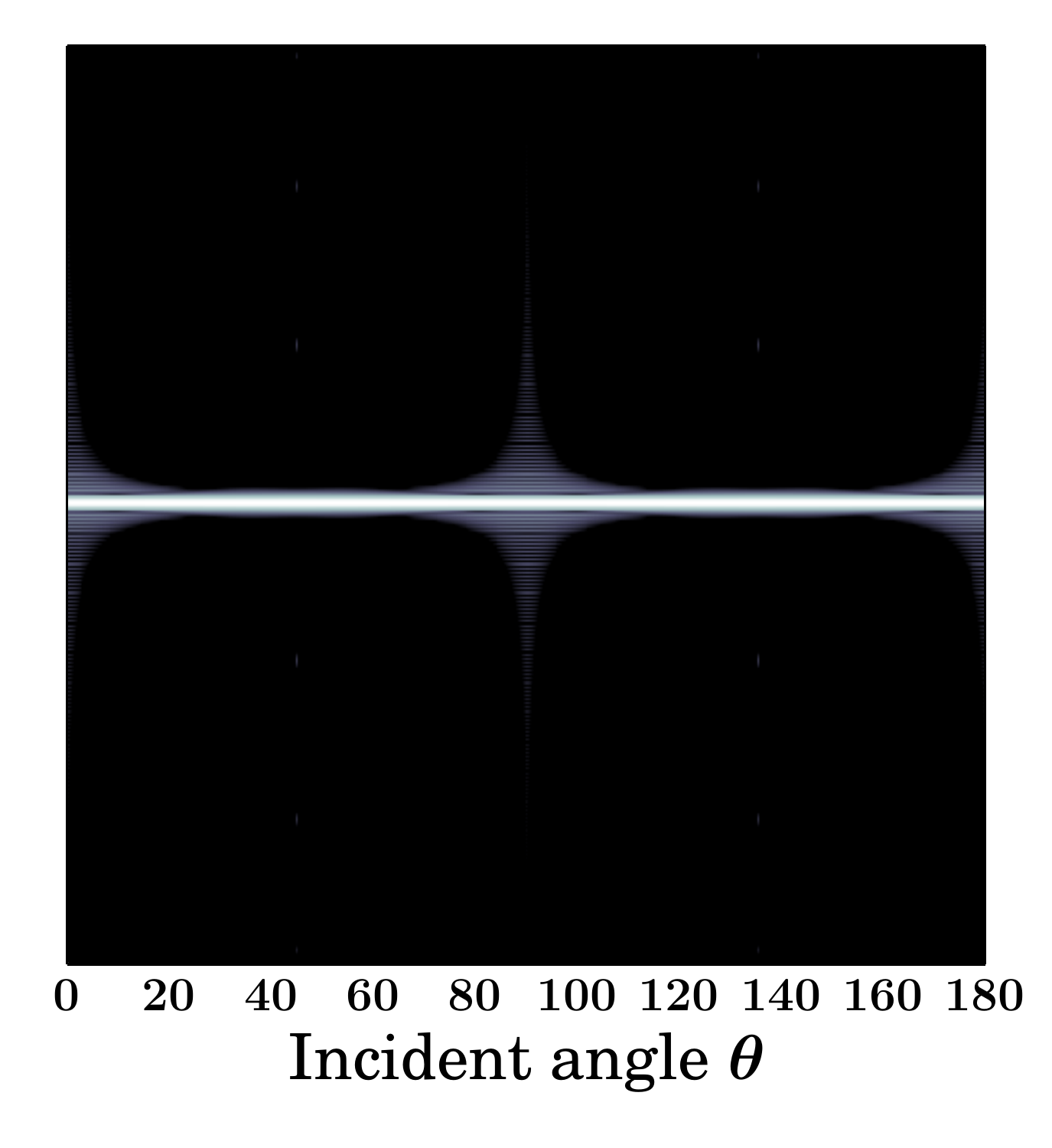}}
	  \hfill
	  \subfloat[Reconstructed data $F_{\Text{dis}}(\bN_{10})$ \label{fig:NumResFF-NoSupp-f}]{\includegraphics[height=0.32\textwidth]{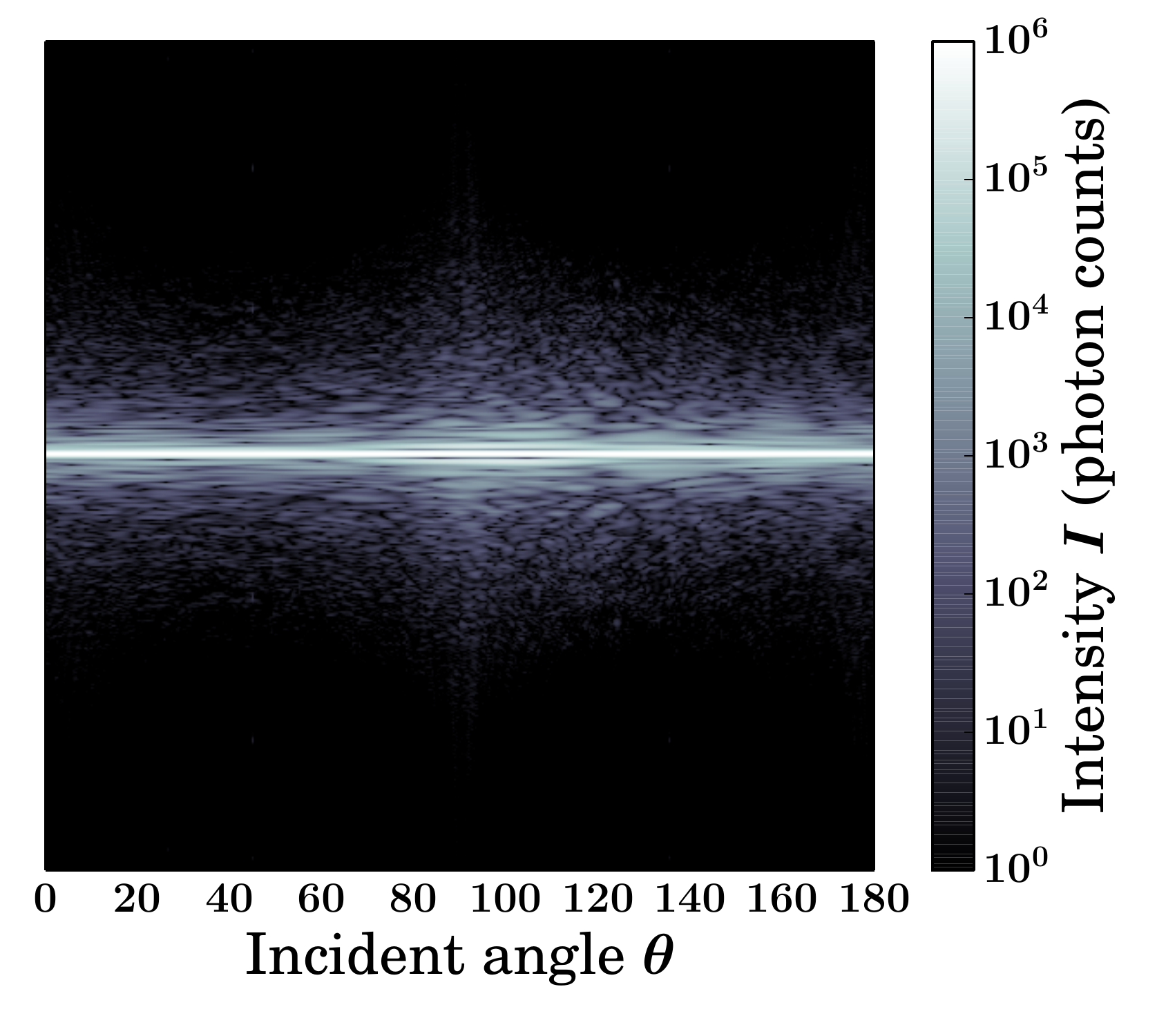}}
	  \caption{2D Far-field reconstruction of a pure phase object $\bN^\dagger$, $\norm{\bN^\dagger} = \pi$ from simulated data with $0.1\,\%$ Poisson noise (total photon counts $\approx 2.4 \cdot 10^9$) without support constraint. Shown: Exact object $\bN^\dagger$, initial guess $\bN_0$ and obtained solution $\bN_{10}$ after 10 Newton iterations, where the reconstruction is found to stagnate (see \figref{fig:NumResFF-Conv}). The relative $L^2$-error of the final iterate $\bN_{10}$ vs. $\bN^\dagger$ is $\approx 51 \, \%$, whereas the residual in the data is only $\approx 0.5 \, \%$. $\Delta x$ is the aspect length of a single pixel. Reconstruction parameters according to \tabref{tab:NumResFFSetup}. \label{fig:NumResFF-NoSupp}} 
     	  \end{figure}
	  
	  Despite the overall poor quality of the achieved numerical result in \figref{fig:NumResFF-NoSupp}, the approximate shape of the cell-phantom is reconstructed sufficiently accurate in order to allow for a support estimate by suitable thresholding. Such iterative support refinements are standard in alternating-projection-type methods for far-field reconstruction such as the Shrinkwrap Algorithm \cite{Fienup1986algorithm}. Although this approach is not straightforward to incorporate into the regularized Newton method considered here, we investigate its potential benefits by repeating the above reconstruction supplemented with two different support constraints:
	  \begin{itemize}
	   \item[(a)] Prescription of the exact cell-shaped support
	   \item[(b)] A ``tight'' rectangular support that precludes translations of the object
	  \end{itemize}
	  Since these constraints already break the translational symmetry, we simply choose bump functions as the initial guess $\bN_0$, assigning the value $\frac 1 2 \norm{\bN^\dagger}$ within the support and zero outside. In order to separate the effect of the support constraint from that of the initial condition, we furthermore perform another reconstruction without support constraint but with the cell-shaped ``support bump'' as the initial guess. All other parameters in these  three supplementary test cases are chosen exactly as for the results shown in \figref{fig:NumResFF-NoSupp}.
	  	\begin{figure}[hbt!]
	 \centering
	 \includegraphics[width=0.7\textwidth]{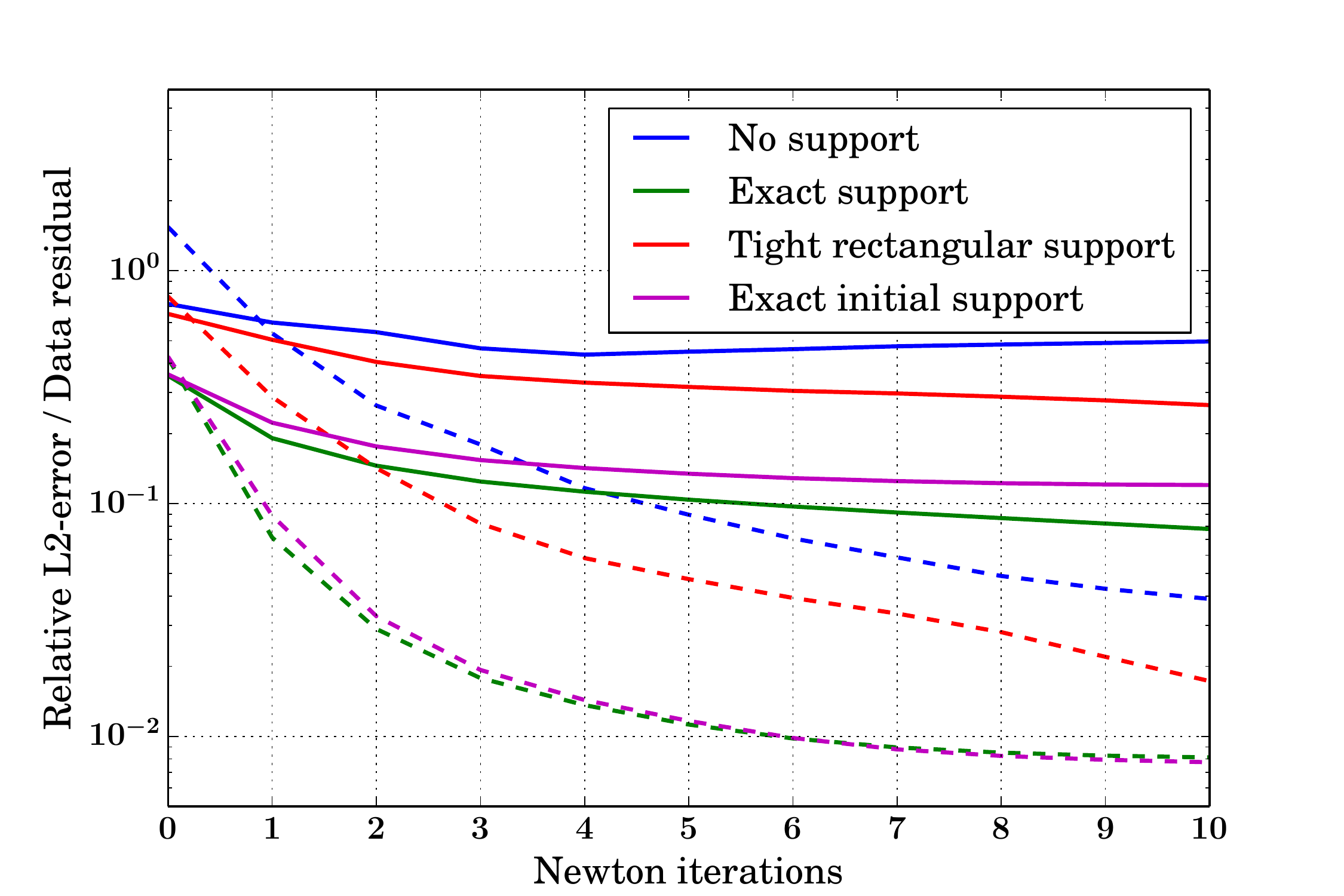}
	 \caption{Convergence of the 2D far-field test cases in \sref{SS:NumResFF-NoReference}. Solid curves show the relative $L^2$-error $\norm{\bN_{k}  - \bN^\dagger}_2 / \norm{\bN^\dagger}_2$ in the reconstructed object $\bN_k$ after $k$ Newton iterations. Dashed lines show the corresponding data residual $\norm{F_{\Text{dis}}( \bN_k ) - \bI^{\Textbf{err}}}_{\mY_{\Text{dis}}}/ \norm{\bI^{\Textbf{err}}}_{\mY_{\Text{dis}}}$ in the error metric induced by the Gramian $\cG_{\mY_{\Text{dis}}}$ in \eqref{eq:GramYDisc}. ``No support'' corresponds to the reconstruction in \figref{fig:NumResFF-NoSupp} which is to found to stagnate after a few iterations. The other curves correspond to the numerical results for the different support constraints and initial conditions  shown in \figref{fig:NumResFF-DiffSupps}. \label{fig:NumResFF-Conv}}
	\end{figure}
	      \begin{figure}[h!]
	  \centering
	  \hfill
	  \subfloat[Reconstructed objects $\bN_{10}$ \label{fig:NumResFF-DiffSupps-a}]{ \includegraphics[height=0.32\textwidth]{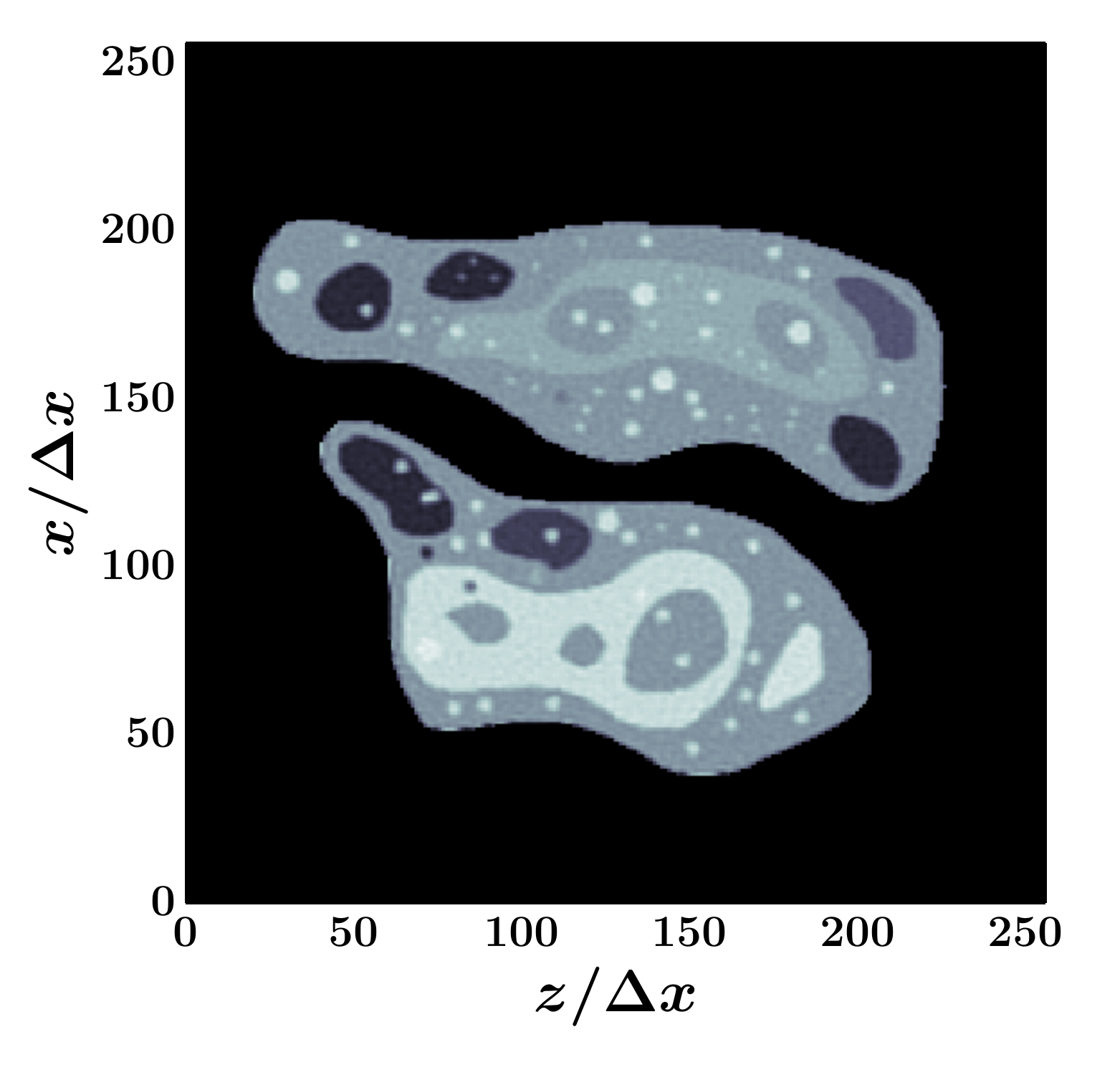} \newline \includegraphics[height=0.32\textwidth]{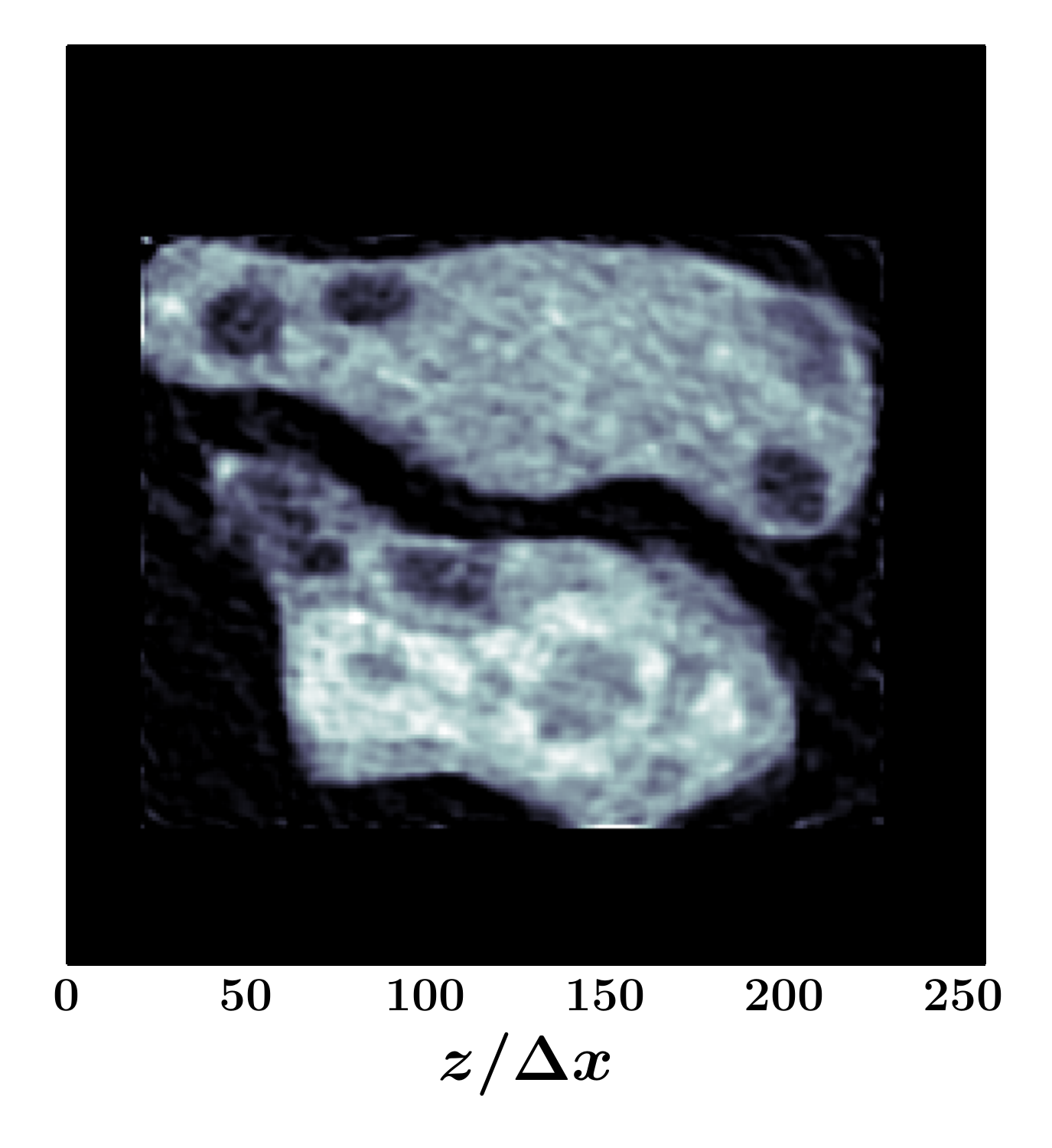}  \includegraphics[height=0.32\textwidth]{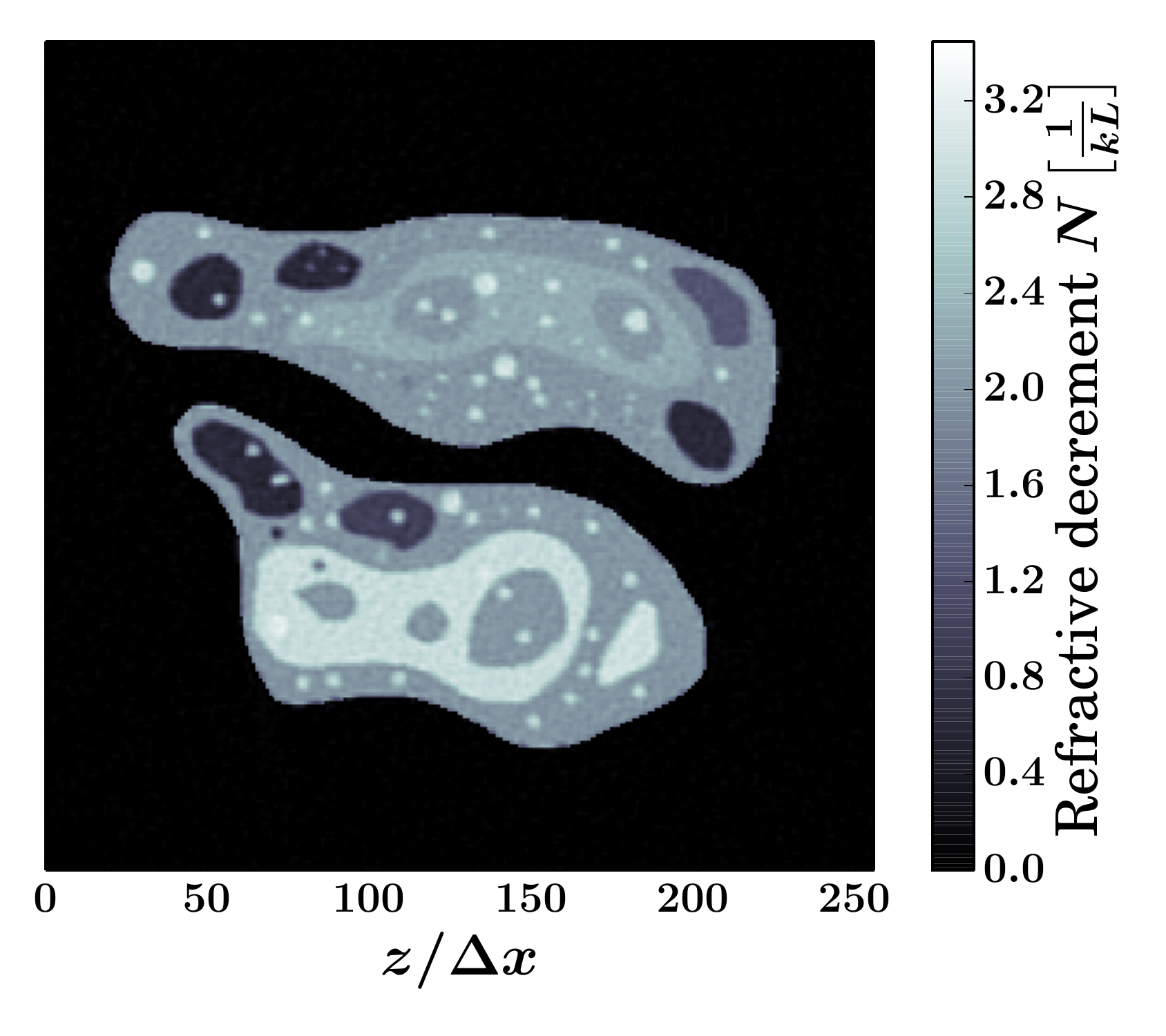} }
	  \hfill
	  \\
	  \hfill
	  \subfloat[Initial guess $\bN_0$ \label{fig:NumResFF-DiffSupps-b}]{\includegraphics[height=0.32\textwidth]{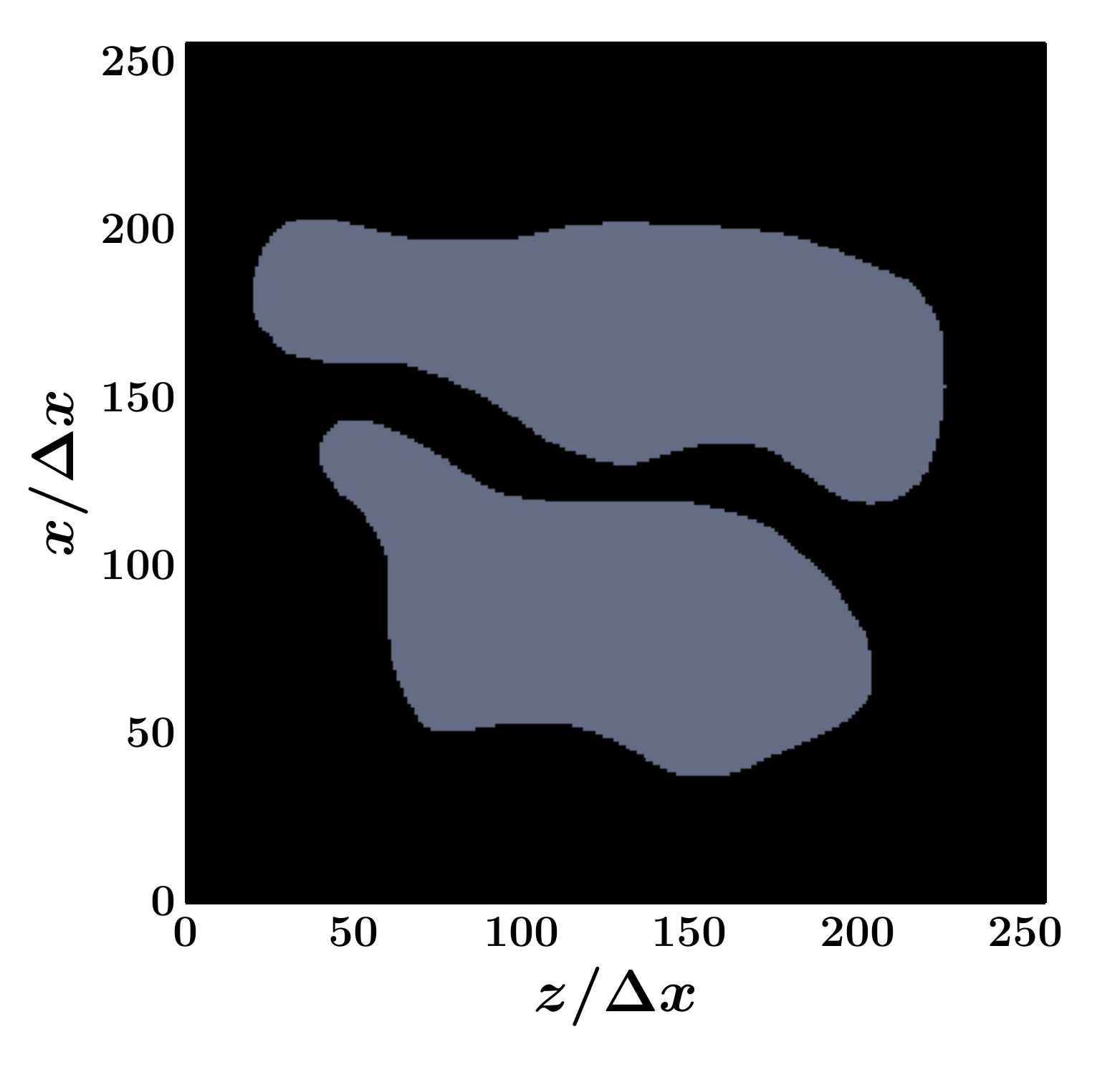} \includegraphics[height=0.32\textwidth]{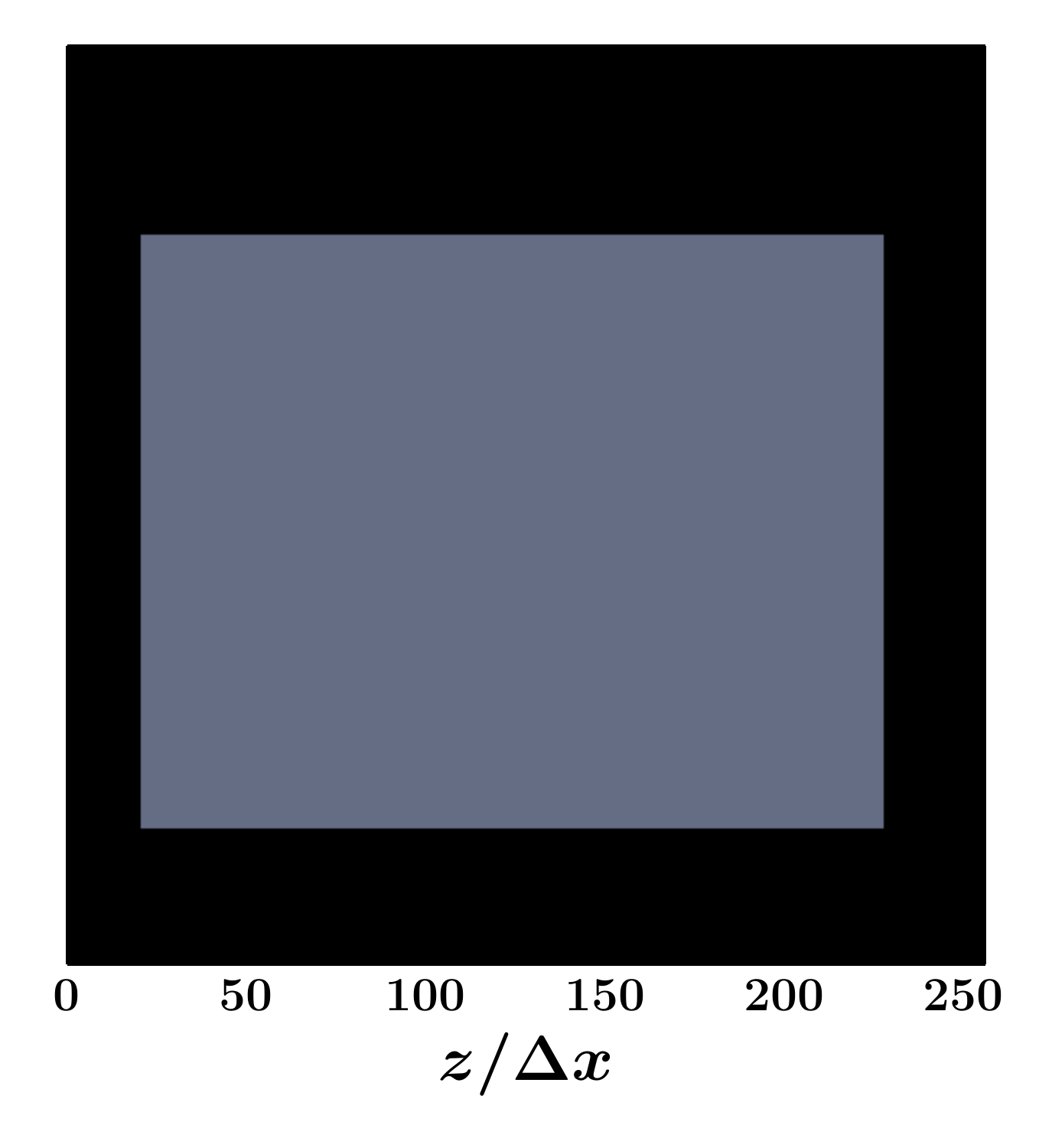} \includegraphics[height=0.32\textwidth]{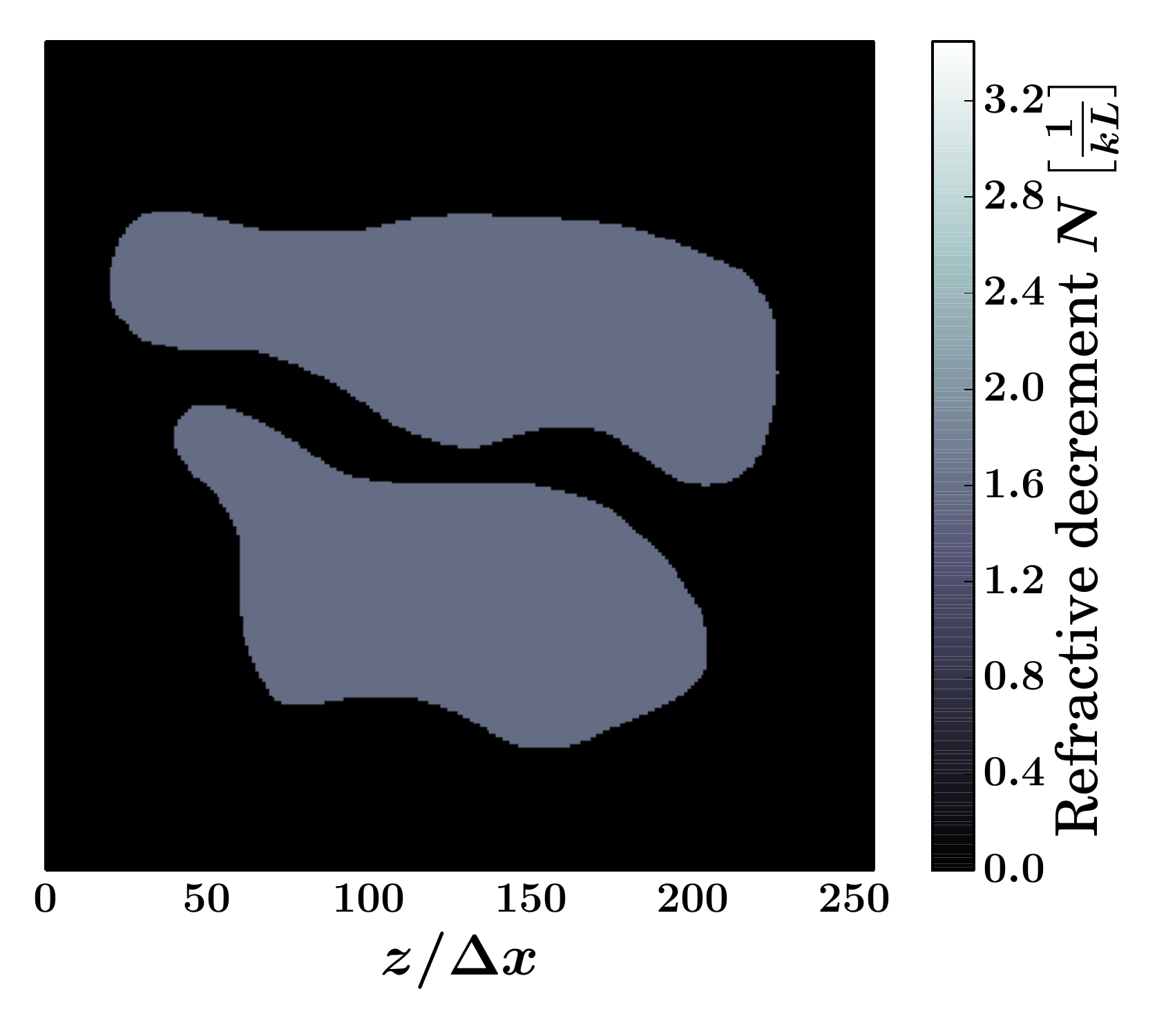} }
	  \\
	  \subfloat[Initial data $F_{\Text{dis}}(\bN_0)$ \label{fig:NumResFF-DiffSupps-c}]{ \includegraphics[height=0.32\textwidth]{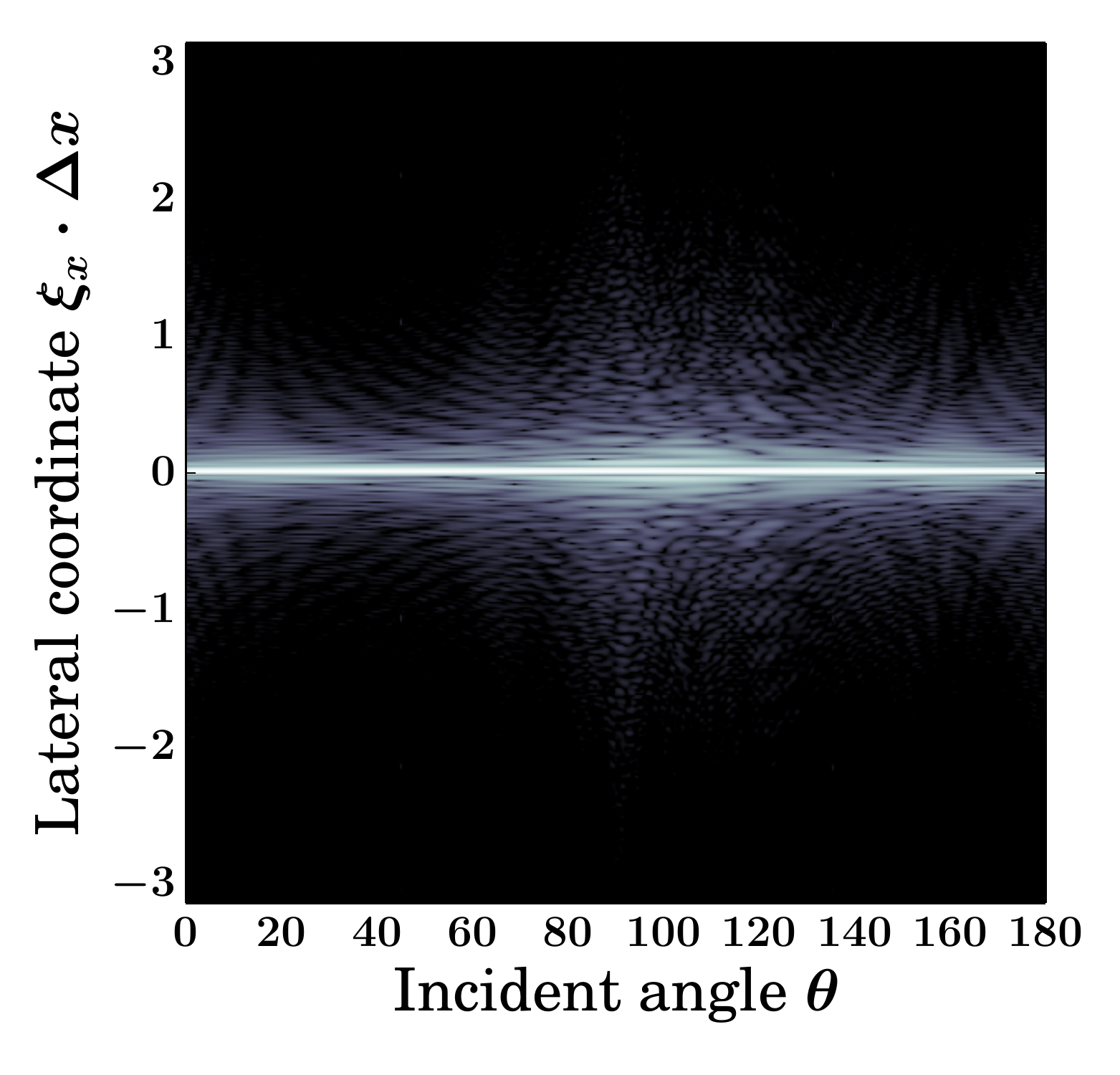} \includegraphics[height=0.32\textwidth]{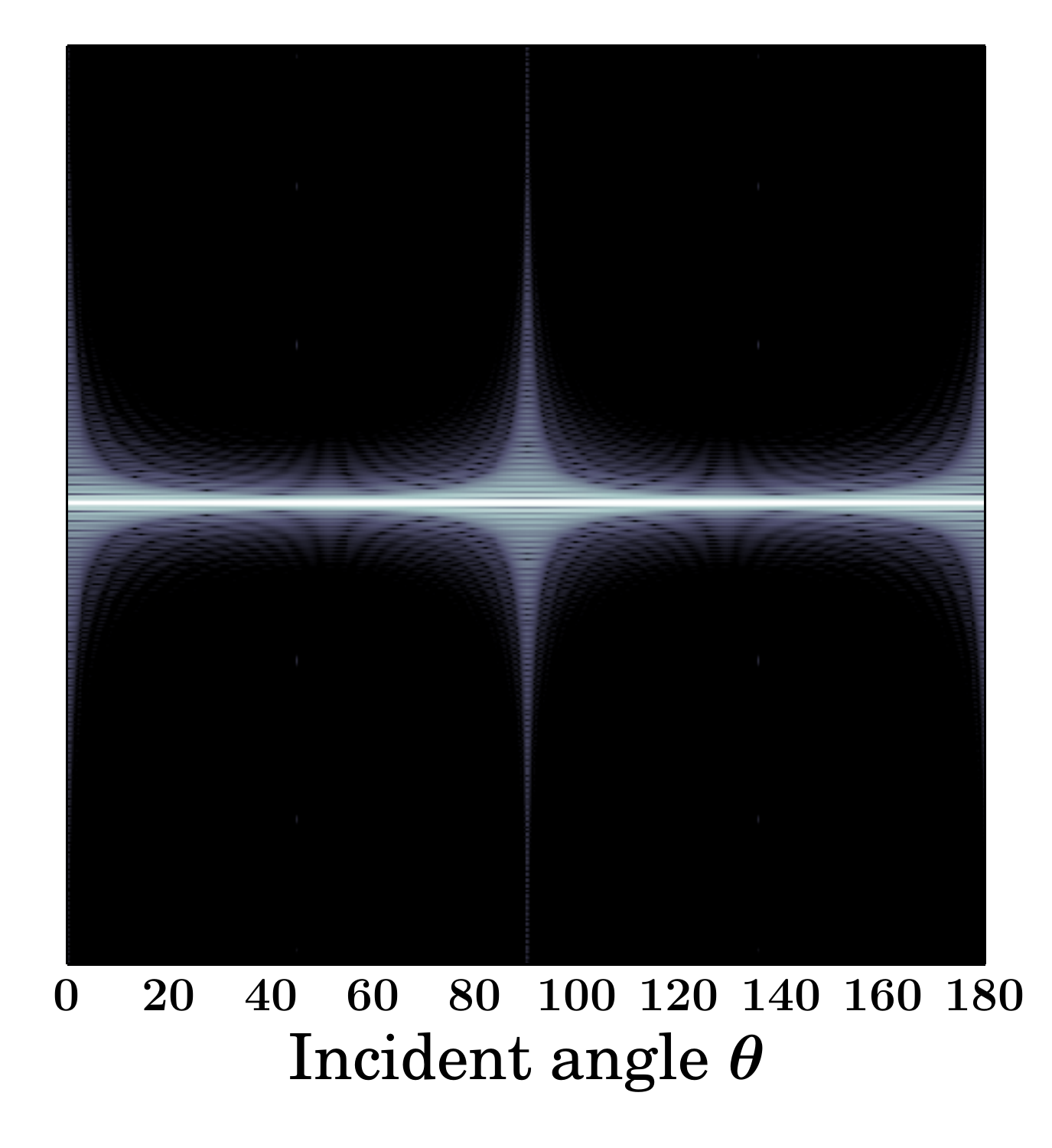} \includegraphics[height=0.32\textwidth]{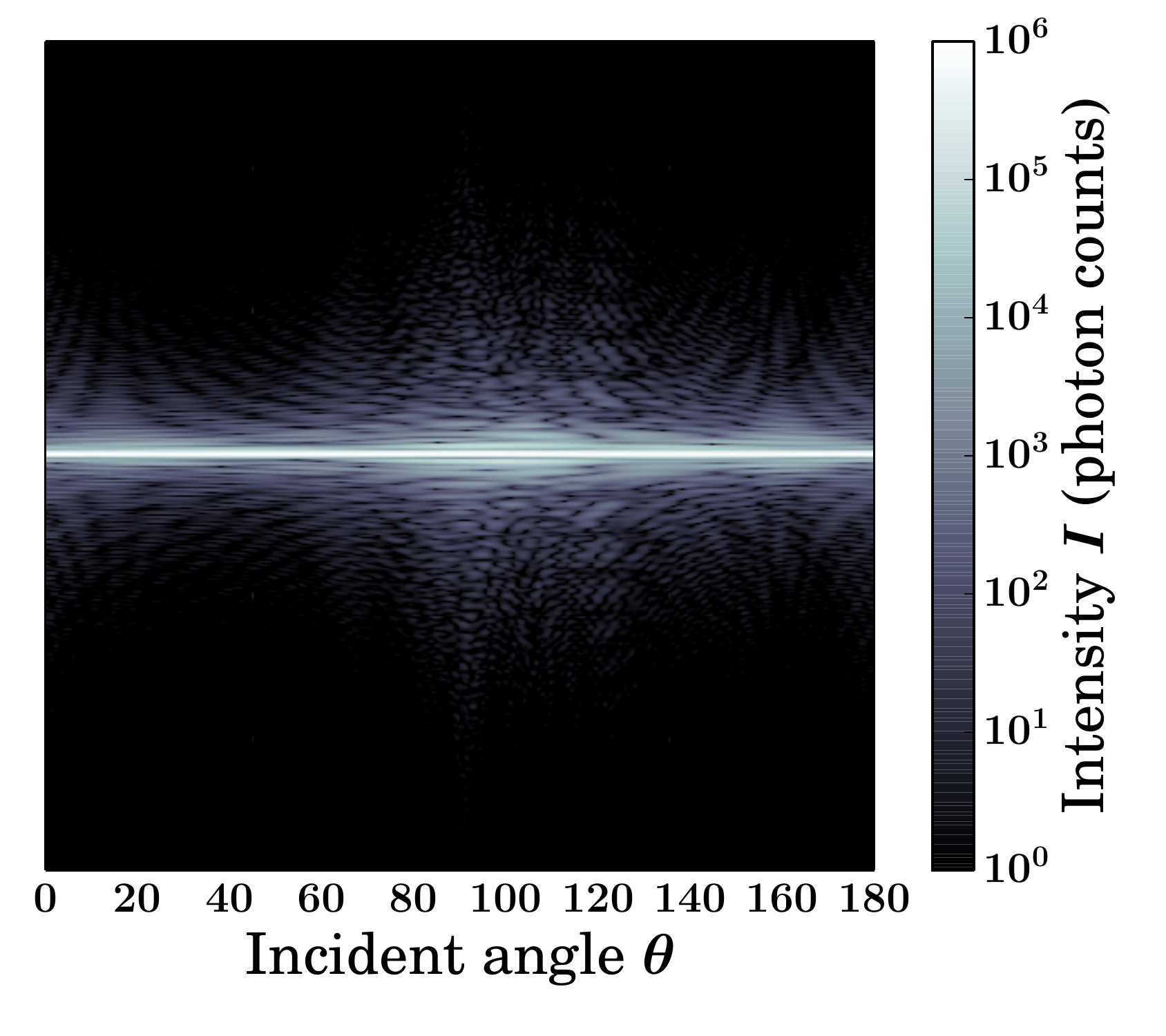}}
	  \caption{2D Far-field reconstructions for different support constraints and initial guesses (other parameters identical to simulation in \figref{fig:NumResFF-NoSupp}). Left: exact support of the cell-shaped phantom as constraint. Center: rectangular support enclosing the object tightly. Right: No support constraint but exact shape as initial guess. \label{fig:NumResFF-DiffSupps}}
	  \end{figure}
	  
	  Convergence of the object reconstructions $\bN_k$ and the corresponding data residual for are plotted in \figref{fig:NumResFF-Conv}. The final iterates after 10 Newton iterations are visualized in \figref{fig:NumResFF-DiffSupps} along with the chosen initial guess and -data. While the reconstruction assuming exact support knowledge comes out apparently artifact-free (left image in \figref{fig:NumResFF-DiffSupps-a}, $\approx 8\,\%$ $L^2$-error), the reconstructed object obtained for the tight but imperfect rectangular support estimate (center, $\approx 27\,\%$ deviation from $\bN^\dagger$) still does not accurately reproduce details of the cell. Yet, it would certainly allow for another refinement of the support estimate compared to a first guess based on the unconstrained computation in \figref{fig:NumResFF-NoSupp}. Note, however, that the intensity data $\bI^{\Textbf{err}}$ considered in the above simulations is still practically ideal. Estimating supports based upon preliminary reconstructions is most likely less stable for realistic data sets. One might argue that the visualized reconstruction for the rectangular support is simply not yet converged, which is confirmed by the red solid curve in \figref{fig:NumResFF-Conv}, so that improvements might be achieved by additional iterations. However, this would be expensive and potentially unstable as the observed convergence is slow and the number of CG-iterations per Newton step  increases rapidly as the regularization parameter further decreased.
	  
	  Notably, the artifacts in the center image of \figref{fig:NumResFF-DiffSupps-a} may not be attributed to trivial ambiguities since there is no finite translation of the object that is consistent with the tight rectangular support constraint. Moreover, comparing the right and the left image in \figref{fig:NumResFF-DiffSupps-a}, we find that the reconstruction \emph{without} support constraint for the cell-shaped initial guess comes out nearly as accurate as in the case where the exact support is incorporated. This indeed suggests that the prescribed support is not the significant factor
	  here but that the reconstruction quality depends most sensitively on the choice of the initial guess $\bN_0$.
	  
	  This can be qualitatively understood by considering the corresponding initial data $F_{\Text{dis}}(\bN_0)$ shown in \figref{fig:NumResFF-DiffSupps-c}: while the far-field data for the cell-shaped support bump already contains many features of the observed data (see \figref{fig:NumResFF-NoSupp-d}), the rectangle as an initial guess gives rise to quite different structures in Fourier space.
	  This is highly relevant to the outcome of the applied Newton-type method as the iterates solve \emph{linearizations} of the forward problem of the form
	  \begin{align}
	   F_\infty( N_{k+1} ) &\approx F_\infty( N_{k} ) + F_\infty'[N_{k}] (N_{k+1}- N_{k}) \nonumber \\
	   &= F_\infty( N_{k} ) + \cc{\cF_{\overline 2} ( G( N_{k} ) )} \cdot \cF_{\overline 2} (T (N_{k+1} - N_{k})) \label{eq:FFFrechetFormal} ,
	  \end{align}
	  where $G,T$ are some operators, compare \eqref{eq:FrechetForwOpFF}. The bilinear structure in \eqref{eq:FFFrechetFormal}
	  implies that the Newton step will be inaccurate in all Fourier frequencies that are underrepresented by $\cF_{\overline 2} (G( N_{k}))$ as the \Frechet derivative practically vanishes in these components, i.e.\ is not stably invertible. Consequently, the initial guess needs to be chosen such that it sufficiently populates Fourier space,  presampling the unknown object's spectrum.
	  On the other hand, overestimated Fourier modes in the initial guess $\bN_0$ will take many iterations to correct
	  as the deviations of the iterates $\bN_{k+1}$ are limited by the regularization term $ \sim \alpha_k \norm{ \bN_{k+1} - \bN_{0} }_2^2$ incorporated into the Newton step \eqref{eq:DiscNewtonStepPCT}.
	  For the latter reason, choosing random noise as $\bN_0$ neither represents a viable alternative. It should be emphasized that the bilinear structure in \eqref{eq:FFFrechetFormal} does not arise from nonlinearity of the object transmission function but from that of the squared modulus. The resulting peculiarities therefore also manifest in far-field reconstructions of \emph{weak} objects, governed by the operator in \eqref{eq:WeakFwOpFF}. We summarize our findings by stating the following result:
	\vspace{1em}
	\begin{res}[Ab Initio Far-Field Tomography by Regularized Newton Methods] \label{res:FFTomoByNewton}
	  	Far-field tomographic reconstructions by the Newton-type \algref{alg:PCT} depend sensitively on the chosen initial guess. If the choice does not sufficiently reflect the structure of the unknown object to be reconstructed, severe artifacts result. In particular, \emph{ab initio} reconstructions without specific a priori knowledge require iteratively updated structural (support-)estimates to improve the initial guess on the fly.
	\end{res}
	\vspace{1em}
	
	  Certainly, it is also possible that the poor reconstruction result for the rectangular support in \figref{fig:NumResFF-DiffSupps} are manifestations of non-trivial phase retrieval ambiguities rather than being only due to issues with the Newton-type approach. However,  note that the considered numerical test cases for far-field phase contrast tomography resemble \emph{two}-dimensional phase retrieval, for which non-trivial ambiguities are known to be ``pathologically rare'' \cite{Barakat1984}. On the other hand, it should be emphasized that independent phase retrieval of the \emph{one}-dimensional diffraction patterns for the different incident angles $\theta$ is expected to suffer from severe non-uniqueness according to the analysis of the 1D phase retrieval problem in \sref{SS:PhaseRetrFF}, compare \resref{res:1DFTPhaseRetrIllposed}. Hence, it can be regarded as a first \emph{proof of concept} for the pursued simultaneous phase- and tomographic reconstruction approach (see \sref{S:CombinedApproach}) that at least the results in \figref{fig:NumResFF-DiffSupps} starting from the exact support shape are artifact-free.

	\end{subsection}

        \begin{subsection}{2D Reconstructions using Reference Signals} \label{SS:NumResFF-WithReference}

	  In \sref{SS:NumResFF-NoReference}, we have found that far-field tomography by the considered regularized Newton-type \algref{alg:PCT} requires a good initial guess incorporating dominant features of the object to be reconstructed, e.g.\ its exact shape. As iteratively updated support estimates may not be implemented in a straightforward manner and would require a large number of Newton iterations, we investigate a different approach inspired by the uniqueness results in \thmref{thm:UniquePhaseFF1D} and \thmref{thm:UniquePhaseFFMD}: if the unknown object is embedded in a known \emph{reference signal}, then the exact support of the superposition is certainly accessible as that of the reference. Moreover, the latter provides a canonical choice for the initial guess incorporating characteristic features of the total object. Two major issues of ab initio far-field tomography identified in the preceding section are thus resolved by the approach.
	  
	  In order to study its potential, we compute far-field reconstructions for the cell phantom superimposed with reference signals $\bN_0$ of different shape: a rectangle, a circle and a non-pointsymmetric bullet-shaped bump as well as a discretized version of the 2D-exponential ramp (see \eqref{eq:DefBaseMD}) motivated by \thmref{thm:UniquePhaseFFMD}. We scale these and the (pure phase object) phantom $\Delta \bN$ to have equal magnitude, such that we obtain total objects
	  \begin{equation} \bN^\dagger = \underbrace{\bN_0}_{\substack{\text{reference} \\ \text{=initial guess}}} + \underbrace{\Delta \bN}_{\substack{\text{unknown} \\ \text{phantom}}}  \MTEXT{of magnitude}   \norm{\bN^\dagger} \in \{ 0.1, \pi \}. \end{equation} 
	  Hence, we investigate both \emph{weak} objects and moderately strong ones. Support constraints are imposed according to the exact support of the reference signal $\bN_0$, which defines the initial guess. In order to ensure a fair comparison, the support sizes are chosen such that an equal number of object pixels has to be reconstructed in each of the test cases.
	  All other parameters of the simulation setup - in particular the error level of $\varepsilon = 10^{-3}$ - are chosen exactly is in the preceding test cases.
	  
	  The resulting reconstructions after 10 Newton iterations are shown in \figref{fig:NumResFF-DiffRefs}. A first surprising aspect to note is the poor quality of the results for the rectangular reference signals (cf. \figref{fig:NumResFF-DiffRefs-a}-d) - even more so as the choice of the exponential ramp shown in \figref{fig:NumResFF-DiffRefs-c}-d is theoretically motivated by the uniqueness result \thmref{thm:UniquePhaseFFMD}. In fact, all of these reconstructions turn out to be subject to considerable artifacts, which are however significantly stronger in the weak object case $\norm{\bN^\dagger} = 0.1$ visualized in the left column of \figref{fig:NumResFF-DiffRefs}.
	  
	  Concerning uniqueness, the exponential ramp on the other hand does show some benefit compared to the constant rectangular reference: in
	  the weak object case, the point symmetry of the latter allows for a manifestation of the \emph{twin-image} of the cell phantom (see \sref{SSS:TrivialAmbiguities} - Trivial Ambiguities), which fades out only slowly in the course of the Newton iterations. Remainders of the point-reflected twin-image can still be identified in \figref{fig:NumResFF-DiffRefs-a} along with other artifacts. Owing to the symmetry breaking by the initial exponential profile, this is not the case in \figref{fig:NumResFF-DiffRefs-c}, i.e.\ the choice of the reference signal eliminates the twin-image ambiguity. 
	  \floatsetup[figure]{style=plain,subcapbesideposition=top}
	        \begin{figure}[H]
	  \centering
	  \sidesubfloat[]{ \includegraphics[height=.34\textwidth]{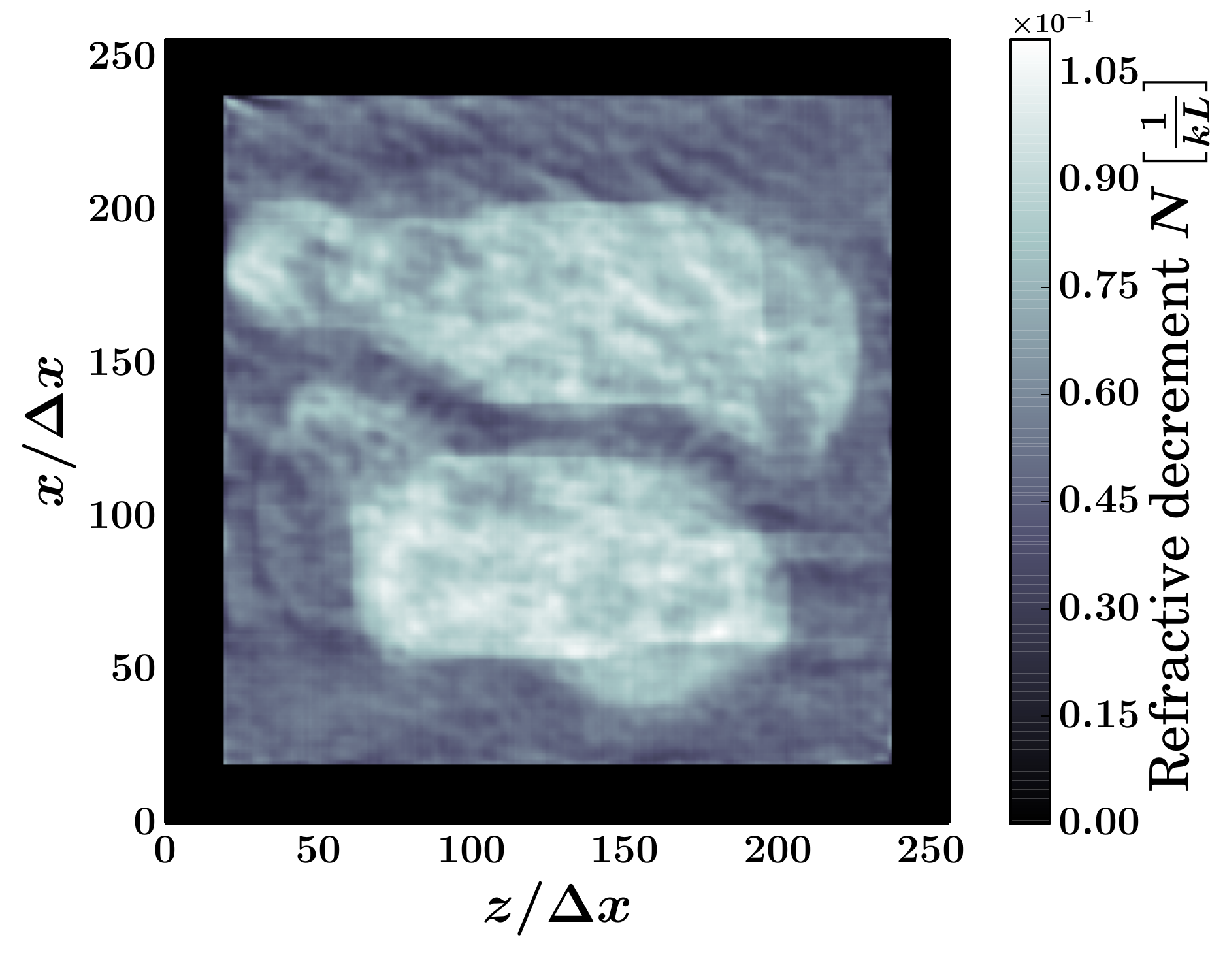} \label{fig:NumResFF-DiffRefs-a} }
	  \hfill
	  \sidesubfloat[]{ \includegraphics[height=.34\textwidth]{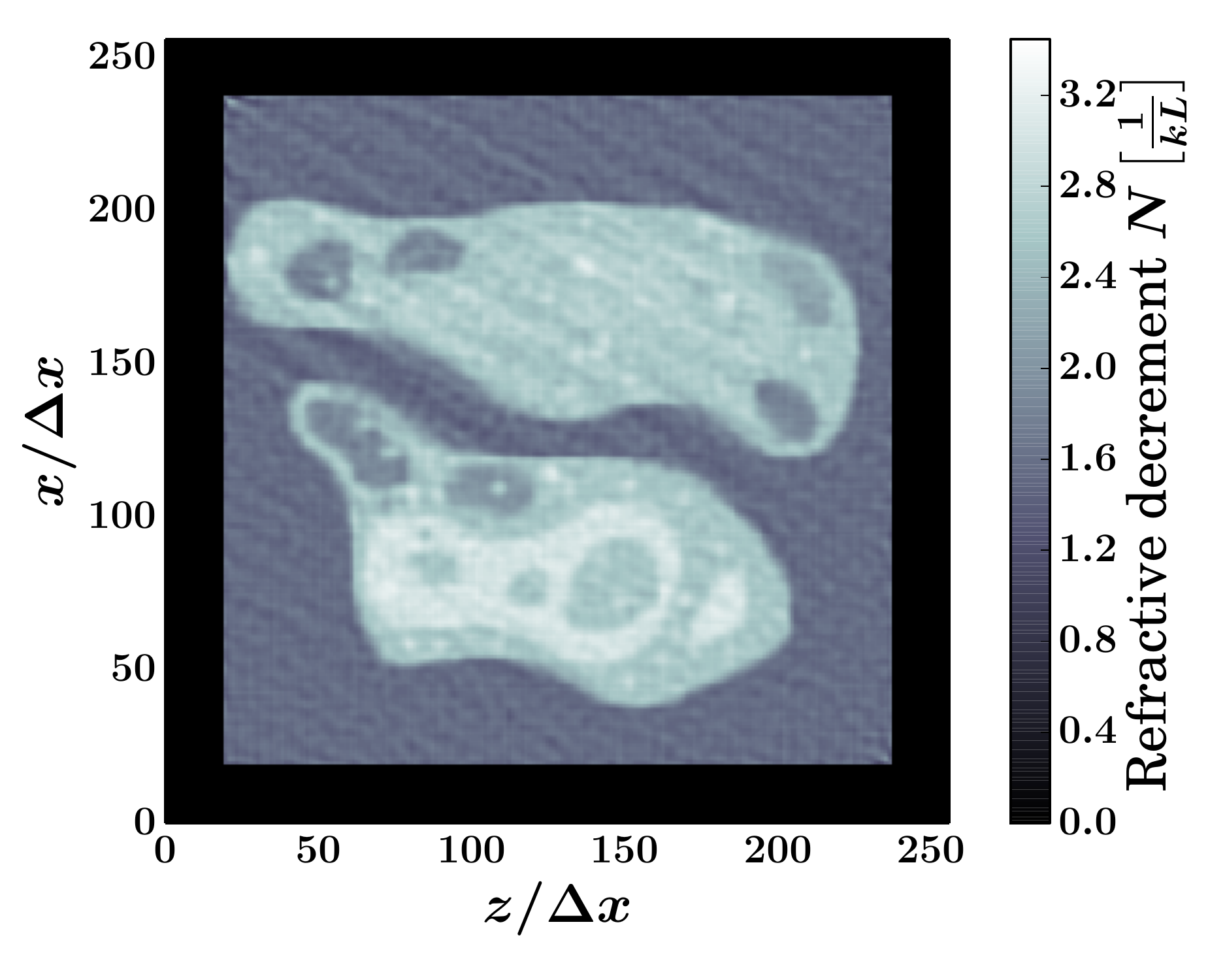} \label{fig:NumResFF-DiffRefs-b} } \vspace{-.5em}
	  \\
	  \sidesubfloat[]{ \includegraphics[height=.34\textwidth]{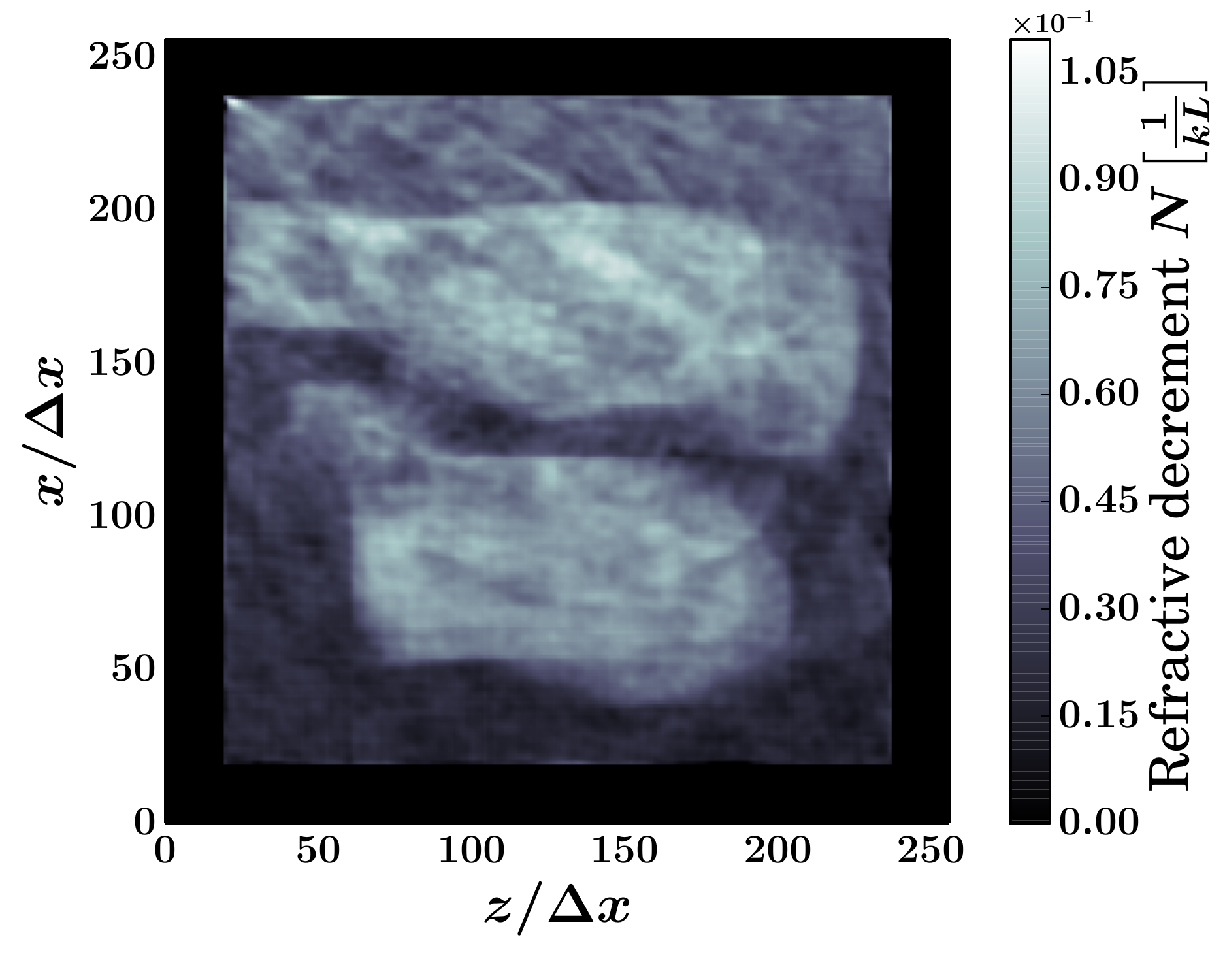} \label{fig:NumResFF-DiffRefs-c} }
	  \hfill
	  \sidesubfloat[]{ \includegraphics[height=.34\textwidth]{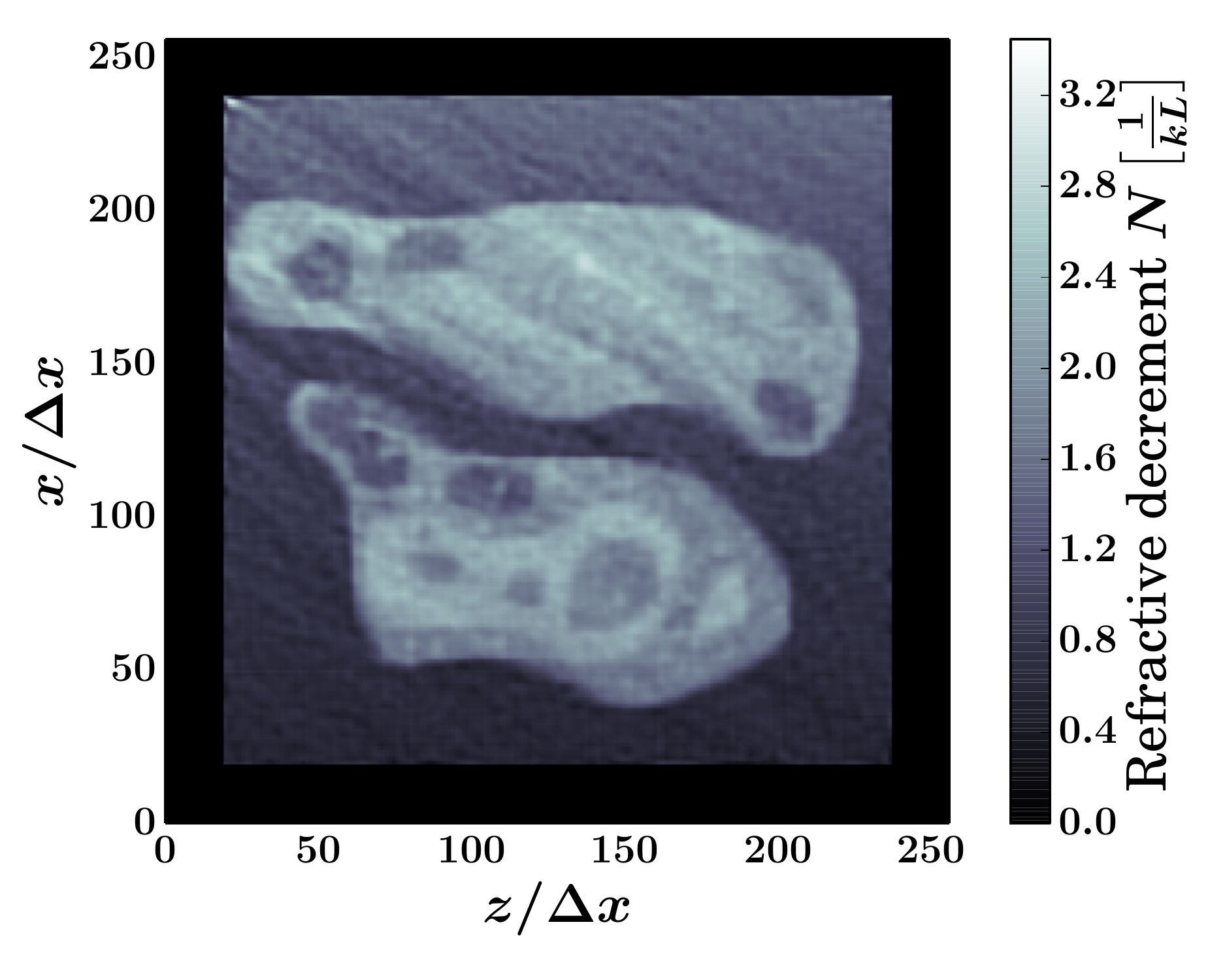} \label{fig:NumResFF-DiffRefs-d} }\vspace{-.5em}
	  \\
	  \sidesubfloat[]{ \includegraphics[height=.34\textwidth]{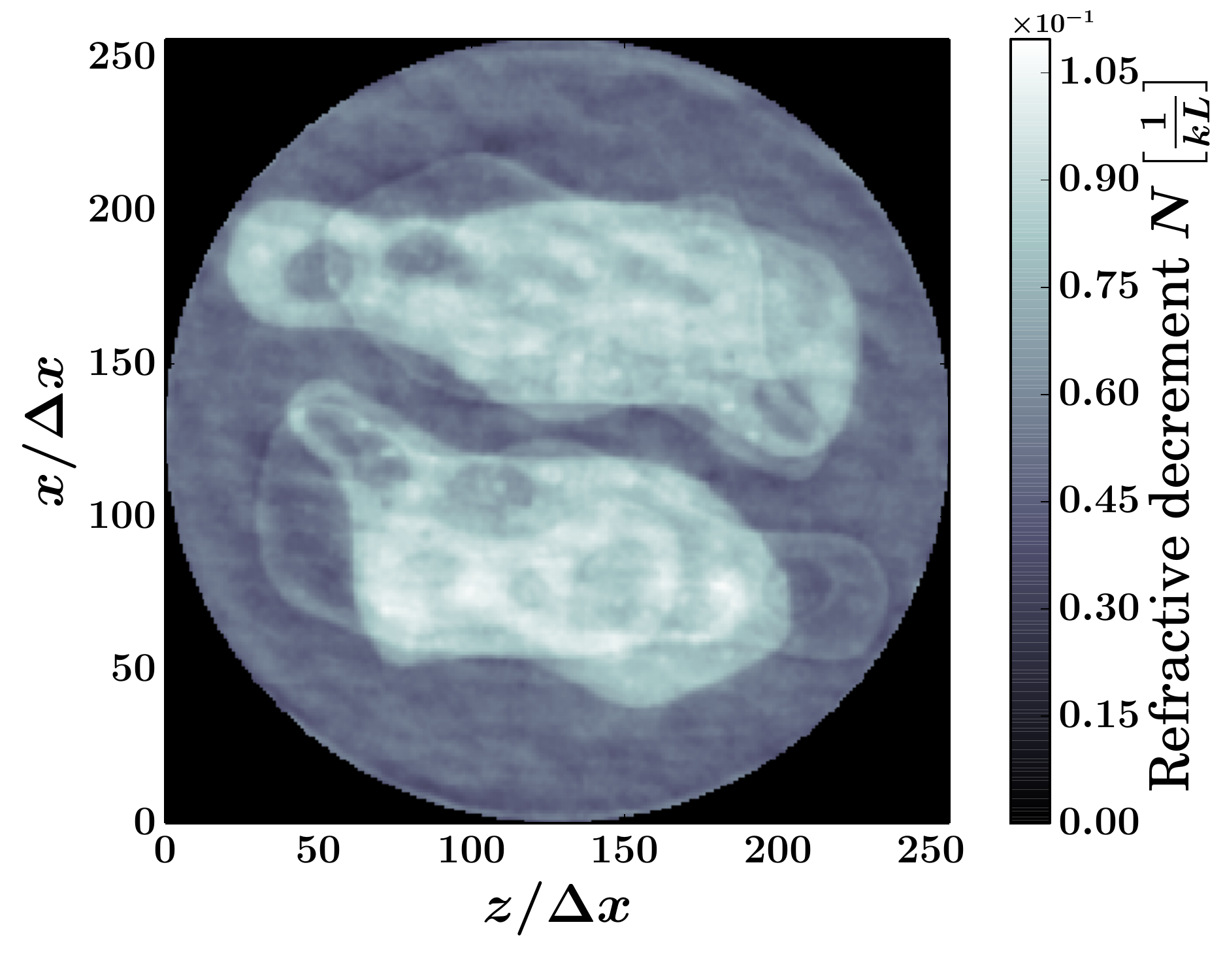} \label{fig:NumResFF-DiffRefs-e} }
	  \hfill
	  \sidesubfloat[]{ \includegraphics[height=.34\textwidth]{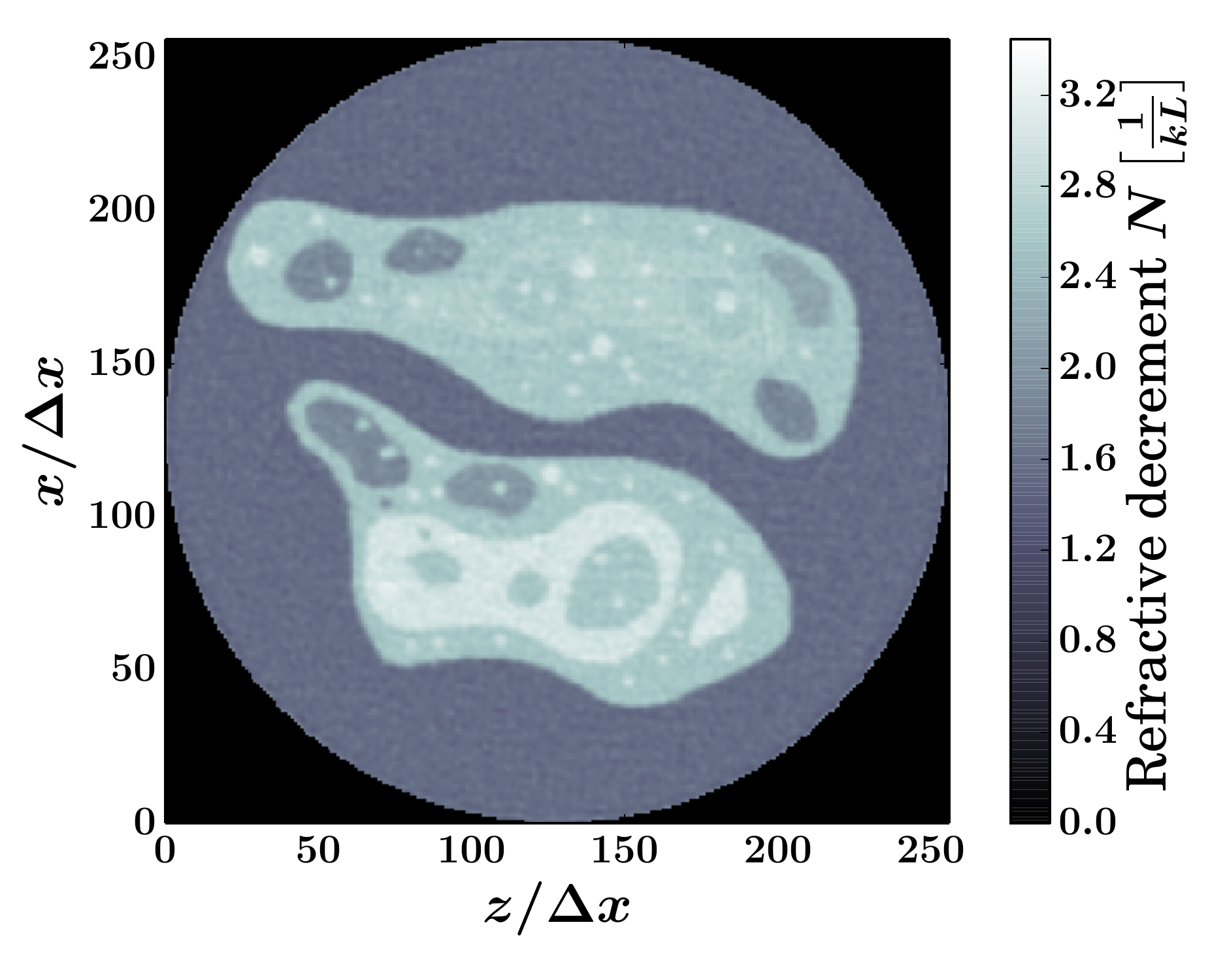} \label{fig:NumResFF-DiffRefs-f} } \vspace{-.5em}
	  \\
	  \sidesubfloat[]{ \includegraphics[height=.34\textwidth]{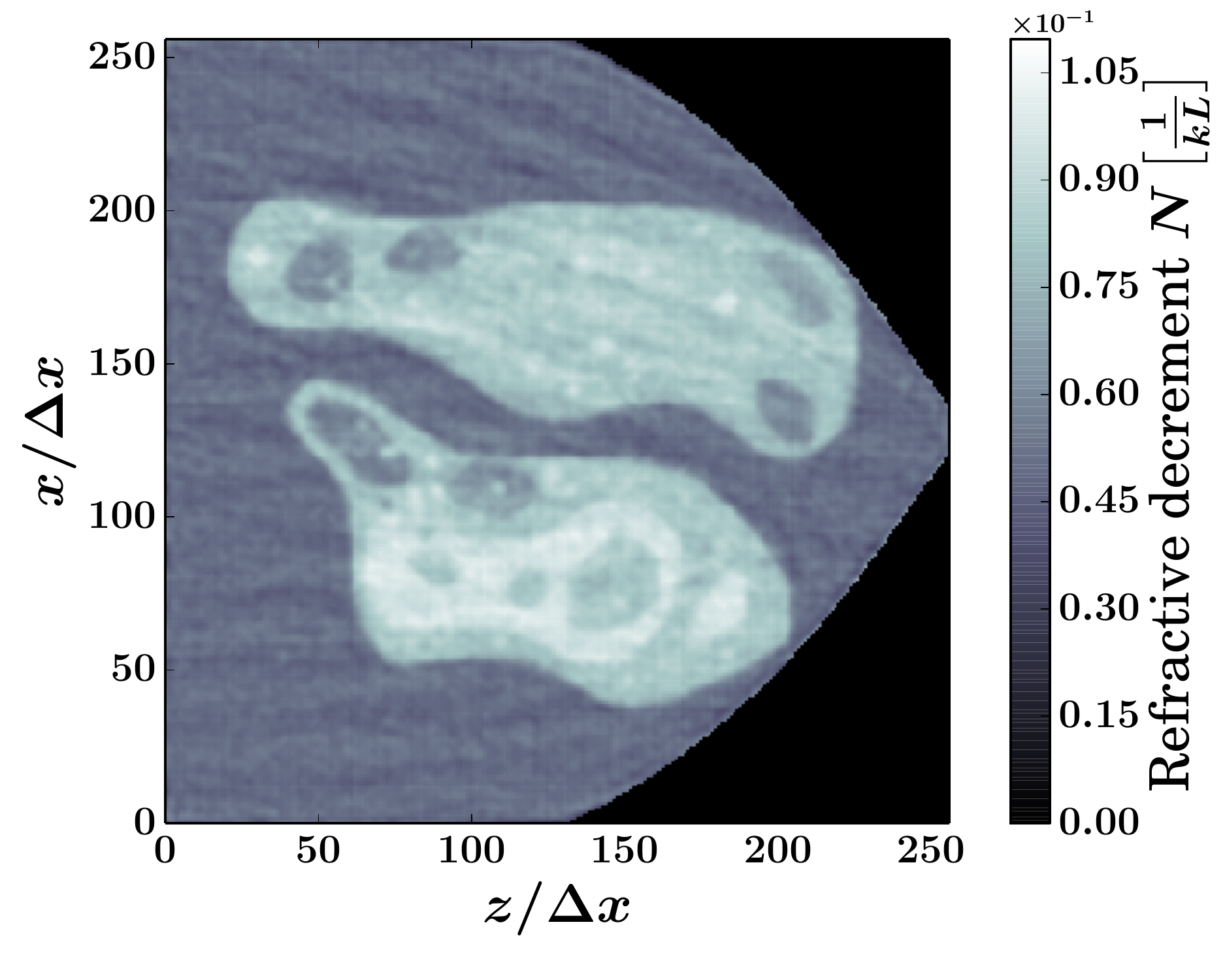} \label{fig:NumResFF-DiffRefs-g} }
	  \hfill
	  \sidesubfloat[]{ \includegraphics[height=.34\textwidth]{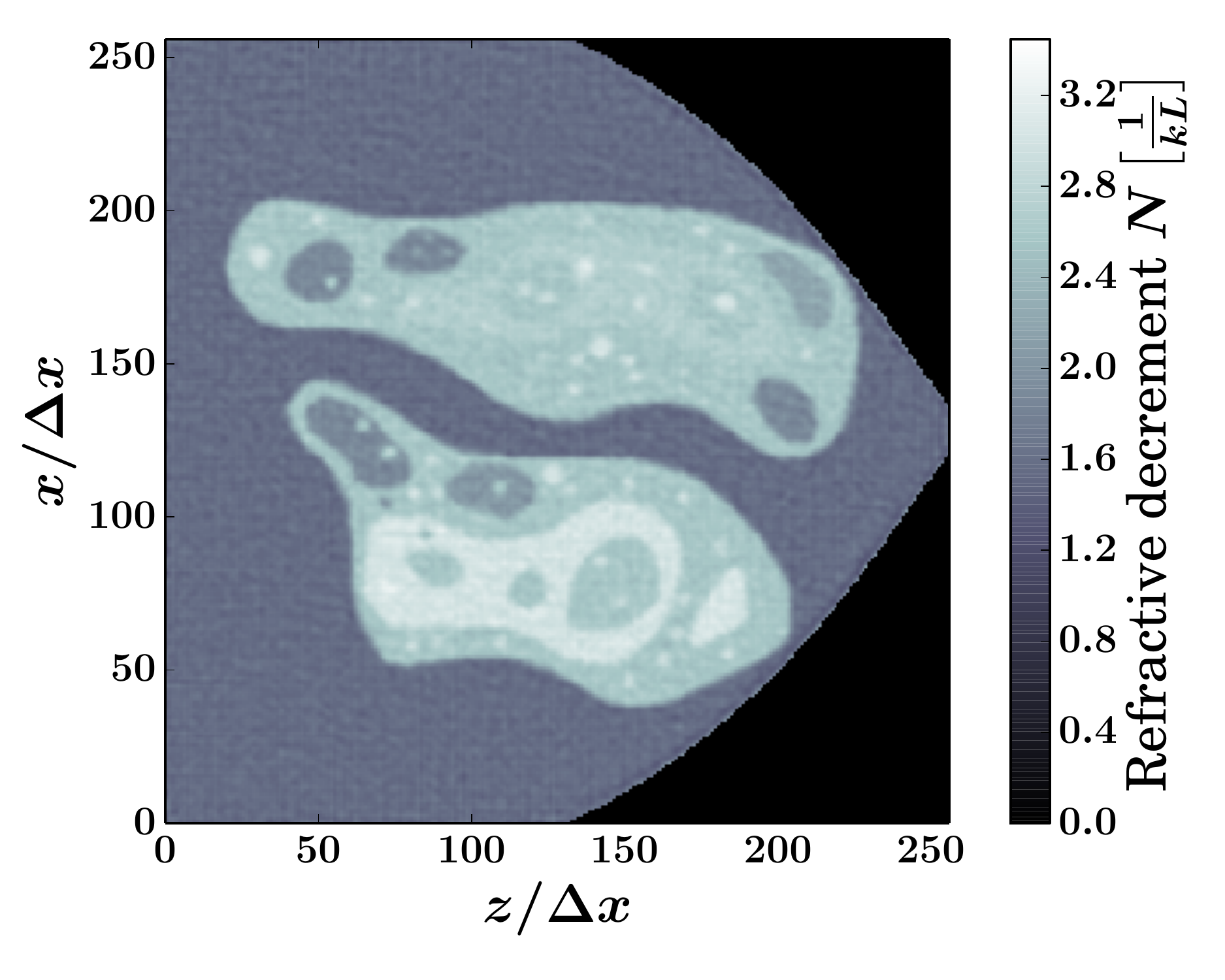} \label{fig:NumResFF-DiffRefs-h} } \vspace{-.5em}
	  \caption{Reconstructed objects $\bN_{10}$ for 2D far-field tomography with different (known) reference signals $\bN_0$ superimposed upon the unknown cell-shaped phantom $\Delta \bN$. $\bN_0$ varied from top to bottom: constant rectangular bump (= scaled indicator function), exponential ramp from \thmref{thm:UniquePhaseFFMD} varying in $x$-direction, constant circular bump, constant bullet-shaped bump. Reference and unknown phantom are scaled to pure phase objects of equal magnitude $\norm{ \Delta \bN } = \norm{ \bN_0 }$ such that $\bN^\dagger = \bN_0 + \Delta \bN$ satisfies $\norm{\bN^\dagger} = 0.1$ (weak object, left column) or $\norm{\bN^\dagger} = \pi$ (moderately strong, right column). Support constraint and initial guess according to $\bN_0$. Other parameters analogous to simulation in \figref{fig:NumResFF-NoSupp}. See \tabref{tab:NumResFFRefErrors} for reconstruction errors. \label{fig:NumResFF-DiffRefs}}
	  \end{figure}
	  
	  Notably, the cell phantom comes out much clearer for the non-rectangular reference signals in \figref{fig:NumResFF-DiffRefs-e}-h. This visual impression is supported by the relative reconstruction errors ${ \norm{ \bN_{10} -  \bN^\dagger }_2 } / { \norm{ \bN^{\dagger} -  \bN_0 }_2}$
	  summarized in \tabref{tab:NumResFFRefErrors}. Here, the known reference component $\bN_0 $ in the object is subtracted since we are interested in the accuracy of the reconstructed \emph{unknown} part. 
	  
	  According to the errors in \tabref{tab:NumResFFRefErrors}, the reconstruction quality depends strongly on the shape of the chosen reference signal. This can be understood by a similar argument as the dependence on the initial guess in \sref{SS:NumResFF-NoReference}: rectangles correspond to sinc-functions in Fourier space, giving rise to a sparse and highly anisotropic presampling of the intensity data, compare \figref{fig:NumResFF-DiffSupps-c}. Accordingly, when starting a reconstruction from a rectangular reference signal as in \figref{fig:NumResFF-DiffRefs-a}-d, only those spatial frequencies representing directions of the edges are strongly pronounced in the initial data. Due to the bilinear structure of the underlying linearization \eqref{eq:FFFrechetFormal}, the Newton method may thus hardly be accurate in the underrepresented bulk Fourier components in between, which are in turn relevant for the non-rectangular phantom. This explanation is supported by the anisotropy of the observed artifacts in \figref{fig:NumResFF-DiffRefs-a}-d. On the contrary, the circular reference signal gives rise to \emph{isotropic} initial far-field data which seems to stabilize the reconstruction considerably, although it retains twin-image symmetry as seen from \figref{fig:NumResFF-DiffRefs-e}. The more complex bullet-shaped reference signal breaks this symmetry according to \figref{fig:NumResFF-DiffRefs-g}, but yields slightly less accurate reconstruction results (see \tabref{tab:NumResFFRefErrors}) as the anisotropy of the support shape once more promotes artifacts along certain directions.
	  
	  Another insight from \figref{fig:NumResFF-DiffRefs} and \tabref{tab:NumResFFRefErrors} is that the reconstruction \emph{improves} as the nonlinearity of the object transmission function (OTF) comes into play for stronger objects - although error levels and regularization are chosen comparably! Indeed, it seems that the assumption of a pure phase object provides a stronger constraint outside the linear regime of the weak object limit (see \sref{SS:SpecialObjects}).

	    \begin{table}[htb!]
	  \centering
	  \begin{tabular}{ccccc} 
	    \toprule
	         &  Rectangle (a,b) & Exp. ramp (c,d) & Circle (e,f) & Bullet (g,h) \\ 
	    \midrule
	      Weak object  &  $30.4 \, \%$ &  $31.7 \, \%$  &  $24.6 \, \%$   &  $18.7 \, \%$ \\
	      Strong object   &  $16.3 \, \%$  &  $19.7 \, \%$  &  $12.1 \, \%$   &  $13.0 \, \%$ \\
	    \bottomrule
	  \end{tabular}
	  \caption{Relative reconstruction errors $\frac{ \norm{ \bN_{10} -  \bN^\dagger }_2 }{ \norm{ \bN^{\dagger} -  \bN_0 }_2}$ w.r.t. the unknown cell $ \Delta \bN = \bN^{\dagger} -  \bN_0$ for the 2D far-field test cases with different reference signals visualized in \figref{fig:NumResFF-DiffRefs} (subfigures a-h). The weak and strong object cases are characterized by  magnitudes $\norm{\bN^\dagger} = 0.1$ and $\norm{\bN^\dagger} = \pi$, respectively.}
	  \label{tab:NumResFFRefErrors}
	\end{table}
	   This conjecture is true at least for the case of twin-image ambiguities, which are visible in \figref{fig:NumResFF-DiffRefs-a} and -\ref{fig:NumResFF-DiffRefs-e} but not in the corresponding reconstructions of moderately strong objects (\figref{fig:NumResFF-DiffRefs-b} and -\ref{fig:NumResFF-DiffRefs-f}).
	  This difference can be understood by studying the underlying symmetry of the far-field intensities (compare \eqref{eq:ForwardOpFF})
	  
	  \begin{equation}
	   I_\infty = | \cF_{\overline 2} ( O_0 ( N ) ) |^2, \label{eq:FFIntensitisReminder} 
	  \end{equation}
	  being the invariance under complex conjugation and lateral reflection of the OTF
	  \begin{equation} O_0(N) = \exp( - \I k \CR( \N ) ) - 1 \mapsto \exp(  \I k \CR( \N )_r ) - 1  = O_0( -N_r) \end{equation}
	  with $N_r(\bx) := \cc{N_r(-\bx)}$, see \sref{SSS:TrivialAmbiguities}. Accordingly, the twin-image $-N_r$ is always an equally valid solution to the tomographic phase reconstruction problem, which is however \emph{negative} for real and positive $N$ and hence suppressed in the considered numerical reconstructions by the positive initial guess. On the other hand, \eqref{eq:FFIntensitisReminder} is also trivially invariant under a change of the sign of $ O_0 ( N )$. In the weak object limit \eqref{eq:WeakObjectCond} where $O_0$ is (approximately) linear, this implies that the \emph{positive} twin-image $N_r$ also constitutes an approximate solution the far-field phase problem since  $O_0( N_r) \approx -O_0(-N_r)$. For stronger objects, the latter symmetry is broken by higher order contributions in the exponential OTF so that the positive twin-image  is suppressed by nonlinearity. This nonlinear symmetry breaking is also the reason why the twin-image is not observed in the ab initio reconstructions in \sref{SS:NumResFF-NoReference}.
	  
	  As a conclusion, we summarize the observations of this section in the form of the following result:
	\vspace{1em}
	\begin{res}[Newton-based Far-Field Tomography with Reference Signals] \label{res:FFTomoWithReference}
		By superposition of known reference signals, far-field tomography via \algref{alg:PCT} may accurately reconstruct unknown pure phase objects without support knowledge. The quality of the result depends strongly on the shape of the reference object, where isotropic support geometries seem preferable to rectangles. Moderately strong objects tend to be reconstructed more stably than weak ones owing to the nonlinearity of the object transmission function, breaking for instance  twin-image symmetry.
	\end{res}        
        \end{subsection}
        
        \newpage
        
        \begin{subsection}{3D Reconstructions from Realistic Data} \label{SS:NumResFF-3DRecon}
        
		Having identified algorithmic peculiarities and -remedies of our Newton-type approach to far-field phase contrast tomography, the final endeavor of this numerical study is to prove the method's applicability to realistic data sets. Most importantly, we investigate the physically relevant case of \emph{three}-dimensional tomography instead of the 2D toy model considered before. As seen in the preceding sections, ab initio reconstructions with \algref{alg:PCT} require a known reference to achieve reasonable accuracy. Unfortunately, experimental tomographic far-field data of such a specific form is not available. For a proof of concept in realistic settings, we therefore design a numerical simulation incorporating principal features of experimental data:
		\vspace{1em}
		\begin{itemize}
		 \item \emph{General phantom} $\Delta \bN = \Delta \bdelta - \I \Delta \bbeta$: Ensemble of randomly shaped ellipsoids, each of constant refractive index $\Delta \bdelta_j, \Delta \bbeta_j$  drawn from a normal distribution:
                 \begin{equation} kL  \Delta \bdelta_j \sim \cN(\mu, \sigma\mu) \MTEXT{and} kL \Delta \bbeta_j \sim c_{\bbeta/\bdelta} \cN(\mu, \sigma\mu ) \end{equation}
		 \item \emph{Reference signal} $\bN_0 = \bN^\dagger -  \Delta \bN$: Uniform sphere (pure phase object)
		 \item \emph{Beam stop}: Centre of diffraction patterns (frequencies $\bxi_{\Text{dis}}$ with $\norm{\bxi_{\Text{dis}}} < \frac{\pi}{30 \Delta x}$) excluded from data fit as usually not measurable (see \figref{fig:FarfieldBeamStop})
		  \item \emph{Missing wedge}: Incident angles limited to $\theta \in [0^{\circ}; 160^{\circ})$, e.g.\ by obstructing instruments in the setup
		  \item \emph{Poisson noise}: According to an average count of $92\,\frac{\Text{photons}}{\Text{pixel}}$ observed in \cite{Chapman2006}
		  \item \emph{Stop rule}: Newton iterations stopped according to the discrepancy principle \eqref{eq:Discrepancy} with $\tau = 1$
		\end{itemize}
		\vspace{1em}

                
                It may seem counter-intuitive in the phantom definition that we fix a coupling constant $c_{\bbeta/\bdelta}$ between refraction and absorption, which is merely perturbed by Gaussian deviations, although we are to simulate \emph{general} objects. Note, however, that this model is more realistic than assumption of entirely uncorrelated $\Delta \bdelta $ and $\Delta \bbeta$ as both parameters should actually correspond to structures of one and the same real-world specimen. The piecewise constant ratio of refraction and absorption corresponds to separated parts of the object that are composed of different materials.
                
                The spherical reference is not only chosen because it provides the best reconstruction results according to \figref{fig:NumResFF-DiffRefs} - it is also favorable from an experimental point of view (provided that accurate manufacturing is feasible): reference objects may be introduced in experimental setups by placing the former in front of the unknown specimen in the beam line. Within the framework of the projection approximation (see \sref{SS:ProjApprox}), the contact image will then look as if the unknown object was enclosed by the reference. The advantage of a uniform sphere is that rotational alignment can be completely omitted by symmetry. In particular, only the unknown sample needs to be rotated for tomographic measurements - other than for a cuboid reference.
	        \begin{figure}[hbt!]
	  \centering
	  \subfloat{\includegraphics[height=.44\textwidth]{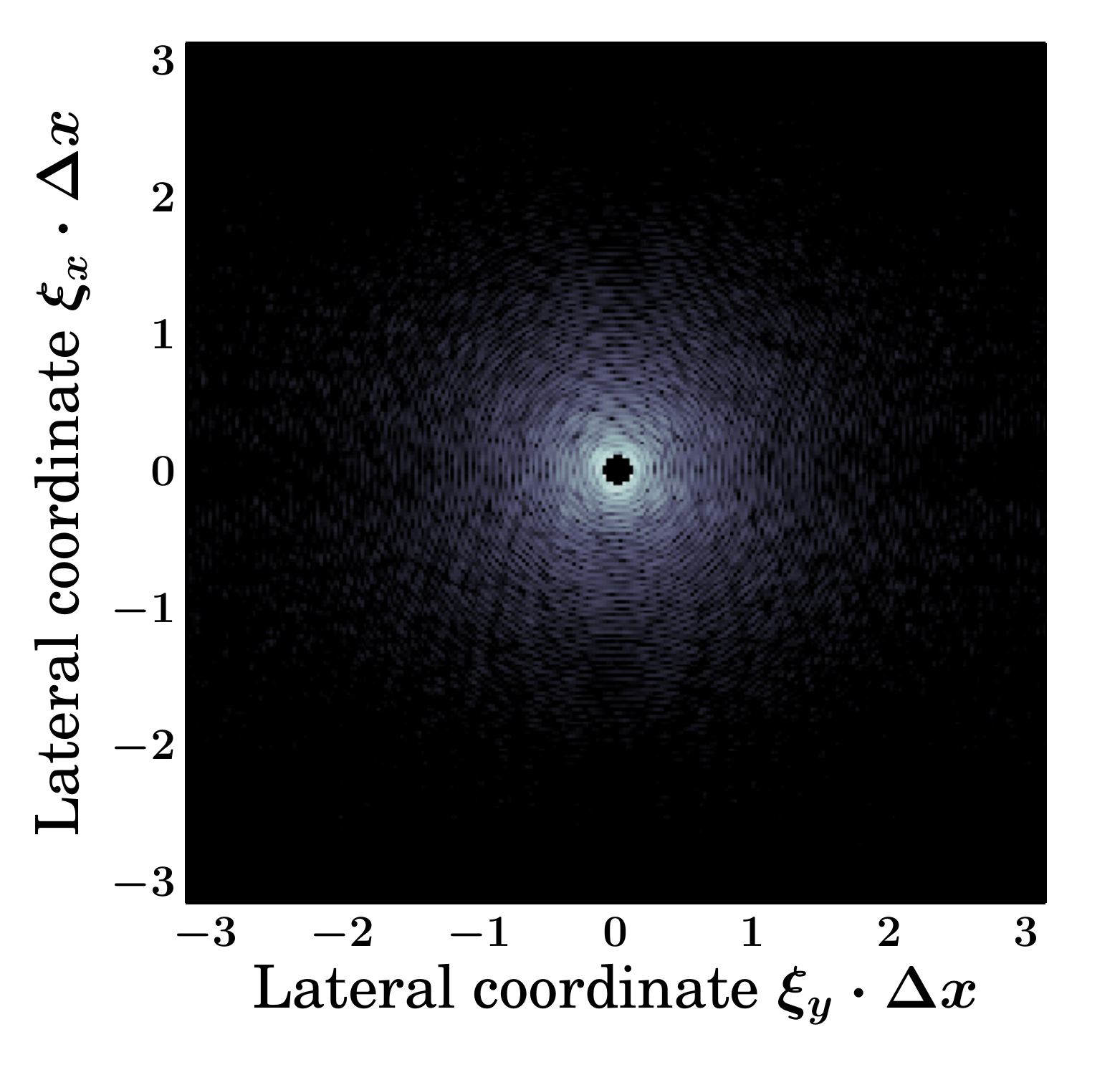} \label{fig:NumResFF-NonIdealData-a} }
	  \hfill
	  \subfloat{\includegraphics[height=.44\textwidth]{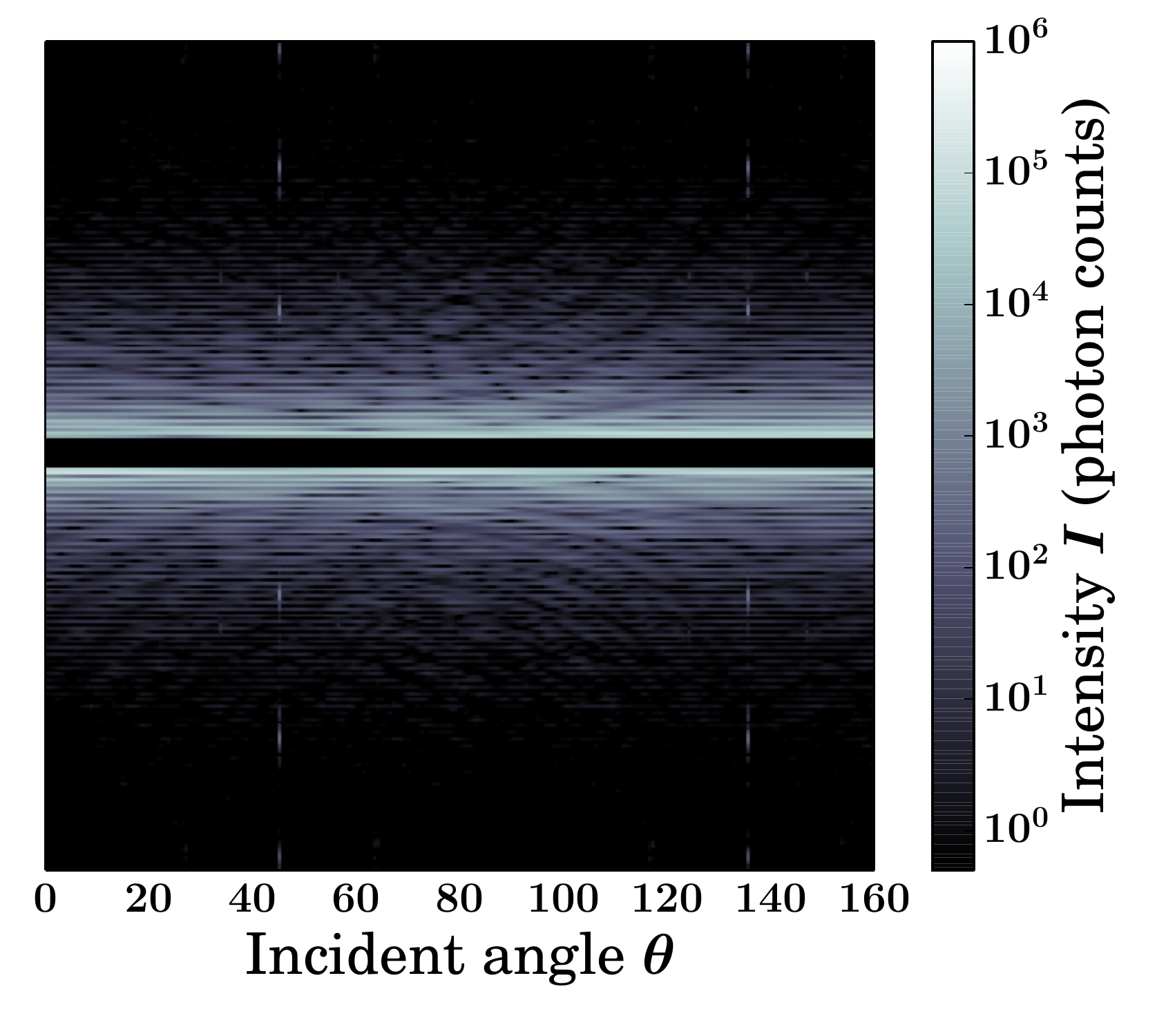} \label{fig:NumResFF-NonIdealData-b} }
	  \caption{Simulated Data $\bI^{\Textbf{err}}$ computed by \eqref{eq:DataGenFF} from the exact object shown in \figref{fig:NumRes-FF3D-exObj}. Left: lateral diffraction pattern for $\theta = 0$, right: far-field ``sinogram'' at $\xi_y = 0$. Beam stop (black shading) and photon count scaling according to experimental measurements in \cite{Chapman2006} (total counts for all angles $\approx 1.5\cdot10^{9}$). Note the missing wedge of $20^\circ$ in the simulated incident angles. \label{fig:NumResFF-NonIdealData}}
	  \end{figure}
	       \begin{table}[htb!]
  \centering
  \begin{tabular}{ccccccc} 
    \toprule
        $\mX_{\Text{dis}}$ &  $\mY_{\Text{dis}}$   & $\cG_{\mY_{\Text{dis}}}$  & $\alpha_0$ & Init. guess $\bN_0$ & Incident $\theta $  & Constraints \\
    \midrule
         $\mC^{128^3}$ & $\mR^{256^3}$   &   $\substack{\text{\scriptsize{\eqref{eq:GramYDisc}}} \\  (0\text{ at beam stop})}$  & $ 10^{16}$ & $\substack{\text{reference sphere} \\  \text{(pure phase obj.)}}$ & $[0^{\circ}; 160^{\circ})$ & \scriptsize{support of $\bN_0$} \\ 
    \bottomrule
  \end{tabular}
  \caption{Simulation parameters for realistic 3D far-field tomography test case via \algref{alg:PCT}. The exact test object $\bN^\dagger \in \mX_{\Text{dis}} = \bN_0 + \Delta \bN$ is a superposition of a purely phase shifting reference sphere plus a random ensemble of ellipsoids $\Delta \bN = \Delta \bdelta - \I \Delta \bbeta$ with $\Delta \bbeta / \Delta \bdelta = 0.1 \pm 30 \, \%$. Results shown in \figref{fig:NumResFF-NonIdealRecon}. Non-specified parameters according to \sref{SS:NumResFF-GenSetup} and \tabref{tab:NumResFFSetup}.}
  \label{tab:NumResFF3DSetup}
\end{table}
	  
          We choose $\sigma = 0.3$ and $c_{\bbeta/\bdelta} = 0.1$ (corresponding to rather strong absorption) and scale $\Delta \bN$ and $\bN_0$ to have equal magnitude such that $\norm{\bN^\dagger} = \norm{ \bN_0 + \Delta \bN } = \pi$ as in \sref{SS:NumResFF-WithReference}. The resulting test object is visualized in \figref{fig:NumRes-FF3D-exObj}. A resolution of $128^3$ voxels is prescribed in object space and the intensity data is simulated on $256^2$
          
          \begin{figure}[H]
	  \centering
	  \subfloat[Exact object $\bN^\dagger = \bdelta^\dagger - \I \bbeta^\dagger$]{\includegraphics[height=.32\textwidth, clip=true, trim = -4cm 0cm -4cm 0cm]{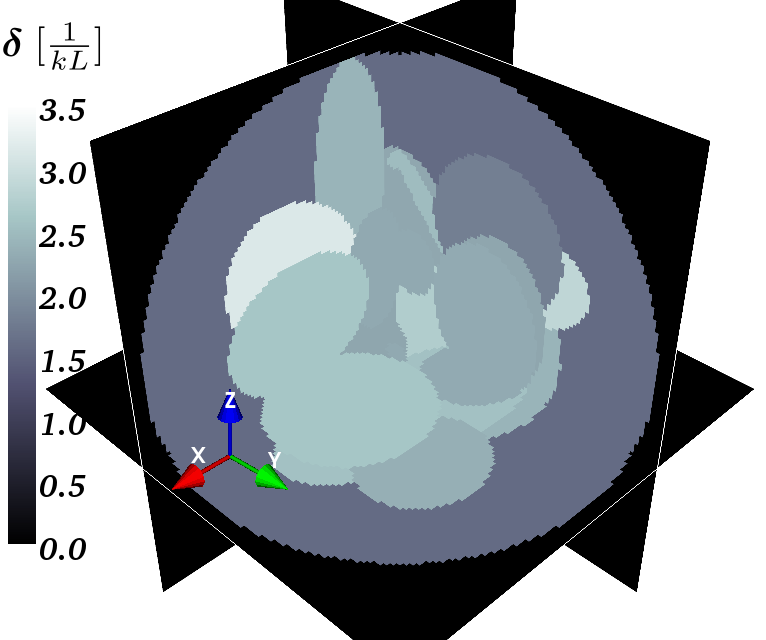} \hfill \includegraphics[height=.32\textwidth, clip=true, trim = -4cm 0cm -4cm 0cm]{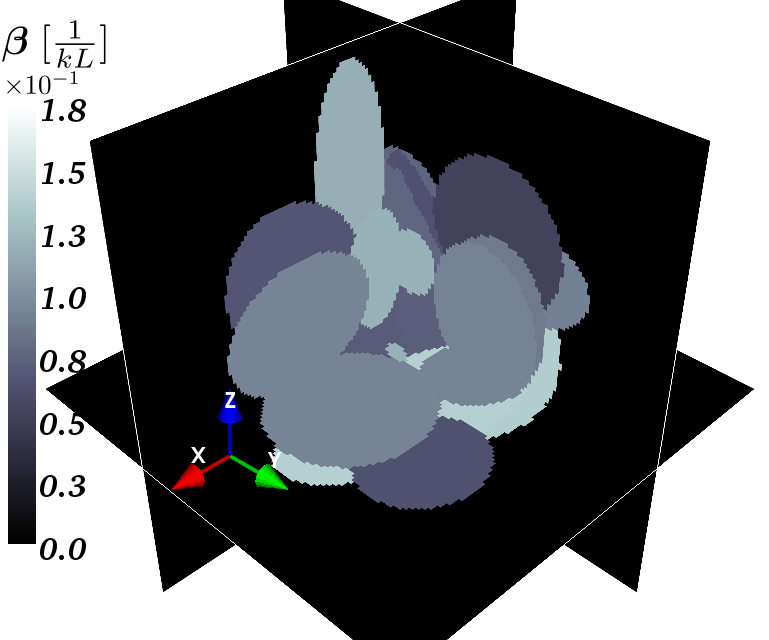} \label{fig:NumRes-FF3D-exObj} }
	  \\
	  \subfloat[Reconstruction $\bN_{k_{\Text{stop}}}^{\Text{phase}} = \bdelta_{k_{\Text{stop}}}^{\Text{phase}}$ with a pure phase object constraint ($\bbeta = 0$)]{\includegraphics[height=.32\textwidth, clip=true, trim = -4cm 0cm -4cm 0cm]{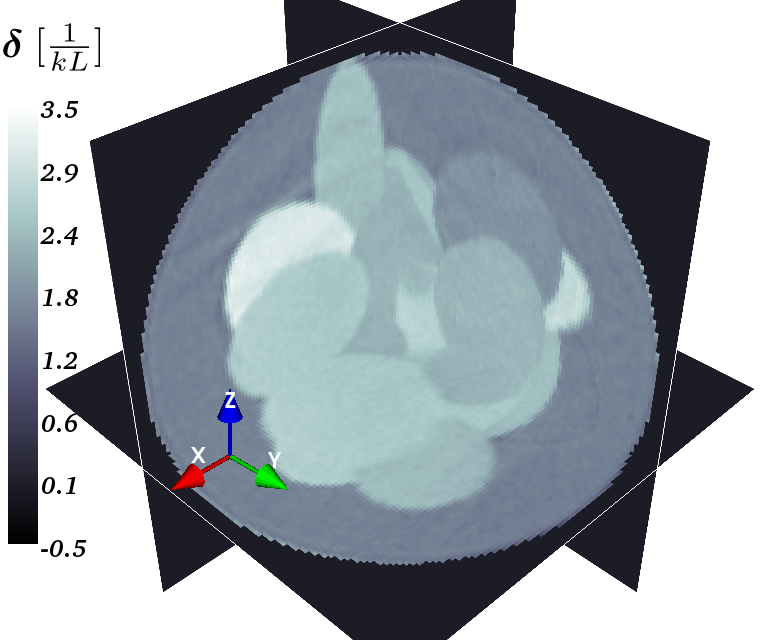} \hfill \includegraphics[height=.32\textwidth, clip=true, trim = -4cm 0cm -4cm 0cm]{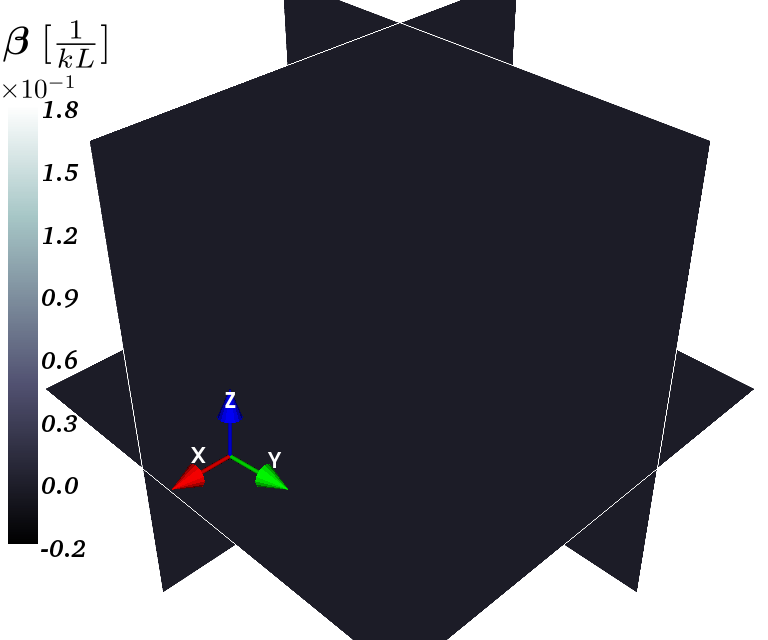} \label{fig:NumRes-FF3D-phaseRec} }
	  \\
	  \subfloat[Reconstruction $\bN_{k_{\Text{stop}}}^{\Text{gen}} = \bdelta_{k_{\Text{stop}}}^{\Text{gen}} - \I \bbeta_{k_{\Text{stop}}}^{\Text{gen}}$ as a general object]{\includegraphics[height=.32\textwidth, clip=true, trim = -4cm 0cm -4cm 0cm]{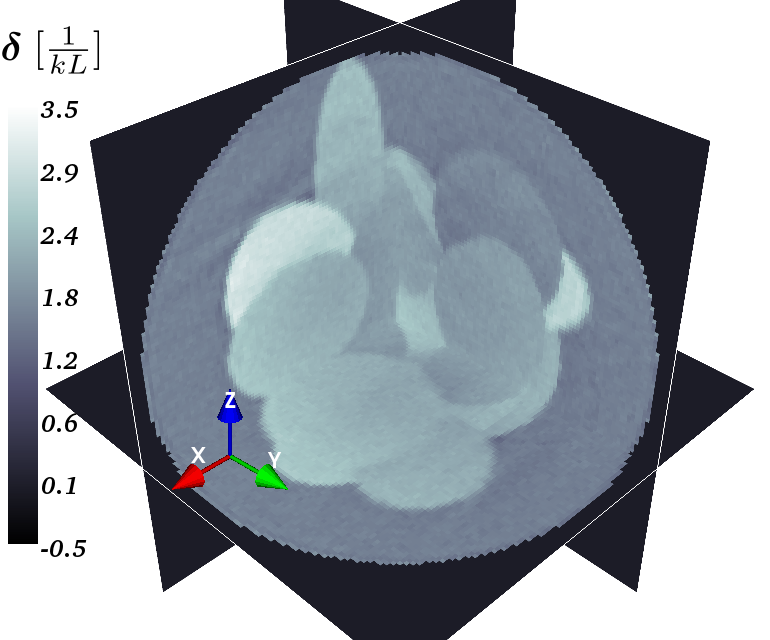} \hfill \includegraphics[height=.32\textwidth, clip=true, trim = -4cm 0cm -4cm 0cm]{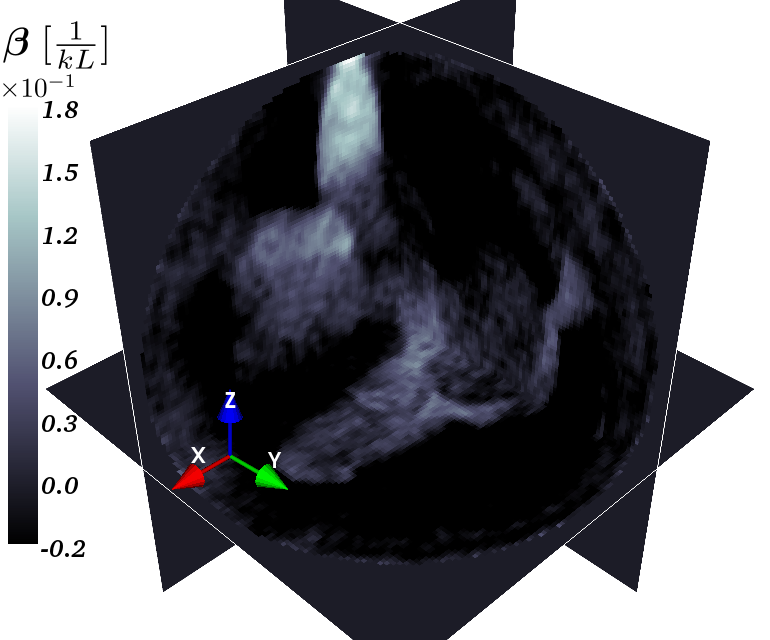} \label{fig:NumRes-FF3D-genRec} }
	  \caption{3D far-field tomography results from the simulated intensity data shown in \figref{fig:NumResFF-NonIdealData}. Reconstructions were carried out both using a (false) pure phase object constraint (b) and assuming a general object with independent refraction $\bdelta$ and absorption $\bbeta$ (c). Stop index $k_{\Text{stop}}$ for \algref{alg:PCT} chosen by discrepancy principle \eqref{eq:Discrepancy} with $\tau = 1$. For simulation parameters, see \tabref{tab:NumResFF3DSetup}. Relative $L^2$-error of $\bN_{k_{\Text{stop}}}^{\Text{phase}}$ and $\bN_{k_{\Text{stop}}}^{\Text{gen}}$ w.r.t. $\bN^\dagger$ is $5.9 \, \%$ and $11.1\,\%$, respectively. \label{fig:NumResFF-NonIdealRecon}}
	  \end{figure}
          \noindent detector  pixels for $256$ equispaced incident angles $\theta \in [0^{\circ}; 160^{\circ})$. Photon counts according to \cite{Chapman2006} are obtained by adjusting the intensity factor in \eqref{eq:DataGenFF}. The resulting data is visualized in \figref{fig:NumResFF-NonIdealData} where the beam stop area is shaded in black. These values are suppressed in the reconstruction by setting the corresponding weights in the Gramian $\cG_{\mY_{\Text{dis}}}$ (cf. \eqref{eq:GramYDisc}) to zero. By this modification, the formula \eqref{eq:NumResFF-alpha0} for the initial regularization parameter no longer applies. Instead we choose $\alpha_0 = 10^{16}$. Customized parameters are summarized in \tabref{tab:NumResFF3DSetup}. The remainder is chosen according as  described in \sref{SS:NumResFF-GenSetup}, see \tabref{tab:NumResFFSetup}.
          
          Two different Newton reconstructions are computed with \algref{alg:PCT}, both starting from the pure phase reference object $\bN_0$ as the initial guess: we reconstruct  $\bN^\dagger = \bdelta^\dagger - \I \bbeta^\dagger$ for once as a general, i.e.\ complex-valued object. Here, we only incorporate the a priori knowledge of the \emph{average} ratio $c_{\bbeta/\bdelta} = 0.1$ between absorption and refraction by adjusting the $L^2$-regularization such that deviations in $\bbeta$ by $c_{\bbeta/\bdelta}$ are punished equally strongly as deviations by 1 in $\bdelta$. This prevents initial overestimation of the absorption, which would obstruct convergence. For comparison, a second reconstruction is computed in which the (false) constraint of a pure phase object is imposed, neglecting the $10 \, \%$ absorption in the unknown phantom.
          
          Instead of prescribing a fixed number of Newton iterations, we apply the discrepancy principle (see \sref{SS:RegPar}) as an implementable stop rule: as the Poisson noise level is uniquely determined by the \emph{exact} intensities (compare \eqref{eq:DataGenFF}, it may be faithfully estimated from the \emph{observed} data $\bI^{\Textbf{err}}$ and is thus accessible from experimental observations. For the reconstruction assuming a pure phase object, the stopping criterion is reached after $k_{\Text{stop}} = 14$ Newton steps corresponding to a total of 268 CG-iterations. In the unconstrained run, convergence is slightly slower taking 15 Newton iterates ($\sim$ 378 CG-steps). The reconstruction results $\bN_{k_{\Text{stop}}}^{\Text{phase}}$ (pure phase constraint) and $\bN_{k_{\Text{stop}}}^{\Text{gen}}$ (general object) are visualized in \figref{fig:NumResFF-NonIdealRecon} as volumetric slice plots along the different coordinate planes.
          
          A first aspect to note is that the reconstructions in the refractive real part $\bdelta^\dagger$ of the object $\bN^\dagger$ come out almost artifact-free (compare left column of \figref{fig:NumResFF-NonIdealRecon}) - in spite of the incompleteness of the data in \figref{fig:NumResFF-NonIdealData} due to the beam stop and missing wedge. In \sref{S:RadonIll-posed} and \sref{S:PhaseRetrieval} it was shown that such missing information may always be recovered by analytic continuation for \emph{exact} and \emph{continuous} data. The quality of the reconstructions now seems to demonstrate that this also works in practice for noisy and discrete data. The only visible traces of this highly ill-posed implicit data completion in the results are weak stripe artifacts, for instance emanating from the edge of the ellipsoid with the maximum $\bdelta$-value in the left images in \figref{fig:NumRes-FF3D-phaseRec} and \ref{fig:NumRes-FF3D-genRec}. These artifacts occur along a characteristic direction with insufficient information coinciding with the axes of the missing projections.
          
          From the right hand slice plots in \figref{fig:NumResFF-NonIdealRecon}, it can be seen that the agreement in the recovered \emph{absorption} $\bbeta_{k_{\Text{stop}}}^{\Text{gen}}$ is poor. Only a few pronounced features of the exact object in \figref{fig:NumRes-FF3D-exObj} can be identified in the reconstruction visualized in \figref{fig:NumRes-FF3D-genRec}. Moreover, the resulting absorption values are even widely negative - with magnitudes partly as large as the maximum $\norm{\bbeta^\dagger} = \frac \pi {20}$ of the exact (non-negative) solution. Recall that positivity constraints may not be implemented in the considered Newton-type method, see \sref{SS:Constraints}. Hence, the attempted reconstruction of $\bbeta^\dagger$ along with the much stronger refractive component $\bdelta^\dagger$ seems to fail, with a final relative $L^2$-error of $\norm{ \bbeta_{k_{\Text{stop}}} - \bbeta^\dagger }_2 /  \norm{ \bbeta^\dagger }_2 = 109 \, \%$. This is perhaps not too surprising because the signal-to-noise-ratio in $\bbeta$ suffers from the much larger contributions of refractive effects to the diffraction patterns, parametrized by $\bdelta$. Yet, note that the considered $\beta$-$\delta$-ratio of $10 \, \%$ in the unknown phantom is already rather large from a physical perspective, compare \sref{SS:XrayRefrIdx}. It would thus be artificial to consider test cases with absorption and refraction of equal magnitude. 
          
          In turn, it should be emphasized that the incorrectly determined absorption values have a considerable negative backlash onto the reconstruction of the \emph{refraction} $\bdelta$: the final error in the case of the unconstrained solution shown in \figref{fig:NumRes-FF3D-genRec} is
          \begin{equation*}
           \norm{ \bdelta_{k_{\Text{stop}}}^{\Text{gen}} - \bdelta^\dagger }_2 /  \norm{ \bdelta^\dagger }_2 = 10.6 \, \%,
          \end{equation*}
          whereas the achieved agreement in the reconstruction assuming a pure phase object is as good as $\norm{ \bdelta_{k_{\Text{stop}}}^{\Text{phase}} - \bdelta^\dagger }_2 /  \norm{ \bdelta^\dagger }_2 = 5.1 \, \%$.
          Accordingly, it seems that the additional ill-posedness arising in the simultaneous recovery of both refraction $\bdelta$ and absorption $\bbeta$ outweighs the systematic errors made by the false neglect of absorption. Indeed, taking into account the incompleteness of the simulated measurements in \figref{fig:NumResFF-NonIdealData} along with these systematic deviations, the reconstruction of the refractive part $\bdelta$ alone appears to be fascinatingly robust in the considered setting. The same can be expected to hold if a general single-material constraint is assumed.
          
	\begin{res}[Newton-based 3D Far-Field Tomography from Realistic Intensity Data] \label{res:FF3DTomoNonideal}
		3D far-field tomography via \algref{alg:PCT} using spherical reference objects and single-material constraints is accurate and robust against both incomplete data due to missing wedges or beam stops  and residual absorption. On the contrary, independent reconstruction of refraction $\delta$ and absorption $\beta$ is too unstable to be competitive in this setting. The discrepancy principle \eqref{eq:Discrepancy} with $\tau = 1$ provides a reasonable stop rule for tomographic far-field data with Poisson noise.
	\end{res}

        \end{subsection}

    \end{section}

        \begin{section}{Near-Field Tomography from Simulated Data}  \label{S:NumResNFSim}
        
        From the far-field case considered in \sref{S:NumResFF}, we now proceed to the discussion of numerical test cases for \emph{near-field} phase contrast tomography. Despite the apparent similarity of the governing forward operators in \eqref{eq:ForwardOpNF} and \eqref{eq:ForwardOpFF}, there is a principal structural difference in the near-field case owing to the holographic superposition of the unscattered probe beam with the propagated contact image. It has been found in \sref{S:PhaseRetrieval} that this naturally arising reference signal completely eliminates non-uniqueness. As will be seen in the following study of numerical reconstructions, this specialty of near-field phase retrieval likewise gives rise to a significantly changed solution behavior of the considered Newton-type \algref{alg:PCT}. In order to allow for a quantitative comparison of this work's simultaneous approach to phase contrast tomography with competing methods, we uniquely consider the physically relevant case of tomographic setups in $m+1 = 3$ spatial dimensions.

         \begin{subsection}{Simulation Setup}  \label{SS:NumResNFSim-Setup}
     
	    We consider 3D objects $\bN \in  \mX_{\Text{dis}}  := \mC^{64 \times 64 \times 64}$ in a cubic domain of relatively low spatial resolution in order to save computation time for the numerical test cases studied herein. As opposed to the far-field setup in \sref{SS:NumResFF-GenSetup}, no support constraints are prescribed in addition to the rough restriction by the  computational domain. Single-material- and pure phase object constraints are imposed as described in \ref{SS:Constraints}. In order to avoid sampling issues, we typically choose a resolution $\bI^{\Textbf{err}} \in \mY_{\Text{dis}} := \mR^{128 \times 128 \times 128}$ in the intensity data via zero-padding by a factor of two in the lateral dimensions (cf. \sref{SS:ZeroPad}), corresponding to measurements by $128 \times 128$ pixel detectors under $128$ tomographic incident angles $\theta \in [0^{\circ}; 180^{\circ})$.
	    
	    As argued in \sref{SS:RecMethodBase}, the statistical errors in near-field intensities may typically be approximated as Gaussian. We account for this in the numerical simulations by computing synthetic data to an exact object $\bN^\dagger \in \mX_{\Text{dis}}$ via
	   \begin{equation}
	    \bI^{\Textbf{err}} =   \underbrace{F_{\Text{dis}}(\bN^\dagger)}_{= \bI^\dagger } + \Textbf{err} \MTEXT{with independent errors} \Textbf{err}_j \sim \cN(0, \sigma_\varepsilon),  \label{eq:NumResNF-DataGen}
	    \end{equation}
	    i.e.\ we supplement the exact data with additive Gaussian white noise. The standard deviation $\sigma_\varepsilon$ of the normal distribution $\cN(0, \sigma_\varepsilon)$ is chosen such that a prescribed noise level $\varepsilon  = \norm{\bI^{\Textbf{err}} - \bI^\dagger } / \norm{\bI^\dagger}$ is obtained. According to the error statistics, $L^2$-data fidelity functionals are implemented by choosing the Gramian $\cG_{\mY_{\Text{dis}}}$ as the identity. For simplicity, we again restrict to plane wave illumination, setting $P =1$ in the discretization $F_{\Text{dis}} \approx F_d$ of the forward operator in \eqref{eq:ForwardOpNF}. Thereby, all intensities are implicitly measured in units of the constant background intensity $I_0$ of the unscattered incident probe beam.
	    
	    As generic test objects $\bN^\dagger = \bdelta^\dagger - \I \bbeta^\dagger$ we choose random ensembles of nested ellipsoids $j$ as in \sref{SS:NumResFF-3DRecon}. We prescribe a mean ratio $c_{\bbeta/\bdelta}$ for the expectation values of the refractive $\bdelta^\dagger_j$ and the absorptive parts $\bbeta^\dagger_j$ ($c_{\bbeta/\bdelta} = 0$ for pure phase objects) and assign values to the ellipsoids drawn from independent normal distributions
                 \begin{equation}    k L \bdelta^\dagger_j \sim \cN(\mu, \sigma\mu) \MTEXT{and} k L  \bbeta^\dagger_j \sim c_{\bbeta/\bdelta} \cN(\mu, \sigma\mu ). \label{eq:NumResNF-ObjGen} \end{equation}
	    Throughout this section, we choose $\sigma  = 0.3$. By scaling the obtained random objects, we prescribe their strength in terms of the norm $\norm{\bN^\dagger} \propto \mu$  defined in \eqref{eq:ObjScaleNorm}.
	    
	    Other than in the far-field case, the physical setup parameters $k$, $L$ and $d$, i.e.\ wavenumber, specimen thickness and propagation distance, are not just relevant for the object magnitude.
	    The second dimensionless problem parameter besides $\norm{\bN^\dagger}$ is given by the Fresnel number $\NF = \frac{k b ^2}{2 \pi d}$, see \sref{SS:FresnelFarfield}, governing the discrete near-field propagator according to \eqref{eq:DiscPropNF}. In the numerical results discussed here, we take the object lengthscale $b$ as the size $\Delta x = L / M_x$ of a single pixel. Prescribing $\NF$, the significance of $k$, $L$ and $d$ reduces to scaling of the coordinate axes.
	    
	    In this section, we consider general Sobolev norm penalty functionals (see \sref{SS:Constraints} - Regularity Constraints), parametrized by a Gramian $\cG_{\mX_{\Text{dis}}}$  of the form \eqref{eq:GramXDisc}. Like in the far-field case the number of CG-iterations for the initial Newton step is taken as a heuristic measure for the initial regularization parameter $\alpha_0 > 0$. A good choice in the considered setup turns out to be 
	    \begin{equation}
	     \alpha_0  = \frac{ \norm{\bI^{\Textbf{err}} -1 }_2^2 }{ \norm{ \bN^\dagger }_{\mX_{\Text{dis}}}^2 }. \label{eq:NumResNF-alpha0}
	    \end{equation}
	    Note that we have to subtract the constant probe beam intensity $1$ from the data in order to obtain a reasonable scaling. We choose $\alpha_{k+1} = \frac 2 3 \alpha_k$, i.e.\ $r_\alpha = \frac 2 3$ just like in the far-field simulations.
	    
	    To explore the potential of our method, we mostly use an optimal ``best stop'' rule choosing the stop index $k_{\Text{stop}}$ such that the final Newton iterate $\bN_{k_{\Text{stop}}}$ minimizes~the $L^2$-reconstruction error
	    \begin{equation}
	     \rho_k := \frac{ \norm{ \bN_k - \bN^\dagger }_2 }{ \norm{ \bN^\dagger }_2 } \label{eq:NumResNF-RecErr}.
 	    \end{equation}
	    \begin{table}[htb!]
	      \centering
	      \begin{tabular}{cccccccccc} 
		\toprule
		  $\bI^{\Textbf{err}}$  &  $\mX_{\Text{dis}}$ &  $\mY_{\Text{dis}}$ & $\cG_{\mX_{\Text{dis}}}$  & $\cG_{\mY_{\Text{dis}}}$  & $\alpha_0$ & $r_\alpha$  & $\bN_0$ & $P$ & Constraints \\
		\midrule
		  \scriptsize{\eqref{eq:NumResNF-DataGen}} & $\mC^{64^3}$ & $\mR^{128^3}$  &  \scriptsize{\eqref{eq:GramXDisc}} &  $\id_{\mY_{\Text{dis}}}$  & $ \frac{ \norm{\bI^{\Textbf{err}}  -1 }_{\mY_{\Text{dis}}}^2 }{   \norm{\bN^\dagger}^2_{\mX_{\Text{dis}}} } $  & $\frac 2 3$  & $ 0 $ & $ 1$ & $\substack{  \text{pure phase obj.} \\ \text{(optional)} } $ \\ 
		\bottomrule
	      \end{tabular}
	      \caption{Setup parameters for the numerical test cases of near-field tomography via \algref{alg:PCT}. The test object $\bN^\dagger = \bdelta^\dagger - \I \bbeta^\dagger \in \mX_{\Text{dis}}$ is taken as a random ensemble of nested ellipsoids with normally distributed values for the refraction $\bdelta^\dagger$ and (optional) absorption $\bbeta^\dagger$ as in \sref{SS:NumResFF-3DRecon}.}
	      \label{tab:NumResNFSetup}
	    \end{table}
	    
	    The simulation parameters assigned to \algref{alg:PCT} for the numerical test cases of near-field tomography, governed by the operator $F_{\Text{dis}} \approx F_d$ in \eqref{eq:ForwardOpNF}, are summarized in \tabref{tab:NumResNFSetup}. From \eqref{eq:FrechetForwOpNF}, it can be seen that the \Frechet derivative $F_{\Text{dis}}'[\bN]$ is non-vanishing for $\bN = 0$ - in contrast to the far-field case, see \sref{SS:NumResFF-GenSetup}. Hence, we may always choose $\bN_0 = 0$ as a canonical initial guess without risking an immediate stagnation of the regularized Newton method.
	    

	 \end{subsection}

         \begin{subsection}{Parametric Study for Pure Phase Objects}  \label{SS:NumResNFSim-ParamStudy}
         
	    As a first step, we investigate the influence of the different problem parameters on the reconstruction by \algref{alg:PCT}. In addition to the Fresnel number $\NF$ and the object magnitude  $\norm{\bN^\dagger}$, smoothing effects by Sobolev space regularizations of different order as well as the impact of a ``missing wedge'' of incident angles are subject to separate parametric studies. For simplicity, these studies are restricted to the important special case of pure phase objects, fixing  $c_{\bbeta/\bdelta}  = 0$ in \eqref{eq:NumResNF-ObjGen}. Parameters which are not specified are chosen according to \tabref{tab:NumResNFSetup} and \sref{SS:NumResNFSim-Setup}.
	   
	    \subsubsection{Fresnel Number}
	    
	    The first parameter to be studied is the Fresnel number. To this end, we compute reconstructions by \algref{alg:PCT} for $\NF \in \{ 0.1, 0.05, 0.02, 0.01, 0.005, 0.002, 0.001 \} $ 
	    for one and the same pure phase object $\bN^\dagger$ of magnitude $\norm{\bN^\dagger} = \pi$ at a fixed noise level $\varepsilon = 3 \%$. For simplicity, we use an $L^2$-regularization term in this study. In order to allow for the physical relevant case of fringes propagating out of the lateral domain, we zero-pad by a factor of 4 but truncate the propagated data to the simulated $128 \times 128$ detector pixels. For a fair comparison, we choose the ideal ``best stop'' rule from \sref{SS:NumResNFSim-Setup}, stopping the Newton iterations $\bN_k$ at a minimum reconstruction error $\rho_k$  given by \eqref{eq:NumResNF-RecErr}. The exact object along with exemplary reconstruction results for $\NF = 0.1$, $0.01$ and $0.001$ are depicted in \figref{fig:NumResNF-FresnelTest-Obj}.
	       \begin{figure}[hbt!]
	  \centering
	  \subfloat[Exact object $\bN^\dagger$]{\includegraphics[height=.34\textwidth, clip=true, trim = 0cm 0cm -4cm 0cm]{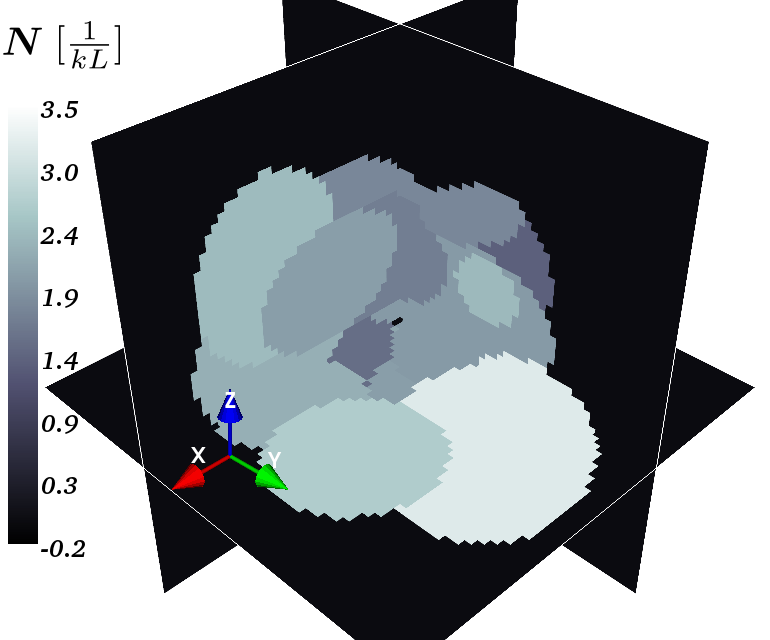} \label{fig:NumResNF-FresnelTest-Obj-ex}  }
	  \hfill
	  \subfloat[Reconstruction $\bN_{k_{\Text{stop}}}$ for $\NF = 0.1$]{\includegraphics[height=.34\textwidth, clip=true, trim = 0cm 0cm -4cm 0cm]{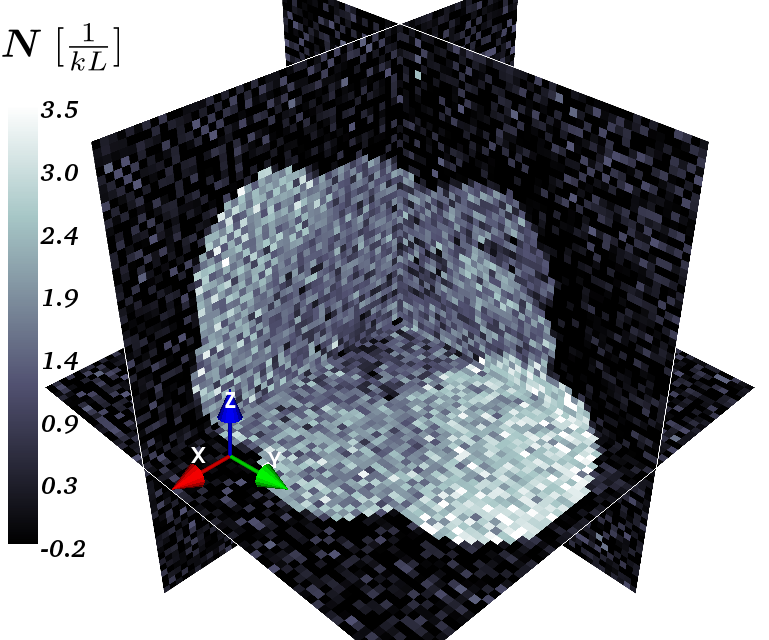} \label{fig:NumResNF-FresnelTest-Obj-NF1}  }
	  \hfill
	  \subfloat[Reconstruction $\bN_{k_{\Text{stop}}}$ for $\NF = 10^{-2}$]{\includegraphics[height=.34\textwidth, clip=true, trim = 0cm 0cm -4cm 0cm]{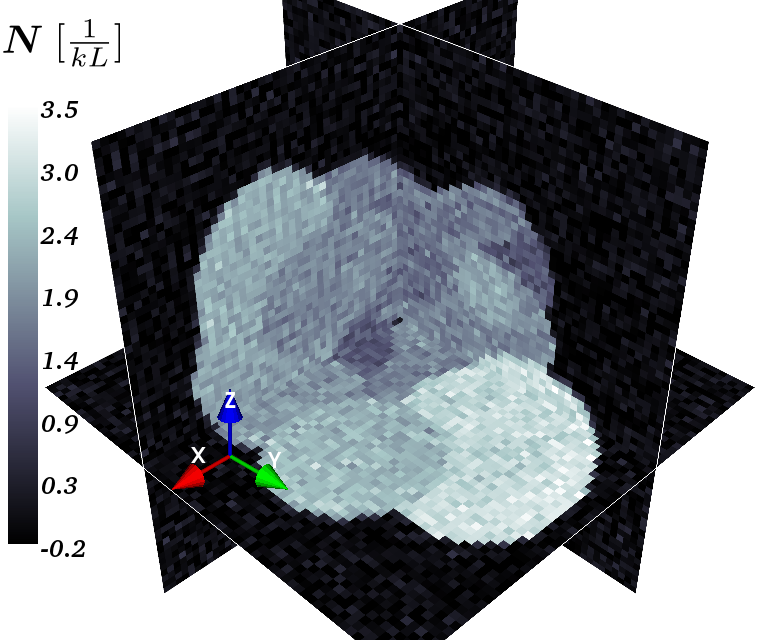} \label{fig:NumResNF-FresnelTest-Obj-NF4}  }
	  \hfill
	  \subfloat[Reconstruction $\bN_{k_{\Text{stop}}}$ for $\NF = 0.001$]{\includegraphics[height=.34\textwidth, clip=true, trim = 0cm 0cm -4cm 0cm]{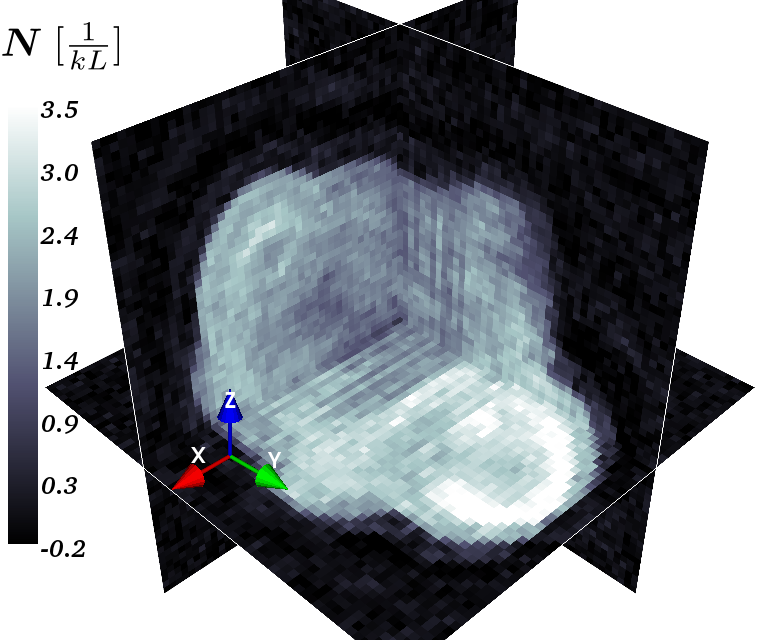} \label{fig:NumResNF-FresnelTest-Obj-NF7}  }
	  \caption{Near-field tomography results by \algref{alg:PCT} for different Fresnel numbers $\NF$. The intensity data $\bI^{\Textbf{err}}$ visualized in \figref{fig:NumResNF-FresnelTest-Data} is taken as the exact data $F_{\Text{dis}}(\bN^\dagger)$ plus  $\varepsilon = 3\,\%$ Gaussian noise. Stop rule: ``best stop'' (minimize $L^2$-error $\rho_{k }$). The noisy reconstruction for $\NF = 0.1$ (b) is a manifestation of low \emph{phase contrast}, whereas the halo- and stripe artifacts in (d) result from fringes leaving  the computational domain. Detailed statistics given in \tabref{tab:NumResNF-FresnelTest}. \label{fig:NumResNF-FresnelTest-Obj}}
	  \end{figure}
     \begin{figure}[hbt!]
	  \centering
	  \subfloat{\includegraphics[height=.3\textwidth]{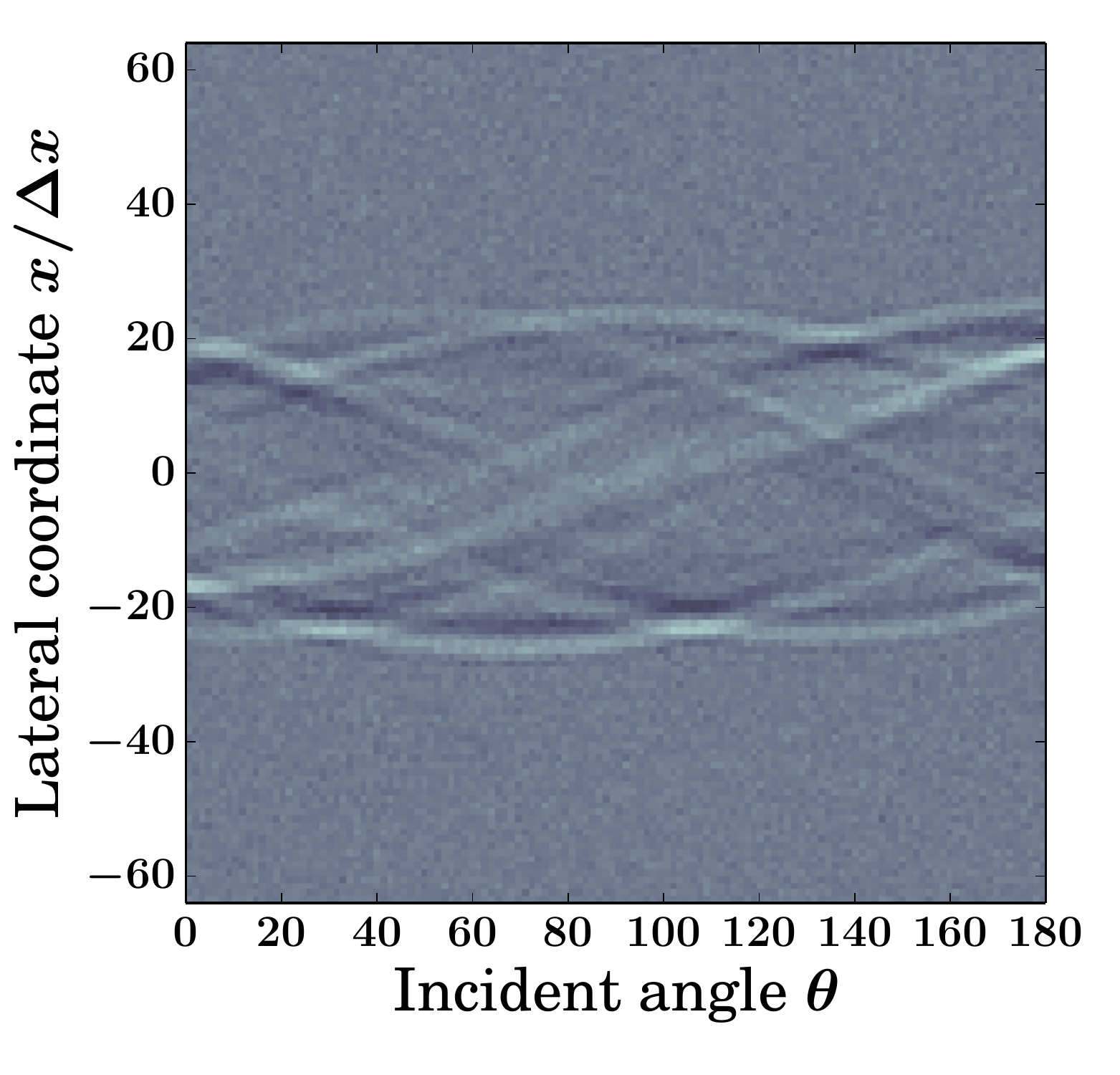} \label{fig:NumResNF-FresnelTest-Data-NF1}  }
	  \hfill
	  \subfloat{\includegraphics[height=.3\textwidth]{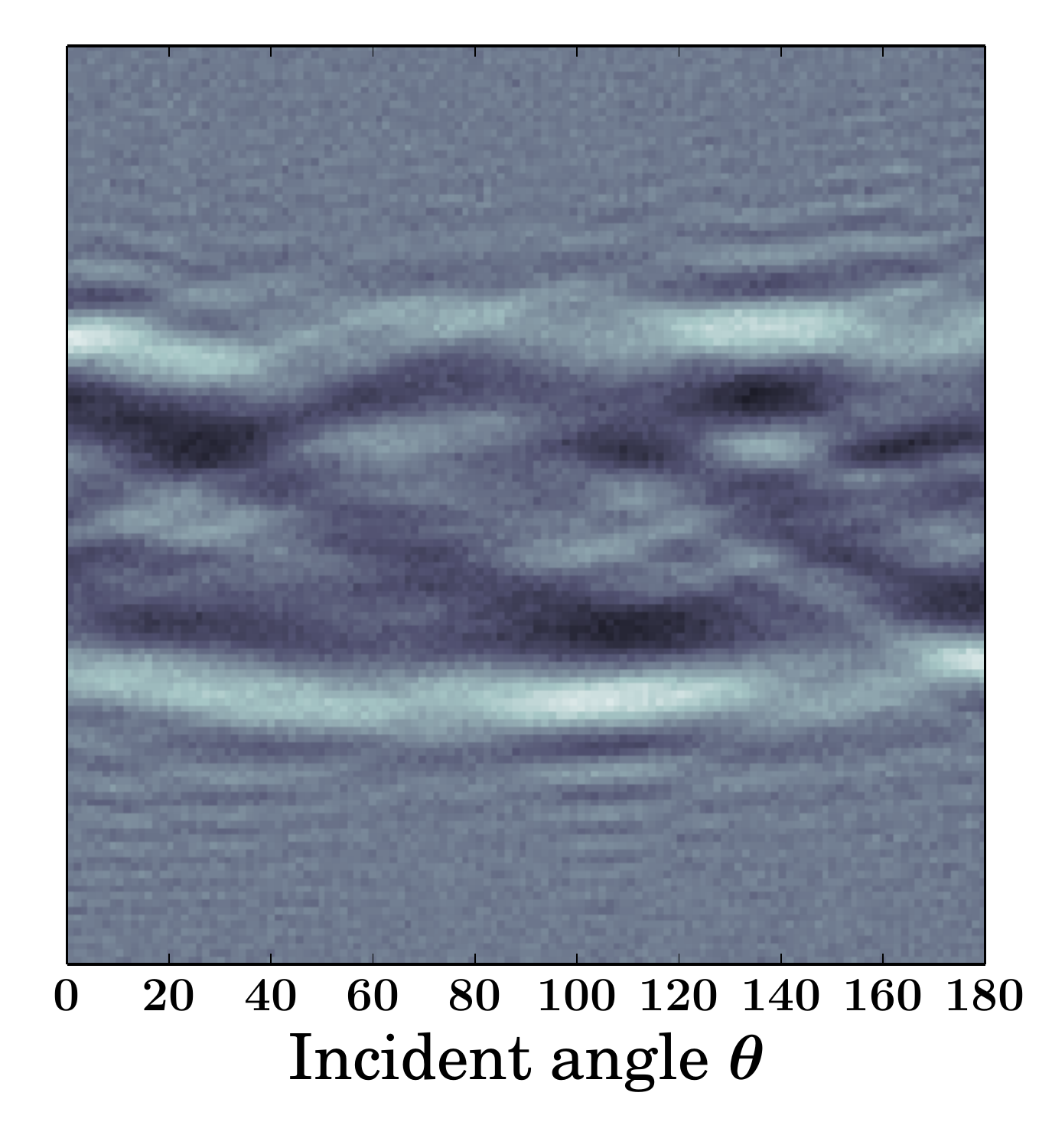} \label{fig:NumResNF-FresnelTest-Data-NF4}  }
	  \hfill
	   \subfloat{\includegraphics[height=.3\textwidth]{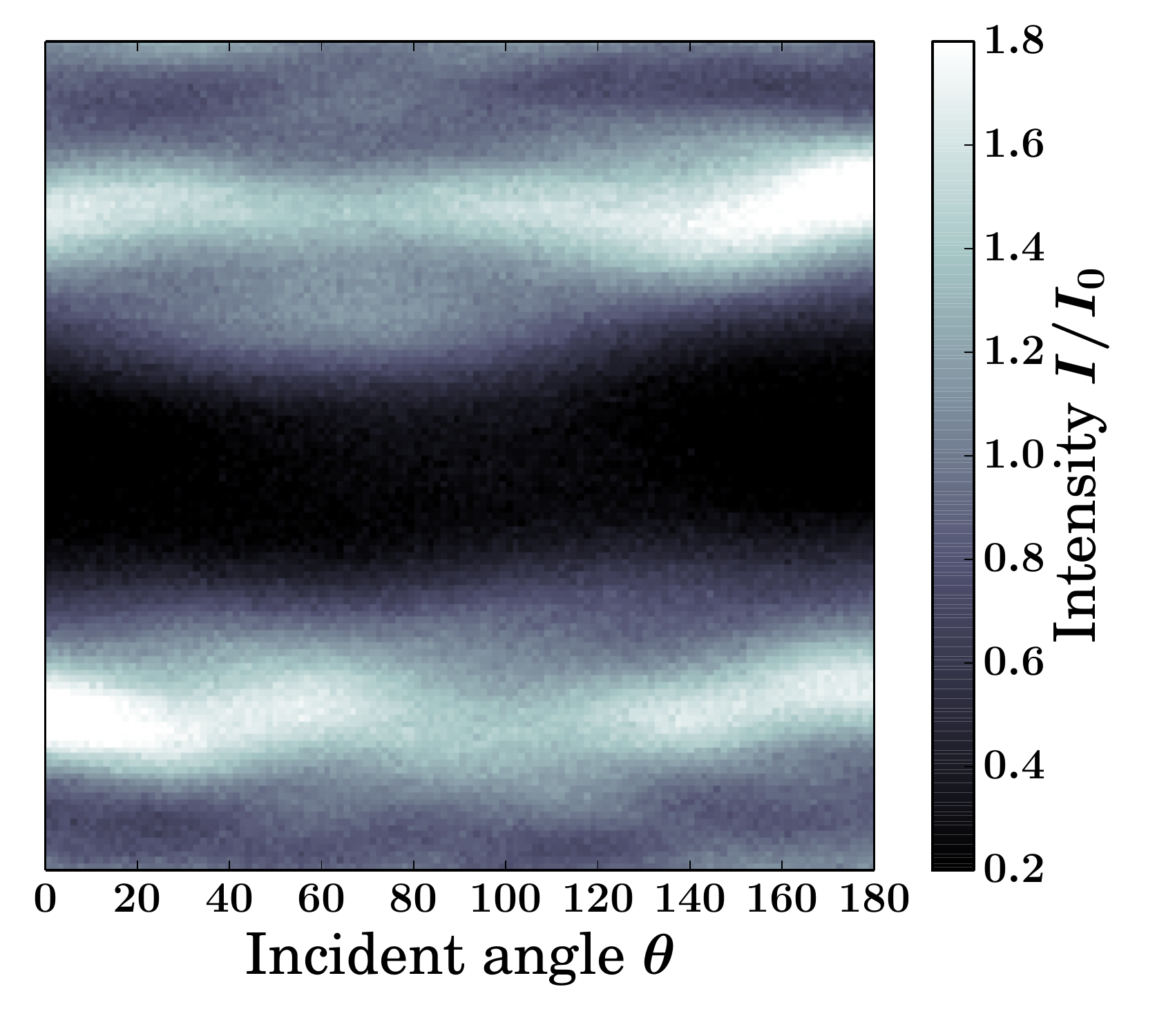} \label{fig:NumResNF-FresnelTest-Data-NF7}  }
	  \caption{Simulated data $\bI^{\Textbf{err}}$ (shown: central propagated holo-sinogram, i.e.\ 2D slice at $y = 0$) corresponding to the object in \figref{fig:NumResNF-FresnelTest-Obj-ex} at different Fresnel numbers $\NF$. Computed using \eqref{eq:NumResNF-DataGen} at noise level $\varepsilon = 3\,\%$. Corresponding exact object and reconstructions shown in \figref{fig:NumResNF-FresnelTest-Obj}. \label{fig:NumResNF-FresnelTest-Data}}
	  \end{figure}
	  
	  The computed solutions show considerable differences in their agreement with the exact object in \figref{fig:NumResNF-FresnelTest-Obj-ex}: while the reconstruction for  $\NF = 0.01$ (\ref{fig:NumResNF-FresnelTest-Obj-NF4}) is surprisingly accurate, showing apparently no artifacts and only slight effects of the moderately high noise level of $3 \%$, the latter has a much stronger impact on the result for $\NF =0.1$ plotted in \figref{fig:NumResNF-FresnelTest-Obj-NF1}. On the other hand, \figref{fig:NumResNF-FresnelTest-Obj-NF7}, reconstructed at a Fresnel number of $\NF = 0.001$, contains spurious halos of the ellipsoids and stripe artifacts whereas noise is less pronounced in this solution.
	  
	  Both effects may be understood by considering the corresponding data visualized in \figref{fig:NumResNF-FresnelTest-Data}. For the small Fresnel number $\NF = 0.1$ only the edges are imprinted in the intensity data, corresponding to low \emph{phase contrast} and thus a bad signal-to-noise-ratio. In the picture of the contrast transfer function in \figref{fig:CTF}, the problem is that the relevant Fourier frequencies of the pure phase object are located too close to the origin on the ascending branch of the oscillating contrast curve. As seen in \figref{fig:NumResNF-FresnelTest-Data}, the largest Fresnel number $\NF = 0.001$ apparently attains the maximum contrast and hence the weakest impact of noise in the reconstruction. However, in this case - physical, not merely numerical - finite domain effects permit propagation of fringes beyond the lateral detector domain, constituting a leak of  object information in the computational or experimental setup. The recovery of the missing information is - although possible in principal by the uniqueness result in \ref{thm:NFUnique} - severely ill-posed and thus gives rise to the characteristic artifacts.
	    \begin{table}[htb!]
	      \centering
	      \begin{tabular}{lccccccc} 
		\toprule
		 Fresnel number $\NF$ & 0.1 & 0.05 & 0.02 & 0.01 & 0.005 & 0.002 & 0.001 \\
		\midrule
		 Reconstruction error $\rho_{k_{\Text{stop}}}$ & 0.58 & 0.41 & 0.27 & 0.22 & 0.39 & 0.26 & 0.26 \\
		 Total CG-iterations & 725 & 541 & 406 & 260 & 61 & 62 & 206 \\
		\bottomrule
	      \end{tabular}
	      \caption{Reconstruction statistics for numerical near-field test cases at different Fresnel numbers $\NF$. Exemplary reconstructed objects $\bN_{k_{\Text{stop}}}$ are shown in \figref{fig:NumResNF-FresnelTest-Obj}. The indicated $L^2$-errors $\rho_{k} = \norm{ \bN_k - \bN^\dagger }_2 / \norm{ \bN^\dagger }_2$ are minimized by the chosen stopping index $k = k_{\Text{stop}}$ (``best stop''). Large numbers of CG-steps for low $\NF$ indicate slow convergence of the Newton \algref{alg:PCT}.}
	      \label{tab:NumResNF-FresnelTest}
	    \end{table}
	  
	  The general tendencies for different $\NF$ identified by visual inspection are confirmed by the quantitative reconstruction statistics summarized in \tabref{tab:NumResNF-FresnelTest}: owing to increasing phase contrast, the reconstruction error initially decreases with $\NF$, before artifacts due to the finite domain begin to corrupt the solution for $\NF > 10^{-2}$ (apparently most significantly for $\NF = 0.005$ in the given example). Moreover, note that the Newton method converges terribly slowly for low Fresnel numbers, giving rise to an excessive number of CG-iterations despite the poor reconstruction quality according to \tabref{tab:NumResNF-FresnelTest}. In this regime, non-iterative methods based on the transport-of-intensity equations (see \ref{S:CombinedApproach}) thus seem preferable. On the other hand, large Fresnel numbers require large field of views in image space, i.e.\ excessive oversampling in the data, which is computationally expensive. 
	  Indeed, our Newton-type method seems to work best if the Fresnel number $\tilde \NF$ based on the \emph{characteristic} lengthscales of the object to be reconstructed is in the order of 1: the typical distance between two edges in \figref{fig:NumResNF-FresnelTest-Obj-ex} is roughly $10$ pixels, yielding $\tilde \NF \approx 100 \NF$, so that the optimum is attained for $\NF \approx 0.01$ - exactly as observed in the numerical simulation.
	 \end{subsection}

	    \subsubsection{Weak and Strong Objects}
	      
	  Next, we study the impact of weak or strong objects onto near-field tomography by our regularized Newton method, as parametrized by the object norm $\norm{\bN^\dagger}$ controlling the (non-)linearity of the object transmission function in \eqref{eq:OTF}. To this end, we scale a single phantom to different magnitudes $\norm{\bN^\dagger} \in \{ \frac \pi {16}, \frac \pi {8}, \frac \pi {4}, \frac \pi {2}, \pi , 2 \pi, 4 \pi, 8 \pi\}$  and compute reconstructions via \algref{alg:PCT} using an $L^2$-regularization term for a fixed data noise level $\varepsilon = 1 \, \%$ and Fresnel number $\NF = 0.01$. As before, we apply the ``best stop'' rule terminating the Newton method at the point where the $L^2$-object error $\rho_k = \norm{ \bN_{k } - \bN^\dagger }_2 /  \norm{ \bN^\dagger }_2$ begins to increase again. However, we prescribe a \emph{minimum} number of six Newton iterations in order to rule out cases where the error initially increases due to  nonlinearity. The numerical results for the different object strengths are summarized in \tabref{tab:NumResNF-NLTest}. Exemplary reconstructed objects for $\norm{\bN^\dagger} \in \{ \frac \pi {8},  2 \pi, 8 \pi \} $ are visualized in \figref{fig:NumResNF-NLTest-Obj}.
	    \begin{table}[htb!]
	      \centering
	      \begin{tabular}{lcccccccc} 
		\toprule
		 Object magnitude $\norm{\bN^\dagger}$ & $ \frac \pi {16} $ & $ \frac \pi {8} $ & $ \frac \pi {4} $ & $ \frac \pi {2} $ & $ \pi $  & $2 \pi$ & $4\pi$ & $8\pi$ \\
		\midrule
		 Reconstruction error $\rho_{k_{\Text{stop}}}$ & 0.65 & 0.47 & 0.31 & 0.20 & 0.13 & 0.10 & 0.31 & 0.99 \\
		 Total CG-iterations & 76 & 156 & 207 & 309 & 365 & 452 & 978 & 115  \\
		\bottomrule
	      \end{tabular}
	      \caption{Reconstruction statistics for numerical near-field test cases for objects $\bN^\dagger$ of different strengths $\norm{\bN^\dagger}$. Values $\norm{\bN^\dagger} \ll 1$ correspond to the weak object limit whereas \emph{phase-wrapping} emerges for $\norm{\bN^\dagger} > 2 \pi$. Computed at Fresnel number $\NF = 0.01$ and data noise level $\varepsilon = 1\,\%$ using $L^2$-regularization. Stop rule: ``best stop''. Further parameters according to \tabref{tab:NumResNFSetup}.
	       Exact phantom and exemplary reconstructed objects $\bN_{k_{\Text{stop}}}$ visualized in \figref{fig:NumResNF-NLTest-Obj}. }
	      \label{tab:NumResNF-NLTest}
	    \end{table}
	       \begin{figure}[hbt!]
	  \centering
	  \subfloat[Exact object $\bN^\dagger$]{\includegraphics[height=.34\textwidth, clip=true, trim = 0cm 0cm -4cm 0cm]{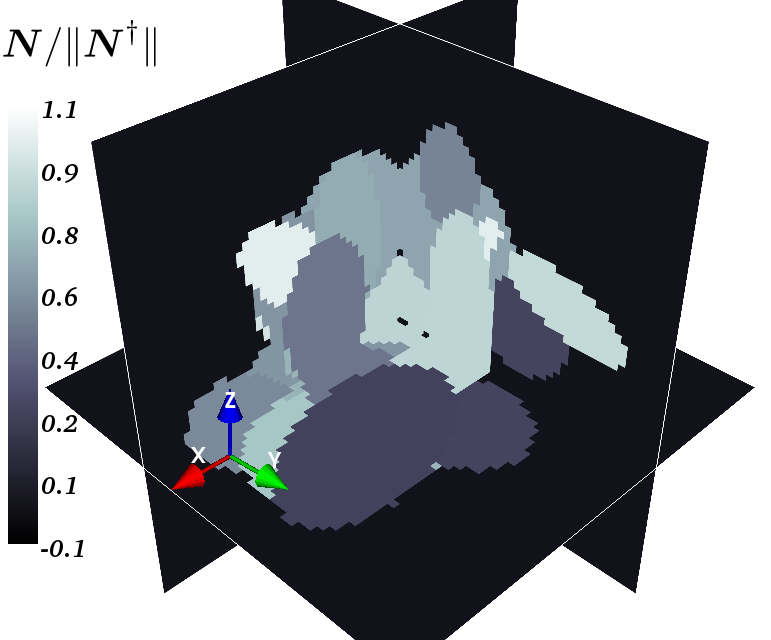} \label{fig:NumResNF-NLTest-Obj-ex}  }
	  \hfill
	  \subfloat[Reconstruction $\bN_{k_{\Text{stop}}}$ for $\norm{\bN^\dagger} = \frac \pi 8$]{\includegraphics[height=.34\textwidth, clip=true, trim = 0cm 0cm -4cm 0cm]{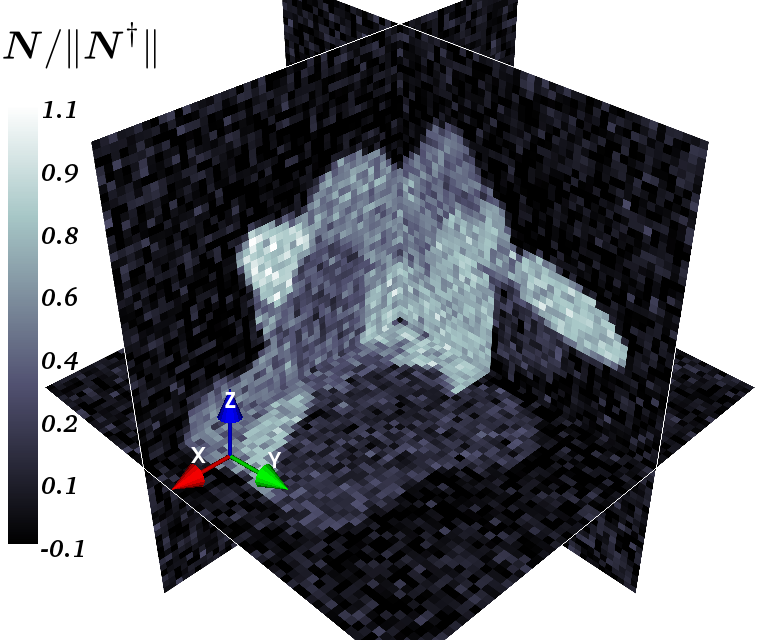} \label{fig:NumResNF-NLTest-Obj-NL2}  }
	  \hfill
	  \subfloat[Reconstruction $\bN_{k_{\Text{stop}}}$ for $\norm{\bN^\dagger} = 2 \pi $]{\includegraphics[height=.34\textwidth, clip=true, trim = 0cm 0cm -4cm 0cm]{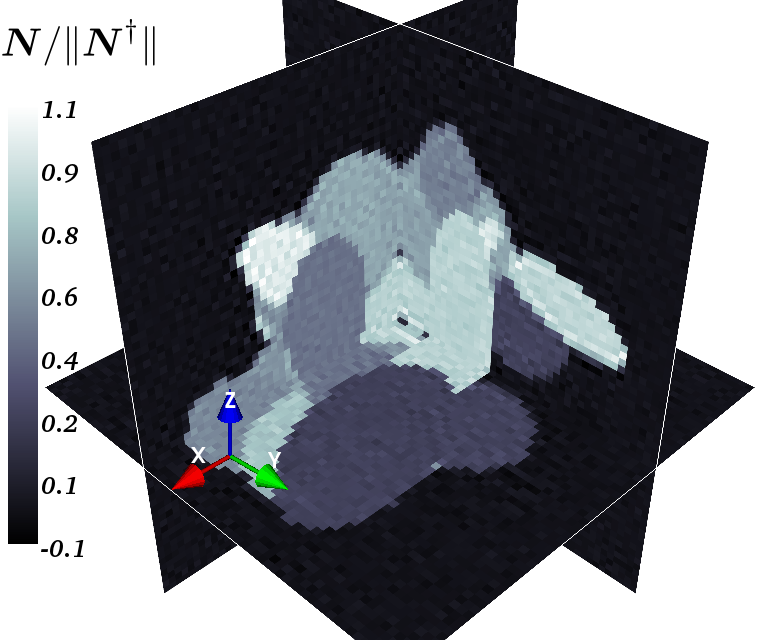} \label{fig:NumResNF-NLTest-Obj-NL6}  }
	  \hfill
	  \subfloat[Reconstruction $\bN_{k_{\Text{stop}}}$ for $\norm{\bN^\dagger} = 8 \pi $]{\includegraphics[height=.34\textwidth, clip=true, trim = 0cm 0cm -4cm 0cm]{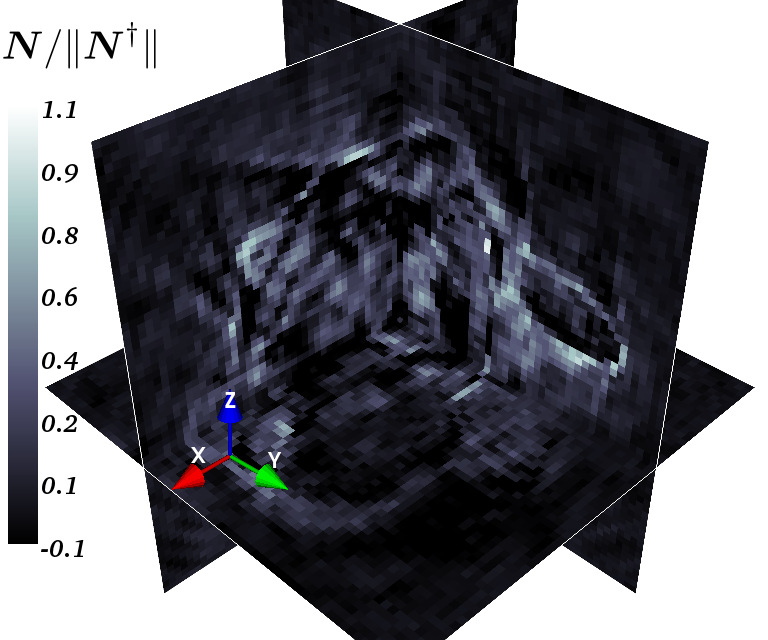} \label{fig:NumResNF-NLTest-Obj-NL8}  }
	  \caption{Exemplary near-field tomography results by \algref{alg:PCT} for different object strengths $\norm{\bN^\dagger}$. Details summarized in \tabref{tab:NumResNF-NLTest}. 
	  The noisy reconstruction in (b) is due to low contrast in the data by weakness of the object, whereas the artifacts in (d) are a manifestation of \emph{phase-wrapping}, see \sref{SS:PhaseWrap}. \label{fig:NumResNF-NLTest-Obj}}
	  \end{figure}
	  
	  It may seem surprising that the reconstruction \emph{improves} for stronger objects, i.e.\ stronger nonlinearity, according to \tabref{tab:NumResNF-NLTest} up to a magnitude of $\norm{\bN^\dagger} = 2 \pi$. However, it can be seen from the weak object limit of the near-field forward operator in \eqref{eq:WeakFwOpNF} that the phase contrast in the intensity data $\bI^{\Textbf{err}}$ is directly proportional to $\norm{\bN^\dagger}$ for $\norm{\bN^\dagger} \ll 1$. Accordingly, the larger reconstruction errors for smaller $\norm{\bN^\dagger} $ may be attributed to a poor signal-to-noise-ratio as the data noise level is identical in all reconstructions. This interpretation is supported by comparison of the examples in \figref{fig:NumResNF-NLTest-Obj-NL2} and \ref{fig:NumResNF-NLTest-Obj-NL6} with the exact object in \figref{fig:NumResNF-NLTest-Obj-ex}, showing no artifacts but different noise levels in the  reconstructions.
	  
	  As argued in \sref{SS:PhaseWrap}, \emph{phase-wrapping} may occur for object strengths $\norm{\bN^\dagger} > 2 \pi$ caused by the periodicity of the object transmission function $\exp( - \I k \CR( \bN ) )$ in the refractive real part of $\bN$. This gives rise to the larger object error of $31\,\%$ for  $\norm{\bN^\dagger} = 4 \pi$, where the phenomenon comes into play, and the total failure of the reconstruction for $\norm{\bN^\dagger} = 8 \pi$ shown in \figref{fig:NumResNF-NLTest-Obj-NL8}. It would certainly be astonishing if the Newton-type \algref{alg:PCT}, based on iterative \emph{linearizations}, could cope with this  severely \emph{nonlinear} effect. On the other hand, it can be regarded as a clear proof of concept that accurate solutions are achieved up to moderately strong objects with $\norm{\bN^\dagger} \sim 2 \pi$ - for which direct \emph{linear} methods based on the contrast transfer function (CTF, see \eqref{eq:CTF} and \sref{S:CombinedApproach}) are likely to fail.
	  
	  However, note that starting \algref{alg:PCT} from the initial guess $\bN_0 = 0$ implicitly computes a regularized CTF-solution in the first iterate, simply because the local linearization coincides with the underlying weak object limit, see \sref{SS:ContrastFormationNF}. In this sense, our Newton-type approach to near-field phase contrast tomography \emph{generalizes} CTF-based methods. The relatedness to these linear techniques suggests that the initial guess is of much lesser significance for the Newton reconstructions than in the far-field case (compare \resref{res:FFTomoByNewton}) - at least if the first CTF-like iterate provides a good approximation, i.e.\ for at most moderately strong objects. This interpretation is supported by the observed robustness of the numerical solutions in the near-field test cases considered so far. It should be emphasized that the difference in the algorithmic behavior arises from the unscattered probe contributions in the governing near-field forward operator \eqref{eq:ForwardOpNF}, providing a natural holographic reference for phase retrieval - as exploited in the uniqueness analysis of \sref{SS:PhaseRetrNF}.
	  
	  A final aspect to note is that, due to the quasi-linearity of near-field phase contrast tomography for weak objects by \eqref{eq:WeakFwOpNF}, Newton-type iterations are practically pointless in this limit, as the \Frechet derivative is almost independent of the current iterate. Hence, a good reconstruction may be achieved by a \emph{single} Newton step provided a suitable choice of the regularization parameter $\alpha_0$. In general, this suggests that the decrease of $\alpha_k$ from one Newton iteration to the next, set by the parameter $r_\alpha$ in \algref{alg:PCT}, should be larger the weaker the object - an adjustment which has been omitted here.

	  \subsubsection{Choice of the Regularization Term}
	  
	  In \sref{SS:Constraints}, it has been motivated that Sobolev $H^s$-norm regularization terms may be applied to suppress noise in the reconstruction, exploiting a priori knowledge on the regularity of the unknown object. This approach is examined in the following. To this end, we reconstruct a given object with \algref{alg:PCT} using Gramians $\cG_{\mX_{\Text{dis}}}$ of the form \eqref{eq:GramXDisc} with different parameters $s \in \{ 0, 0.25, 0.5, 0.75, 1\}$, corresponding to different degrees of smoothing ($s = 0$ corresponds to the $L^2$-penalty considered so far). We choose $\NF = 0.01$, an exact object $\bN^\dagger$ with $\norm{\bN^\dagger} = \pi$ and a moderately high data noise level of $\varepsilon = 3 \, \%$, stopping the Newton iterations according to the ideal ``best stop'' criterion, see \sref{SS:NumResNFSim-Setup} and previous test cases. The resulting reconstruction errors $\rho_{k_{\Text{stop}}}$ for the different regularizations are summarized in \tabref{tab:NumResNF-NLTest} along with the number of CG-iterations required to reach the optimum. Exemplary reconstructions for $s \in \{ 0, 0.5, 1 \}$ are shown in \figref{fig:NumResNF-RegTest-Obj}.
	    \begin{table}[htb!]
	      \centering
	      \begin{tabular}{lccccc} 
		\toprule
		 Sobolev exponent $s$ & $ 0 $ & $ 0.25 $ & $ 0.5 $ & $0.75$ & $1$ \\
		\midrule
		 Reconstruction error $\rho_{k_{\Text{stop}}}$ & $ 0.22 $ & $ 0.16 $ & $ 0.14 $ & $0.13$ & $0.13$ \\
		 Total CG-iterations & $ 244 $ & $ 252 $ & $ 273 $ & $438$ & $698$ \\
		\bottomrule
	      \end{tabular}
	      \caption{Reconstruction statistics for numerical near-field test cases using Sobolev regularization of different order $s$, parametrized by Gramians $\cG_{\mX_{\Text{dis}}}$ of the form \eqref{eq:GramXDisc}. Computed for a single object $\bN^\dagger$ of magnitude $\norm{\bN^\dagger} = \pi$ at Fresnel number $\NF = 0.01$ and noise level $\varepsilon = 3\,\%$. Stop rule: ``best stop''. Further parameters according to \tabref{tab:NumResNFSetup}. Exemplary reconstructed objects $\bN_{k_{\Text{stop}}}$ visualized in \figref{fig:NumResNF-NLTest-Obj}.}
	      \label{tab:NumResNF-RegTest}
	    \end{table}
	       \begin{figure}[hbt!]
	  \centering
	  \subfloat[Exact object $\bN^\dagger$]{\includegraphics[height=.34\textwidth]{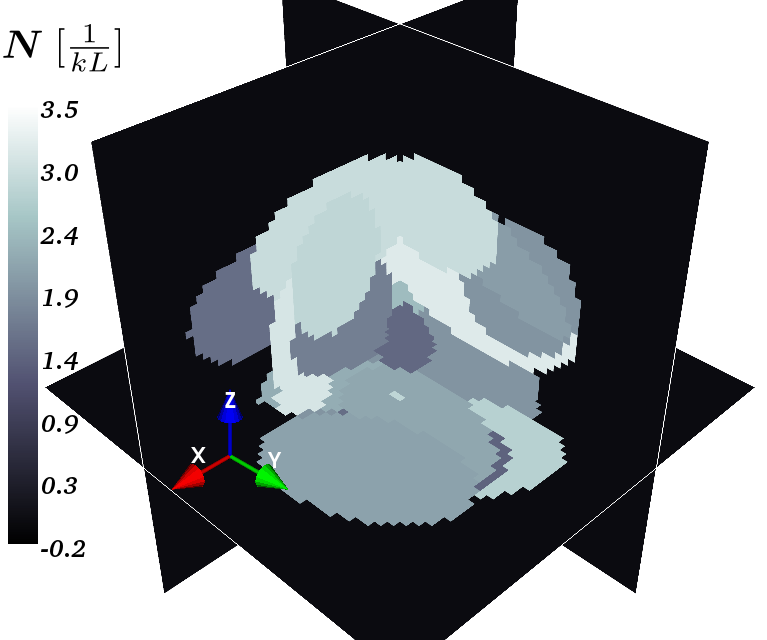} \label{fig:NumResNF-RegTest-Obj-ex}  }
	  \hfill
	  \subfloat[Reconstruction $\bN_{k_{\Text{stop}}}$ for $s  = 0$]{\includegraphics[height=.34\textwidth]{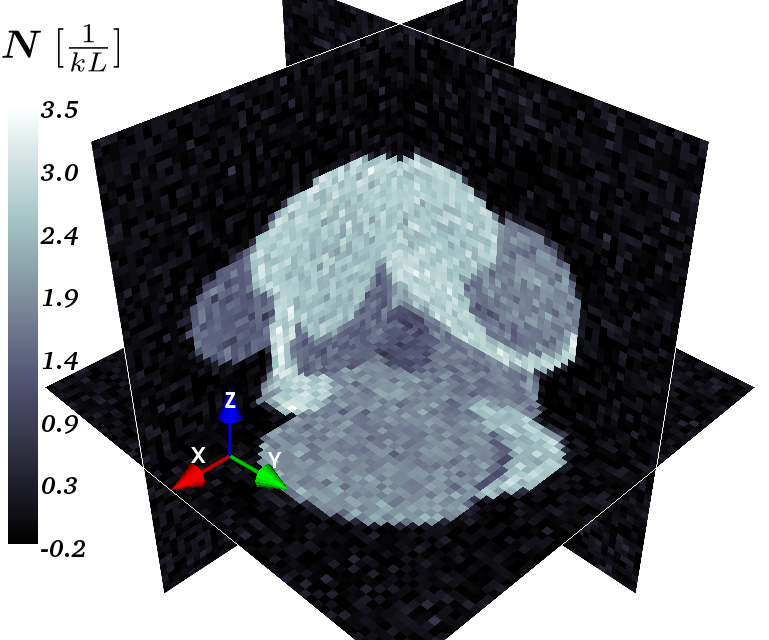} \label{fig:NumResNF-RegTest-Obj-s1}  }
	  \hfill
	  \subfloat[Reconstruction $\bN_{k_{\Text{stop}}}$ for $s  = 0.5$]{\includegraphics[height=.34\textwidth]{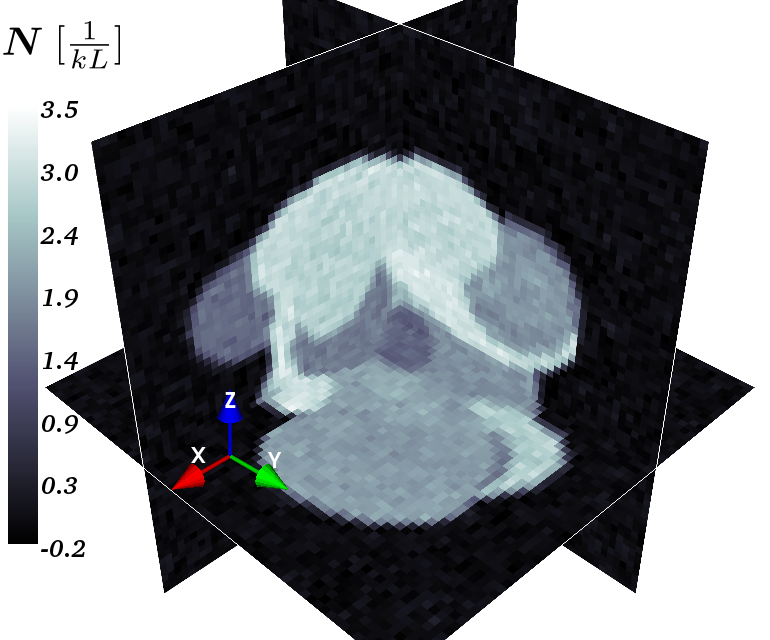} \label{fig:NumResNF-RegTest-Obj-s3}  }
	  \hfill
	  \subfloat[Reconstruction $\bN_{k_{\Text{stop}}}$ for $s  = 1$]{\includegraphics[height=.34\textwidth]{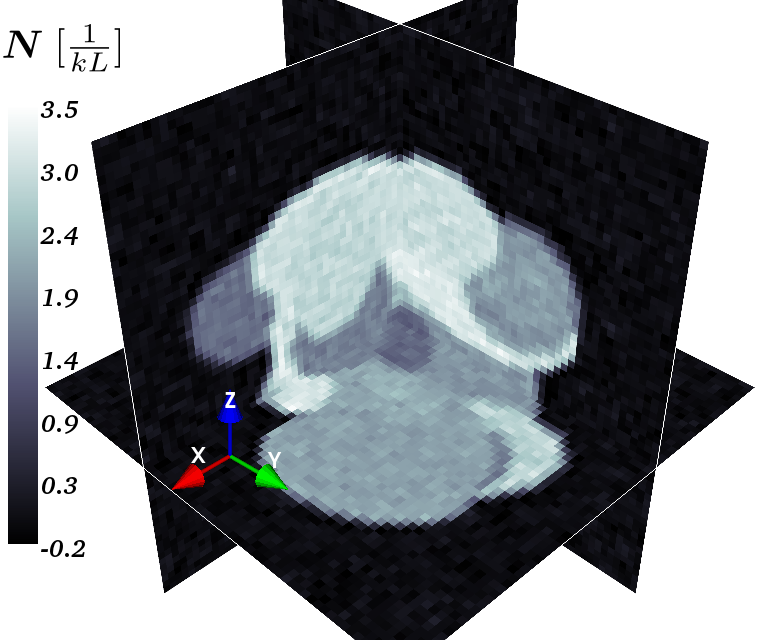} \label{fig:NumResNF-RegTest-Obj-s5}  }
	  \hfill
	  \caption{Exemplary near-field tomography results by \algref{alg:PCT} for Sobolev regularization of different order $s$. Details summarized in \tabref{tab:NumResNF-NLTest}. 
	  The choice $s= 0$ represents standard $L^2$-regularization. Choosing $s>0$ prescribes a certain regularity (see \sref{S:Sobolev}) and thus suppresses noise in the reconstructed objects, possibly at the expense of sharpness and performance. \label{fig:NumResNF-RegTest-Obj}}
	  \end{figure}

	  Visual comparison of the reconstructions in \figref{fig:NumResNF-RegTest-Obj-s1}-d with the exact object (\figref{fig:NumResNF-RegTest-Obj-ex}) seems to confirm the noise-suppressing effect of the Sobolev regularization. The error statistics in \tabref{tab:NumResNF-RegTest} indeed show a quantitative improvement of the results by almost a factor of two compared to the case $s = 0$ representing $L^2$-regularization. Notably, however, the reduction of noise comes at the expense of a slight blur of the edges in the reconstructed objects - just like in the (computationally much less expensive) case of postprocessing with a Gaussian filter. Moreover, \tabref{tab:NumResNF-RegTest} reveals that the numerical result hardly improves if the Sobolev exponent is increased beyond $s = 0.5$, whereas the number of required CG-iterations and thus the computational costs is almost three times greater for $s =1$. The reason is that the underlying assumption of a $H^1$-Sobolev regularity (compare \sref{S:Sobolev}) is overestimated for the realistic ``edgy'' objects simulated here, which causes slow convergence of the regularized Newton method.
	  
	  Indeed, the optimal compromise between noise suppression and performance seems to be given by the choice of $H^{  0.5 }$ regularization terms. Even if not noise but systematic errors are predominant in the data, the resulting damping of high Fourier frequencies may promote robustness of the Newton algorithm by imposing a hierarchical reconstruction from coarser to finer lengthscales.
	  
	  
	  \subsubsection{Convergence Rates and Effects of the Missing Wedge}
	  
	  As in \sref{SS:NumResFF-3DRecon} for the case of far-field tomography, we investigate the effect of incomplete data due to a missing wedge in the recorded incident angles $\theta \in [0^\circ; \theta_{\Text{rec}})$. In this final parametric study, a moderately strong object $\norm{\bN^\dagger} = \pi$ is reconstructed at a Fresnel number of $\NF = 0.01$ for different ranges of recorded incident angles $\theta_{\Text{rec}} \in \{ 180^\circ, 165^\circ, 150^\circ, 135^\circ, 120^\circ \} $ at different data noise levels $0.1 \, \% \leq \varepsilon \leq 10 \, \%$. Note that the  intensity data  $\bI^{\Textbf{err}}$ is always simulated for a fixed number of 128 incident angles, regardless of the size of the missing wedge. Motivated by the results of the preceding paragraph, Sobolev $H^s$-regularization of order $s = 0.5 $ is examined in comparison to computations using the standard $L^2$-penalty. Moreover, the Newton iterations are stopped according to the realistically implementable \emph{discrepancy principle} with $\tau =1$ \eqref{eq:Discrepancy}. Results by this stop rule are compared to reconstruction errors obtained for the optimal yet artificial ``best stop'' criterion applied before. The observed convergence rates for the different incident angle ranges, regularizations and stop rules are plotted in \figref{fig:NumResNF-WedgeTest-Conv}. Exemplary reconstructions for a noise level of $\varepsilon = 0.5 \, \%$ and $H^{0.5 }$-regularization are shown in \figref{fig:NumResNF-WedgeTest-Obj}.
	  	  	\begin{figure}[hbt!]
	 \centering
	  \subfloat[Sobolev $H^{0.5}$-regularization]{\includegraphics[width=.49\textwidth]{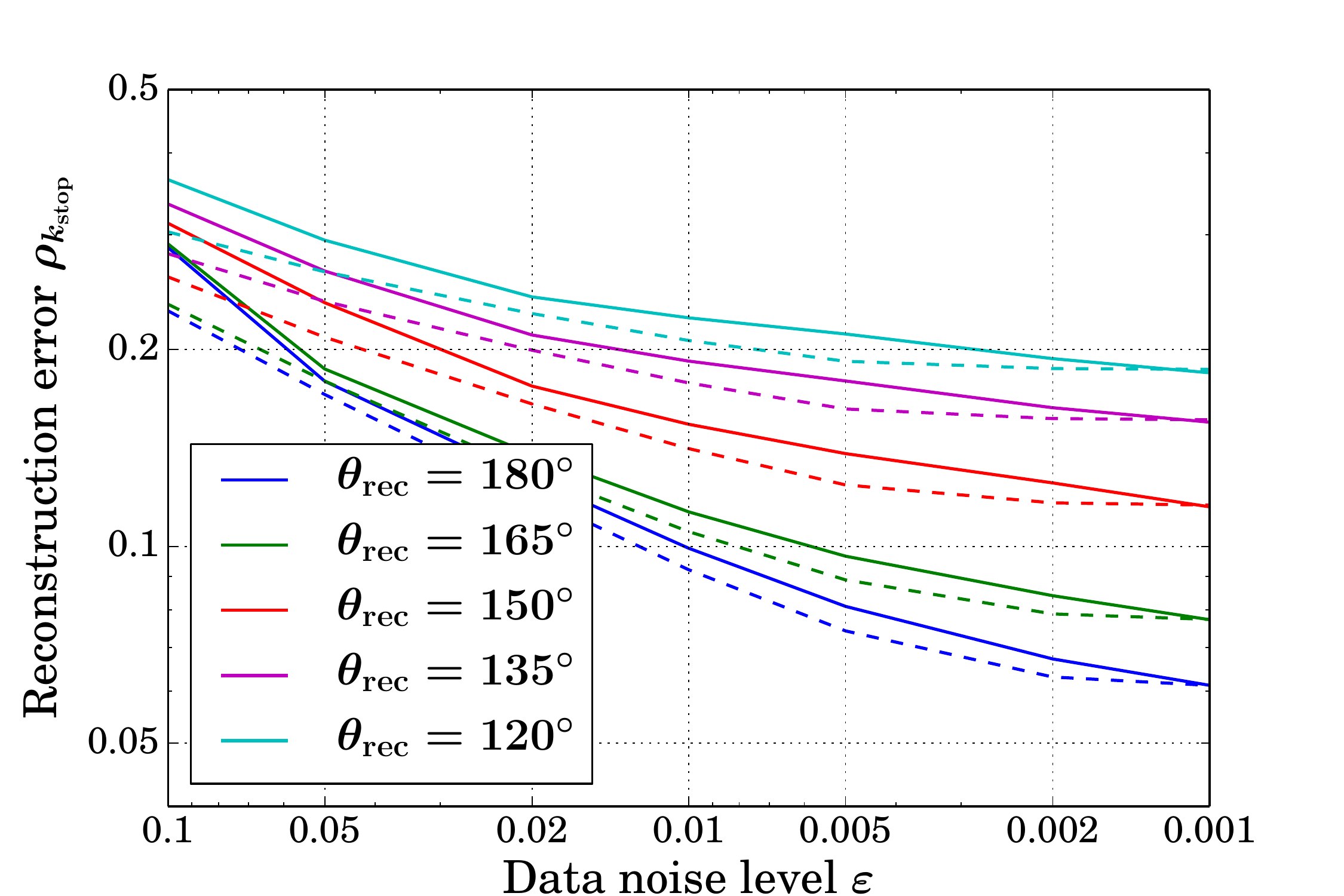}  \label{fig:NumResNF-WedgeTest-Conv-Hs}  }
	  \hfill
	  \subfloat[$L^2$-regularization]{\includegraphics[width=.49\textwidth]{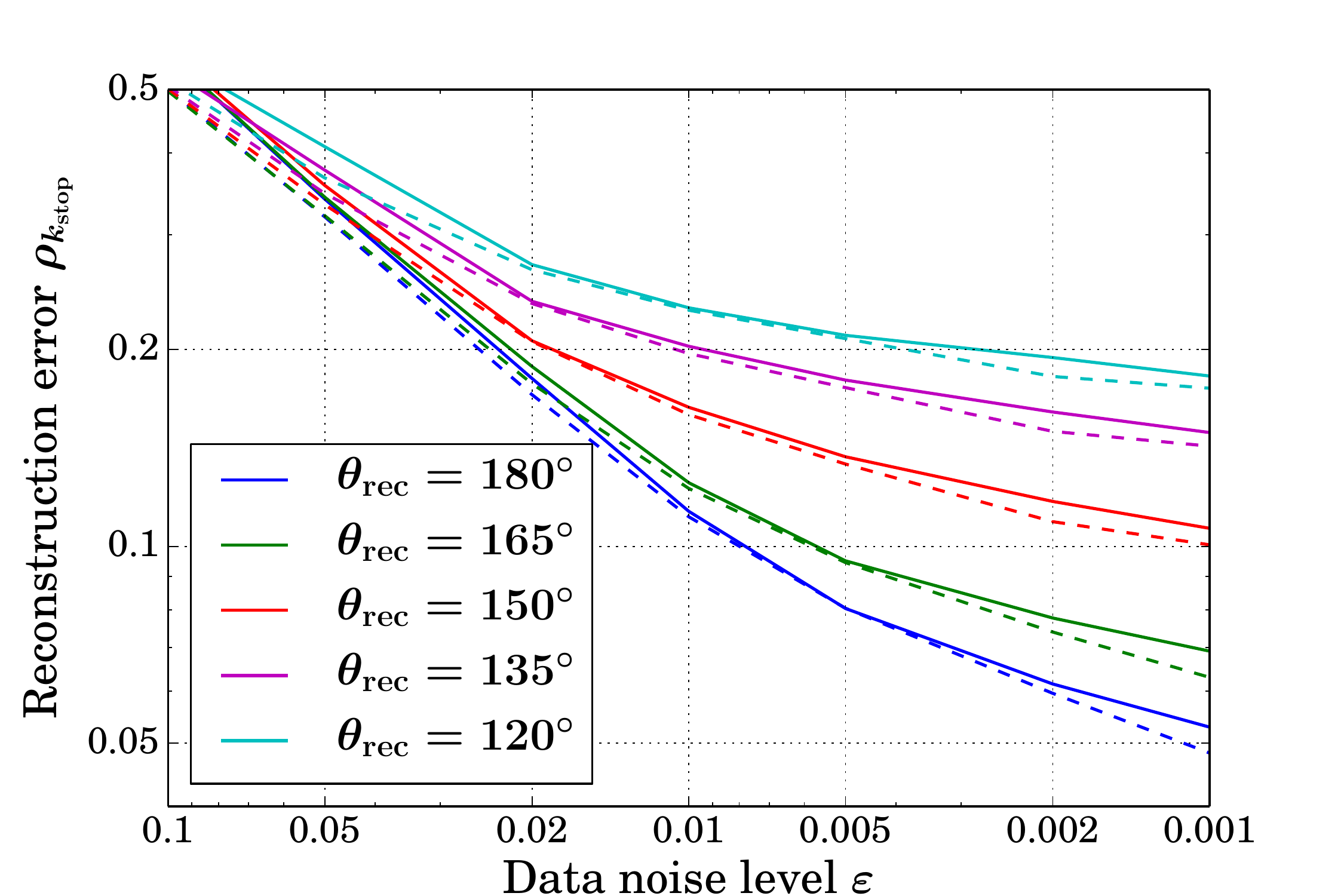}  \label{fig:NumResNF-WedgeTest-Conv-L2}  }
	 \caption{Numerical convergence rates for near-field tomography by \algref{alg:PCT} with the (Gaussian) data error level $\varepsilon$ for different maximum incident angles $\theta_{\Text{rec}}$ in the simulated intensity data $\bI^{\Text{err}}$ (missing wedge: $180^{\circ} - \theta_{\Text{rec}}$). Parameters: $\NF = 0.01$, $\norm{\bN^\dagger} = \pi$. Regularization: $H^{0.5 }$- (a) and $L^2$-penalty terms (b). Stop rule: discrepancy principle (solid lines) vs. ``best stop'' (dashed). Exemplary reconstructions and exact phantom visualized in \figref{fig:NumResNF-WedgeTest-Obj}.} \label{fig:NumResNF-WedgeTest-Conv}
	\end{figure}
	  
	  The curves in the log-log-plots in \figref{fig:NumResNF-WedgeTest-Conv} show no global algebraic convergence rates for any of the considered regularizations and stop rules - even in the case $\theta_{\Text{rec}} = 180^\circ$, i.e.\ without a missing wedge. The observed error decay for data noise levels below $1 \, \%$ is indeed rather slow where the $L^2$-regularization turns out to perform better in this low noise regime: from $\varepsilon = 1 \,\%$ to $\varepsilon = 0.1 \,\%$, the numerical reconstruction error merely reduces by a factor of $\approx 2$ in the $L^2$-case and even less for Sobolev regularization (compare blue curves in \figref{fig:NumResNF-WedgeTest-Conv-Hs} and \figref{fig:NumResNF-WedgeTest-Conv-L2}) This is a manifestation of the \emph{ill-posedness} of the considered phase contrast tomography problem, involving not only a mildly ill-posed Radon inversion (see \sref{S:RadonIll-posed}) but also (supposedly) more severely ill-posed phase reconstructions. However, the subalgebraic convergence may also partly be attributed to numerical difficulties as the regularization parameter for small $\varepsilon \sim 10^{-3}$ is iteratively reduced so strongly that up to $\sim 100$ CG-iterations are required for the final Newton steps, indicating a very ill-conditioned problem. This excessive step number increase is even more pronounced for the $H^{0.5}$-penalty term.
	  
	  Yet, it is noteworthy that Sobolev $H^{0.5}$-regularization supplemented with the non-ideal discrepancy principle outperforms even ``best-stopped'' $L^2$-regularization up to moderate noise levels $\varepsilon \gtrsim 0.01$ for all of the considered missing wedges in \figref{fig:NumResNF-WedgeTest-Conv}. Indeed, the good agreement between the dashed and solid lines in \figref{fig:NumResNF-WedgeTest-Conv-Hs} and \figref{fig:NumResNF-WedgeTest-Conv-L2} shows that the chosen discrepancy principle with $\tau = 1$ provides a \emph{quasi-optimal} stopping rule in the considered problem setup.
	       \begin{figure}[hbt!]
	  \centering
	  \subfloat[Exact object $\bN^\dagger$]{\includegraphics[height=.34\textwidth]{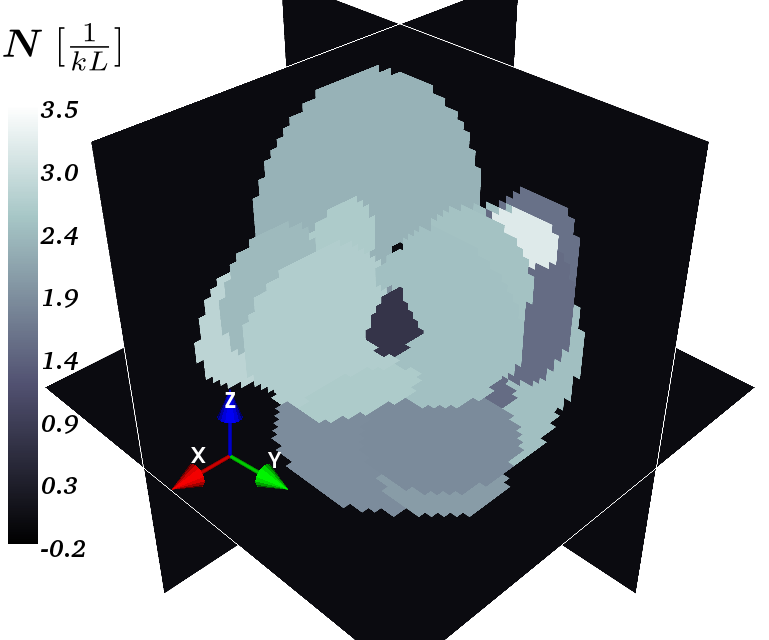} \label{fig:NumResNF-WedgeTest-Obj-ex}  }
	  \hfill
	  \subfloat[Reconstruction  for $\theta_{\Text{rec}} = 165^\circ$]{\includegraphics[height=.34\textwidth]{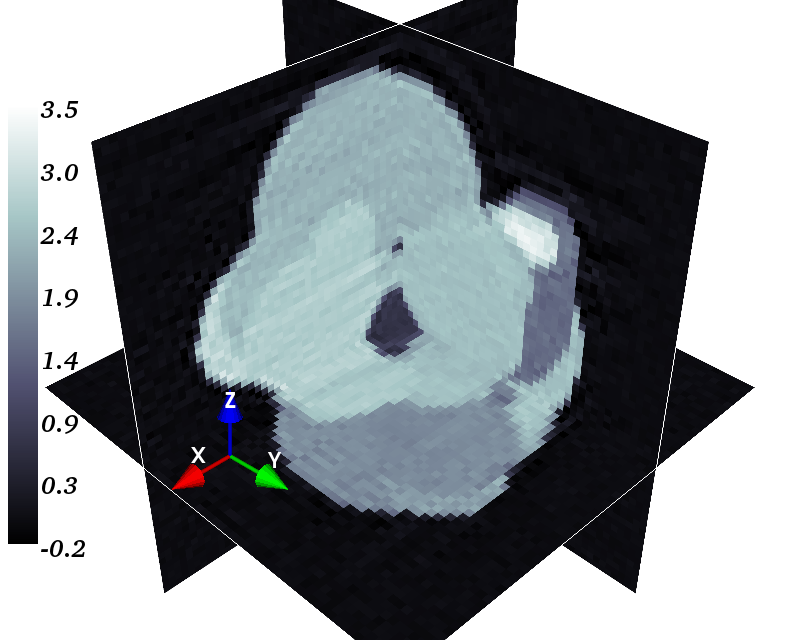} \label{fig:NumResNF-WedgeTest-Obj-w2}  }
	  \hfill
	  \subfloat[Reconstruction  for $\theta_{\Text{rec}} = 150^\circ$]{\includegraphics[height=.34\textwidth]{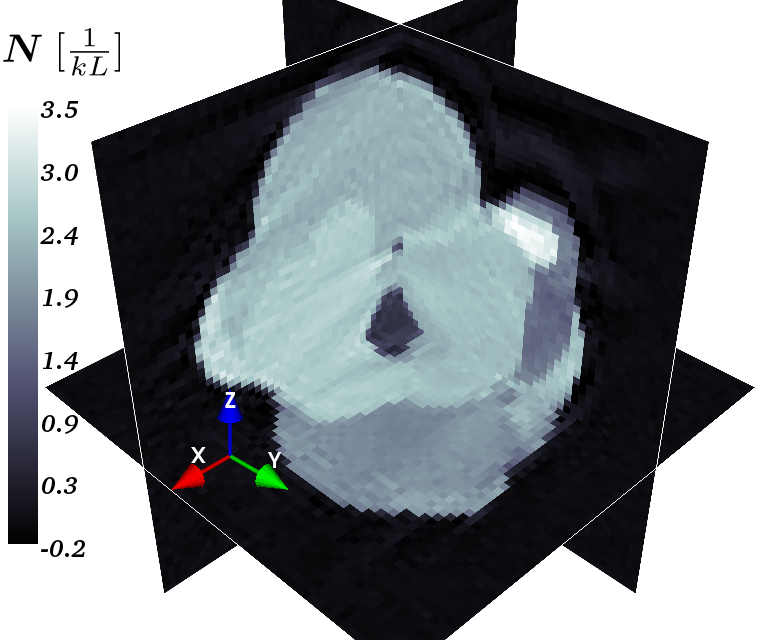} \label{fig:NumResNF-WedgeTest-Obj-w3}  }
	  \hfill
	  \subfloat[Reconstruction  for $\theta_{\Text{rec}} = 120^\circ$]{\includegraphics[height=.34\textwidth]{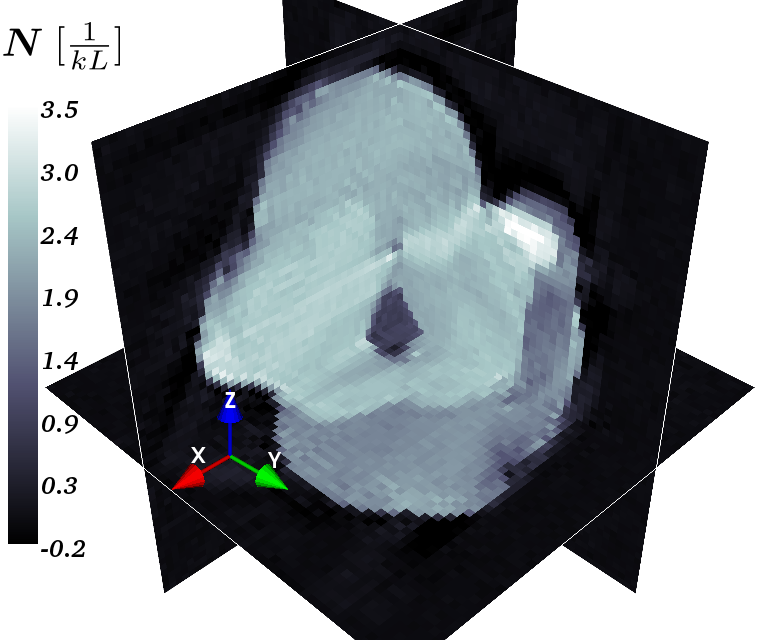} \label{fig:NumResNF-WedgeTest-Obj-w5}  }
	  \hfill
	  \caption{Exemplary near-field tomography results by \algref{alg:PCT} for different ranges of recorded incident angles $\theta \in [0^\circ; \theta_{\Text{rec}})$ in the simulated intensity data, representing \emph{missing wedges} of $15^\circ$ (b), $30^\circ$ (c) and  $60^\circ$ (d), respectively. Noise level: $\varepsilon = 0.5 \, \%$. Regularization: $H^{0.5}$-penalty. Stop rule: discrepancy principle with $\tau = 1$. For details and quantitative results, see \figref{fig:NumResNF-WedgeTest-Conv}.
	  \label{fig:NumResNF-WedgeTest-Obj}}
	  \end{figure}
	  
	  On the other hand, the uniformity of the curves in \figref{fig:NumResNF-WedgeTest-Conv} for different $\theta_{\Text{rec}}$ reveals a certain robustness of the Newton reconstruction against data incompleteness: even in the case $\theta_{\Text{rec}} = 120^\circ$ of a missing wedge of $60^\circ$ covering a third of the relevant tomographic incident angles, a moderate noise level of $\varepsilon = 1 \, \%$ still allows a numerical solution within an error of around $20 \, \%$. This interpretation is supported by the exemplary reconstructions visualized in \figref{fig:NumResNF-WedgeTest-Conv}: like in the far-field example (\figref{fig:NumResNF-GenRecon}), these remain widely accurate for incomplete intensity data apart from slight halo- and stripe artifacts emanating preferably along the axes of the omitted tomographic projections. As the implicit completion of the data corresponds to analytic continuation in Fourier space according to the Fourier Slice \thmref{thm:FourierSlice} - an operation which is highly susceptible to noise - one might have expected considerably worse convergence rates and artifacts. As this is not observed, the regularization seems successfully damp out instabilities of the reconstruction due to incomplete tomographic data.
	
	We conclude this section by summarizing the results of our parametric study of near-field phase contrast tomography for pure phase objects:
	 \vspace{1em}
	\begin{res}[Newton-based Near-Field Tomography of Pure Phase Objects] \label{res:NFTomoByNewton}
	  	Pure phase objects  in near-field tomography are stably reconstructed by \algref{alg:PCT} without incorporating further a priori knowledge as support constraints or via the initial guess. The method is accurate up to moderately strong, i.e.\ non-phase wrapping objects and works best if the sample's characteristic lengthscales correspond to Fresnel numbers in the order of 1. Sobolev norm regularization terms allow noise suppression in the numerical solution. The reconstruction is robust against incomplete data due to a ``missing wedge'' in the tomographic projections and shows stable, yet subalgebraic convergence with the data noise level if the Newton iterations are stopped according to the discrepancy principle, which is observed to be quasi-optimal.
	\end{res}
	\vspace{1em}

	  \begin{subsection}{Reconstruction of General Objects}  \label{SS:NumResNFSim-GenObjects}
	  
	    In the preceding section, we have seen that our regularized Newton approach permits stable and accurate near-field tomography of pure phase objects. Notably,
	    our uniqueness statement for phase contrast tomography (\cref{cor:UniqueNFPCT}) in principal also allows for the reconstruction of \emph{general}, complex-valued objects $\bN^\dagger =   \bdelta^\dagger - \I \bbeta^\dagger$ composed of multiple materials with different ratios $\beta/\delta$ of absorption and refraction. However, it has been seen in \sref{SS:NumResFF-3DRecon} for the far-field case that reconstructing the absorption $\bbeta$ as independent degrees of freedom renders the numerical solution by \algref{alg:PCT} less stable. In this section, we provide a qualitative proof of concept that \emph{near-field} phase contrast tomography of general objects is feasible with the Newton-type method developed in \chapref{C:NumMeth}.
	    
	               \begin{figure}[hbt!]
	  \centering
	  \subfloat[Exact object $\bN^\dagger = \bdelta^\dagger - \I \bbeta^\dagger$]{\includegraphics[height=.32\textwidth, clip=true, trim = -4cm 0cm -4cm 0cm]{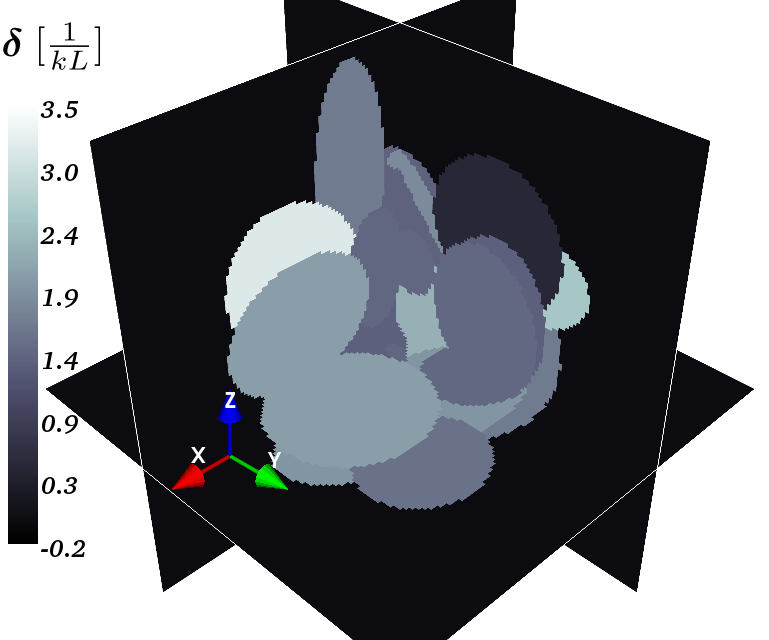} \hfill \includegraphics[height=.32\textwidth, clip=true, trim = -4cm 0cm -4cm 0cm]{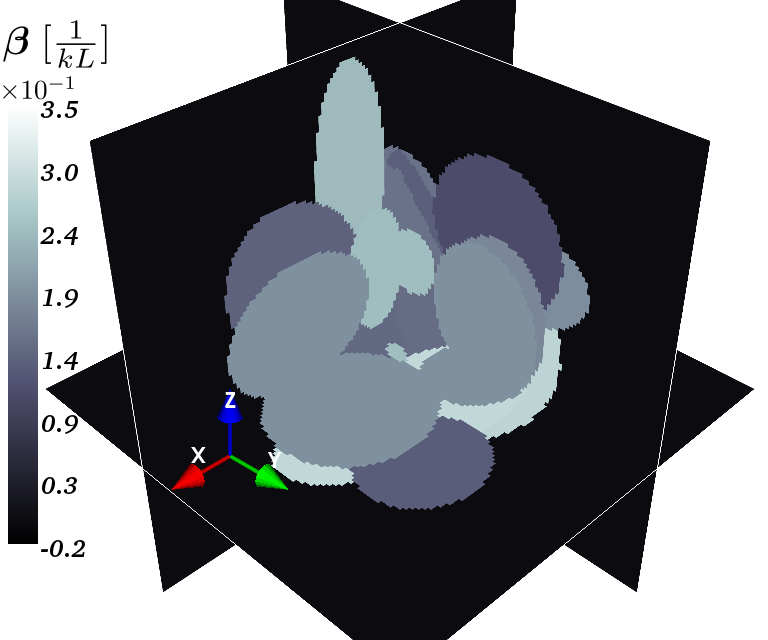} \label{fig:NumRes-NF3D-exObj} }
	  \\
	  \subfloat[Reconstruction $\bN_{k_{\Text{stop}}}^{\Text{gen}} =   \bdelta_{k_{\Text{stop}}}^{\Text{gen}} - \I \bbeta_{k_{\Text{stop}}}^{\Text{gen}}$ as a general object (rel. $L^2$-error: $\rho_{k_{\Text{stop}}}^{\Text{gen}}= 11.9 \, \%$)]{\includegraphics[height=.32\textwidth, clip=true, trim = -4cm 0cm -4cm 0cm]{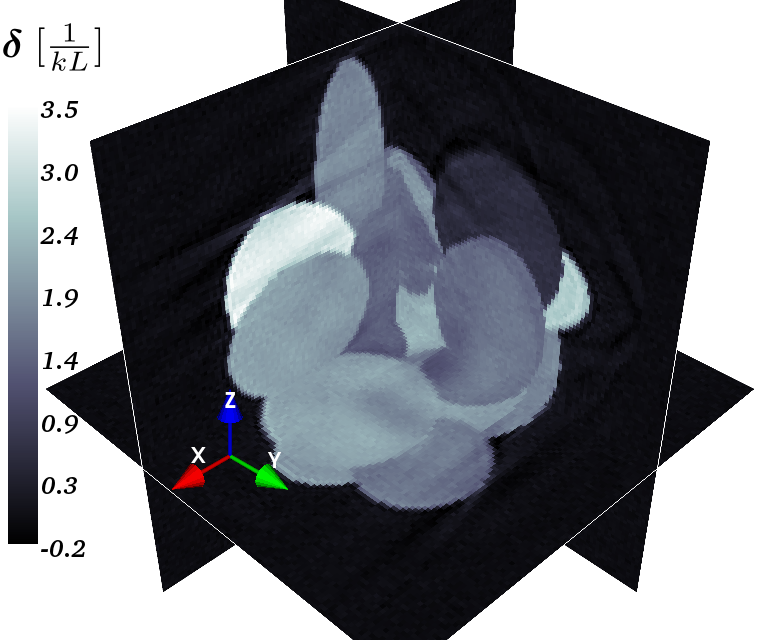} \hfill \includegraphics[height=.32\textwidth, clip=true, trim = -4cm 0cm -4cm 0cm]{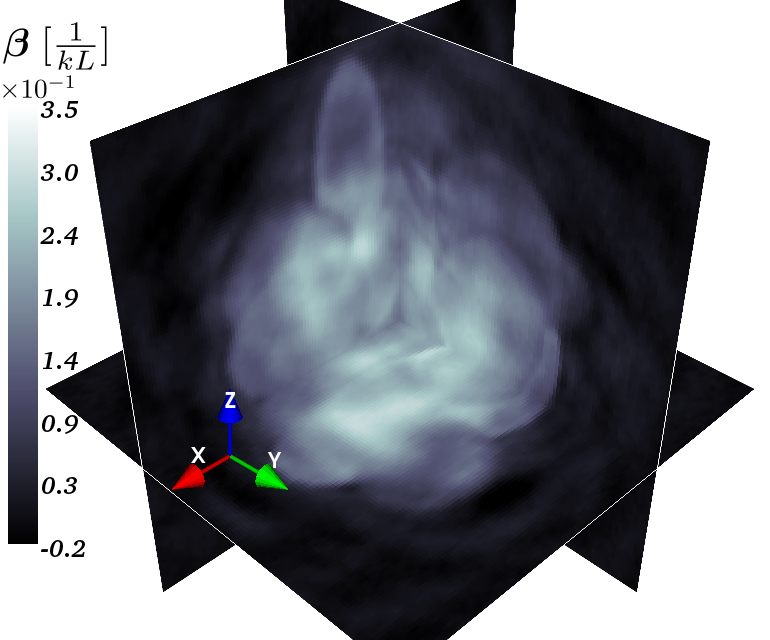} \label{fig:NumRes-NF3D-genRec} }
	  \\
	  \subfloat[Reconstruction $\bN_{k_{\Text{stop}}}^{\Text{single}} = \bdelta_{k_{\Text{stop}}}^{\Text{single}} - \I \bbeta_{k_{\Text{stop}}}^{\Text{single}}$ as a single-material object $\bbeta =  \bdelta / 10$ ($\rho_{k_{\Text{stop}}}^{\Text{single}} = 12.3 \, \%$)]{\includegraphics[height=.32\textwidth, clip=true, trim = -4cm 0cm -4cm 0cm]{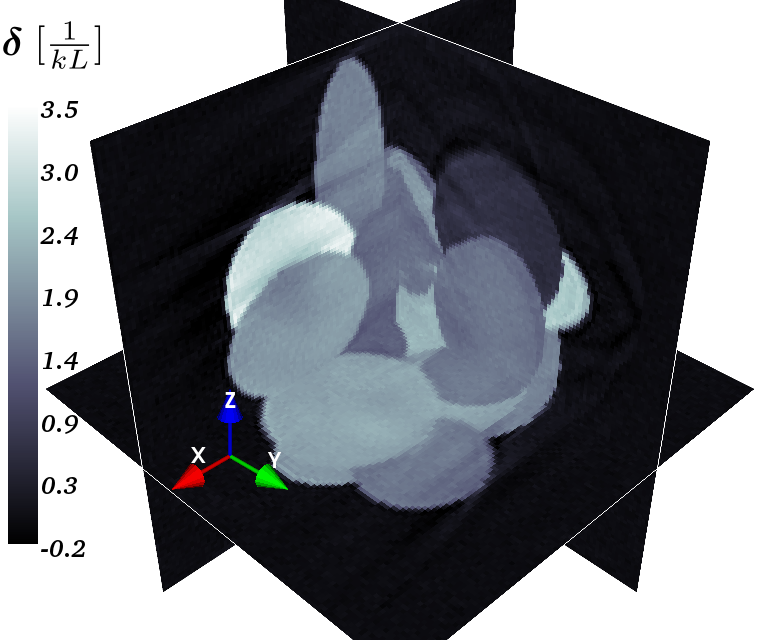} \hfill \includegraphics[height=.32\textwidth, clip=true, trim = -4cm 0cm -4cm 0cm]{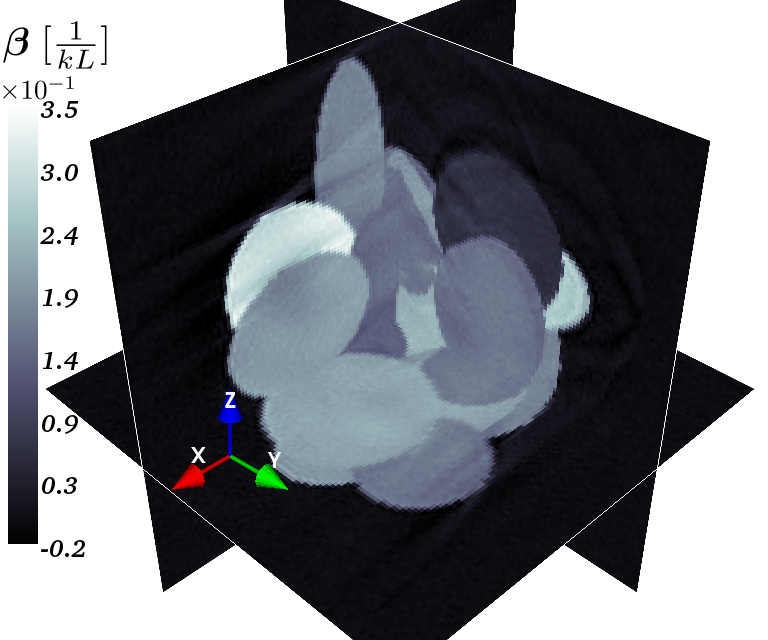} \label{fig:NumRes-NF3D-singleRec} }
	  \caption{Near-field tomography results for a general object $\bN^\dagger = \bdelta^\dagger - \I \bbeta^\dagger \in \mC^{128^3}$ by \algref{alg:PCT} with non-coupled absorption and refraction of mean ratio $c_{\bbeta/ \bdelta} = 0.1$. Simulated intensities $\bI^{\Textbf{err}} = F_{\Text{dis}}( \bN^\dagger ) + \Textbf{err} \in \mR^{256^3}$ contain $\varepsilon = 0.5 \, \%$ noise and a missing wedge of $20^\circ$. Fresnel number: $\NF = 0.0025$. Object magnitude: $\norm{\bN^\dagger} = \pi$. Regularization by $H^{0.5}$-penalty term. Stopping index $k_{\text{stop}}$ chosen by discrepancy principle \eqref{eq:Discrepancy} with $\tau = 1$. Further parameters according to \tabref{tab:NumResNFSetup} } \label{fig:NumResNF-GenRecon}
	  \end{figure}
	    
	    For a direct comparison, we use the same moderately strong ellipsoid object $\bN^\dagger = \bdelta^\dagger - \I \bbeta^\dagger$, $\norm{\bN^\dagger} = \pi$ with mean absorption-refraction-ratio $c_{\bbeta/ \bdelta} = 0.1$ as in  \sref{SS:NumResFF-3DRecon}. However, no reference signal, special initial guess or support constraint is assumed as neither is needed for near-field reconstructions. Note that the resolutions in object and image space $\mX_{\Text{dis}} = \mC^{128^3}, \mY_{\Text{dis}} = \mR^{256^3}$ in this numerical example deviate from \tabref{tab:NumResNFSetup}, whereas all the other simulation parameters are chosen exactly as described in \sref{SS:NumResNFSim-Setup}. We consider the optimal regime of our method according to \resref{res:NFTomoByNewton}, choosing $\NF = 0.0025$ such that the characteristic object lengthscales correspond to Fresnel numbers $\tilde \NF \sim  1$. Moreover, Sobolev $H^{0.5}$-regularization is used and the discrepancy principle \eqref{eq:Discrepancy} with $\tau = 1$ is applied as a stop rule. A data noise level of $\varepsilon = 0.5 \, \%$ is prescribed, comparable to the Poisson errors in \sref{SS:NumResFF-3DRecon}. Likewise the missing wedge $20^\circ$ of is retained, simulating intensity data $\bI^{\Textbf{err}} = F_{\Text{dis}}( \bN^\dagger ) + \Textbf{err}$ for incident angles $\theta \in [0^\circ; 160^\circ)$. For comparison, we compute a second reconstruction, assuming a false but optimally chosen single-material constraint with fixed $\beta$-$\delta$-coupling parameter $c_{\bbeta/ \bdelta} = 0.1$ in \algref{alg:PCT}. For the independent recovery of $\bdelta^\dagger$ and $\bbeta^\dagger$, the \emph{average} ratio $c_{\bbeta/ \bdelta}$ is accounted for  in the regularization by weighting deviations in $\bbeta$ with $1/c_{\bbeta/ \bdelta}$ as outlined in \sref{SS:NumResFF-3DRecon}.
	    
	    The resulting objects $\bN_{k_{\Text{stop}}}^{\Text{gen}}$ and $\bN_{k_{\Text{stop}}}^{\Text{single}}$ from both the general- and the single-material reconstruction are visualized in \figref{fig:NumResNF-GenRecon}. The number of required Newton(CG)-iterations for the discrepancy principle to terminate the reconstructions are $17(748)$ for the general object and $16(658)$ in the case of the single-material constraint, respectively. Hence, we find that the additional computational effort due to the independent recovery of the absorption is minor. Likewise, the reconstruction results in \figref{fig:NumRes-NF3D-genRec} and \ref{fig:NumRes-NF3D-singleRec} show neither a clear advantage nor drawback of the general approach: the refractive part $\bdelta^\dagger$ is equally well resolved by either solution, merely perturbed by weak artifacts arising from the missing wedge, the relative $L^2$-errors of the final iterates $\bdelta_{k_{\Text{stop}}}$ being $11.2 \, \%$ (general object) and $11.9 \, \%$ (single-material constraint). Other than in the far-field case (see \sref{SS:NumResFF-3DRecon}), we thus observe that the simultaneous reconstruction of $\bbeta^\dagger$ as an independent parameter  at least has no negative backlash onto the numerical recovery of refractive phase shifts.
	    
	    On the other hand, the reconstructed absorption in \figref{fig:NumRes-NF3D-genRec} is unfortunately highly inaccurate with an error of $\norm{\bbeta_{k_{\Text{stop}}}^{\Text{gen}} - \bbeta^\dagger}_2 / \norm{\bbeta^\dagger }_2 = 40.8 \, \%$. The latter is even larger than the deviation of $32.2 \, \%$ for the single-material approximation $\bbeta_{k_{\Text{stop}}}^{\Text{single}}$. By visual comparison of the right hand column in \figref{fig:NumRes-NF3D-genRec} and \ref{fig:NumRes-NF3D-exObj}, we find that the absorptions of the individual ellipsoids tend to be distributed in a qualitatively correct manner. However, the obtained reconstruction is subject to strong low-frequency halo-like artifacts. A possible explanation for this phenomenon might indeed be poor \emph{phase} contrast in these low-frequencies (cf. \sref{SS:ContrastFormationNF}), leading to errors in the reconstructed refraction $\bdelta_{k_{\Text{stop}}}^{\Text{gen}}$ which can be identified as slight wafting inhomogeneities of the ellipsoids in \figref{fig:NumRes-NF3D-genRec} (left plot). This incorrect wobbling background, being possibly only of relative magnitude $ \sim 5 \, \%$ in $\bdelta$,  may manifest significantly in the simultaneously reconstructed $\bbeta$ due to the much smaller magnitude of the latter (here: $\bbeta \sim \bdelta/10$). According to this interpretation, ill-posedness in the recovery of the refraction around the corresponding low-phase-contrast-zeros of the CTF (see \figref{fig:CTF}) negatively affects the reconstruction of the absorption - and vice verser.
	    
	    Based on the poor reconstruction of $\bbeta^\dagger$ in \figref{fig:NumResNF-GenRecon}, one may even come to the conclusion that an independent numerical solution for this parameter by our regularized Newton method is never sufficiently accurate to provide additional information. This would constitute a major practical restriction of the theoretical uniqueness statement in \cref{cor:UniqueNFPCT}. In order to disprove this, we substitute the absorptive part $\bbeta^\dagger$ of the object in \figref{fig:NumRes-NF3D-exObj} by scaled tomographic images of $\Text{Ta}_2\Text{O}_5$-coated nano-porous glass, experimentally observed in \cite{Holler2014PorousTomoDataset}. Thereby, we obtain an exact test object $\bN^\dagger$ with \emph{uncorrelated} $\bdelta^\dagger$ and $\bbeta^\dagger$ except for their common support. Apart from this modification, the simulation setup is completely retained. In particular, the refractive and absorptive parts are scaled such that  $\norm{\bN^\dagger} = \pi$ and an average $\beta$-$\delta$-ratio of $c_{\bbeta/ \bdelta} = 0.1$ is obtained as above. Reconstruction results by \algref{alg:PCT} for this modified object are shown in \figref{fig:NumResNF-GenRecon-Mod}.	    
 	    \begin{figure}[hbt!]
	    \centering
	    \subfloat[Exact object $\bN^\dagger = \bdelta^\dagger - \I \bbeta^\dagger$]{\includegraphics[height=.32\textwidth, clip=true, trim = -4cm 0cm -4cm 0cm]{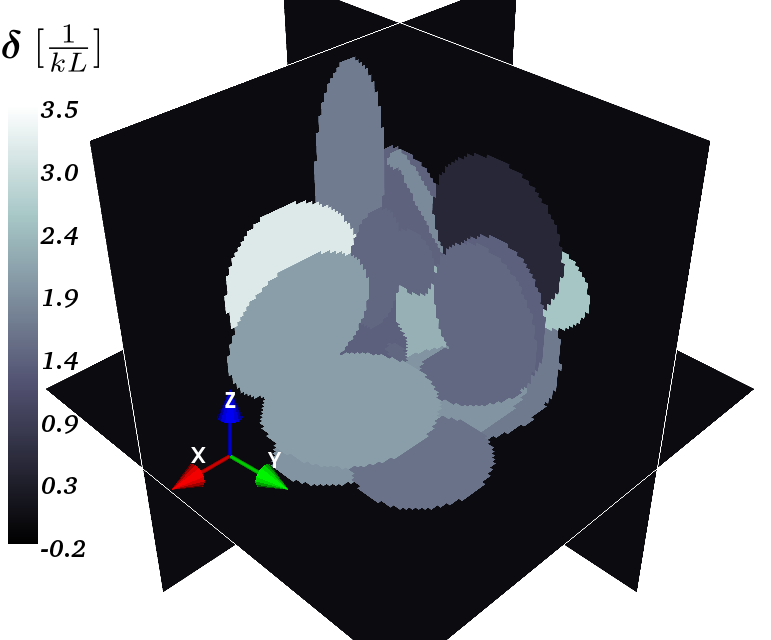} \hfill \includegraphics[height=.32\textwidth, clip=true, trim = -4cm 0cm -4cm 0cm]{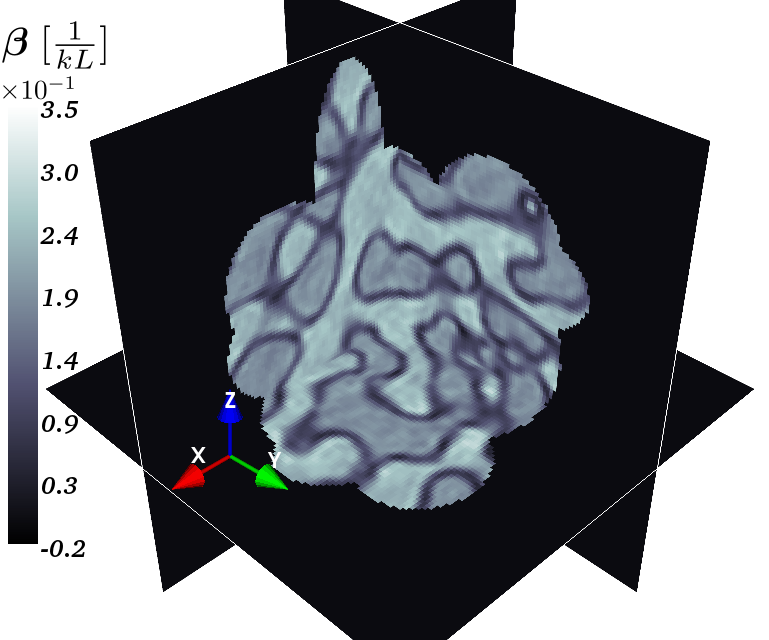} \label{fig:NumResNF-GenRecon-Mod-exObj} }
	    \\
	    \subfloat[Reconstruction $\bN_{k_{\Text{stop}}} =   \bdelta_{k_{\Text{stop}}} - \I \bbeta_{k_{\Text{stop}}}$ (relative $L^2$-error: $\rho_{k_{\Text{stop}}}= 11.8 \, \%$)]{\includegraphics[height=.32\textwidth, clip=true, trim = -4cm 0cm -4cm 0cm]{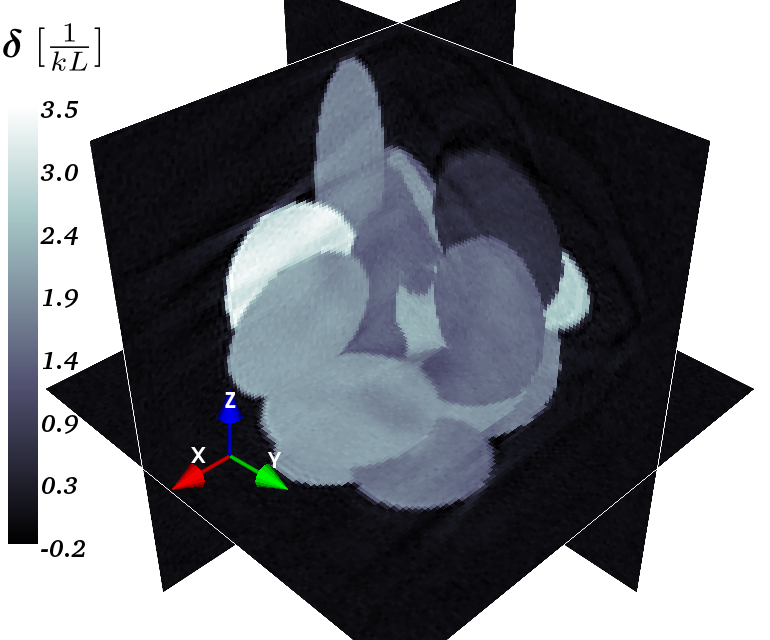} \hfill \includegraphics[height=.32\textwidth, clip=true, trim = -4cm 0cm -4cm 0cm]{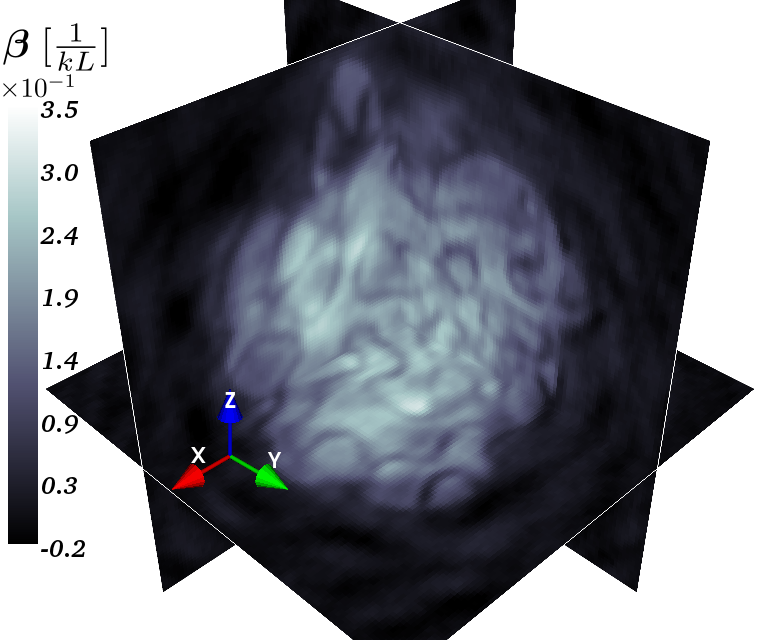} \label{fig:NumResNF-GenRecon-Mod-rec} }
	    \caption{Reproduction of \figref{fig:NumRes-NF3D-exObj}-b with a modified absorption-object $\bbeta^\dagger$, given by scaled tomographic images of nano-porous glass from \cite{Holler2014PorousTomoDataset}. Apart from the low-frequency artifacts in $\bbeta_{k_{\Text{stop}}}$, features in $\delta^\dagger$ and $\bbeta^\dagger$ are reconstructed qualitatively accurately and well-separated between the different components.} \label{fig:NumResNF-GenRecon-Mod}
	    \end{figure}
	    
	    As in \figref{fig:NumResNF-GenRecon}, dominant low-frequency artifacts are observed in the recovered absorption $\bbeta_{k_{\Text{stop}}}$, whereas the solution is accurate within the refractive part $\bdelta^\dagger$ up to a relative $L^2$-error of $\approx 11 \, \%$. Notably, however, the \emph{qualitative} structure of the absorption-object is accurately resolved (compare \figref{fig:NumResNF-GenRecon-Mod}, right column): the edges bounding the different glass- or air-filled segments are located correctly and may be identified clearly behind the ``veil'' of low-frequency errors. Moreover, it should be emphasized that features of $\bdelta^\dagger$ or $\bbeta^\dagger$ do not falsely manifest in the reconstruction of the other parameter: neither are there visible traces of porous structures in the left plot of \figref{fig:NumResNF-GenRecon-Mod-rec} nor spurious ellipsoid boundaries within the plotted slices of the reconstructed absorption $\bbeta_{k_{\Text{stop}}}$. Hence, refraction and absorption of the unknown object are indeed reconstructed in a cleanly separated way.
	    
	    It seems to be a general tendency that the reconstruction of general objects $\bN^\dagger = \bdelta^\dagger - \I \bbeta^\dagger$ via \algref{alg:PCT} works better the less correlated the refractive and absorptive part. Unfortunately, reality is closer to the simulation in \figref{fig:NumResNF-GenRecon} with almost perfectly correlated structures in $\delta$ and $\beta$ of merely a variable ratio $\beta/\delta$. Note that the obtained reconstructions are still significantly better than in the far-field simulation shown in \ref{fig:NumRes-FF3D-genRec} where non-physical \emph{negative} absorption values have been obtained in the numerical solution $\bbeta_{k_{\Text{stop}}}$. Indeed, the results of this section suggest that near-field tomography of general objects by regularized Newton-type methods is not hopeless but may even become quantitatively correct if the observed low-frequency artifacts can be suppressed. A promising approach could be to impose a \emph{loose} coupling of $\delta$ and $\beta$ by prescribing a maximum ratio $\beta/\delta$. This would certainly suppress the halo-like artifacts in $\bbeta_{k_{\Text{stop}}}$ outside the object support in \figref{fig:NumRes-NF3D-genRec} and \ref{fig:NumResNF-GenRecon-Mod-rec} as there are no comparable structures reconstructed in the refractive part $\bdelta_{k_{\Text{stop}}}$. Note, however, that this constraint is \emph{non-smooth} just like the prescription of positivity discussed in \sref{SS:Constraints} and may thus only be incorporated in future semismooth generalizations of this work's Newton-type approach.

	\begin{res}[Newton-based Near-Field Tomography of General Objects] \label{res:GenObjByNewton}
	  	Near-field tomography via \algref{alg:PCT} is applicable to general objects $N^\dagger = \delta^\dagger - \I \beta^\dagger$, simultaneously recovering refraction $\delta$ and absorption $\beta$. For moderately strong absorption $\beta^\dagger/\delta^\dagger \sim \frac 1 {10}$, no negative effects onto the recovery of $\delta^\dagger$ are observed
	  	compared to reconstructions assuming approximate single-material constraints.
	  	The reconstruction of $\beta^\dagger$ is found to be susceptible to low-frequency artifacts, which might be suppressed by loosely coupling refraction and absorption via non-smooth penalty terms.
	\end{res}

	 \end{subsection}

         \begin{subsection}{Evaluation of the Simultaneous Approach}  \label{SS:NumResNFSim-EvalSimApproach}
         
	    As outlined in \sref{S:CombinedApproach}, our iterative Newton-type approach to phase contrast tomography is motivated by two principal objectives: for once, the aim is to overcome the limitations of direct methods, such as CTF- or transport-of-intensity-based techniques, to the regimes of validity of the underlying linearizations. This has been achieved according to the results of \sref{SS:NumResNFSim-ParamStudy} where applicability of the regularized Newton method to both moderately strong objects and small Fresnel numbers has been demonstrated. Yet, this would also be possible with less computational effort by restricting to Newton-based phase retrieval, supplemented with an a posteriori Radon inversion by a method of our choice. 
	    A second motivation for our method, however, is to exploit the consistency conditions (see \thmref{thm:HelgasonLudwig}) between the diffraction patterns for different incident angles via \emph{simultaneous} tomographic- and phase reconstruction by a single algorithm. As argued in \sref{S:CombinedApproach}, this could render phase retrieval more stable and accurate, which we aim to verify in this section.
	    
	    To this end, we compare our simultaneous approach by \algref{alg:PCT} to the results obtained by applying an analogous regularized Newton-type method to the phase retrieval problem $\bI^{\Textbf{err}} =  F_{\Text{dis}}^{\Text{phase}}  ( \CR(\bN^\dagger) ) + \Textbf{err}$, which is constructed simply by omitting the Radon transform in the governing operators. Accordingly,  only the tomographic projections $\CR(\bN^\dagger)$ of the exact object $\bN^\dagger \in \mX_{\Text{dis}} = \mC^{64^3}$ are reconstructed in the latter approach where the phase recovery is computed completely independently for the different incident angles $\theta$. For comparison, we furthermore solve the phase retrieval problems using a CTF-based method by Peter Cloetens \cite{Cloetens1999}, Matthias Bartels \cite{BartelsDiss} and Martin Krenkel \cite{Krenkel2014BCAandCTF}, directly inverting the contrast transfer function \eqref{eq:CTF} with a cut-off around the zeros in \figref{fig:CTF}.
	    
	    In order to allow for a fair comparison with the iteratively regularized methods, the cut-off parameter of the CTF-method if optimized to obtain a minimum error in each reconstruction. Moreover, we restrict to pure phase objects of moderately small magnitude $\norm{\bN^\dagger} = 1$ so that nonlinearity of the problem does not obstruct usage of the CTF too severely. On the other hand, we are also interested in driving the competing Newton-type methods to their optimum in order to explore their \emph{principal} potential. We therefore choose the ``best stop'' rule to terminate the Newton methods and Sobolev $H^s$-regularization terms with $s = 0.5$ for the simultaneous recovery of $\bN^\dagger$ and $s = 1$ in the independent phase reconstruction of the projections $\CR(\bN^\dagger)$, see \sref{SS:NumResNFSim-ParamStudy}. The distinction in $s$ takes into account a higher regularity of the latter due to the mild smoothing effect of the Radon transform discussed in \sref{SS:RadonContiInv}.
	    
	    Our numerical comparison is based on three different random test objects $\bN^\dagger$ of magnitude $\norm{\bN^\dagger} = 1$ of the form described in \sref{SS:NumResNFSim-Setup} at a Fresnel number of $\NF = 0.01$, each reconstructed at various data noise levels $10 \, \% \geq \varepsilon  \geq  0.1 \, \%$. Furthermore, we a impose cylindrical support enclosed within the $64^3$ voxel cube $\mX_{\Text{dis}}$ for the objects $\bN^{\dagger}$. By rotational symmetry, all projections $\CR(\bN^{\dagger})$ are then contained within a $64 \times 64$ pixel domain, which is taken as the support for the Newton- and CTF phase reconstructions. The remaining parameters are chosen according to \tabref{tab:NumResNFSetup}. However, we reduce the number of incident angles $\theta \in [0^{\circ}; 180^\circ)$ (no missing wedge) to 64, i.e.\ consider intensity data $\bI^{\Textbf{err}} = F_{\Text{dis}}(\bN^\dagger) + \Textbf{err} \in \mY_{\Text{dis}} = \mR^{64 \times 128 \times 128 }$, in order not to simulate unrealistically dense sinograms which would supposedly provide a strong advantage for our simultaneous approach.
	    \begin{figure}[hbt!]
	    \centering
	       \subfloat[Exemplary projections for $\varepsilon = 1 \,\%$. From left to right: exact solution vs. reconstructions using the simultaneous Newton approach, Newton-based phase retrieval and CTF-inversion.]{\includegraphics[height = .22\textwidth]{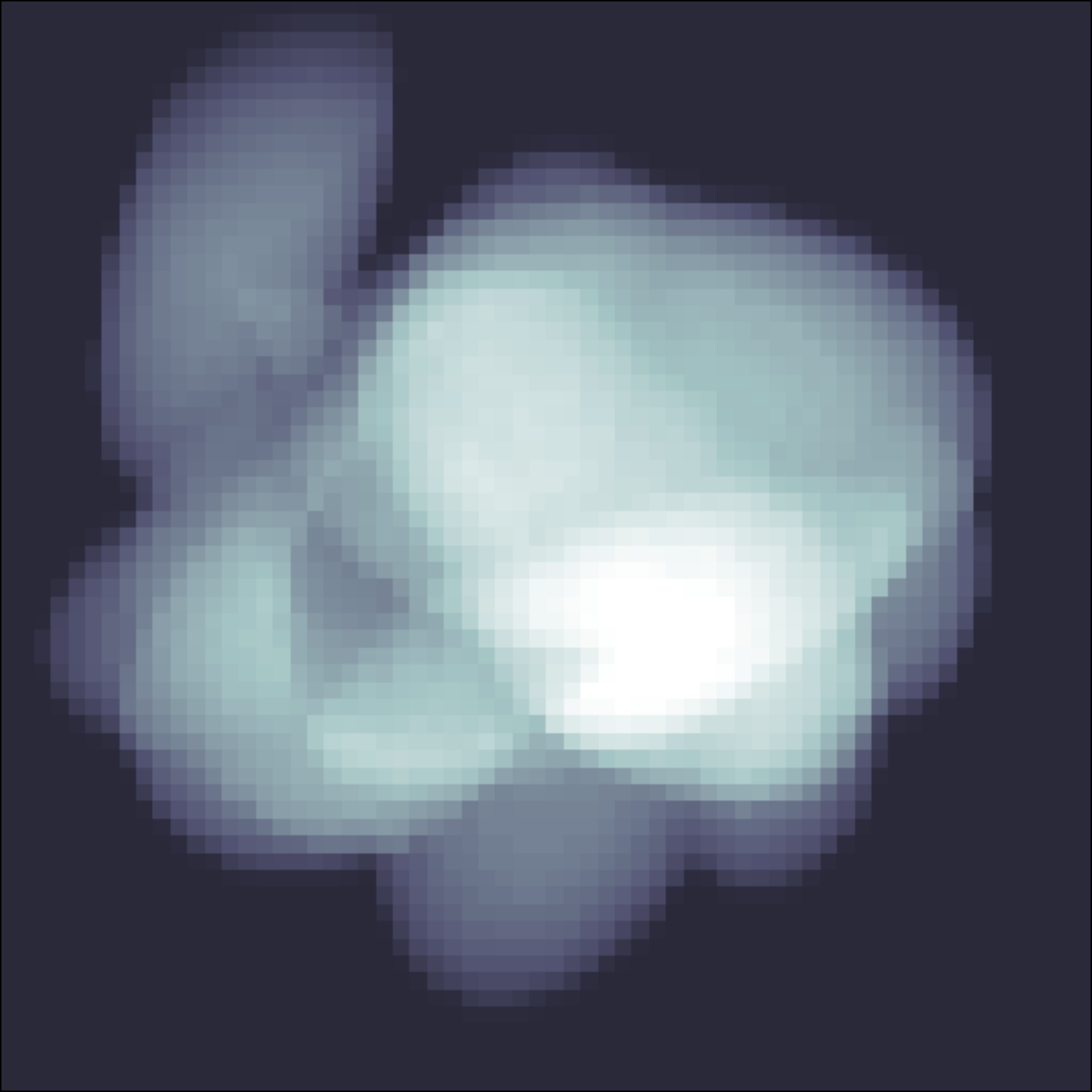}  \includegraphics[height = .22\textwidth]{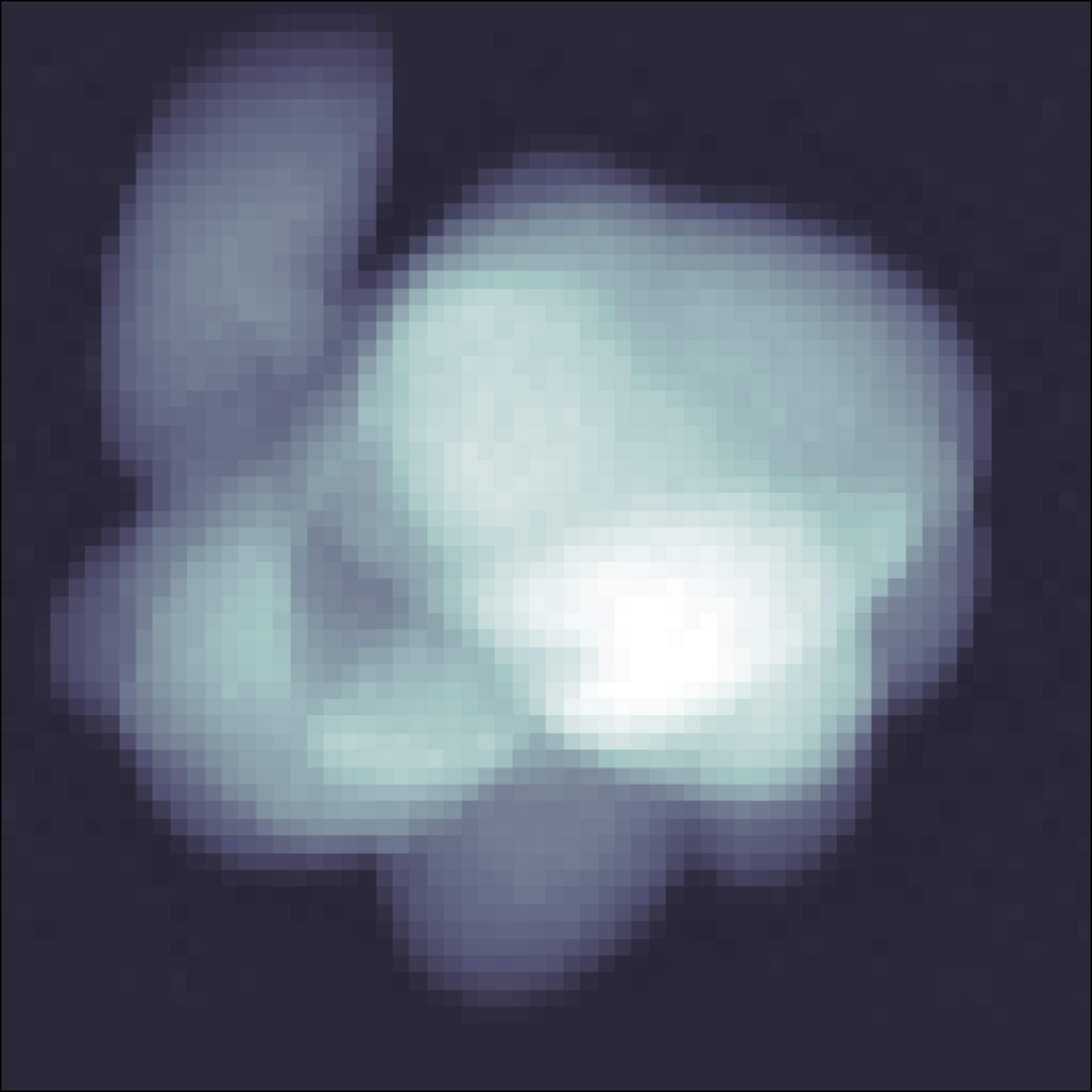} \includegraphics[height = .22\textwidth]{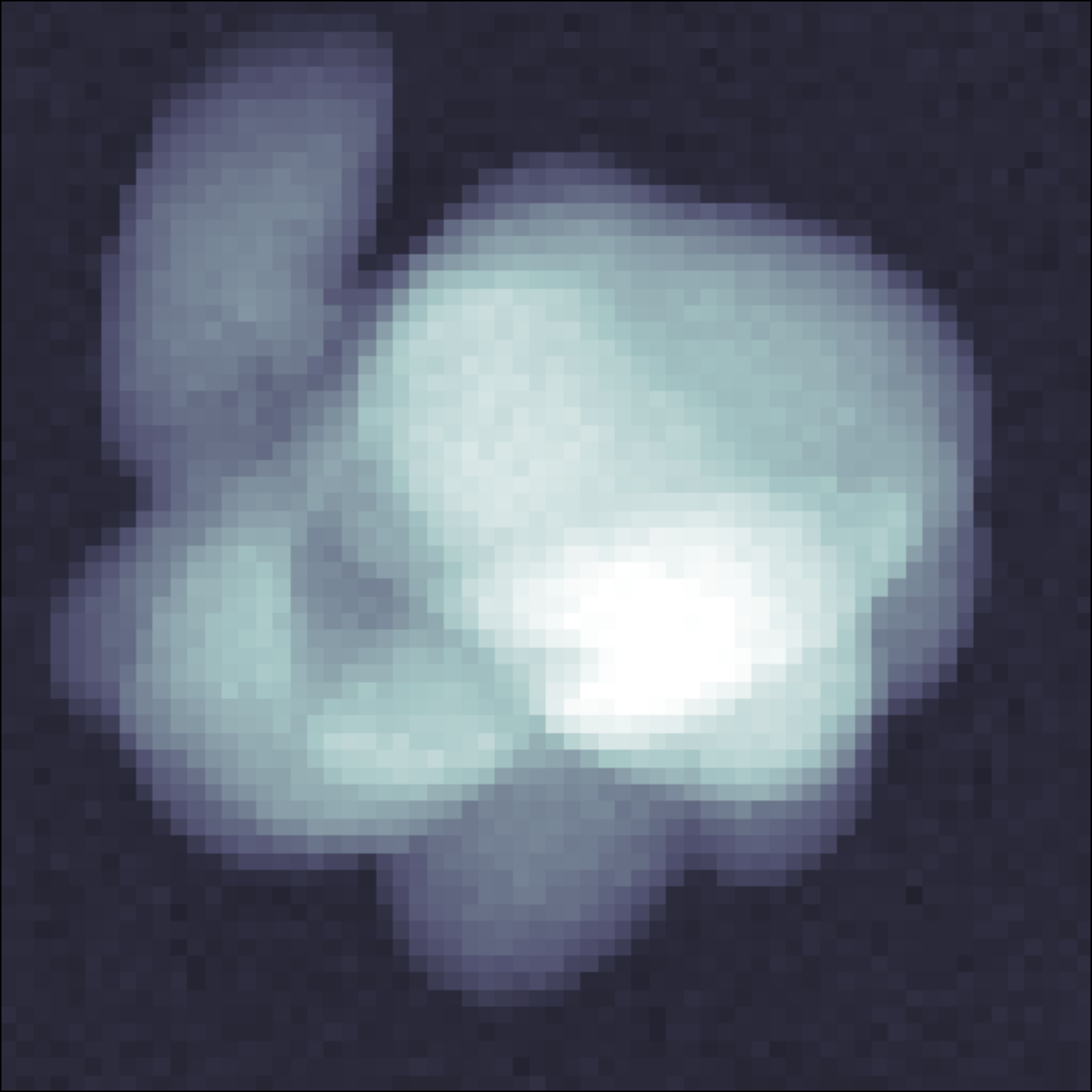} \includegraphics[height = .22\textwidth]{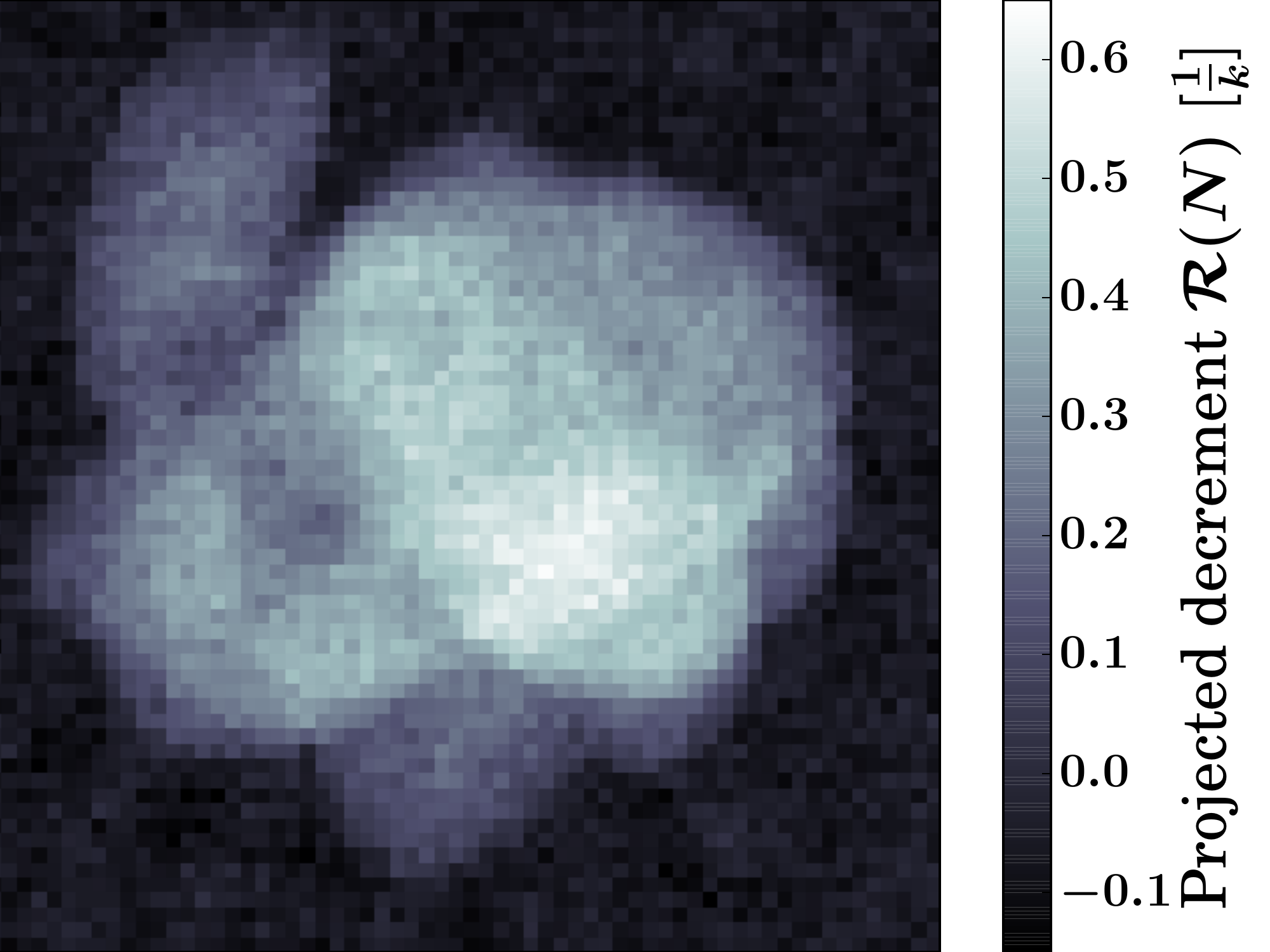} \label{fig:NumResNF-ConsTest-PurePhase-ExplProjs}  } \\
	      \subfloat[Relative errors in the projections $\CR(\bN^\dagger)$\vspace{-.5em}]{\includegraphics[width=.49\textwidth]{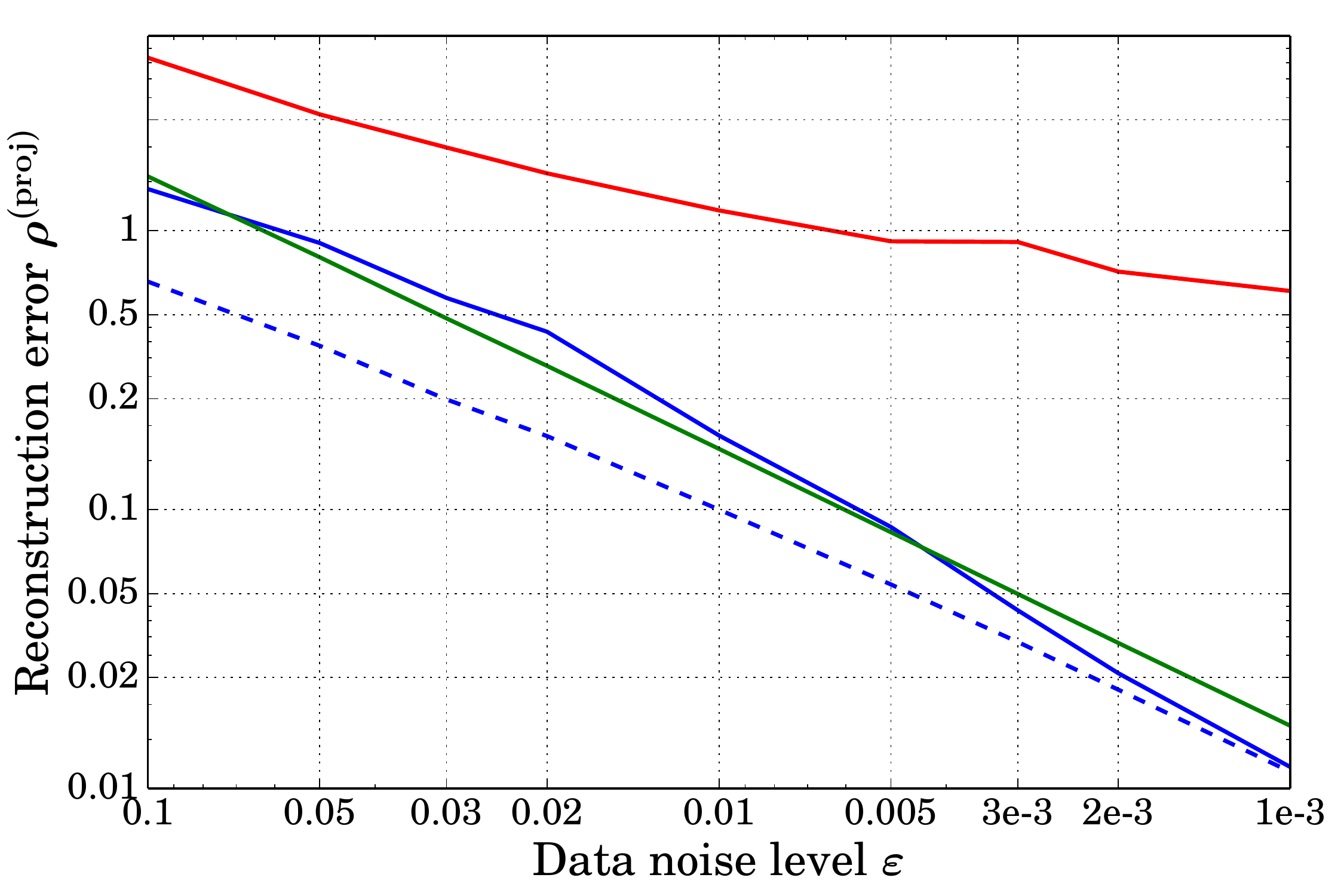} \label{fig:NumResNF-ConsTest-PurePhase-errs-proj} }
	      \hfill
	      \subfloat[Relative errors in the object $\bN^\dagger$\vspace{-.5em}]{\includegraphics[width=.49\textwidth]{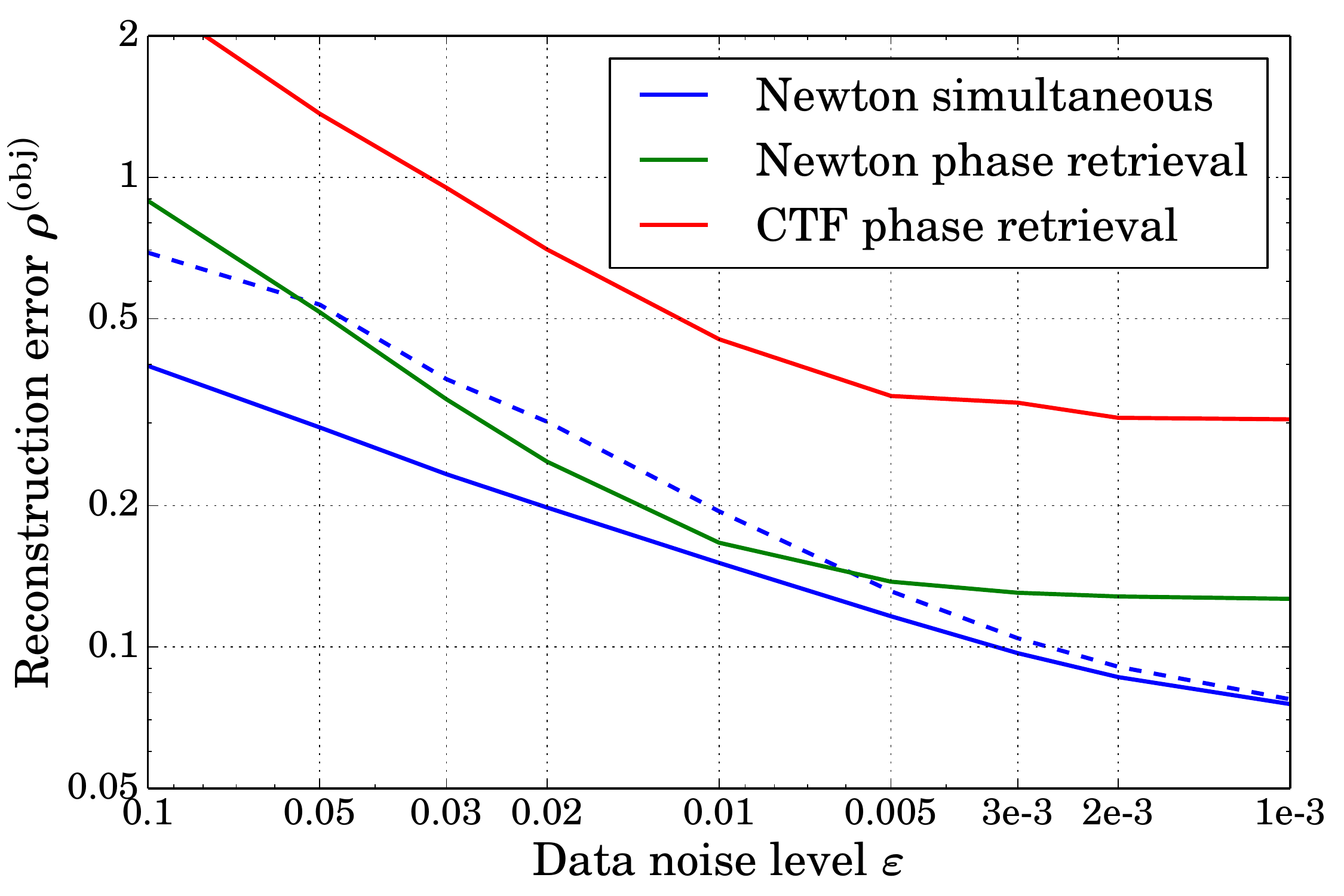}  \label{fig:NumResNF-ConsTest-PurePhase-errs-obj}  }
	      \vspace{-.5em}
	    \caption{ Comparison of different reconstruction methods for pure phase objects $\bN^\dagger \in \mR^{64^3}$ of magnitude $\norm{\bN^\dagger} =1$ at Fresnel number $\NF = 0.01$ reconstructed from simulated intensity data $\bI^{\Textbf{err}} \in \mR^{64 \times 128 \times 128 }$ with different $L^2$-data noise levels $\varepsilon$. ``Newton simultaneous'' represents \algref{alg:PCT} as used in the previous numerical examples. ``Newton phase retrieval'' restricts to the reconstruction of the projections $\CR(\bN^\dagger)$ by solving the 2D phase retrieval problems $\bI^{\Textbf{err}} =  F_{\Text{dis}}^{\Text{phase}}  ( \CR(\bN) )$ in the lateral coordinates independently for all 64 incident angles $\theta \in [0^{\circ};180^\circ)$. ``CTF phase retrieval'': 2D phase reconstructions by inversion of the CTF \eqref{eq:CTF} with optimally chosen truncation around the zeros. The relative $L^2$-errors in (b) and (c) are averaged over three random test objects. The tomographic reconstruction for the latter two methods in (c) is computed by standard filtered backprojection. Stop rule: ``best stop''. Solid vs. dashed blue line: optimal stopping according to the error in $\bN^\dagger$ or $\CR(\bN^\dagger)$, respectively. For details, see text.}
	    \label{fig:NumResNF-ConsTest-PurePhase}
	    \end{figure}
	    
	    For objectivity, we compare the achieved $L^2$-error in the reconstructed projections $\CR(\bN )_{ \Text{rec}}$
	    \begin{equation}
	      \rho^{\Text{(proj)}} := \frac{ \norm{ \CR(\bN )_{ \Text{rec}} -  \CR(\bN^\dagger) }_2 }{ \norm{  \CR(\bN^\dagger) }_2 }
	      \end{equation}
	    because this rules out any bias induced by the choice of the Radon inversion method for the CTF- and Newton-based phase retrieval. In the case of the simultaneous Newton method, $\CR(\bN )_{ \Text{rec}}$ is obtained by applying the forward Radon transform $\CR$ to the reconstructed object $\bN_{k_{\Text{stop}}}$. For completeness, however, we also consider the resulting errors in the object $\rho^{\Text{(obj)}} := \norm{ \bN_{ \Text{rec}} -  \bN^\dagger }_2 / \norm{  \bN^\dagger  }_2$, obtained by applying filtered backprojection to the reconstructed projections $\CR(\bN )_{ \Text{rec}}$ using \textsc{Octave's} \texttt{iradon} routine with default parameters.
	    
	    The reconstruction errors $\rho^{\Text{(proj)}}$ and $\rho^{\Text{(obj)}}$ averaged over the three test objects are depicted in \figref{fig:NumResNF-ConsTest-PurePhase-errs-proj} and \ref{fig:NumResNF-ConsTest-PurePhase-errs-obj}, respectively. Exemplary projections obtained from the different methods are shown in \figref{fig:NumResNF-ConsTest-PurePhase-ExplProjs}. A first aspect to note, both from the visual impression and quantitatively, is that the CTF-based phase reconstruction turns out not to be competitive in the considered setting due to its apparently high sensitivity to noise and its \emph{linearity}. Indeed, the reconstruction errors (red lines in \figref{fig:NumResNF-ConsTest-PurePhase-errs-proj}-c) always exceed those achieved by the competing methods by factors of more than two and quasi stagnate for small noise levels. The latter effect may be attributed to the linearization underlying to the CTF, for which the slight nonlinearity of the considered setting with a chosen object magnitude of $\norm{\bN^\dagger} = 1$ gives rise to \emph{systematic} errors. Notably, however, the non-iterative CTF-based reconstruction is also by far the computationally least expensive method.
	    
	    On the other hand, comparison of the blue and green solid curves in \figref{fig:NumResNF-ConsTest-PurePhase-errs-proj} suggests that Newton-based phase retrieval performs equally well as the simultaneous Newton approach - and even slightly better for $\varepsilon \approx 2\,\%$. Yet, note that the corresponding tomographic reconstructions of the object are still significantly more accurate for the latter according to \figref{fig:NumResNF-ConsTest-PurePhase-errs-obj}. It should be furthermore emphasized that the simultaneous Newton method was stopped at minimum error $\rho^{\Text{(obj)}}$ with respect to the \emph{object} $\bN^\dagger$ and not in the \emph{projections} $\CR(\bN^\dagger)$. If the stop rule is adjusted to optimize the projection error $\rho^{\Text{(proj)}}$, then the achieved agreement with the exact solution $\CR(\bN^\dagger)$ improves compared to independent phase retrieval by factors of $1.5\ldots 2$, see blue dashed line in \figref{fig:NumResNF-ConsTest-PurePhase-errs-proj}. Accordingly, near-field phase retrieval benefits from the pursued simultaneous approach even in the case of pure phase objects reconstructed from quasi-ideal data.
	    
	    The observed superiority  is certainly to be attributed to the exploited mutual consistency of the reconstructed projections. Yet, the improvements according to \figref{fig:NumResNF-ConsTest-PurePhase} might be considered as not significant enough to justify the additional computational effort associated with incorporating the Radon transform in the Newton iterations (see \sref{SS:RadonImplementation}). For this reason, we compare the competing methods in two further settings:
	    \begin{itemize}
	     \item[$\boldsymbol 1$] Truncated, non-oversampled holograms due to a limited field of view
	     \item[$\boldsymbol 2$] General, refracting and absorbing objects $\bN^\dagger = \bdelta^\dagger - \I \bbeta^\dagger$ as studied in \sref{SS:NumResNFSim-GenObjects}
	    \end{itemize}
	    
	    \subsubsection{Truncated Holograms}
	   	    
	    So far we have always considered data that is oversampled by a factor of two in the lateral coordinates, corresponding to many more degrees of freedom in the intensities $\bI^{\Textbf{err}} \in \mY_{\Text{dis}} = \mR^{K_\theta \times 128 \times 128}$ than in the object $\bN^\dagger \in \mX_{\Text{dis}} = \mC^{64}$ to be reconstructed. Thereby, it was furthermore ensured that no fringes containing object information could leave the field of view - an ideal setting that may not always be realized in experimental setups and which is therefore relaxed in the following: other than in the preceding computations, the simulated holograms $\bI^{\Textbf{err}}$ are symmetrically truncated to match exactly the $64 \times 64$-sized projections, i.e.\ we reconstruct from data in the image space $\bI^{\Textbf{err}} \in \mY_{\Text{dis}} = \mR^{64 \times 64 \times 64}$.
	    
	    The results for the different reconstruction methods are plotted in \figref{fig:NumResNF-ConsTest-LowSamp} in analogous manner as in \figref{fig:NumResNF-ConsTest-PurePhase}. Notably, the CTF-reconstruction performs even worse than in the preceding test cases, stagnating at a projection error of $ \rho^{\Text{(proj)}} \approx 38 \, \%$ up to small noise levels according to the red curve in \figref{fig:NumResNF-ConsTest-LowSamp-errs-proj}. The corresponding projections (see rightmost image in \figref{fig:NumResNF-ConsTest-LowSamp-ExplProjs}) show characteristic stripe artifacts and spurious negative halos outside the support caused by the missing high-frequency data fringes in the truncated holograms, which are incorrectly completed by simple padding operations.
	    \begin{figure}[hbt!]
	    \centering
	       \subfloat[Exemplary projections for $\varepsilon = 1 \,\%$. From left to right: exact solution vs. reconstructions using a simultaneous Newton approach, Newton-based phase retrieval and CTF-inversion.]{\includegraphics[height = .22\textwidth]{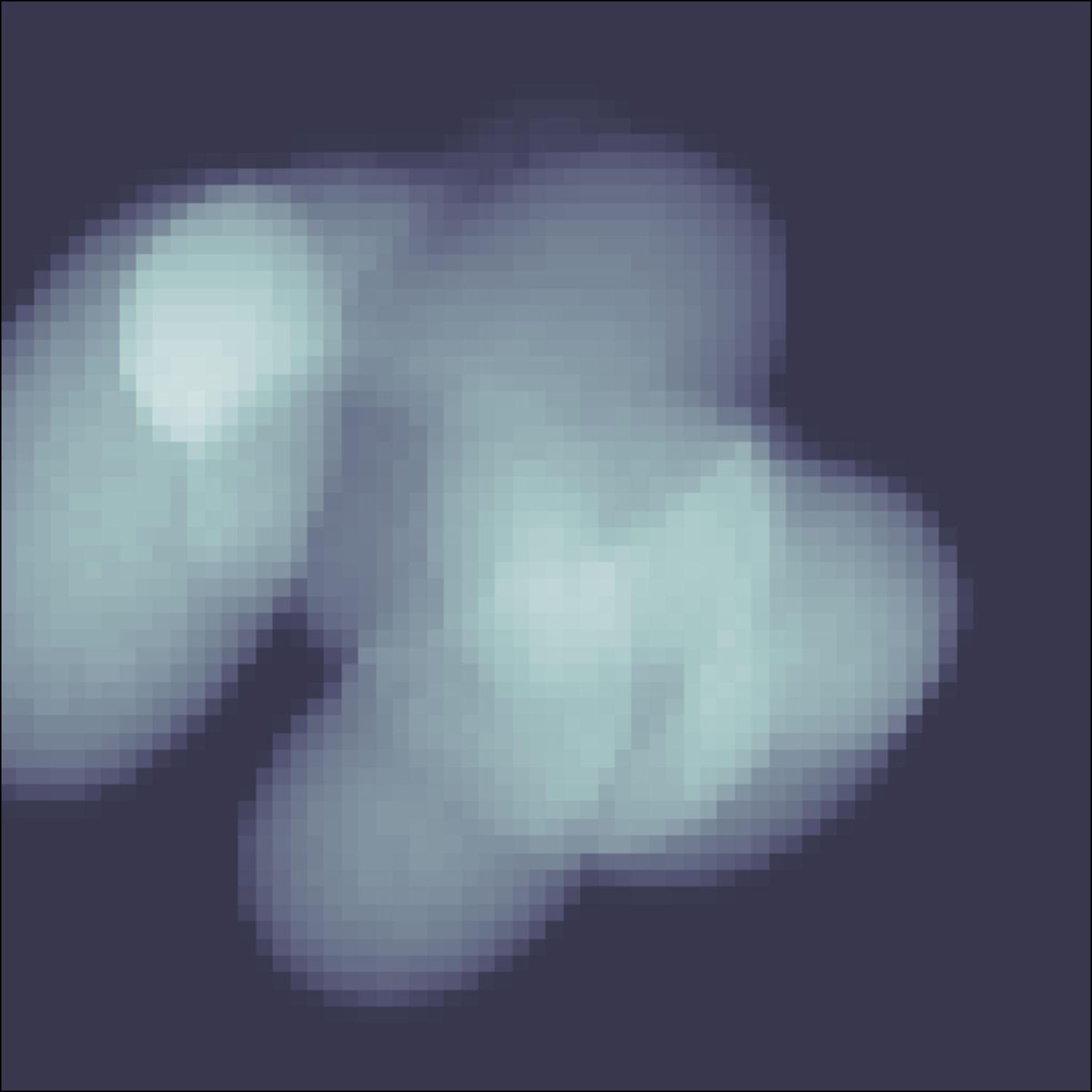}  \includegraphics[height = .22\textwidth]{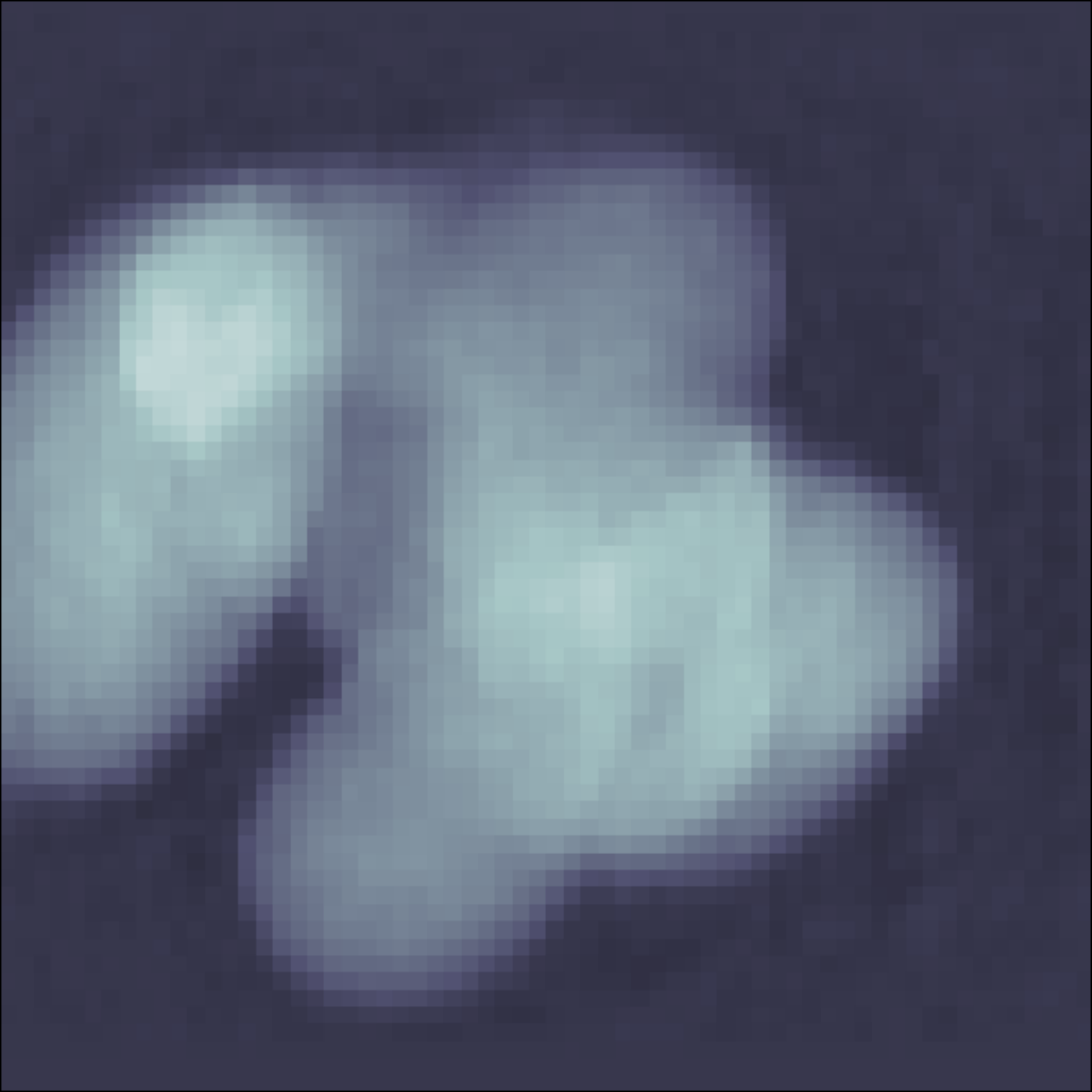} \includegraphics[height = .22\textwidth]{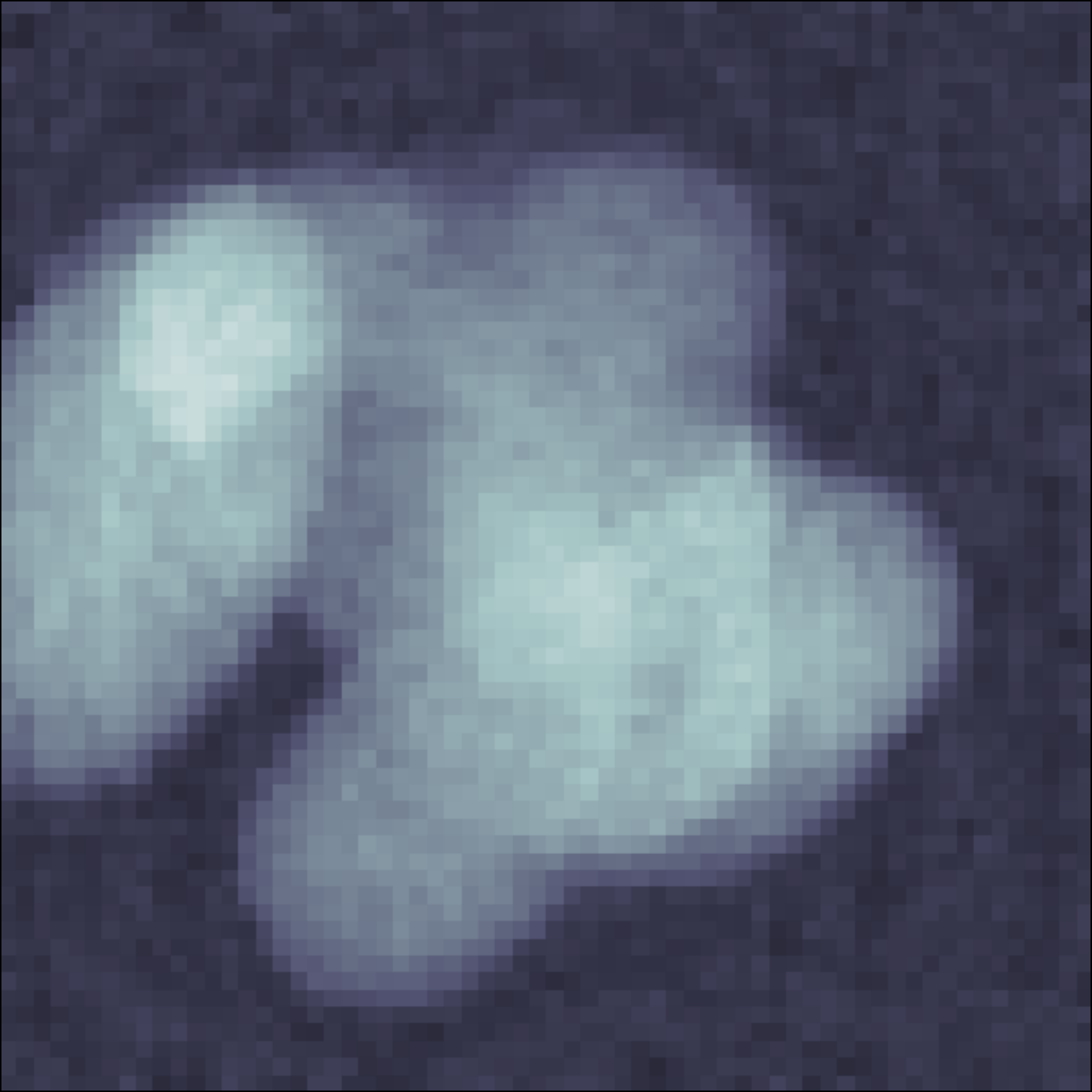} \includegraphics[height = .22\textwidth]{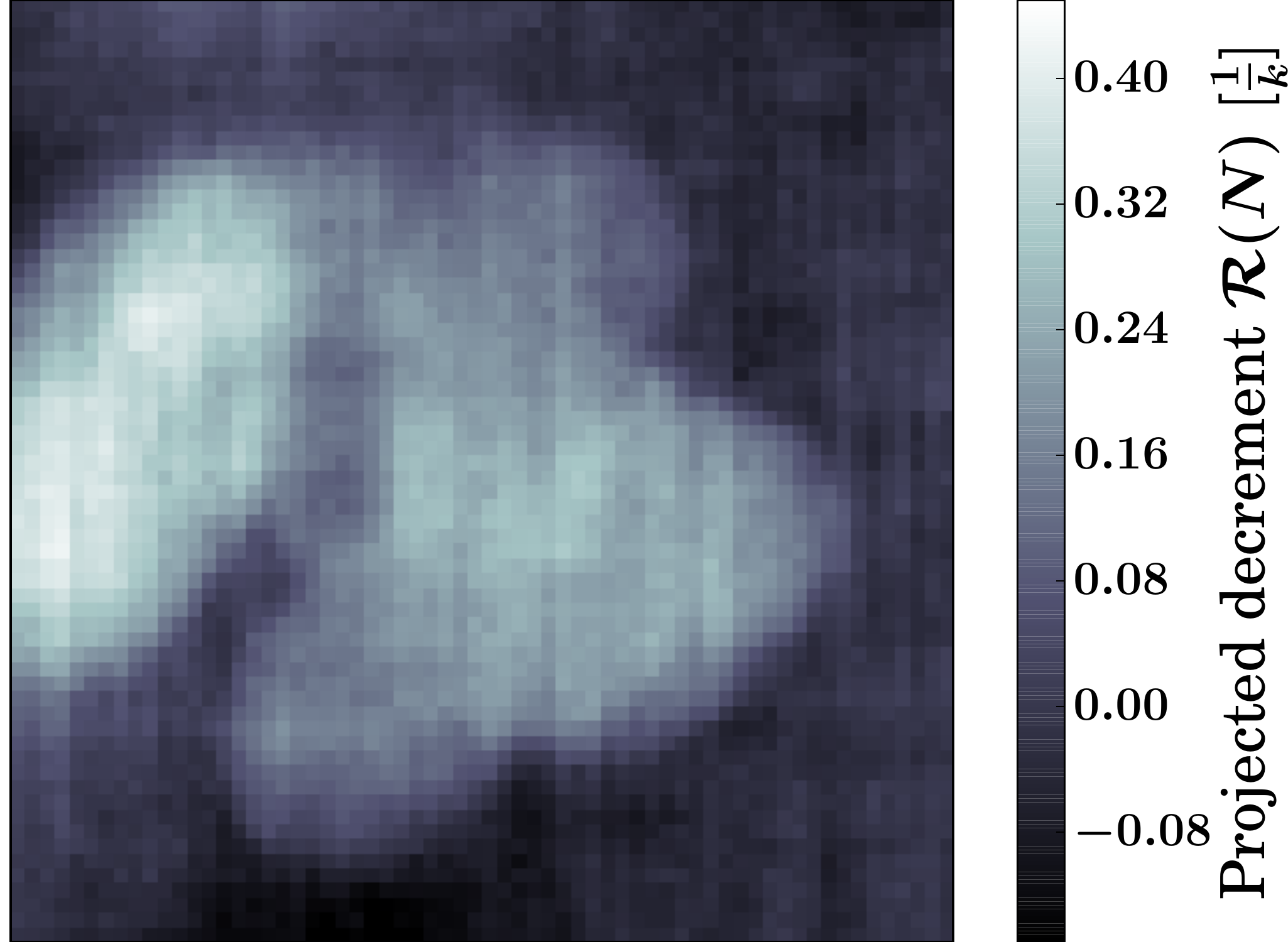} \label{fig:NumResNF-ConsTest-LowSamp-ExplProjs}  } \\
		\vspace{.5em}
	      \subfloat[Relative errors in the projections $\CR(\bN^\dagger)$]{\includegraphics[width=.49\textwidth]{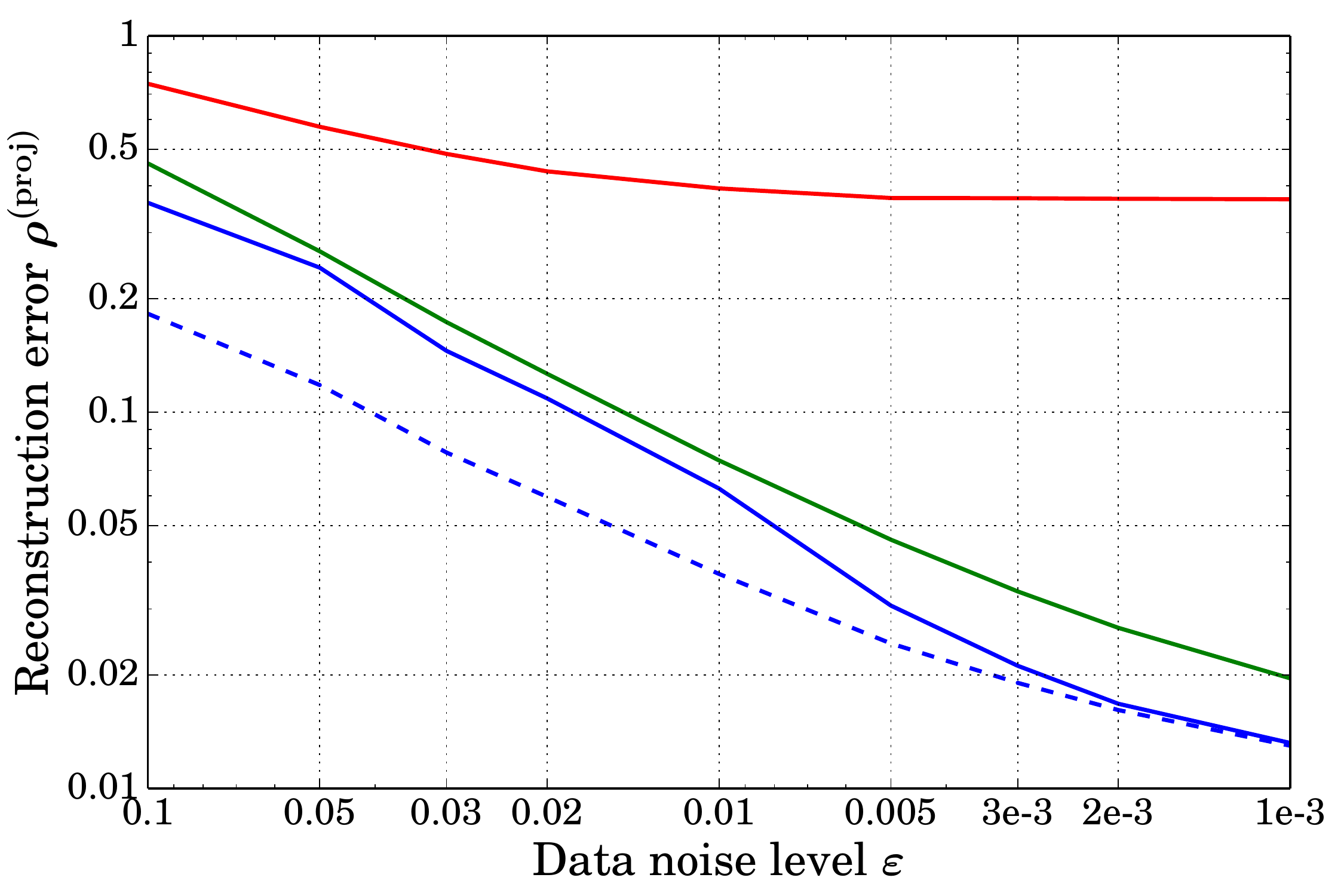} \label{fig:NumResNF-ConsTest-LowSamp-errs-proj} }
	      \hfill
	      \subfloat[Relative errors in the object $\bN^\dagger$]{\includegraphics[width=.49\textwidth]{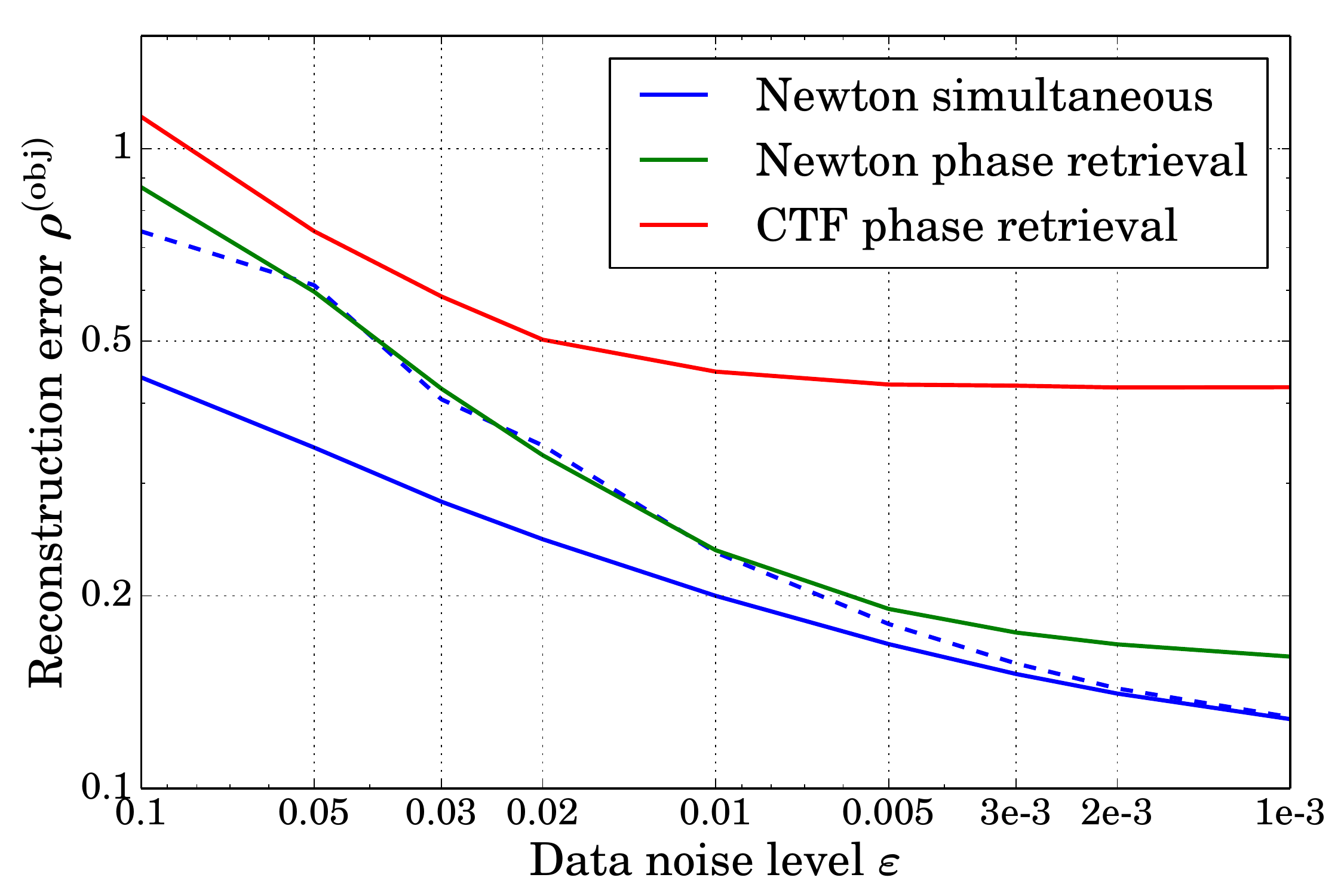}  \label{fig:NumResNF-ConsTest-LowSamp-errs-obj}  }
	    \caption{Comparison of different reconstruction methods for truncated hologram data $\bI^{\Textbf{err}} \in  \mR^{64 \times 64 \times 64}$. Subfigures and simulation parameters analogous to \figref{fig:NumResNF-ConsTest-PurePhase}. The setting corresponds to no oversampling in the intensity data and requires implicit completion of the holograms due to fringes propagating out of the computational field of view. Apparently, this is achieved most stably and accurately by the simultaneous Newton approach, incorporating tomographic correlations in the required inference of the missing data.} \label{fig:NumResNF-ConsTest-LowSamp}
	    \end{figure}
	    
	    On the contrary, the iteratively regularized Newton methods demonstrate their potential by inferring the missing fringes via implicit analytical continuation of the available data in the reconstruction. As can be seen from the blue and green curves in \figref{fig:NumResNF-ConsTest-LowSamp-errs-proj}-c, this works quite well although the achieved reconstruction errors are still significantly higher than for the case of ideal holograms plotted in \figref{fig:NumResNF-ConsTest-PurePhase}. However, note that the simultaneous Newton method outperforms the Newton-based phase retrieval more significantly in the considered truncated hologram setup. This can be seen both from the convergence rates in \figref{fig:NumResNF-ConsTest-LowSamp-errs-proj} and visually from the exemplary projections in \figref{fig:NumResNF-ConsTest-LowSamp-ExplProjs}: the Newton-based phase reconstruction shows high frequency artifacts, corresponding to Fourier modes which are most likely to be represented by fringes leaving the truncated field of view. Apparently, the completion of these partially missing frequencies is much more stable for the simultaneous Newton approach. This could have been anticipated because the implicit inference is performed from a much broader data basis owing to the incorporated consistency correlations between the different holograms.
	    
	    Hence, we find that our regularized Newton method for simultaneous Radon inversion and phase retrieval is particularly beneficial for incomplete holograms, which constitute a frequently encountered imperfection of realistic experimental data. 
	    \newpage
	    
	    \subsubsection{General Objects}

	    We now consider the case of general objects $\bN^\dagger = \bdelta^\dagger - \I \bbeta^\dagger$. It has been seen in \sref{SS:NumResNFSim-GenObjects} that an independent reconstruction of refraction $\bdelta^\dagger$ and absorption $\bbeta^\dagger$ is in general feasible with \algref{alg:PCT}, yet numerically cumbersome, in accordance with the near-field uniqueness result of \sref{SS:PhaseRetrNF}.  Here, we compare the performance in this setting for the competing phase retrieval methods.
	    \begin{figure}[hbt!]
	    \centering
	       \subfloat[Exemplary projected absorption for $\varepsilon = 1 \,\%$. From left to right: exact solution vs. reconstructions using the simultaneous Newton approach, Newton-based phase retrieval and CTF-inversion.]{\includegraphics[height = .22\textwidth]{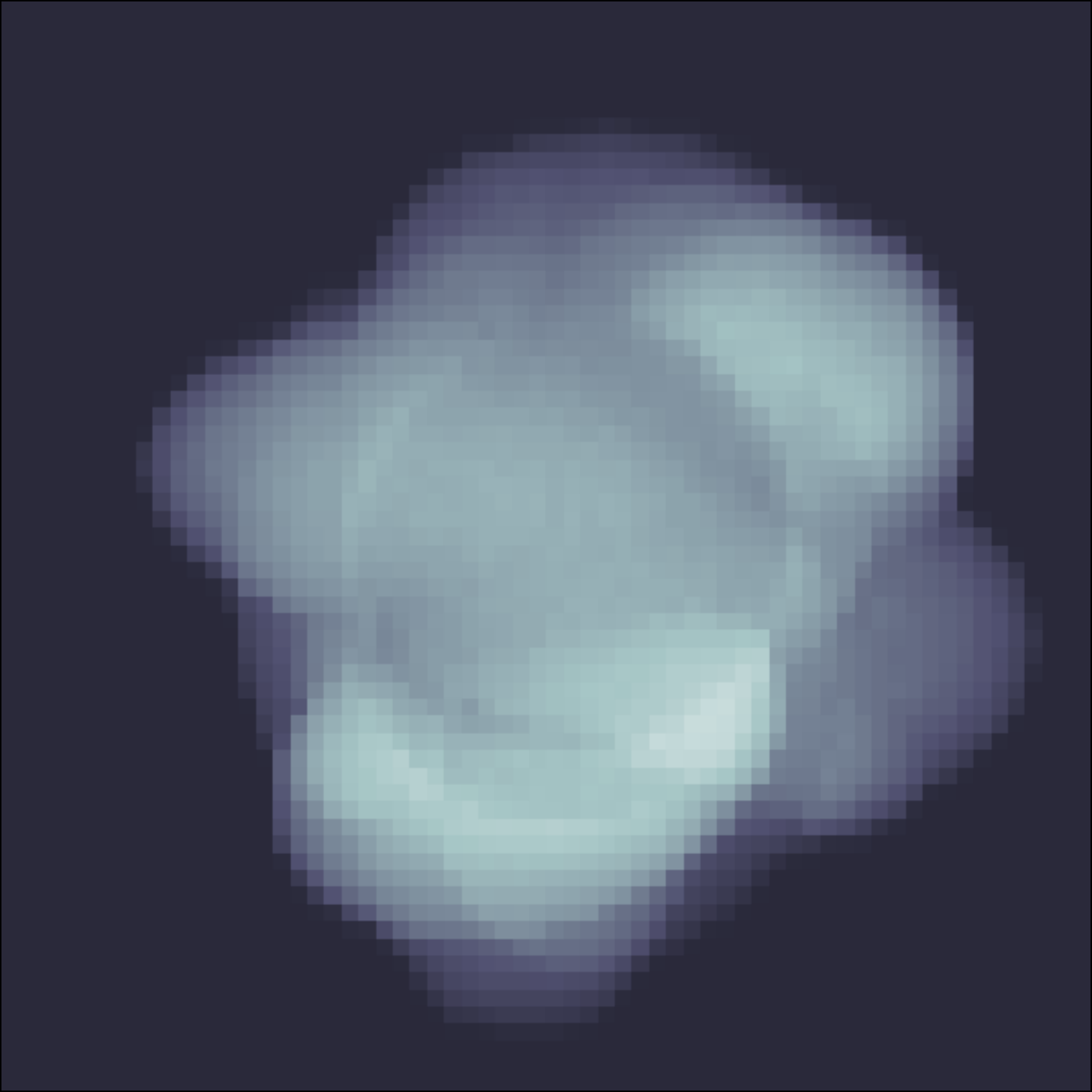}  \includegraphics[height = .22\textwidth]{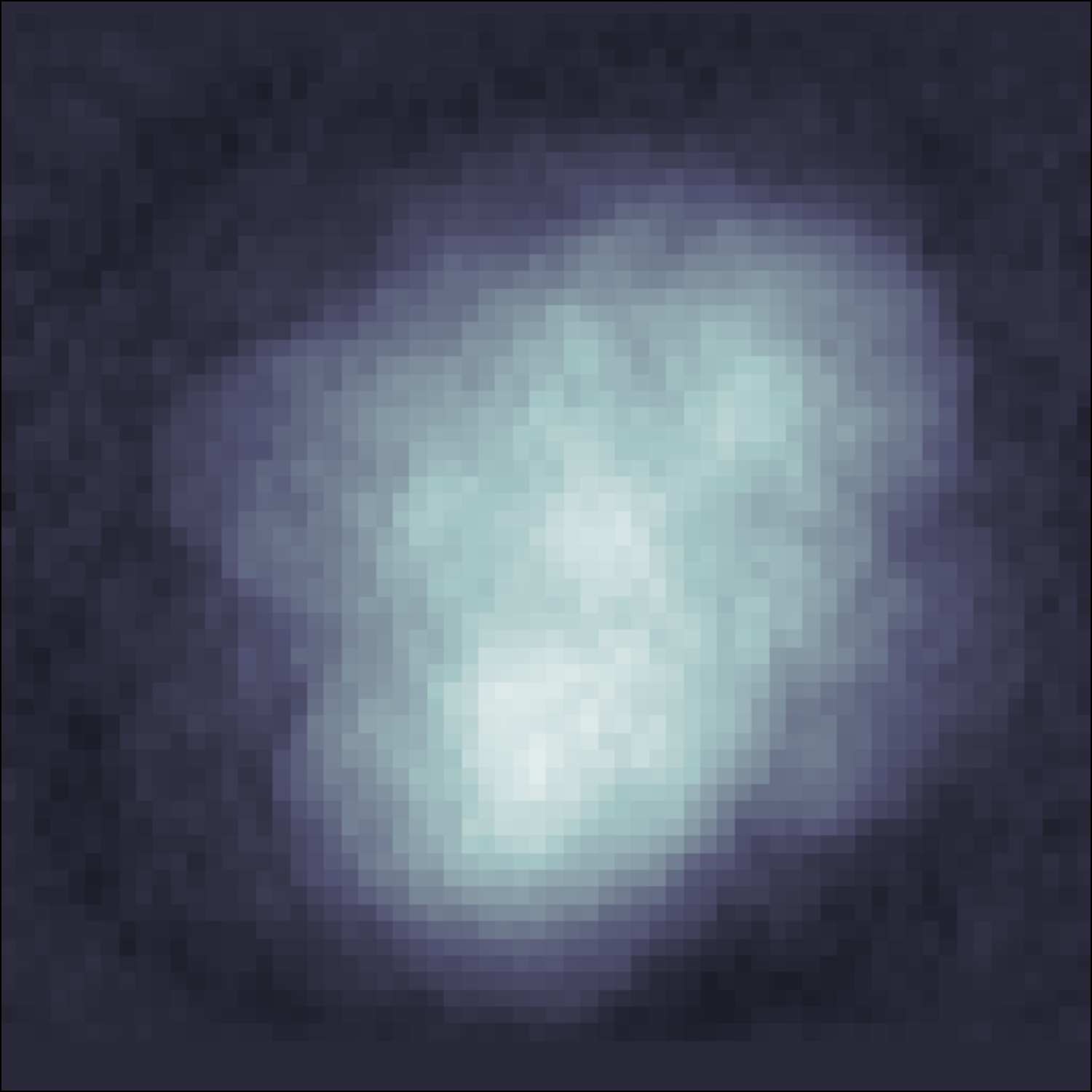} \includegraphics[height = .22\textwidth]{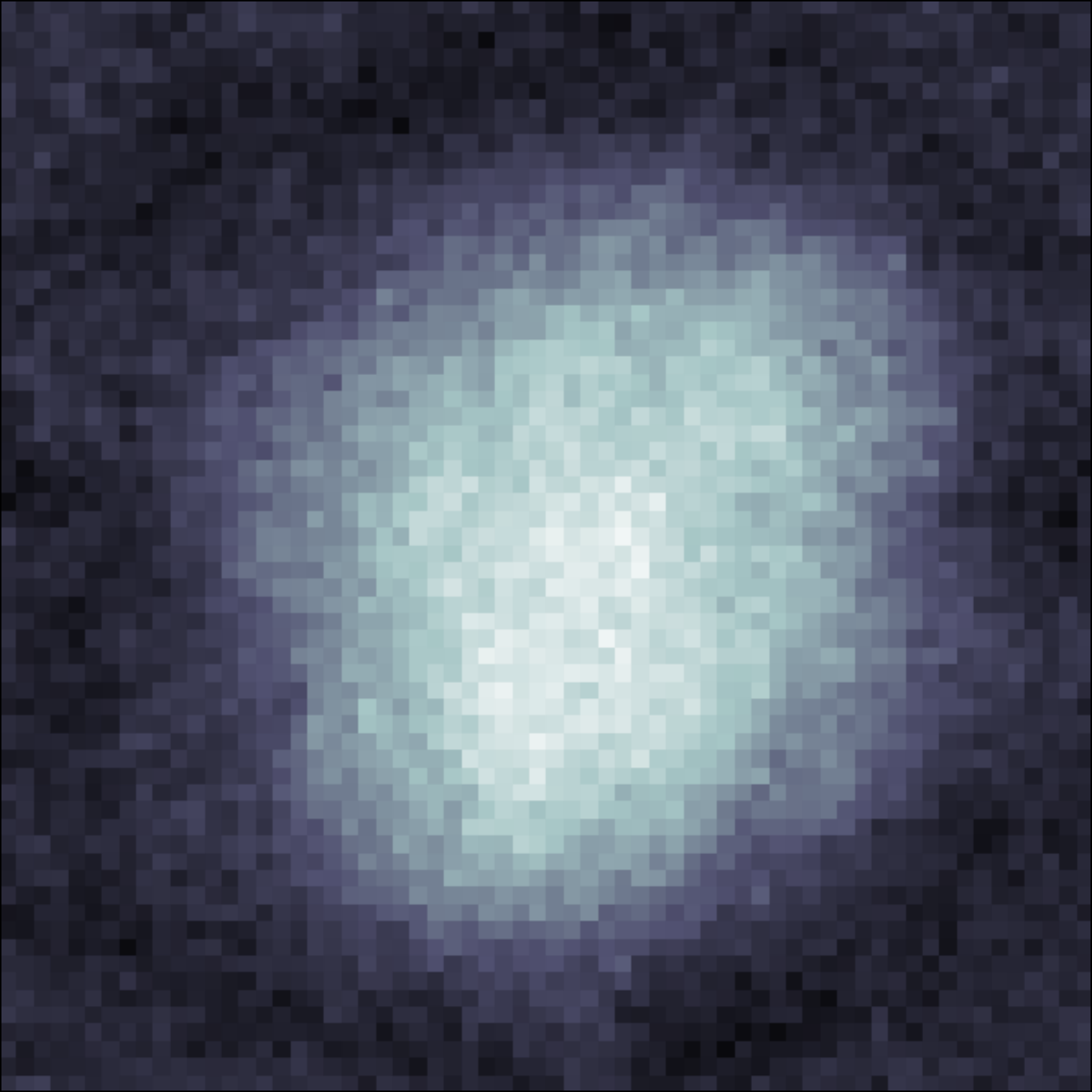} \includegraphics[height = .22\textwidth]{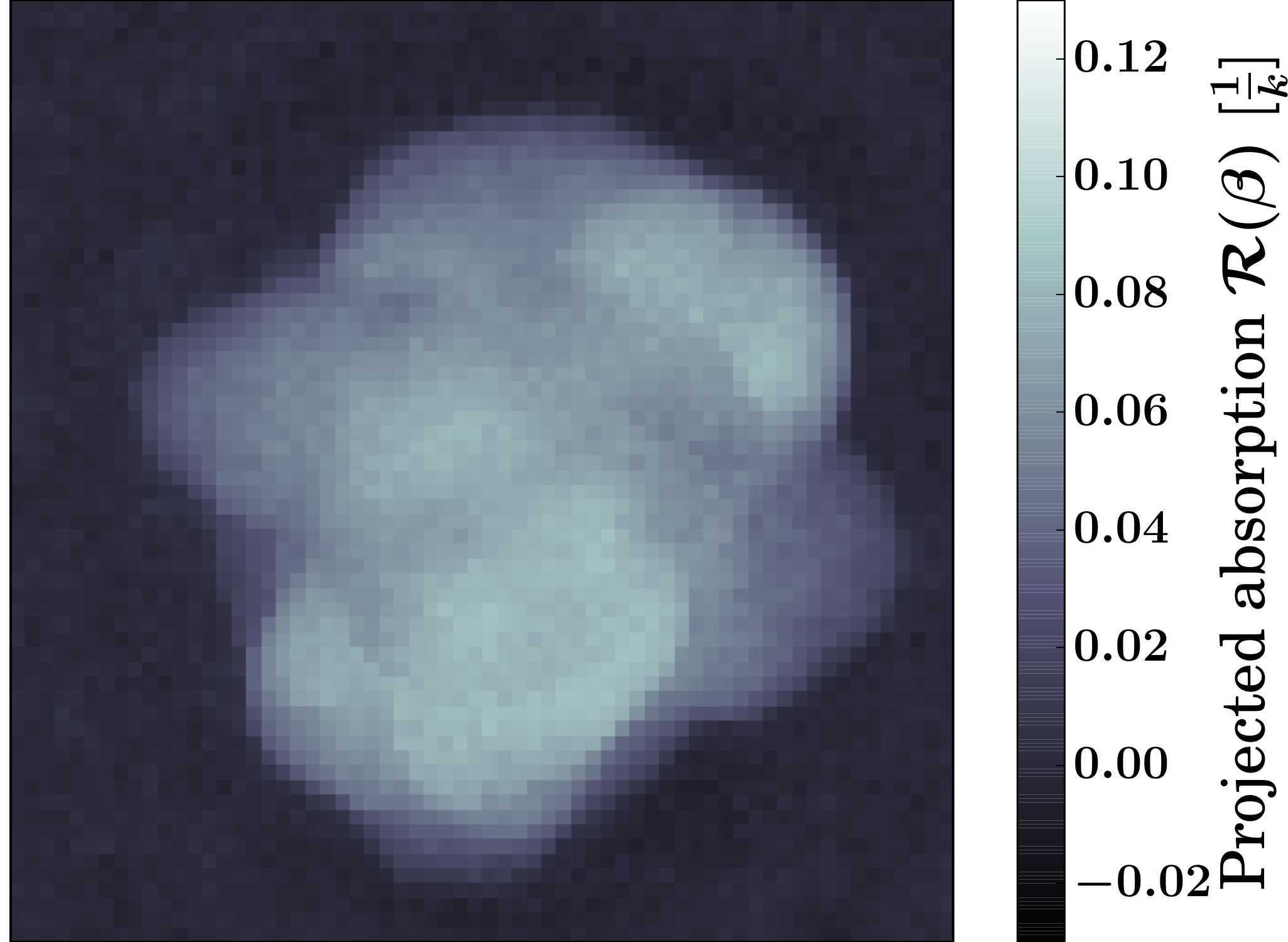} \label{fig:NumResNF-ConsTest-GenObj-ExplProjs}  } \\
	  		\vspace{.5em}
	      \subfloat[Relative errors in the projections $\CR(\bN^\dagger)$]{\includegraphics[width=.49\textwidth]{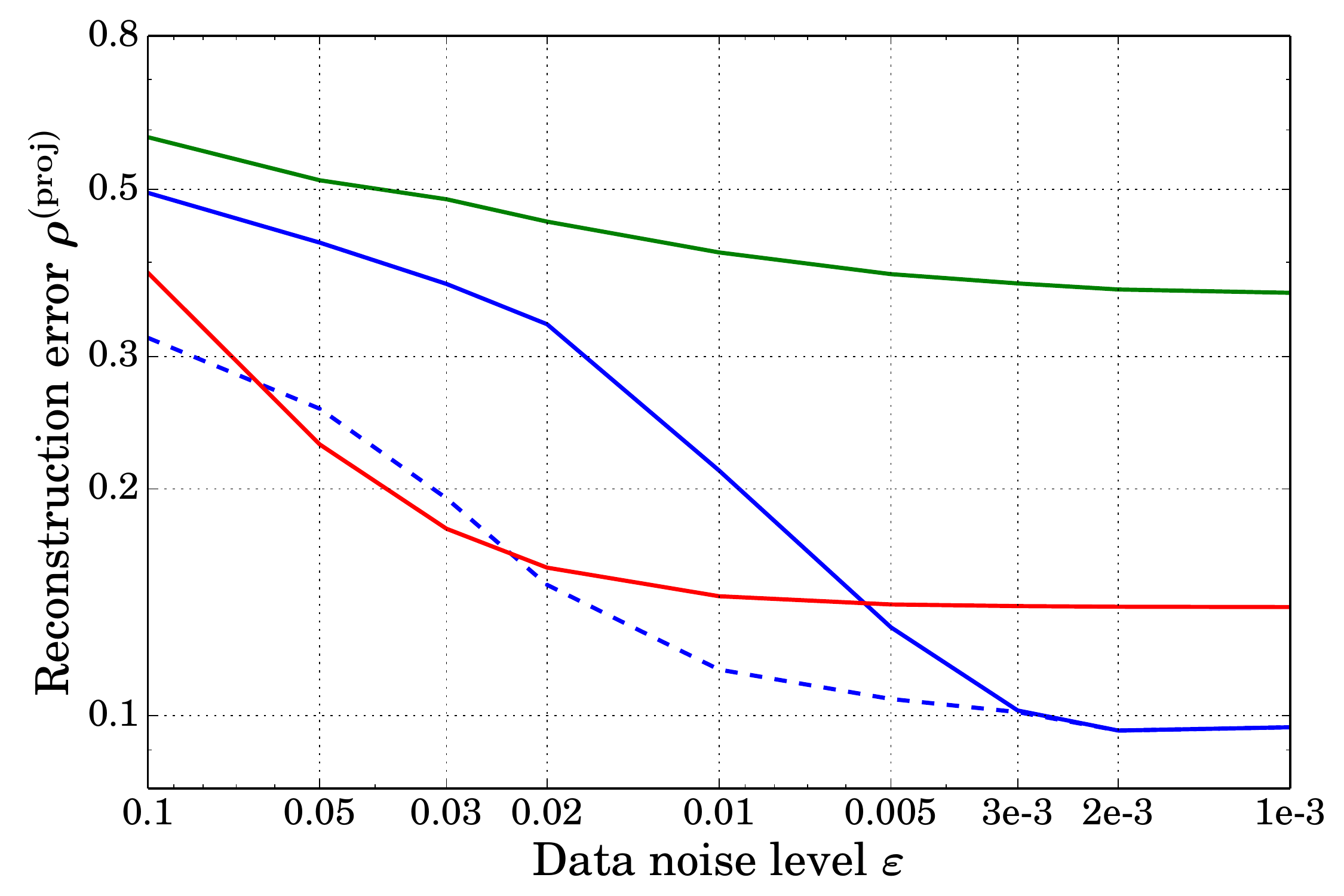} \label{fig:NumResNF-ConsTest-GenObj-errs-proj} }
	      \hfill
	      \subfloat[Relative errors in the object $\bN^\dagger$]{\includegraphics[width=.49\textwidth]{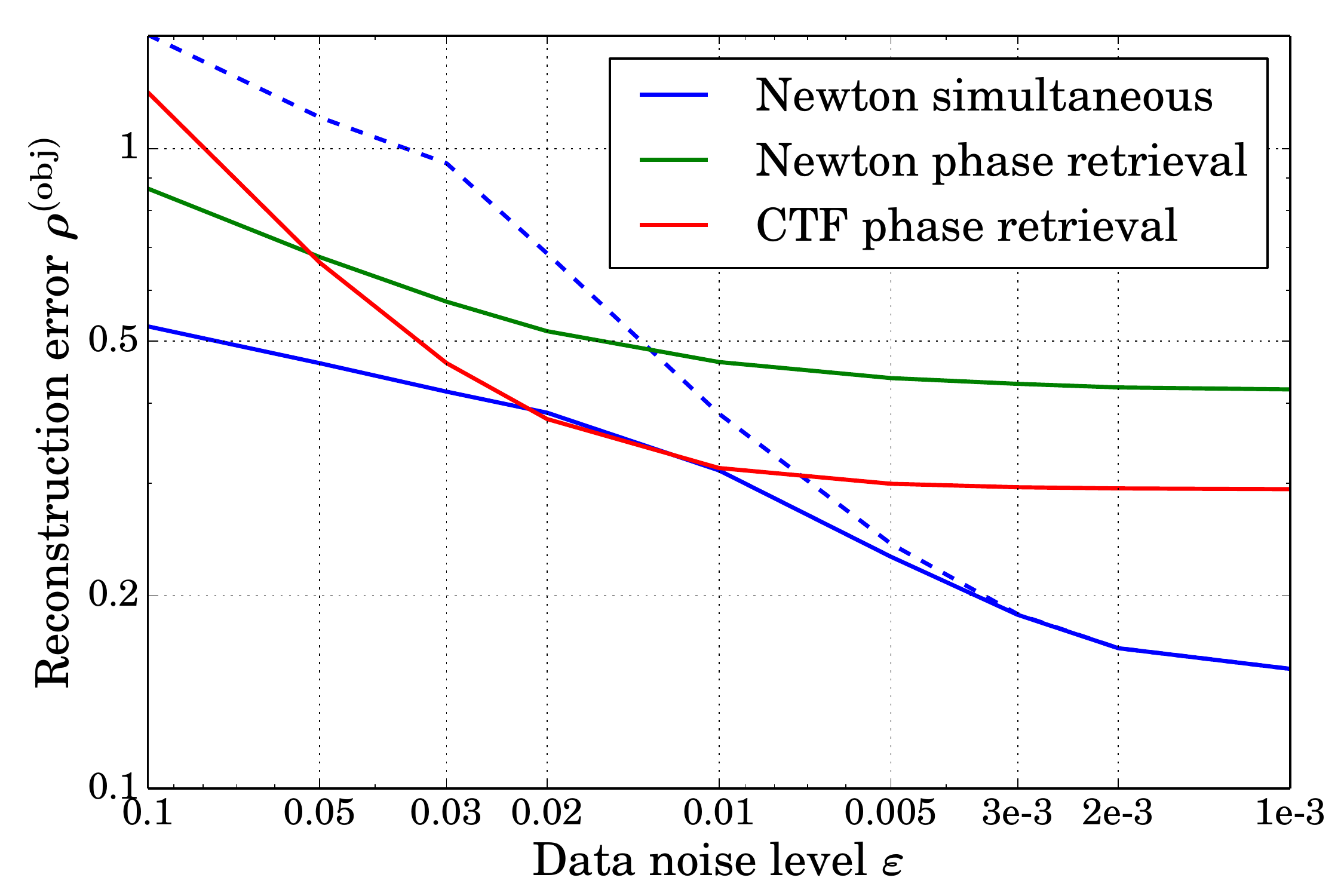}  \label{fig:NumResNF-ConsTest-GenObj-errs-obj}  }
	    \caption{Comparison of different reconstruction methods for general objects $\bN^\dagger = \bdelta^\dagger - \I \bbeta^\dagger$. Subfigures and simulation parameters analogous to results in \figref{fig:NumResNF-ConsTest-PurePhase-errs-proj} where the projections show the recovered \emph{absorption} $\bbeta_{k_{\Text{stop}}}$. In the CTF-reconstruction, a (false) single-material constraint with the approximate $\beta$-$\delta$-ratio $c_{\bbeta /\bdelta } = 0.2$ has to be imposed. In the Newton methods, $\bdelta^\dagger$ and $\bbeta^\dagger$ are reconstructed as independent parameters. Here, the simultaneous approach seems to allow for considerably more accurate reconstructions.} \label{fig:NumResNF-ConsTest-GenObj}
	    \end{figure}
	    
	   To this end, we consider random ellipsoid objects $\bN^\dagger = \bdelta^\dagger - \I \bbeta^\dagger$ with $\norm{\bN^\dagger} = 1$ and a large $\beta$-$\delta$-ratio $c_{\bbeta /\bdelta } = 0.2$ in order to permit accurate recovery of the absorption $\bbeta^\dagger$. The CTF-based method may not solve for the two components independently. Hence, we have to assume a false single-material constraint fixing $\bbeta^\dagger = c_{\bbeta /\bdelta } \bdelta^\dagger$ in these reconstructions. For the iterative Newton methods, the regularization in $\bdelta, \bbeta$ is adjusted to match the average ratio $c_{\bbeta /\bdelta }$ as in \sref{SS:NumResNFSim-GenObjects}. In this test case, we reconstruct once more from the complete laterally oversampled holograms $\bI^{\Textbf{err}} \in \mY_{\Text{dis}} = \mR^{64 \times 128 \times 128}$. The remaining setup parameters are retained.
	   
	   Results for the different methods based on three different  objects are visualized in \figref{fig:NumResNF-ConsTest-GenObj}. The plots suggest that the CTF-reconstruction performs surprisingly well in this setting. This is due to two aspects: for once, only the parameter $\bdelta^\dagger$ has to be reconstructed in this case owing to the fixed coupling $\bbeta^\dagger = c_{\bbeta /\bdelta } \bdelta^\dagger$. On the other hand, this coupling also stabilizes the CTF-inversion as poor phase contrast at small Fourier frequencies (compare \figref{fig:CTF}) is balanced by the absorptive part. However, note that the imposed single-material constraint gives rise to systematic errors for the considered general objects which cause the stagnation of the convergence for small noise levels $\varepsilon$ as observed for the red curves in \figref{fig:NumResNF-ConsTest-GenObj-errs-proj}-c.
	   
	   More surprisingly, stagnation at large reconstruction errors is also found in the case of the (non-simultaneous) Newton-based phase retrieval (green curves). Indeed, the independent recovery of $\bdelta^\dagger$ and $\bbeta^\dagger$ widely seems to fail for this method as confirmed by the noisy and hardly defined projection shown in \figref{fig:NumResNF-ConsTest-GenObj-ExplProjs}. On the contrary, the simultaneous Newton approach converges up to errors of $\approx 10 \, \% $ and  $ 15 \, \%$ in the projections $\CR(\bN^\dagger)$ and the object $\bN^\dagger$, respectively. This is well below the magnitude of the absorption  $\bbeta^\dagger\sim 0.2 \bdelta^\dagger$ so that the latter must be recovered at least roughly. This is confirmed by the exemplary projection in \figref{fig:NumResNF-ConsTest-GenObj-ExplProjs}, showing characteristic structures of the exact object apart from the low-frequency halo that has already observed in \sref{SS:NumResNFSim-GenObjects}. Notably, the simultaneous method even outperforms the very stable CTF reconstruction at low noise levels although the latter is somewhat close to its ideal setting with a moderately weak object $\norm{\bN^\dagger}=1$ of only slightly varying $\delta$-$\beta$-ratio (see \eqref{eq:NumResNF-ObjGen}).
	   
	   Although the results obtained by simultaneous Radon inversion and phase retrieval via \algref{alg:PCT} are far from perfect, the implicit incorporation of consistency thus seems to provide a promising ingredient for the reconstruction of general objects. We summarize further findings of this section:
	\begin{res}[Simultaneous Radon Inversion and Phase Retrieval] \label{res:Simultaneous}
	  	The simultaneous tomographic- and phase reconstruction implemented in \algref{alg:PCT} yields quantitative improvements compared to (non-tomographic) near-field phase retrieval with regularized Newton methods and CTF-based techniques by incorporating consistency of the diffraction patterns. The advantage ist most pronounced for laterally truncated, i.e.\ incomplete holograms and for general, refracting and absorbing samples to be reconstructed.
	\end{res}

	 \end{subsection}

    \end{section}

            \begin{section}{Near-Field Tomography from  Experimental Data}  \label{S:NumResNFReal}
            
            As a final numerical example of this chapter, we study near-field tomography via \algref{alg:PCT} from an experimentally recorded data set. The principal aim is to demonstrate applicability of our regularized Newton-type approach to realistic measurements containing unknown statistical and systematic errors.

            \subsection{Reconstruction Setup} \label{SS:NumResNFReal-Setup}
            
            The considered data set has been measured at the P10 beamline of the third generation synchrotron light source PETRAIII at the DESY facilities (Deutsches Elektron Synchrotron, Hamburg) using a GINIX setup (see \cite{Kalbfleisch2011GINIX}, \cite[sec. 4.3]{BartelsDiss}). The specimen is a colloidal crystal of  $415\,\unit{nm}$ diameter beads of the polymer polystyrene ($(\text{C}_8 \text{H}_8)_m$) on a silicon nitride (SiN) membrane of which holograms have been recorded under $K_\theta = 249$ incident angles $\theta \in [0^\circ; 172.85^\circ]$ at an X-ray wavelength $\lambda = \unit[0.157]{nm}$ with an exposure time of one second each. The missing wedge of $\approx 6^\circ$ is due to the experimental constraint that the SiN-membrane needs to be penetrated by the incident radiation at a sufficiently sharp angle in order to avoid systematic errors due to reflections on or within the planar layer. The diffraction patterns are resolved by $K_x = K_y = 1024$ equidistant quadratic detector pixels (no astigmatism).
            
             In the experimental setup, the specimen is illuminated by a cone beam as sketched in \figref{fig:SetupExp} being located $d_1 = \unit[22.8]{mm}$ behind the focal point. The distance between the focus and the detector is $d_2 = \unit[5.085]{m}$. By the Fresnel scaling theorem (see for instance \cite[Appendix B] {PaganinXRay}, \cite{Pogany1997noninterferometric}), this scattering setup may be approximated by an effective parallel-beam geometry which is characterized by the \emph{magnification} $M= d_2 / d_1$, the effective detector distance $d_{\text{eff}} = (d_2 - d_1 ) / M = \unit[22.7]{mm}$  and -pixel size $\Delta x _{\text{eff}} = \Delta x /  M = \unit[29.3]{nm}$. From these parameters, the numerically relevant Fresnel number $\NF = d_{\text{eff}}^2 / ( \lambda \Delta x _{\text{eff}} )$ is determined. Note that the \emph{real} detector pixel size is $\Delta x = \unit[6.54]{\mu m}$ in the given example and thus much larger than the nanoscale structures of the specimen in question. Hence, without the natural magnification associated with the cone-beam setup, the sample could not be resolved.
 	        \begin{figure}[hbt!]
	  \centering
	  \subfloat{\includegraphics[height=.44\textwidth]{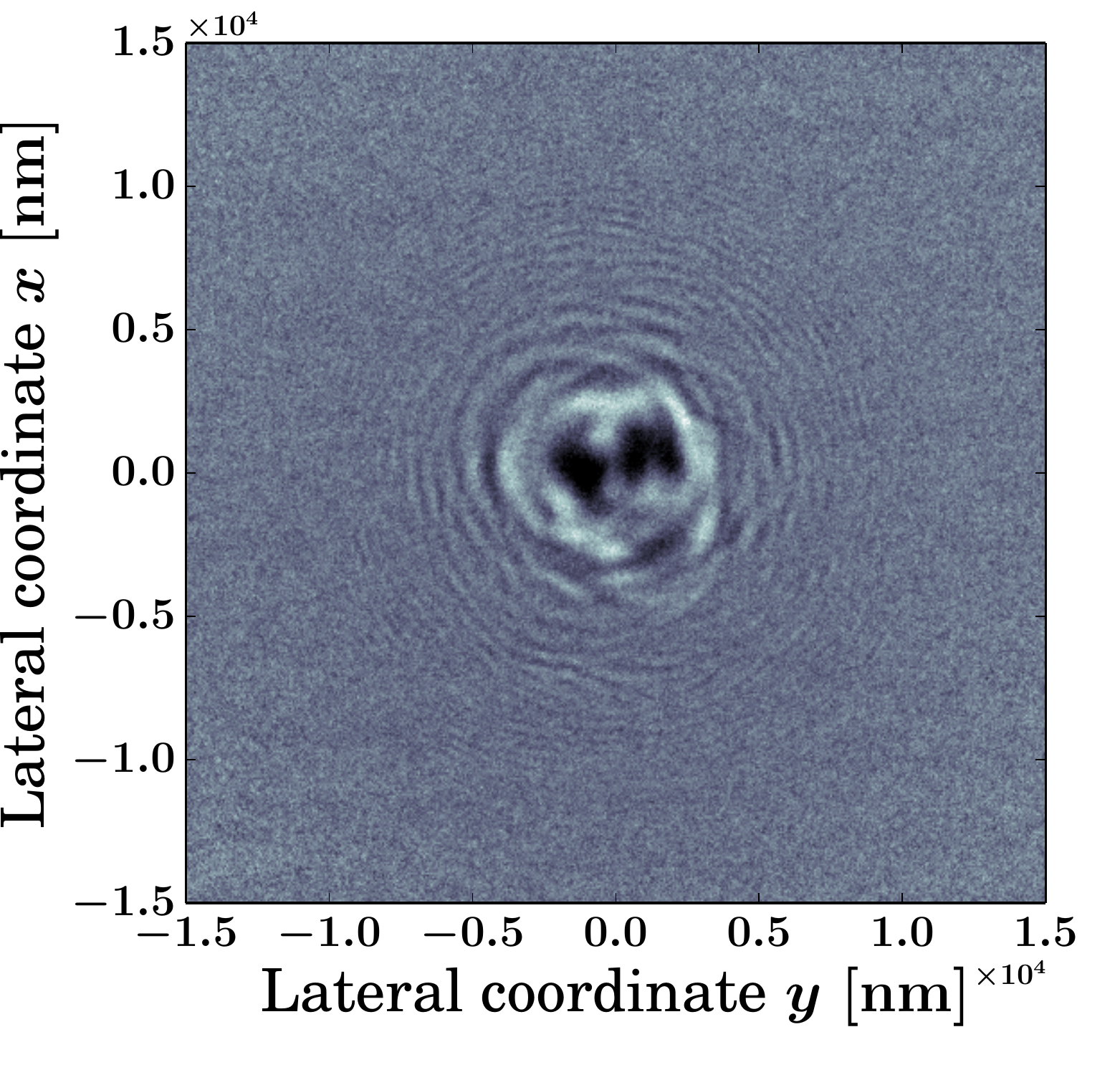} \label{fig:NumResExp-Data-a} }
	  \hfill
	  \subfloat{\includegraphics[height=.44\textwidth]{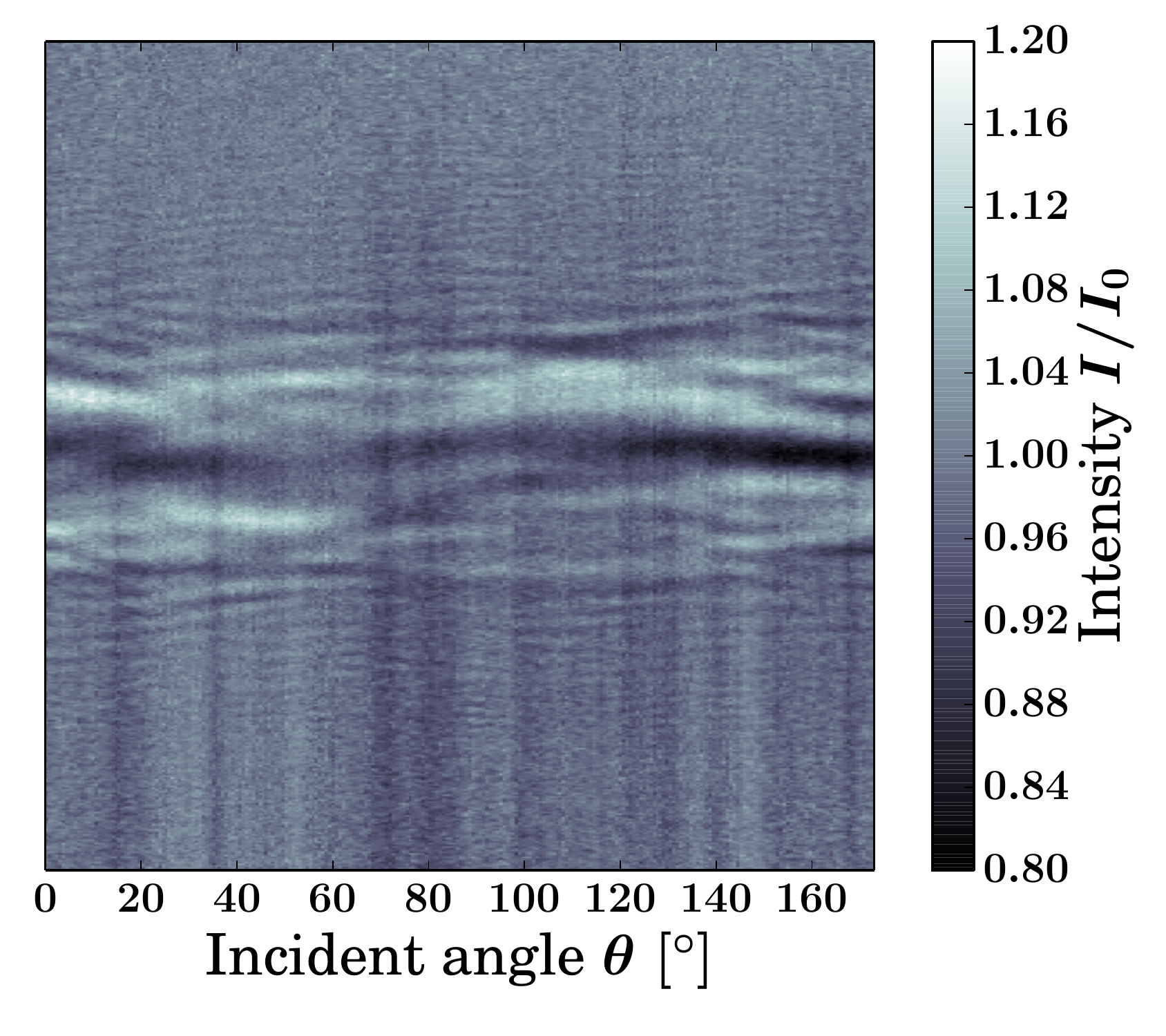} \label{fig:NumResExp-Data-b} }
	  \caption{Experimentally observed near-field tomography data $\bI^{\Textbf{err}}$ for a colloidal crystal of $415\,\unit{nm}$ polystyrene beads (flat-field corrected). Measured on GINIX setup \cite{Kalbfleisch2011GINIX} at P10, PETRAIII, DESY. Fresnel number: $\NF = 2.41\cdot 10^{-4}$, Effective pixel size: $\Delta x _{\text{eff}}   = \unit[29.3]{nm}$. Left: lateral hologram for an incident angle $\theta = 0^\circ$, right: holographic ``sinogram'' slice at $y = 0$. Note the missing wedge of approximately $6^\circ$ in the recorded incident angles and residual variations in the background intensity due to imperfect flat-field correction. \label{fig:NumResExp-Data}}
	  \end{figure}
 	    
 	    The measured holograms have been aligned to correct translational shifts due to vibrations or drifts of the specimen during the measurement. Moreover, the studied data set has been \emph{flat-field-corrected} in preprocessing by division of the diffraction patterns by the \emph{empty-beam image}, i.e.\ the intensities recorded without a scattering object in the beam line. Thereby, the resulting data is made to approximately equal the hypothetical holograms recorded under ideal plane wave illumination, corresponding to a constant probe $P = 1$. This approximation turns out to be accurate if the real probe beam varies on larger lengthscales than the specimen, see  \cite{Hagemann2014EmptyBeam} for details. The preprocessed data set $\bI^{\Textbf{err}} \in \mY_{\Text{dis}} = \mR^{249 \times 1024 \times 1024}$ used for the reconstruction is visualized in \figref{fig:NumResExp-Data}. 	    
	  
	  Polystyrene gives rise to negligible absorption $\beta  \sim 10^{-3} \delta$ for the considered incident hard X-rays of energy $\unit[7.9]{keV}$ according to \cite{Cloetens1999}. Hence, a pure phase object constraint is assumed in the reconstruction. Moreover, it can be inferred from the hologram data in \figref{fig:NumResExp-Data} that refracting matter is present only near the center of the field of view. We exploit this by assuming the sample to be located in a central cube of $256^3$ voxels, i.e.\ we choose a discrete object domain $\bN = \bdelta \in \mX_{\Text{dis}} = \mR^{256^3}$, corresponding to a rough support constraint. The projections $\CR( \bN )$ are symmetrically zero-padded as described in \sref{SS:ZeroPad} to match the $1024\times 1024$ lateral resolution of the intensity data $\bI^{\Textbf{err}}$.
	  
 	 \begin{table}[htb!]
	  \centering
	  \begin{tabular}{cccccccc} 
	    \toprule
		$\mX_{\Text{dis}}$ &  $\mY_{\Text{dis}}$  & $\NF\cdot 10^{4}$ & Reg. term & $\alpha_0$ & $k_{\Text{stop}}$ & Incident $\theta $  & Constraints \\
	    \midrule
		$\mR^{256^3}$ & $\mR^{249 \cdot 1024^2}$ & $2.41$ & $H^{0.5}$ & $ 10^{10}$ & $ 12$   & $[0^{\circ}; 172.85^\circ]$ & \scriptsize{pure phase obj.} \\ 
	    \bottomrule
	  \end{tabular}
	  \caption{Setup parameters for near-field tomography via \algref{alg:PCT} for the experimental data set in \figref{fig:NumResExp-Data}. Non-specified parameters according to \sref{SS:NumResNFSim-Setup} and \tabref{tab:NumResNFSetup}. Physical parameters: $\lambda = \unit[0.157]{nm}$ (wave length), $\Delta x _{\text{eff}}  = \unit[29.3]{nm}$ (effective pixel size), $d _{\text{eff}}  = \unit[22.7]{mm}$ (effective detector distance).}
	  \label{tab:NumResExp-SetupParams}
	\end{table}
	  
	  Despite the flat-field-correction, the holograms contain considerable variations of the background intensity by up to $\pm 10 \, \%$ manifesting as the stripes in the sinogram in \figref{fig:NumResExp-Data}. These systematic errors 
	  render the discrepancy principle hardly applicable as a stop rule, because this would require a very accurate estimate of the data error. The heuristic choice for the initial regularization parameter $\alpha_0$ given in \eqref{eq:NumResNF-alpha0} may in principal be approximated by a guess for the $L^2$-norm $\norm{\bN^\dagger}_2$ of the exact object, based on its material composition and spatial extent. However, the systematic data errors render this once more inaccurate. Hence, $\alpha_0$ is determined by trial and error using the heuristic criterion of 5-10 initial CG-iterations as an indicator of an adequate regularization (cf. \sref{SS:NumResFF-GenSetup}). A Sobolev $H^{0.5}$-regularization term is used for noise suppression as motivated in \sref{SS:NumResNFSim-ParamStudy}. The numerical reconstruction via \algref{alg:PCT} is stopped after 12 Newton steps according to empirical observations for moderately noisy data.
	  
	  The setup parameters for near-field tomography from the experimental data set in \figref{fig:NumResExp-Data} are summarized in \tabref{tab:NumResExp-SetupParams}.

          \subsection{Reconstruction Results} \label{SS:NumResNFReal-Results}
          
          The reconstruction via \algref{alg:PCT} terminates after a total number of 305 CG-iterations, taking approximately two hours on a workstation with an 8-core Intel Xeon CPU E5-2609 at $\unit[2.40]{GHz}$ with 256 Gigabytes main memory. Hence, our Newton-based approach is numerically feasible with relatively little effort for the given data set. Yet, note that the computation time would increase dramatically by the complexity of the performance-critical Radon transform (see \sref{SS:RadonImplementation}) if we had chosen a larger object domain than $256^3$ voxels, i.e.\ a weaker support constraint.
          
          Notably, the residual $\norm{F_{\Text{dis}}(\bN_k) - \bI^{\Textbf{err}}}_{\mY_{\Text{dis}}}$, measuring the agreement of the data corresponding to the current Newton iterate $\bN_k$ and the observed intensities, reduces by a factor of less than $20 \, \%$ over the whole reconstruction starting from the initial guess $\bN_0  = 0$. From the sixth to the final iterate $\bN_{12}$ the reduction is even no more than $0.5 \, \%$. One might thus come to the conclusion that the reconstruction fails due to stagnation. However, comparing the reconstructed intensities $\bI_{12} = F_{\Text{dis}}(\bN_{12})$ visualized in \figref{fig:NumResExp-DataRecon} to the experimental measurements in \figref{fig:NumResExp-Data} suggests a different interpretation: while the holographic fringes in the data are apparently well-fitted, the observed systematic errors by variations of the background intensity and data noise seem to be effectively filtered out in the reconstruction. Accordingly, the large final data residual is indeed not a sign of failure but of the robustness of the regularized Newton method. Yet, it unfortunately precludes usage of the discrepancy principle as a stop rule since these data errors may hardly be estimated a priori. 
 	        \begin{figure}[hbt!]
	  \centering
	  \subfloat{\includegraphics[height=.44\textwidth]{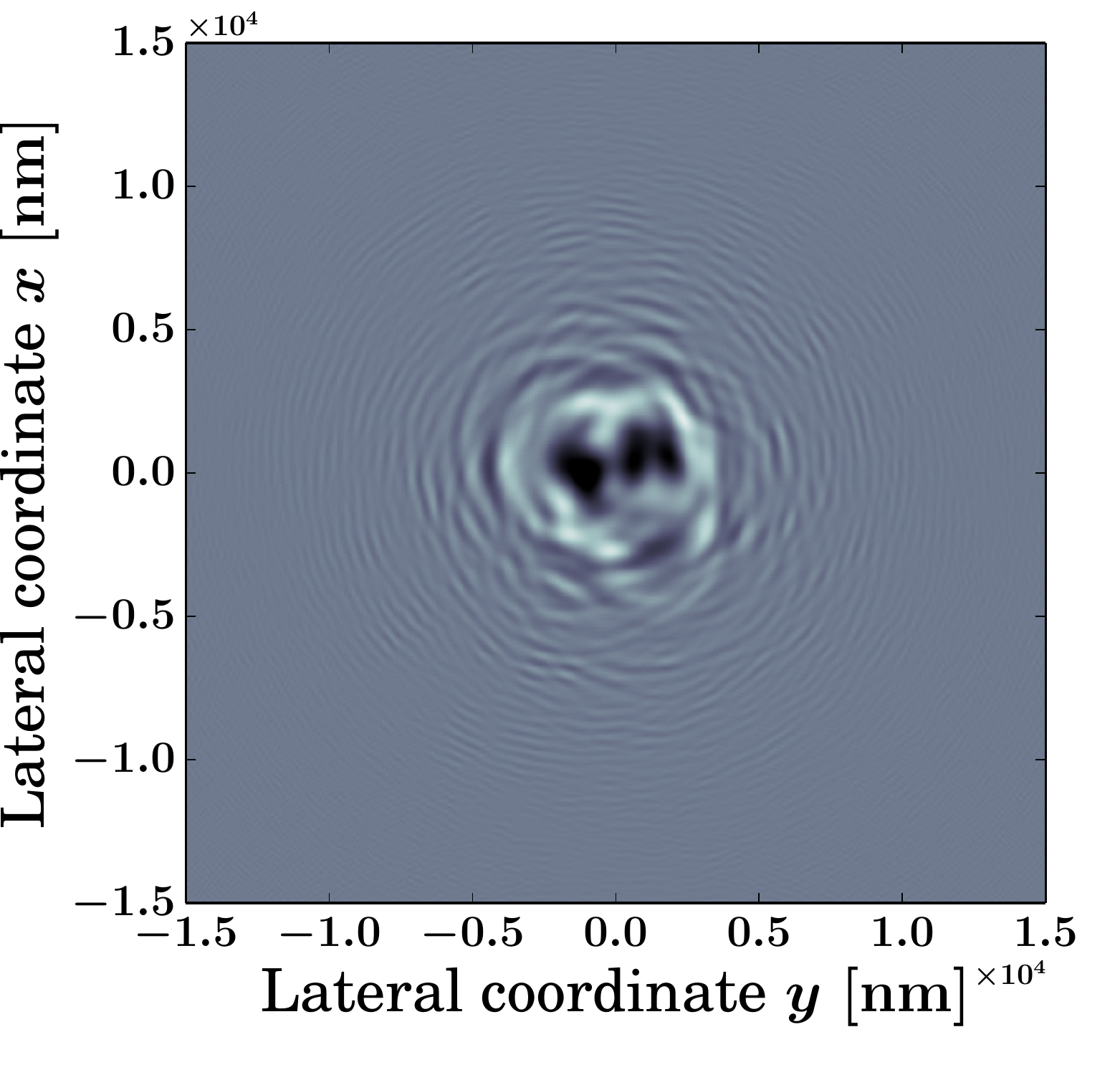} \label{fig:NumResExp-DataRecon-a} }
	  \hfill
	  \subfloat{\includegraphics[height=.44\textwidth]{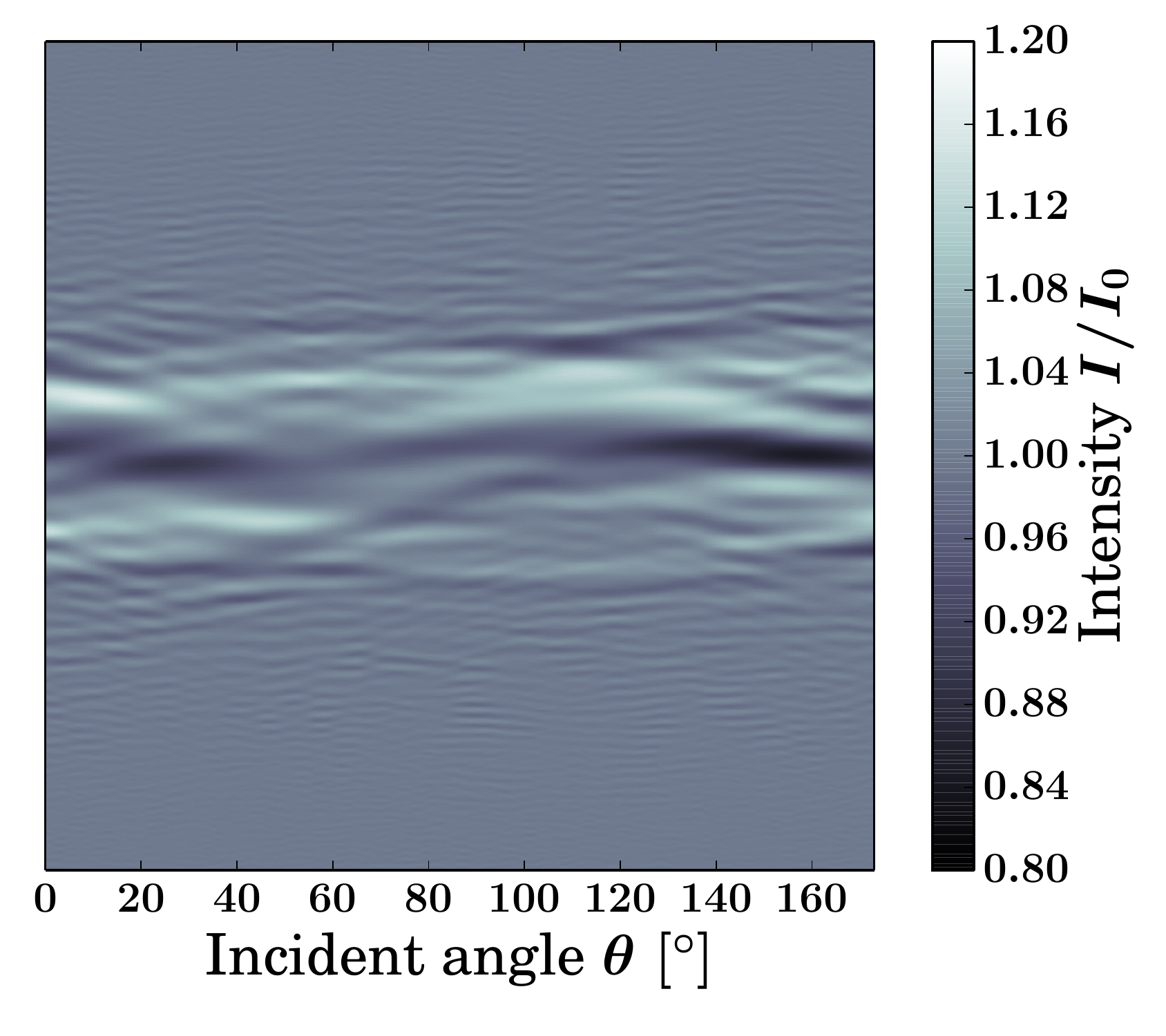} \label{fig:NumResExp-DataRecon-b} }
	  \caption{Reconstructed intensity data $\bI = F_{\Text{dis}}(\bN_{k_{\Text{stop}}})$ for the final Newton iterate ($k_{\Text{stop}} = 12$) in the reconstruction of the near-field polystyrene bead data set via \algref{alg:PCT}. Same 2D-hologram and -sinogram slices shown as for the measured data in  \figref{fig:NumResExp-Data}. Note that the holographic fringes are reproduced well whereas inhomogeneities in the background intensity and noise are apparently filtered out by the regularized Newton method. \label{fig:NumResExp-DataRecon}}
	  \end{figure}
	  
	  The central two-dimensional slice of the reconstructed object $\bN_{12} = \delta$ in the $x$-$z$-plane is shown in \figref{fig:NumResExp-ObjRecon}. For the colloidal crystal of polystyrene beads, we would expect a binary distribution of the refractive decrement $\delta \in \{0, \delta_{(\text{C}_8 \text{H}_8)_m} \}$ separated into uniform spheres of $(\text{C}_8 \text{H}_8)_m$ and vacuum. This expectation is qualitatively confirmed by \figref{fig:NumResExp-ObjRecon} up to noise and slight inhomogeneities. Different diameters of the approximately circular spots are due to the fact that merely a 2D slice plot is depicted, showing intersections of the spheres at different latitudes. From the inset plot of $\delta$ along the red line in \figref{fig:NumResExp-ObjRecon}, it can be seen that the diameter $\approx \unit[400]{nm}$ of the intersected spot roughly matches the bead size. Moreover, by measuring the length of the peak flanks in the cross section plot, giving an approximate full width at half maximum, we obtain an estimated resolution of
	  \begin{equation}
	   a_{\Text{obs}} \approx \unit[130]{nm}.
	  \end{equation}
	  This is significantly larger than the theoretical limits given by the effective pixel size $a_{\Text{pixel}} = \Delta x _{\Text{eff}} \approx \unit[29.3]{nm}$ and the regime of validity of the \emph{projection approximation}. For an object diameter of $L = 256\cdot \Delta x _{\Text{eff}}$ and the wavenumber $k = \frac{2\pi }{ \unit[0.157]{nm} }$, the latter bounds the resolution by $ a_{\Text{proj}} \gtrsim  \unit[10]{nm} $ according to \eqref{eq:ProjApproxBound}. Hence, the observed limitations must be of practical nature, for instance caused by 
	  the finite size of the ideally point-like nanofocus in \figref{fig:SetupExp}, limited coherence of the radiation
	  or systematic errors related to the flat-field correction and to the approximation by an effective parallel-beam geometry. 
	  Likewise, imperfect alignment of the holograms, correcting translations of the specimen, may cause a resolution-limiting blur in the data - in addition to errors resulting from the Newton reconstruction, of course. 
	  \begin{figure}[hbt!]
	  \centering
	  \hfill
	  \subfloat{\includegraphics[height=.44\textwidth]{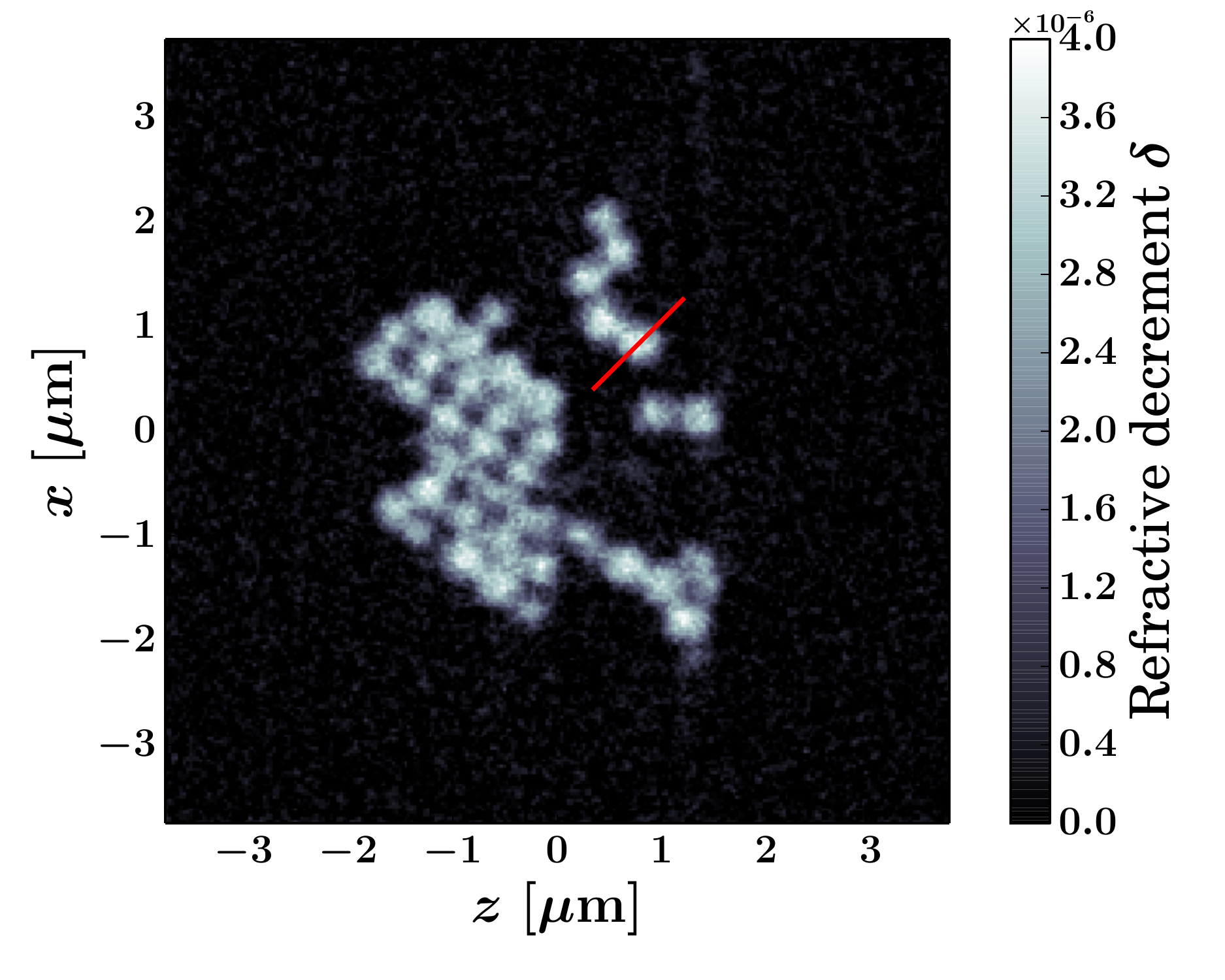} \label{fig:NumResExp-ObjRecon-a} }
	  \subfloat{ \includegraphics[width=.38\textwidth]{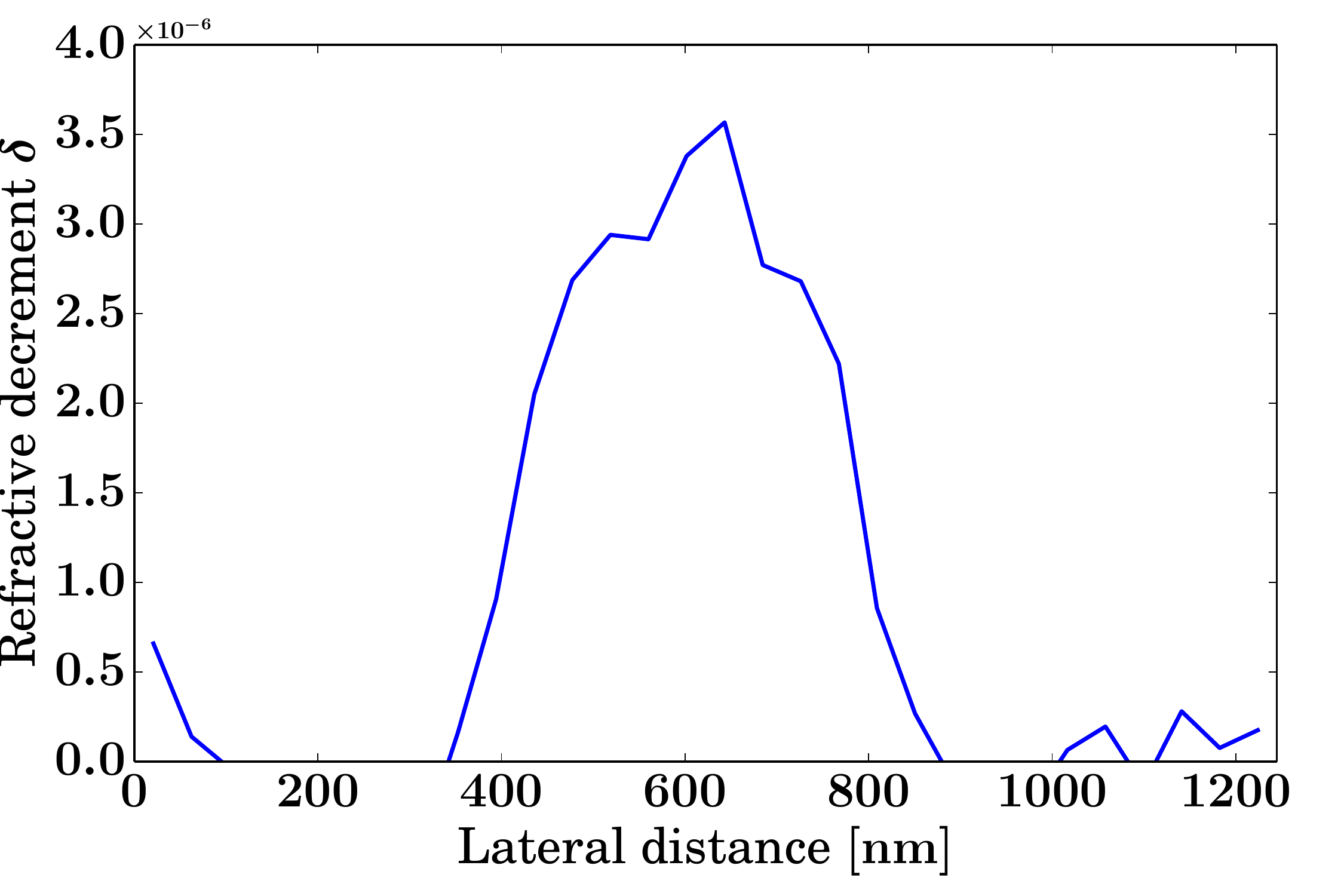} \label{fig:NumResExp-ObjRecon-b} }
	  \hfill
	  \caption{Central $x$-$z$-slice at $y= 0$ of the reconstructed object $\bN_{12} = \delta$ (left figure) obtained from the application of \algref{alg:PCT} to the data set in \figref{fig:NumResExp-Data}. The reconstructed data is shown in \figref{fig:NumResExp-DataRecon}. According to the cross section plot (right figure) along the red line, the expected binary refractive decrement $\delta \in \{0, \delta_{(\text{C}_8 \text{H}_8)_m} \}$ of the $\unit[415]{nm}$ spheres is well resolved up to moderate noise and blurry edges caused by the limited resolution. The peak value of the blue curve matches the theoretical prediction \eqref{eq:NumresExp-deltaTheo} for $\delta_{(\text{C}_8 \text{H}_8)_m}$ within an error of $\approx 5 \, \%$. \label{fig:NumResExp-ObjRecon}}
	  \end{figure}
	  
	  For a closer investigation of the resolution of the binary refractive decrement, we compute a histogram of the reconstructed $\delta$-values on the $256^3$ voxels. The result is shown in \figref{fig:NumResExp-Histo}. In addition to a strong peak around zero corresponding the background values associated with noise, a second local maximum is found as anticipated for a binary object, yet with a relatively wide peak. By estimating the maximum and its width via a local Gaussian fit visualized in \figref{fig:NumResExp-Histo}, we obtain for the material-specific refractive decrement  of polystyrene
	  \begin{equation}
	    \delta_{(\text{C}_8 \text{H}_8)_m, \Text{obs}} = (2.4 \pm 0.8) \cdot 10^{-6}. \label{eq:NumresExp-deltaObs}
	  \end{equation}
 	        \begin{figure}[hbt!]
	  \centering
	   \includegraphics[width = .8\textwidth]{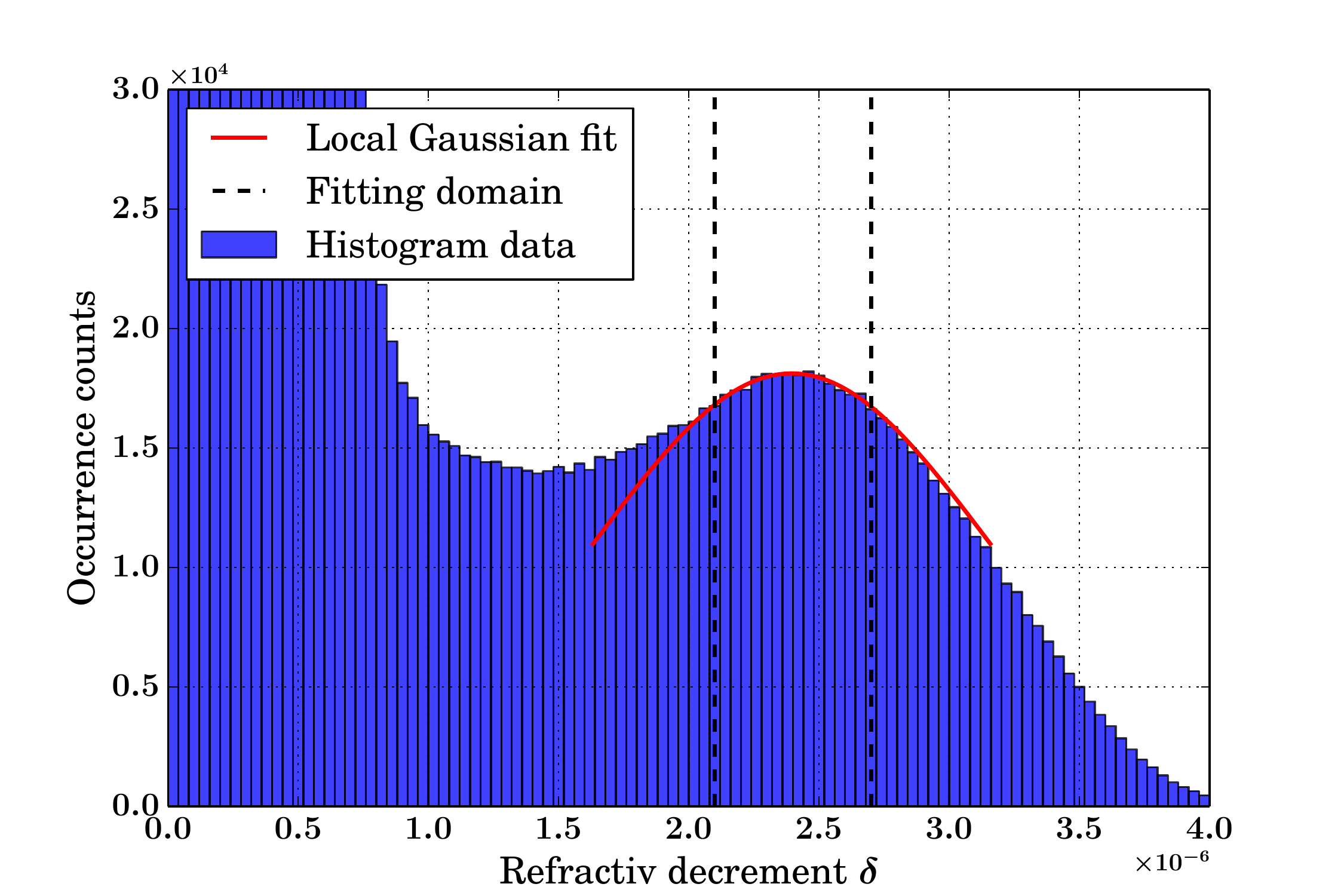} 
	  \caption{Histogram of the reconstructed refractive decrement $\delta = \bN_{12}$ of the polystyrene colloidal crystal sample. The large number of counts around the origin is due to background noise in the Newton reconstruction, whereas the second peak around $\delta \approx 2.4$ is associated with the resolved polystyrene beads visualized in \figref{fig:NumResExp-ObjRecon}. The Gaussian fit around the maximum yields an estimate $\delta_{(\text{C}_8 \text{H}_8)_m, \Text{obs}} = (2.4 \pm 0.8) \cdot 10^{-6}$ for the material's refractive decrement. \label{fig:NumResExp-Histo}}
	  \end{figure}
	  
	  Polystyrene is composed of an equal number of carbon (six neutrons and protons) and hydrogen atoms (one proton). Hence, one gram of $(\text{C}_8 \text{H}_8)_m$ contains $\approx \frac 7 {13} \cdot 6.02 \cdot 10^{23}$ electrons. Taking into account the mass density of $\unit[1.05]{\frac g {cm^3}}$ of the polystyrene nano-beads \cite{PolyBeads}, we obtain the electron density and thereby a theoretical prediction for the refractive decrement according to \eqref{eq:deltaVSrho}: 
	  \begin{equation}
	    \delta_{(\text{C}_8 \text{H}_8)_m, \Text{theo}} \approx 3.76\cdot 10^{-6} \label{eq:NumresExp-deltaTheo}.
	  \end{equation}
	  This value is about $50 \, \%$ larger than the empirical one in \eqref{eq:NumresExp-deltaObs}, deviating by more than the standard deviation of the Gaussian fit. On the other hand, note that  \eqref{eq:NumresExp-deltaTheo} is in good agreement with the upper edge of the non-negligible histogram counts in \figref{fig:NumResExp-Histo} and with the peak value of the intersected spot in \figref{fig:NumResExp-ObjRecon}. This suggests that the computed Newton reconstruction is nevertheless quantitatively correct in principal. Indeed, comparing the resolution $a_{\Text{obs}} \approx \unit[130]{nm}$ to the sphere diameter of $\unit[415]{nm}$, it becomes clear that the smeared out edges of the reconstructed polystyrene beads occupy a significantly larger volume fraction of the object domain than the peak values in their interior. The corresponding $\delta$-values thus give rise to more counts in the histogram. Accordingly, it is a characteristic value of the blurry transitional regions marking the colloids' interfaces that has been fitted in \figref{fig:NumResExp-Histo} - as is confirmed by the cross section plot in \figref{fig:NumResExp-ObjRecon}.
	  
	  This suggests to use $\delta_{(\text{C}_8 \text{H}_8)_m, \Text{obs}}$ as a threshold value defining the boundaries of the individual spheres. \figref{fig:NumResExp-Render} shows the corresponding three-dimensional contour surface computed from the numerically reconstructed $256^3$-voxel object $\bN_{12}$, yielding a 3D-rendering of the observed colloidal crystal. The uniform spherical shapes are by and large well resolved except for spurious transition pieces between neighboring beads, which arise from overlapping blur at the interfaces. At any rate, the obtained result is sufficiently accurate for the principal endeavor of the tomographic experiment: to determine the crystalline structure of the colloidal sample.
	  \begin{figure}[hbt!]
	  \centering
	   \includegraphics[width = .6\textwidth]{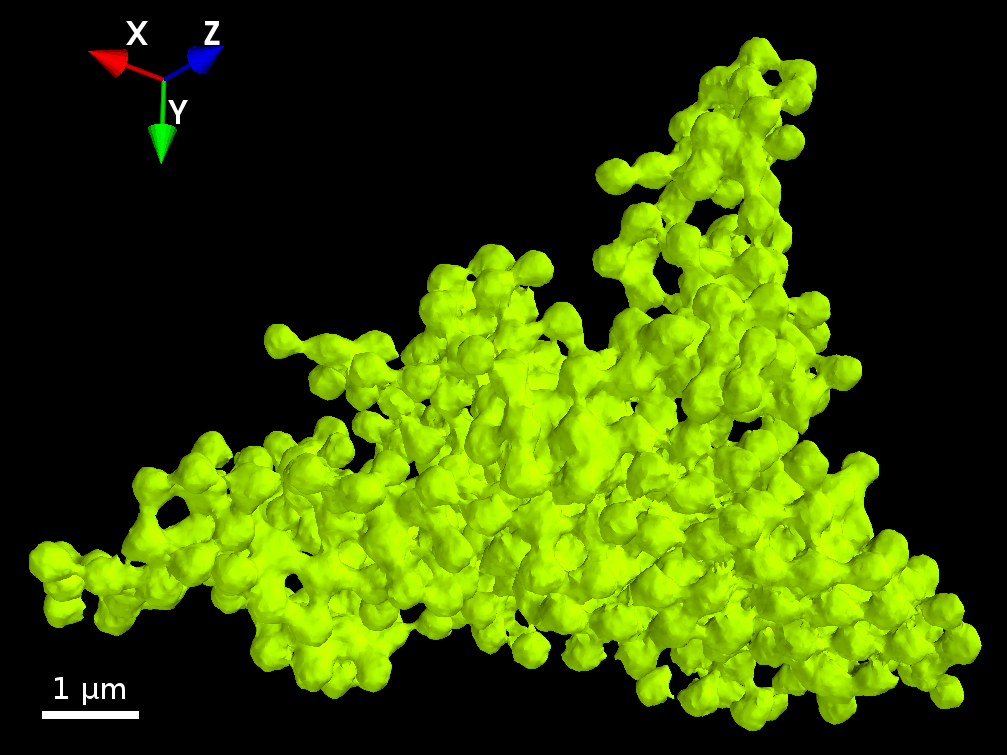} 
	  \caption{Contour plot of the reconstructed refractive decrement $\delta = \bN_{12}$ showing the three-dimensional structure of the observed colloidal crystal, using the approximated peak value in \figref{fig:NumResExp-Histo} given by \eqref{eq:NumresExp-deltaTheo} as a threshold. Spherical shapes and uniform sizes are well resolved up to transition tubes between neighboring colloids, arising from overlapping blurry regions around the theoretically sharp interfaces, which are due to the limited resolution. \label{fig:NumResExp-Render}}
	  \end{figure}
	  
	  The findings of this section's application of our regularized Newton method to an experimental near-field data set are summarized in the form of a final result:
	\begin{res}[Newton-based Near-Field Tomography from Experimental Data] \ \\ \label{res:ExpDataByNewton}
	  \algref{alg:PCT} permits quantitative near-field tomography of non-absorbing nanoscale specimen. The regularized Newton method is robust against realistic data noise and moderate variations of the background intensities, which however rule out the discrepancy principle as a stop rule. Resolution improvements might be achieved by numerically correcting for residual shifts of the holograms in the reconstruction.
	\end{res}

	  \end{section}

\end{chapter}


\begin{chapter}{Summary and Conclusions}\label{C:Concl}

%

	In this work, we have studied the problem of propagation-based X-ray phase contrast tomography and designed regularized Newton methods for numerical reconstructions. The overall aim is the recovery of the spatially varying refractive index $n = 1- \delta + \I \beta$ of an unknown specimen - for example a biological cell - from diffraction patterns recorded under illumination with coherent X-rays at different incident angles. An exemplary setup for such measurements is sketched \figref{fig:Setup}.
	
	In \chapref{C:PhysProb}, a physical model has been developed for the tomographic imaging problem based on the paraxial Helmholtz equation, adopting a geometrical optics description of the radiation-matter interaction by the \emph{projection approximation}. Thereby, nonlinear forward operators $F_d, F_\infty$ have been obtained, which map the sample information $N = 1-n$ onto the expected tomographic intensity data to be detected in the near-field or far-field, i.e.\ at moderate or large distances between the sample and the detector. This work's principal \emph{inverse problem} of reconstructing the specimen structure amounts to inverting these maps. The near-field or far-field propagation of the scattered wave field onto the detector, incorporated into the model in the form of the Fresnel propagator and the Fourier transform, respectively, has been shown to yield \emph{phase contrast}: by interference, not only X-ray absorption $\sim \beta$ manifests itself in measurable intensities but also the refractive phase shifts that are imprinted upon the transmitted radiation according to the parameter $\delta$.
	
	The encoding of specimen structure in the observable data has been further analyzed in \chapref{C:Analysis} within the derived mathematical formulation of phase contrast tomography. 
	While the forward operators have been proven to be \Frechet differentiable in \sref{S:WellDefFrechet}, i.e.\ well-posed, the corresponding inverse problem turns out to be \emph{ill-posed} in a number of different aspects: \emph{Radon inversion}, representing tomographic reconstruction (cf. \sref{S:Radon}, \sref{S:RadonIll-posed}), is not only discontinuous as it unboundedly amplifies measurement errors in large Fourier frequencies, but also imposes strong \emph{consistency conditions} between the observed diffraction patterns under different incident angles. This practically rules out existence of exact solutions for noisy intensity data.
	Additionally, \emph{phase-wrapping} prevents the unique recovery of strong objects which induce refractive phase shifts by more than a wavelength.
	
	The major source of ill-posedness in the reconstruction  $F_{\ast}(N) \mapsto N$, however, is given by the involved \emph{phase retrieval problem} analyzed in \sref{S:PhaseRetrieval}, induced by the characteristic loss of phase information in the  detection of the scattered wave field. In the far-field case, where phase retrieval corresponds to the recovery of a signal from the squared modulus (no phase information) of its Fourier transform, it is observed that non-uniquely reconstructible objects in general exist - even if additional priori constraints on support, regularity and real-valuedness are assumed. By proving the Theorems \ref{thm:UniquePhaseFF1D} and  \ref{thm:UniquePhaseFFMD}, we have demonstrated that these ambiguities may be overcome by superimposing a suitable known \emph{reference signal} upon the unknown object. The obtained result in \thmref{thm:UniquePhaseFFMD} for dimensions $\geq 2$ turns out to be significantly stronger than the 1D analogue in accordance with the general tendency  that higher dimensionality facilitates phase retrieval \cite{Fienup19782DPhaseRetrFeasible,Barakat1984}.
	
	In near-field phase contrast imaging, the unscattered part of the incident beam constitutes a \emph{natural} holographic reference for the  phase shifts and absorption induced upon the traversing X-ray wave field, i.e.\ for the sample's \emph{contact image}. As shown in \thmref{thm:NFUnique}, this leads to the startling conclusion that \emph{any} compactly supported complex-valued contact image may be recovered uniquely from near-field intensities for a suitable illumination, e.g.\ by plane waves or a Gaussian beam. Together with \cref{cor:UniqueNFPCT}, stating uniqueness of near-field phase contrast tomography for non-phase-wrapping compact specimen $N = \delta - \I \beta$, this uniqueness theorem constitutes the present work's principal theoretical result, submitted as the manuscript \cite{maretzke2014uniqueness}. In fact, the statement that arbitrary refracting and absorbing objects (modulo phase-wrapping) may be recovered from near-field intensities measured at only \emph{one} detector distance is unprecedented in its generality and has even been commonly argued to be untrue \cite{Jonas2004TwoMeasUniquePhaseRetr,Nugent2007TwoPlanesPhaseVortex,Burvall2011TwoPlanes}.
	
	Motivated by these theoretical results, regularized Newton methods \cite{Bakushinskii1992IRNM} for reconstructions in phase contrast tomography have been developed in \chapref{C:NumMeth}. This algorithmic choice both accounts for the nonlinearity and ill-posedness of the inverse problem and exploits the \Frechet differentiability of the forward operators, iteratively solving regularized local linearizations of the reconstruction problem.
	As a benefit, no \emph{global} linearizations with limited regimes of validity have to be incorporated in the approach. This renders it significantly more general than direct methods inverting the contrast transfer function (CTF) \cite{Cloetens1999,Cloetens1999Diss,BartelsDiss} or techniques based on the transport-of-intensity-equations (TIE) \cite{Reed1983TIEPhaseRetr,Nugent1996Propagationbased,Wilkins1996,Paganin1998}, which are restricted to weak objects or small propagation distances, respectively.
	Another crucial feature of the chosen reconstruction method is that phase retrieval and Radon inversion are performed \emph{simultaneously} as the forward operators are inverted as a whole. Thereby, the aforementioned consistency conditions between the tomographic projections are incorporated in the phase reconstruction, promising improved stability and accuracy.
	In the discretization of our regularized Newton algorithm for phase contrast tomography, constructed by sampling object- and intensity data on finite voxel- and pixel-grids, the iterations amount to the solution of a self-adjoint positive-definite linear problem. The latter are solved efficiently by the conjugate gradient (CG) method. However, the required evaluations of the discrete Radon transform in each CG-iteration renders the constructed \algref{alg:PCT} computationally expensive.
	
	In \chapref{C:NumRes}, numerical reconstruction results have been discussed, obtained by an implementation of \algref{alg:PCT} in \textsc{Matlab}/\textsc{Octave} \cite{Octave}. Like in the uniqueness theory, significant differences are observed concerning the numerical solution behavior in the far-field and the near-field imaging case. In the former setting treated in \sref{S:NumResFF}, ab initio reconstructions even of pure phase objects $N  = \delta$ turn out to be practically impossible as the quality of the achieved solutions is found to depend strongly on the choice of the initial guess.
	In addition to  latent phase retrieval ambiguities, this effect is attributed to the dominant quadratic nonlinearity of the far-field forward operator $F_\infty$, which is only poorly approximated by the linearizations in the Newton iterations if the initial guess is far from the exact object. A remedy is once more found in superimposing known reference signals, providing a canonical choice for a support constraint and the initial guess. Using spherical or general non-rectangular reference objects, robust artifact-free 3D reconstructions are achieved in the conducted numerical simulations of realistic far-field setups, including non-vanishing absorption as well as incomplete intensity data due to a central beam stop and a missing wedge of incident angles.
	
	Although reference objects are implementable in principal, their necessity for Newton-based far-field tomography constitutes a considerable constraint in the design of experimental setups. An alternative would be to supplement our approach with iteratively updated support estimates as in the Shrinkwrap Algorithm \cite{Fienup1986algorithm}, constructing a suitable initial guess \emph{on the fly} in some sense. Yet, the supposedly large number of iterations associated with this trial-and-error strategy would render it computationally expensive. Reconstruction algorithms based on convex optimization  such as \emph{Relaxed Averaged Alternating Reflections} (RAAR) \cite{Luke2005RAAR} thus seem generally better suited for far-field imaging owing to their greater flexibility~and robustness to the quadratic nonlinearity. Nevertheless, regularized Newton iterations could still be applied to improve initial reconstructions obtained by other methods.
	
	In the near-field case studied in \sref{S:NumResNFSim}, on the other hand, numerical reconstructions of simulated pure phase objects turn out to always stably converge up to noise level, except for the highly nonlinear problem of phase-wrapping occurring for strong objects.  While the latter may not be overcome, this work's regularized Newton approach turns out to be applicable to a wide range of Fresnel numbers, i.e.\ propagation distances, and objects inducing phase shifts up to the order of one wavelength. The optimal near-field regime for the method indeed seems to be given by such moderately strong objects of characteristic lengthscales corresponding to Fresnel numbers in the order of one or less. Fortunately, this implies that our approach exactly fills the gap in which neither CTF- nor TIE-based methods are reasonably applicable. In contrast to these, it furthermore allows for an independent reconstruction of refraction $\delta$ and absorption $\beta$ according to the numerical proof of concept in \sref{SS:NumResNFSim-GenObjects}.
	
	The benefits of our simultaneous approach to Radon inversion and phase retrieval have been evaluated in \sref{SS:NumResNFSim-EvalSimApproach}. By comparison to a regularized Newton method performing separate phase reconstructions for all incident angles, it is confirmed that the exploitation of tomographic consistency greatly improves the reconstruction - especially if $\delta$ and $\beta$ are to be recovered as independent parameters. This observation is in good agreement with results in \cite{Ruhlandt2014} obtained by the alternating-projection-type IRP algorithm. Another setting where the simultaneous approach turns out to be particularly beneficial is when the recorded holograms do not contain all fringes encoding object information due to a limited field of view. Here, the implicit data completion seems to be  stabilized by the incorporated consistency. This might motivate an adaption of our method to \emph{region of interest tomography} where the detection typically only captures a small section of  a much smaller object. Finally, it should be emphasized that also the non-tomographic Newton method performs significantly better than CTF-based reconstructions considered for comparison.
	
	Hence, we may conclude in general that regularized Newton methods are a promising approach to near-field phase contrast imaging and -tomography. Particular benefits are given by their applicability for moderately strong objects and a large bandwidth of Fresnel numbers, in addition to their robustness against noise via the choice of the regularization term (compare \sref{SS:NumResNFSim-ParamStudy}) as well as to systematic data errors. The latter has been observed in the successful application of our reconstruction method to experimental tomographic near-field data of a colloidal crystal of polystyrene-nanobeads in \sref{S:NumResNFReal}. The mathematical reason for the excellent performance of Newton methods in the near-field case seems to be once more related to the structure of the forward operator arising from contributions of the unscattered probe beam: by superposition with the latter, the imprint of the scattering object $N = \delta - \I \beta$ in the measured intensities is always \emph{linear} to leading order - different from the far-field setting. In this sense, the near-field imaging problem is only \emph{weakly} nonlinear up to moderately strong objects so that the Newton iterations may easily find their way along the predominantly linear dependence.
	
	On the other hand, the reconstruction method presented in this work may still be improved and extended in various ways. For once, it may easily generalized to incorporate intensities measured at \emph{multiple} propagation distances as used e.g.\ in CTF-reconstructions, see for instance \cite{Krenkel2014BCAandCTF}. From the promising numerical results obtained for a \emph{single} distance, it may be inferred that already two measurements are likely to permit a quantitatively accurate independent recovery of both absorption $\beta$ and phase shifts $\delta$. More data might even enable simultaneous recovery of the illumination function $P$, which is often unknown in experiments. At any rate, the reconstructions may benefit from \emph{positivity constraints} as physics dictates non-negative values for $\beta$ and $\delta$. Such may be incorporated in a generalization of the present approach by \emph{semismooth Newton methods} \cite{Hintermuller2002SSNewton,Griesse2008Semismooth}. Similarly, the latter would allow for the prescription a maximum ratio $\beta/\delta$ of say $\frac 1 {10}$ or $\frac 1 {100}$, which enforces a physically reasonable coupling between the two parameters as any absorbing matter to be reconstructed is necessarily also refracting. These constraints might significantly reduce the halo-artifacts observed in the simultaneous reconstruction of $\delta $ and $\beta$ and thereby indeed permit an accurate recovery of both parameters from intensity data at a \emph{single} propagation distance.
	
	The latter is possible in principal according to our near-field uniqueness result. However, in order to tell whether practically relevant or merely  a mathematical curiosity, the statement has to be supplemented with \emph{stability estimates} bounding the reconstruction error. A good starting point for such an analysis is to investigate whether \thmref{thm:NFUnique} remains valid in the weak object limit, i.e. for a CTF-like linearization in the contact image, because showing stability based on the linear case would simplify matters considerably. These questions are subject to future work.
	
	In any event, simultaneous Radon inversion and phase retrieval is likely to provide a considerable stabilization of the reconstructions according to the results of this work. In order to retain numerically feasibility also for discrete objects of $1024^3$ voxels or more, however, significant efficiency gains have to be achieved. For example, such could be obtained via a divide-and-conquer modification of our regularized  Newton algorithm, reconstructing only with respect to a small set of (neighboring) incident angles in each Newton iteration in the spirit of the \emph{Algebraic Reconstruction Technique} \cite{Kaczmarz1937ART,Gordon1970ART}. By reducing memory requirements, this would also allow for a massively parallel implementation on graphic cards. Non-simultaneous Newton-based phase retrieval, on the other hand, is likely to provide a numerically efficient and accurate substitute for CTF-based methods already in near future.
	
	To conclude, the present work has thus revealed that regularized Newton methods yield promising - not to say excellent - results in (near-field) phase contrast tomography for both simulated and experimental data - even though the potential of the approach has by far not been fully exploited yet.

\end{chapter}


\appendix

\begin{chapter}{Mathematical Background}\label{C:MathPre}

In this appendix, some basic notions and central results from functional analysis are reviewed providing a mathematical background for the main chapters. In particular, we introduce the Fourier- and Radon transforms which yield a mathematical description of tomographic imaging. No claim of completeness is raised for the given overview, which is mainly based on the books \cite{Werner2007FA,Hoermander,Natterer,Evans}. Whenever proofs are omitted, these can be found in the references.

\begin{section}{Operators and Adjoints} \label{S:OpsAdjoints}

As shown in \chapref{C:PhysProb}, image reconstruction in phase contrast tomography amounts to solving an equation of the form
\begin{equation}
 F(f) = g \label{eq:GeneralOpEq}
\end{equation}
for a map $F :  \mX \supset  U \to \mY $ between topological spaces $\mX $ and $\mY$, i.e.\ to finding its \emph{inverse}
\begin{equation}
 F^{-1}: F(U) \to \mX; \; g \mapsto f \label{eq:GeneralInverse}
\end{equation}
In the following, we introduce some notions from functional analysis providing a framework for the analysis of such general \emph{inverse problems}. In the remainder of this section, let $\mX$, $\mY$ be real or complex \emph{Banach-} or \emph{Hilbert Spaces}. For theoretical background on these, refer to \cite{Werner2007FA}. We begin by introducing \emph{linear} operators:
	\vspace{1em}
	\begin{df1}[Bounded Operators and their Adjoints \text{\cite[pp. 238 f.]{Werner2007FA}}]\label{def:OpsAdjoints}
		We call a linear map $T: \mX \to \mY $ a \emph{bounded linear operator} and write $\sL(\mX, \mY)$ if
		\begin{equation}
		\norm{T} := \sup_{x \in \mX \setminus \{0\}} \frac{ \norm{T x }_\mY }{ \norm{x}_\mX } < \infty. \label{eq:defOpNorm}
		\end{equation}
		In  Hilbert spaces  $\mX$ and $\mY$ with inner products $\ip{\cdot}{\cdot}_\mX$ and $\ip{\cdot}{\cdot}_\mY$, any $T \in \sL(\mX, \mY)$ has a unique \emph{adjoint operator} $T^\ast \in \sL(\mY, \mX)$, defined by
		\begin{equation}
		\ip{T x }{y}_\mY = \ip{x }{T^\ast y}_\mX  \MTEXT{for all} x \in \mX, y \in \mY. \label{eq:defAdj}
		\end{equation}
		$T$ is called an \emph{isometry} if  $\norm{T x }_\mY = \norm{ x }_\mX$ for all $x \in \mX$. A surjective isometry is denoted as a \emph{unitary operator} and is characterized by $T^{-1} = T^\ast.$
	\end{df1}
	\vspace{1em}
Note that a linear operator $T: \mX\to\mY$ is continuous if and only if it is bounded \cite[p. 45]{Werner2007FA}. By definition, linear combinations and compositions of bounded linear operators are again bounded. The adjoint operator introduced in \eqref{def:OpsAdjoints} is characterized by the following properties:
	\vspace{1em}
	\begin{th1}[Properties of the Adjoint \text{\cite[pp. 238 f.]{Werner2007FA}}]\label{thm:AdjointProps}
		For Hilbert spaces $\mX$, $\mY$, $\mW$, let $S,T \in \sL(\mX,\mY)$, $R \in \sL(\mY,\mW)$ and $\lambda, \mu \in \mK  \in \{\mR, \mC \}$. Then
		\begin{itemize}
		 \item[\Text{(a)}] $(\lambda S+\mu T)^\ast = \cc \lambda S^\ast + \cc \mu T^\ast$
		 \item[\Text{(b)}] $(RS)^\ast =  S^\ast R^\ast$
		 \item[\Text{(c)}] $\norm{S^\ast} = \norm{S} $
		  \item[\Text{(d)}] $ S^{\ast \ast} = S $.
		\end{itemize}
	\end{th1}
	\vspace{.5em}

 Any \emph{complex} Banach space $\mX $ can be turned into a \emph{real} Banach space containing the same elements by restricting scalar multiplication to reals. We denote this real analogue by  $\mX_{\mR}$. Any $T \in \sL(\mX , \mY )$ can be interpreted as an operator in $\sL(\mX_{\mR}, \mY_{\mR})$ within the framework of this identification. If $\mX$ is a complex Hilbert space, $\mX_{\mR}$ equipped with the inner product
\begin{equation}
 \ip{x}{y}_{\mX_{\mR}} = \Re(\ip{x}{y}_{\mX }) \MTEXT{for all} x,y \in \mX_{\mR} = \mX .
\end{equation}
 becomes a real Hilbert space. Moreover, the adjoint $T^\ast$ of a bounded operator $T: \mX \to \mY$ is retained under the transition $\mX \mapsto \mX_{\mR}$,  $\mY \mapsto \mY_{\mR}$ by definition.
 
 The subsequent examples illustrate the above definitions, characterizing certain operators that are needed in the sequel of this work. The considered $L^p$-spaces are introduced in \sref{S:SchwartzLp}.
\vspace{1em}
	\begin{ex1}[Adjoint Operators] \ \\ \label{ex:Adjoints}
	\begin{itemize}
	 \item[\Text{(a)}] \vspace{-1em} For a closed subspace $\mV \subset \mX$ of a Hilbert space $\mX$, the canonical embedding
	 \begin{equation*}\iota: \mV \to \mX; \; v \mapsto v \end{equation*}
	 defines a bounded linear operator. Its adjoint is given by the \emph{orthogonal projection} $\cP: \mX \to \mV$ onto $\mV$ \Text{(see \cite[pp. 226 f.]{Werner2007FA} for definition)}.
	 
	 \item[\Text{(b)}] Let $\Omega \subset \mR^m$, $\varphi : \Omega \to \mR$ measurable and $\Phi(x) := \exp(\I \varphi(x))$ for all $ x \in \Omega$. Define the pointwise multiplication operator
	 \begin{equation*}\cM_{\Phi} : \Lp{2}{\Omega} \to \Lp{2}{\Omega}; \; f \mapsto \Phi \cdot f  \end{equation*}
	 Then $\cM_{\Phi}$ is unitary with $ \cM_{\Phi}^{-1} = \cM_{\Phi}^\ast =  \cM_{\cc \Phi}$.

	 
	  \item[\Text{(c)}] For $\Omega \subset \mR^m$, $1 \leq p \leq \infty$, the pointwise real part of a function $f$, given by $\Re(f)(x) = \Re(f(x))$, defines an bounded $\mR$-linear operator
	 \begin{equation*}
	  \Re: \Lp{p}\Omega _{\mR} \to \Lp{p} \Omega _{\mR}; \; f \mapsto \Re(f).
	 \end{equation*}
	 Moreover, this operator is self-adjoint for $p = 2$.
	\end{itemize}
	\end{ex1}
	\begin{pf} \ \\ 
	\begin{itemize}
	 \item[\Text{(a)}] For all $v \in \mV$, $x \in \mX$, we have
	 \begin{equation*}\ip{v }{\iota^\ast (x)}_{ \mV } = \ip{\iota(v) }{x}_{ \mX }   \stackrel{v \in \mV}= \ip{v }{ \cP x}_{ \mV } + \underbrace{ \ip{ \iota(v) }{ (x - \cP x)}_{ \mX } }_{= 0}   \end{equation*}
	 since $x - \cP x$ is in the orthogonal complement of $\mV$.
	 
	 \item[\Text{(b)}] Isometry and surjectivity follow from the relations  
	  \begin{equation*} \norm{ \cM_{\Phi} f }_{\Lp 2 \Omega}^2  = \int_{\Omega} \underbrace{|\exp(\I \varphi )|^2}_{= 1} \cdot |f|^2  \; \D x =  \norm{  f }_{\Lp 2 \Omega}^2,  \end{equation*}
	  and $\cM_{\Phi}(\cM_{\cc \Phi} f) = \exp(\I \varphi) \cdot \exp( - \I \varphi ) \cdot f  = f$ for all $f \in \Lp 2 \Omega$. The latter furthermore implies $\cM_{\cc \Phi} = \cM_{\Phi}^{-1} = \cM_{\Phi}^{\ast}$.

	 
	  \item[\Text{(c)}] Boundedness with $\norm{\Re} = 1$ follows from the estimate for all  $f \in \Lp{p}\Omega _{\mR}$
	  \begin{equation*} \norm{ \Re(f) }_{\Lp{p}\Omega _{\mR}}^2 \leq \norm{ \Re(f) }_{\Lp{p}\Omega _{\mR}}^2 + \norm{ \Im(f) }_{\Lp{p}\Omega _{\mR}}^2 = \norm{ f }_{\Lp{p}\Omega _{\mR}}^2. \end{equation*}
	  In addition, we have for all $f,g \in \Lp{2}\Omega _{\mR}$
	   \begin{align*}\pushQED{\qed} 
	     \ip{f }{\Re^\ast (g)}_{\Lp{2}\Omega _{\mR}} &= \ip{\Re(f) }{g}_{\Lp{2}\Omega _{\mR}} = \Re\left(\ip{\Re(f) }{\Re(g) + \I \Im(g)}_{\Lp{2}\Omega }\right) \\
	  &=   \ip{\Re(f) }{\Re(g) }_{\Lp{2}\Omega }  = \Re\left(\ip{\Re(f) + \I \Im(f) }{\Re(g)}_{\Lp{2}\Omega }\right) \\
	  &=  \ip{f }{\Re(g)}_{\Lp{2}\Omega _{\mR}}.  \qedhere \popQED
	  \end{align*}
	 \end{itemize}
	 	  \renewcommand{\qedsymbol}{}	 
	\end{pf}

\end{section}

\begin{section}{\Frechet Derivatives} \label{S:Frechet}

For the inverse problem of the form \eqref{eq:GeneralOpEq} studied in this work, the operator $F$ is \emph{nonlinear}. Its solution is sought by Newton-type methods introduced in \chapref{C:NumMeth}, iteratively solving linearized versions of \eqref{eq:GeneralOpEq}. To this end, we need a notion of differentiability in the sense of a local best linear approximation. This is established by the concept of \Frechet differentiability:
	\vspace{1em}
	\begin{df1}[\Frechet Derivative \text{\cite[p. 123]{Hutson2005}}]\label{def:FrDer}
		Let $\mX$, $\mY$ be Banach spaces, $U \subset \mX$ open and $F: U \to \mY$. Then $F$ is called \emph{(Fr\'echet) differentiable} at $x_0 \in U$ if there exists an operator $F'[x_0] \in \sL(\mX, \mY) $ such that
		\begin{equation}
		 \lim_{\norm{h}_\mX \to 0} \frac{ \norm{F(x_0 + h) - F(x_0) - F'[x_0]h}_\mY }{\norm{h}_\mX} = 0 \label{eq:defFrechet}
		\end{equation}
		In this case $F'[x_0]$ is called the \emph{\Frechet derivative} of $F$ at $x_0$. $F$ is called \emph{(Fr\'echet) differentiable}, if it is differentiable for all $x_0 \in U$.
	\end{df1}
	\vspace{1em}
	Note that \Frechet differentiability implies in particular \begin{equation*}\lim_{\norm{h}_\mX \to 0}  F(x_0 + h) - F(x_0) = 0 \end{equation*} and thus continuity of $F$ by definition. In the following, we summarize further
	properties, which show its analogy to derivatives of functions in $\mR$:
	\vspace{1em}
	\begin{th1}[Properties of the \Frechet Derivative \text{\cite[p. 125]{Hutson2005}, \cite[pp. 120 f.]{Werner2007FA}}]\label{thm:FrechetProps}
		Let $\mX$, $\mY$, $\mW$ be Banach spaces, $U \subset \mX$ and $V \subset \mY$ open and $F,H: U \to \mY$, $G: V \to \mW$  such that $F(U) \subset V$. Then the following holds true:
		\begin{itemize}
		 \item[\Text{(a)}] \emph{(sum rule)} If $F,H$ are differentiable at $x_0 \in \mX$ then $\lambda F + \mu H$ is differentiable at $x_0 \in \mX$ for all $\mu ,\lambda \in \mR$ and
		 \begin{equation}
		  (\lambda F + \mu H)'[x_0] = \lambda F'[x_0] + \mu H'[x_0] \label{eq:SumRule}
		 \end{equation}
		 \item[\Text{(b)}] \emph{(chain rule)} If $F$ is differentiable at $x_0 \in U$ and $G$ is differentiable at $F(x_0) \in V$, then $G \circ F$ is differentiable at $x_0$ and 
		 \begin{equation}
		    (G \circ F)'[x_0] = G'[F(x_0)] \circ F'[x_0] \label{eq:ChainRule}
		 \end{equation}
		  \item[\Text{(c)}] \emph{(product rule)} A bounded bilinear map $b: \mX \times \mY \to \mW$ is differentiable with
		   \begin{equation}
		  b'[(x_0,y_0)](h_x,h_y) = b(x_0,h_y)  + b(h_x,y_0) \label{eq:ProdRule}
		 \end{equation}
		\item[\Text{(d)}] \emph{(constant maps)} If $F$ is \emph{constant}, i.e.\ $F(x) = c$ for all $x \in U$ and some $c \in \mY$, then $F$ is differentiable with
		 \begin{equation}
		  F'[x_0] = 0 \MTEXT{for all} x_0 \in U \label{eq:ConstOpFrechet}
		 \end{equation}
		 \item[\Text{(e)}] \emph{(linear maps)} For $U = \mX$ and $F$ \emph{linear}, $F$ is differentiable if and only if $F$ is bounded with
		 \begin{equation}
		  F'[x_0] = F \MTEXT{for all} x_0 \in \mX \label{eq:LinOpFrechet}
		 \end{equation}
		 \item[\Text{(f)}] \emph{(extrema)} If $\mY = \mR$ and $F$ is \Frechet differentiable with a local extremum at $x_0 \in U$, then $F'[x_0] =0$
		\end{itemize}
	\end{th1}
	\vspace{1em}
	\newpage
	We conclude this section by explicitly computing the \Frechet derivatives of some simple nonlinear operators that are needed in this work:
	\vspace{1em}
	\begin{ex1}[\Frechet Derivatives of special Operators] \label{ex:Frechet}  \ \\
	 \begin{itemize}
	  \item[\Text{(a)}] \vspace{-1em} For $\mX$ a real Hilbert space
	  $F: \mX \to \mR; \; x \mapsto \norm{x}_\mX^2$ 
	  is \Frechet differentiable with
	  \begin{equation*}
	   F'[x_0]h =  2 \ip{x_0}{h}_\mX \MTEXT{for all} x_0, h \in \mX.
	  \end{equation*}
	  \item[\Text{(b)}] For $\Omega \subset \mR^m$ measurable and $1 \leq p \leq \infty$, the pointwise squared modulus
	  	   $F: \Lp{2p}{\Omega} \to \Lp{p}{\Omega}; \; f \mapsto |f|^2 $ 
	  is \Frechet differentiable with 
	  \begin{equation*}
	   F'[f_0]h = 2 \Re ( \cc{f_0} \cdot h ) \MTEXT{for all} f_0 , h \in \Lp{2p}{\Omega}.
	  \end{equation*}
	  \item[\Text{(c)}] For $\Omega \subset \mR^m$ measurable, the pointwise exponential
	  	  	   \begin{equation*} F: L^\infty ( \Omega ) \to L^\infty ( \Omega );\; f \mapsto [x \mapsto \exp(f(x))]   \end{equation*}
	  is \Frechet differentiable with derivative
	  \begin{equation*}
	   F'[f_0]h = \exp(f) \cdot h \MTEXT{for all} f_0 , h \in \Lp{\infty}{\Omega}
	  \end{equation*}
	 \end{itemize}
	\end{ex1}
	\begin{pf} \ \\
	 \begin{itemize}
	  \item[\Text{(a)}]
	   \vspace{-1em} Using $\norm{x}_\mX^2 = \ip{x}{x}_\mX$ and bilinearity and symmetry of $\ip{\cdot}{\cdot}_\mX$, we obtain for all $x_0, h \in \mX$
	  \begin{equation*}
	   F(x_0 + h) -F(x_0) = \ip{x_0 + h}{x_0 + h}_\mX- \ip{x_0}{x_0}_\mX = 2 \ip{x_0}{ h}_\mX +  \norm{h}_\mX^2.
	  \end{equation*}
	  \item[\Text{(b)}] Invoking the relations $|f|^2 = \cc f \cdot f $ and $2\Re(f) =f +  \cc f $ yields
	  \begin{equation*}
	   F(f_0 + h) -F(f_0) = \cc{(f_0 + h)}\cdot (f_0 + h) - \cc{f_0} \cdot f_0   = 2 \Re( \cc{f_0}\cdot h) + |h|^2
	  \end{equation*}
	  for all $f_0, h \in \Lp{2p}{\Omega}$. Consequently,
	  \begin{equation*}
	   \norm{F(f_0 + h) -F(f_0) - 2 \Re ( \cc{f_0} \cdot h )}_{\Lp{p}{\Omega}} = \norm{|h|^2}_{L^p(\Omega)} =    \norm{h}_{L^{2p}(\Omega)}^2.
	  \end{equation*}
	  \item[\Text{(c)}] For $f_0,h \in L^\infty ( \Omega )$, $\norm{\cdot} := \norm{\cdot}_{L^\infty ( \Omega ) }$, this follows from the estimate
	   \begin{align*}\pushQED{\qed} 
	   \norm{ \exp(f_0 + h) - \exp(f_0) - \exp(f_0) \cdot h } &=  \norm{ \exp(f_0)\cdot \left(  \exp(h) -  1 - h \right) } \\
	    \leq \norm{ \exp(f_0) }  \left\|  \sum_{k=2}^\infty \frac{ h^k  }{k!} \right\|  &\leq  \norm{ \exp(f_0) } \norm{h}^2 \left(  \sum_{k=0}^\infty \frac{ \norm{h}^k  }{(k+2)!} \right) \\
	     &\leq \norm{ \exp(f_0) } \exp\left( \norm{h} \right) \norm{h}^2.  \qedhere \popQED
	  \end{align*}
	 \end{itemize}
	 	  \renewcommand{\qedsymbol}{}
	\end{pf}
	\vspace{1em}

%
%
%
%
 
\end{section}

\begin{section}{Function- and Distribution Spaces} \label{S:SchwartzLp}

\subsection{Lebesgue $L^p$-Spaces} \label{SS:Lp}

For an open subset $\Omega \subset \mR^m$, let $\Ck k \Omega, \, k \in \mN_0 \cup\{\infty\}$ denote the space of all $k$-times continuously differentiable functions. We further define
\begin{equation}
 \Ckc k \Omega := \{ \phi \in \Ck k \Omega :  \supp(\phi) \subset \Omega \text{ compact} \}  \label{eq:DefCc}
\end{equation}
as the $\Ck k \Omega$-functions with \emph{compact support} $\supp(\phi):= \closure{ \{ \bx \in \mR^m: \phi(\bx ) \neq 0 \} }$.

Moreover, let $\Lp{p} \Omega$ for $1 \leq p \leq \infty$ denote the Banach spaces of all measurable functions $f: \Omega \to \mC$ such that
\begin{equation}
\infty > \norm{f}_{\Lp{p}{\Omega}} = \begin{cases}
                              \left( \int_{\Omega} |f(\bx)|^p \; \D x \right)^{\frac 1 p } &\text{for }p< \infty \\
                              \inf_{g = f}\sup_{\bx\in \Omega} |g(\bx)|  &\text{for }p =  \infty
                             \end{cases},                              
 \label{eq:DefLpNorm}
\end{equation}
with the usual identification $f=g$ if and only if $f(\bx) = g(\bx)$ for almost all $\bx \in \Omega$. See \cite[sec. I.1]{Werner2007FA} for details.
 
 Recall that the \emph{dual space} $\Lp{p} \Omega'$ of all linear and continuous functionals on $\Lp{p} \Omega$ can be identified with $\Lp{q} \Omega$ for $q = (1- \frac 1 p)^{-1}$ by the isometric anti-isomorphism \cite[p. 60]{Werner2007FA}
 \begin{equation}
 T: \Lp{q} \Omega \to \Lp{p} \Omega'; \; (T \phi)(f) := \int_{\Omega} \cc{\phi} f \; \D x  \label{eq:LpLqDual}
 \end{equation}
$\Lp{p}{\mR^m}$-functions with support in $\Omega$ may furthermore be canonically identified with elements in $\Lp{p} \Omega$ via the embedding
\begin{equation}
\iota: \Lp{p} \Omega \hookrightarrow \Lp{p}{\mR^m}; \; \iota(f)(\bx) = \begin{cases}
                              f(\bx)  &\text{for }\bx \in \Omega \\
                              0 &\text{for } \bx \notin \Omega
                             \end{cases}. \label{eq:LpLpREmbed}
\end{equation}
This identification is frequently made implicitly, for instance whenever Fourier- or Radon transforms (see \sref{S:FT} and \sref{S:Radon}) are evaluated on $\Lp{p} \Omega$. Furthermore, we use the inclusion $\Lp{q} \Omega \subset \Lp{p} \Omega$ for $ p < q  $, valid on \emph{bounded} domains $\Omega$:
\vspace{1em}
\begin{th1}[$L^p$-Embeddings on Bounded Domains] \label{thm:LpBoundedEmbed}
 Let $1 \leq p < q \leq \infty$ and $\Omega \subset \mR^m$ with finite measure $\mu(\Omega) < \infty$. Then $\Lp{q} \Omega \subset \Lp{p} \Omega$ and the embedding
 \begin{equation*}
  \iota : \Lp{q} \Omega  \hookrightarrow \Lp{p} \Omega
 \end{equation*}
 is continuous with norm $\norm{\iota} \leq \mu(\Omega)^{\frac 1 p - \frac 1 q }$.
\end{th1}
\begin{pf}
 Let $f \in \Lp{q} \Omega$. For $q < \infty$ the statement follows from \emph{Jensen's inequality} which reads for measurable $g: \Omega \to \mR_{\geq 0}$ and convex $\varphi: \mR_{\geq 0} \to \mR$ \cite[p. 152]{KlenkeWTheorie}
 \begin{equation*}
  \varphi\left(\frac{1}{\mu(\Omega) } \int_{\Omega} g \; \D x  \right) \leq \frac{1}{\mu(\Omega) } \int_{\Omega} \varphi \circ g \; \D x.
 \end{equation*}
 Setting $g := |f|^p$, $\varphi: x \mapsto x^{\frac q p }$, this yields
 \begin{equation*}
  \mu(\Omega)^{-\frac q p } \norm{f}_{\Lp{p} \Omega}^{q} \leq \mu(\Omega)^{-1} \norm{f}_{\Lp{q} \Omega}^{q} < \infty
 \end{equation*}
 which implies $ f \in \Lp{p} \Omega$ and $\norm{\iota} \leq \mu(\Omega)^{\frac 1 p - \frac 1 q }$ by \defref{def:OpsAdjoints}. For $q = \infty$, the result is obtained by estimating the integrand in \eqref{eq:DefLpNorm} by $\norm{f}_{\Lp{q} \Omega} = \sup_{\bx \in \Omega} |f(x)|$:
  \begin{equation*}\pushQED{\qed} 
  \norm{f}_{\Lp{p} \Omega} = \left( \int_{\Omega} |f(\bx)|^p \; \D x \right)^{\frac 1 p } \leq \norm{f}_{\Lp{\infty} \Omega} \left( \int_{\Omega} \D x \right)^{\frac 1 p } \stackrel{q = \infty}= \mu(\Omega)^{\frac 1 p - \frac 1 q } \norm{f}_{\Lp{q} \Omega}. \qedhere \popQED
 \end{equation*}
 \renewcommand{\qedsymbol}{}
\end{pf}

\subsection{Schwartz Spaces} \label{SS:Schwartz}

A convenient space to study the Fourier- and Radon transforms is the \emph{Schwartz space} $\sS(\mR^m)$, given by all smooth, rapidly decaying functions, i.e.\ by all $\phi \in \Cinf{\mR^m}$ such that for all multi-indices $\alpha, \beta \in \mN_0^m$
	\begin{equation}
	 \sup_{\bx \in \mR^m} | \bx^\beta \partial^\alpha \phi (\bx) | < \infty. \label{eq:DefSchwartz}
	\end{equation}
The semi-norms on the left hand side of \eqref{eq:DefSchwartz} induce a topology on $\sS(\mR^m)$ which turns it into a \Frechet space \cite[p. 160]{Hoermander}, i.e.\ a locally convex complete metric space that is a little more general than a Banach space \cite[p. 464]{Werner2007FA}.

From the definitions it is clear that $\Cinfc{\mR^m} \subset \sS(\mR^m) \subset \Lp{p}{\mR^m}$. An important result is that these inclusions are \emph{dense} (see \cite[p. 28]{Werner2007FA} for the definition) under certain conditions, meaning that for instance elements in $\Lp{2}{\mR^m}$ may be approximated arbitrarily well by $\sS(\mR^m)$-functions in $L^2$-norm:
\vspace{1em}
\begin{th1}[Dense Inclusions \text{\cite[p. 163]{Hoermander}, \cite[p. 209]{Werner2007FA}}] \label{thm:DenseIncl}
 Let $m \in \mN$, $ 1 \leq p < q \leq \infty$ and $\Omega \subset \mR^m$ open. Then
 \begin{itemize}
  \item[\Text{(a)}] $\Cinfc{\mR^m} \subset \sS(\mR^m) \subset \Lp{p}{\mR^m}$ where the inclusions are dense for $p < \infty$
    \item[\Text{(b)}] $\Cinfc{\Omega}  \subset \Lp{p}{\Omega}$ and denseness holds whenever $p < \infty$ or $\Omega$ is bounded
   \item[\Text{(c)}] If $\Omega \subset \mR^m$ is bounded, then $\Lp{q}{\Omega} \subset \Lp{p}{\Omega}$ and the inclusion is dense
 \end{itemize}
\end{th1}
\vspace{1em}

\subsection{Tempered Distributions} \label{SS:TempDistri}
Let $\sS'(\mR^m)$ denote that \emph{dual space} of $\sS(\mR^m)$, i.e.\ the space of all continuous linear functionals $\phi: \sS(\mR^m) \to \mC$.
By virtue of the identification \eqref{eq:LpLqDual}, $L^q$-functions may be regarded as elements in $\sS'(\mR^m)$ where \thmref{thm:DenseIncl} implies
\begin{equation}
 \Lp{q}{\mR^m} \subset \sS'(\mR^m) \MTEXT{for all} 1 \leq q \leq \infty. \label{eq:LpinSdash}
\end{equation}
The elements in $\sS'(\mR^m)$ are called \emph{tempered distributions} and of much lesser regularity than $\Lp{q}{\mR^m}$-functions. This can be seen from the fact that the maps
\begin{align}
\partial^\alpha: &\sS'(\mR^m) \to \sS'(\mR^m); \; (\partial^\alpha u)(f) = (-1)^{|\alpha|} T( \partial^\alpha  f) \label{eq:DistriDeriv} \\
X^\alpha: &\sS'(\mR^m) \to \sS'(\mR^m); \; (X^\alpha u)(f) =   T( \bx^\alpha \cdot  f) \label{eq:DistriPolyMulti}
\end{align}
are well-defined for all $\alpha  \in  \mN_0^m$ \cite[p. 437]{Werner2007FA}, which implies in particular that any $\phi \in \Lp{q}{\mR^m}$ has derivatives $\partial^\alpha \phi$ in $\sS'(\mR^m)$ for arbitrary $\alpha \in \mN_0^m$.

The degree of singularity of a distribution $u \in \sS'(\mR^m)$ is expressed by its \emph{order}, defined as the minimum $  N \in \mN_0$ such that for some $C>0$ \cite[pp. 33 f.]{Hoermander}
\begin{equation}
 |u(\phi)| \leq C \sum_{\alpha \in \mN_0^m: |\alpha| \leq N} \sup_{\bx \in \mR^m} |\partial^\alpha \phi(\bx)| \MTEXT{for all} \phi \in \Ckc{\infty}{\mR^m}. \label{eq:DefDistrOrder}
\end{equation}
For $u \in \sS'(\mR^m)$ and $V \subset \mR^m$ open, we write $u_{|V} = 0 $ iff $u(\phi) = 0$ for all $\phi \in \sS(\mR^m)$ with $\supp(\phi) \subset V$. Generalizing the definition for functions via \eqref{eq:LpLqDual}, the \emph{support} of a distribution $u \in \sS'(\mR^m)$ may be defined by \cite[p. 41]{Hoermander}
\begin{equation}
\mR^m \setminus   \supp(u) = \bigcup_{\substack{ V \subset \mR^m\text{ open}: \\ u_{|V} = 0  }} V. \label{eq:SuppDistri}
\end{equation}
The subspace of \emph{compactly supported} $u \in \sS(\mR^m)$ is denoted by $\Sdashc{m}$.

As an illustration of the above definitions, we consider a derivative of the \emph{Dirac delta} $\delta_0: \sS(\mR^m) \to \mC; \; \phi \mapsto \phi(0)$:
\vspace{1em}
	\begin{ex1}[Derivatives of the Dirac Delta] \label{ex:DiracDelta}
	     For $\alpha \in \mN_0^m$, let $u := \partial^{\alpha}\delta_0$. Then $u \in \Sdashc{m}$ has compact support $\supp(\partial^\alpha \delta_0) = \{ 0 \}$ and is of order $|\alpha|$.
	\end{ex1}
	\begin{pf}
	Linearity and well-definedness on $\sS(\mR^m)$ follow from \eqref{eq:DistriDeriv}. Furthermore, $u$ is continuous and of order $|\alpha|$ by the estimate
	 	    \begin{equation*}
	    |u(\phi)|  =  |\delta_0(\partial^\alpha \phi )|  =  | \partial^\alpha \phi(0 ) | \leq  \sup_{\bx \in \mR^m} |\partial^\alpha  \phi(\bx)| \MTEXT{for all} \phi \in \Ckc{\infty}{\mR^m}.
	    \end{equation*}
	    Moreover, $u$ has compact support $\{ 0 \}$ as $ \delta_0(\partial^\alpha \phi ) $ vanishes for all $\phi \in \sS(\mR^m)$ for which $0$ is not in the support.
	\end{pf}
\vspace{1em}

\end{section}

\begin{section}{The Fourier Transform} \label{S:FT}

	In the following, we review the properties of the Fourier transform, being a central tool in tomographic imaging. 
	\vspace{1em}
	\begin{df1}[Fourier Transform \text{\cite[p. 212]{Werner2007FA}}]\label{def:FT}
		The $m$-dimensional \emph{Fourier transform} of a function $f \in L^1(\mR^m)$ is defined by
		\begin{equation}
		 	\cF(f)(\bxi):= \hat f (\bxi) = (2\pi)^{-\frac m 2}  \int_{\mR^m} \E^{-\I \bxi \cdot \bx} f(\bx) \; \D x \MTEXT{for all} \bxi \in \mR^m. \label{eq:FT}
		\end{equation}
		If $\hat f \in L^1(\mR^m)$, then $f$ is given by \emph{Fourier's inversion formula}
		\begin{equation}
		 	 f (\bx)  = (2\pi)^{-\frac m 2} \int_{\mR^m} \E^{\I \bxi \cdot \bx} \hat f (\bxi) \; \D \xi  \MTEXT{for all} \bx \in \mR^m. \label{eq:iFT}
		\end{equation}
 		For indices $1\leq k \leq m$, let $\cF_{\underline k} (f)$ and $\cF_{\overline k} (f)$ denote the Fourier transform of $f$ with respect to the first $k$ or last $n+1-k$ arguments, respectively.
	\end{df1}
		\vspace{1em}
	\begin{th1}[Boundedness of the Fourier Transform \text{\cite[p. 212]{Werner2007FA}}]\label{thm:FTBounded}
		The Fourier transform defines a bounded linear operator $\cF: \Lp 1 {\mR^m} \to \Lp \infty {\mR^m}$ with norm $\norm{\cF} = (2\pi)^{-\frac m 2}$. Moreover, $\cF(f)$ is continuous for all $f \in \Lp 1 {\mR^m}$.
	\end{th1}
	\vspace{1em}
	
	Due to its significance in imaging, we further recall the definition of the \emph{convolution} $f \ast g \in \Lp{\infty}{\mR^m}$ of two functions $f,g \in \Lp{2}{\mR^m}$, given by \cite[p. 339]{Werner2007FA}
		\begin{equation}
		(f \ast g)(\bx) := \int_{\mR^m} f(\bx-\by)g(\by) \; \D y \label{eq:DefConv}.
		\end{equation}
	We study the properties of the Fourier transform in the Schwartz space $\sS(\mR^m)$ (see \sref{S:SchwartzLp}), which is closed under differentiation, multiplication and convolution:
	\vspace{1em}
	\begin{th1}[Properties of the FT \text{\cite[pp. 161-163]{Hoermander}, \cite[p. 189]{Evans}, \cite[p. 31]{Strichartz}}]\label{thm:PropFT}
		The Fourier transform defines an isomorphism  $\cF: \sS(\mR^m) \to \sS(\mR^m)$ with inverse given by \eqref{eq:iFT}. Moreover, we have for all $f,g \in \sS(\mR^m), \bxi \in \mR^m$, multi-indices $\alpha \in \mN_0^m$ and translations $\tau_{\ba}: h \mapsto (\bx \mapsto h(\bx+\ba))$ by $\ba \in \mR^m$
		\begin{subequations} \label{eq:PropFT}
		\begin{align}
		  \ip{f}{g}_{\Lp 2 {\mR^m}} &= \ip{\cF(f)}{\cF(g)}_{\Lp 2 {\mR^m}}              \label{eq:Parseval} \\
		  (2\pi)^{- \frac m 2} \cF(f \ast g) &= \cF(f) \cdot \cF ( g)           \label{eq:ConvThm} \\
		  (2\pi)^{\frac m 2}\cF(f \cdot g) &= \cF(f) \ast \cF ( g)              \label{eq:iConvThm} \\
		  \cF(\partial^\alpha f )(\bxi)   &= (\I \bxi)^\alpha \cF( f )(\bxi)    \label{eq:FTDeriv} \\
		  \cF(\bx^\alpha  f )(\bxi) &= (\I \partial)^\alpha \cF( f )(\bxi)      \label{eq:iFTDeriv} \\
		  \cF(\tau_{\ba} f)(\bxi) &= \E^{\I \bxi \cdot \ba}   \cF( f )(\bxi)    \label{eq:FTTrans} \\
		  \cF(\E^{\I \ba \cdot \bx}   f)(\bxi) &= \tau_{-\ba} \cF( f )(\bxi)    \label{eq:iFTTrans} \\
		   \cF(\cc{f})(\bxi) &= \cc{ \cF( f )(-\bxi) }   \label{eq:FTReflect} 
		\end{align}
		\end{subequations}
	\end{th1}
	\vspace{1em}
	Relation \eqref{eq:Parseval} is known as \emph{Parseval's formula} and states that $\cF$ is isometric with respect to the inner product in $L^2(\mR^m)$. By denseness of $\sS(\mR^m) \subset \Lp{2}{\mR^m}$ (see \thmref{thm:DenseIncl}) and \defref{def:OpsAdjoints}, this yields the following result:
		\vspace{1em}
	\begin{cor1}[Fourier Transform on $L^2(\mR^m)$ \text{\cite[p. 218]{Werner2007FA}}] \label{cor:FourierL2}
	 $\cF$ has a unique extension to a \emph{unitary operator}
	 	\begin{equation*}
	  \cF: L^2(\mR^m) \to L^2(\mR^m).
	\end{equation*}
	\end{cor1}
	\vspace{1em}
	Beyond \cref{cor:FourierL2}, the Fourier transform may even be extended to tempered distributions. In fact, the map
	\begin{equation}
	  \cF: \sS'(\mR^m) \to \sS'(\mR^m); \; (\cF u)(\phi) = u( \cF \phi ) \label{eq:DefFTSdash}
	\end{equation}
	defines an isomorphism and the relations \eqref{eq:PropFT} remain valid in a distributional sense. See \cite[pp. 164 ff.]{Hoermander} for details. By the inclusion $\Lp{p}{\mR^m} \subset \sS'(\mR^m)$, the properties in \thmref{thm:PropFT} generalize in particular to suitable $L^p$-spaces.
	
	Realistic specimens in tomographic applications are of bounded spatial extent, corresponding to functions or - more generally - distributions of compact support. 
	As the Fourier transform is a bijection on  $\sS'(\mR^m)$, general tempered distributions are mapped onto arbitrarily singular objects - like the Dirac delta considered in \exref{ex:DiracDelta}. It is thus surprising that the Fourier transform of any \emph{compactly supported} $u \in \Sdashc{m}$ is represented by an \emph{entire function} via the identification \eqref{eq:LpLqDual}, i.e.\ by $\sC^\infty$-functions with a globally convergent Taylor series in $\mC^m$. This entire function representation of $\cF(\Sdashc{m})$ is adopted throughout this work:
		\vspace{1em}
  	\begin{th1}[Paley-Wiener-Schwartz Theorem \text{\cite[p. 181]{Hoermander}}]\label{thm:PaleyW}
  		Let $K\subset \mR^m$ be compact and convex. Then, for any $u \in \Sdashc{m}$ of order $N \in \mN_0$ and support $\supp(u) \subset K$, $\cF(u) $ has an extension to an entire function $\hat u: \mC^m \to \mC$ and there exists a constant $C> 0$ such that
  		\begin{equation}
 		 |\hat u(\bxi)| \leq C ( 1+ \norm{\bxi})^N \exp\left( \sup_{\bx \in K} \Im(\bxi ) \cdot \bx \right) \MTEXT{for all} \bxi \in \mC^m. \label{eq:PaleyWBound}
  		\end{equation}
  		Conversely, any entire function $\hat u$  satisfying \eqref{eq:PaleyWBound}  is the complex extension of the Fourier transform of a distribution $u$ of order $\leq N$ and support in $K$.
  	\end{th1}
  		\vspace{1em}
 	\begin{ex1}[Fourier Transform of the Dirac delta] \label{ex:FTDiracDelta}
	     Let $u$ as in \Text{\exref{ex:DiracDelta}}. Then $\cF(u)(\bxi) = (2\pi)^{- \frac n 2 } (\I \bxi)^\alpha $ for all $\bxi \in \mR^m$ by \eqref{eq:DefFTSdash} and \eqref{eq:FTDeriv}, i.e.\ $\cF(u)$ is polynomial and thus entire. By the estimate \eqref{eq:PaleyWBound}, its purely algebraic growth behavior is a manifestation of its support $\supp(u) = \{0 \}$ and order $| \alpha|$.
	\end{ex1}
  	\vspace{1em}
  	Beyond the identification of compactly supported distributions and entire functions, it should be noted that the estimate \eqref{eq:PaleyWBound} relates regularity and support shape in real space to algebraic and exponential growth behavior in Fourier space. \thmref{thm:PaleyW} is the principal tool in the uniqueness analysis of \chapref{C:Analysis}.

\end{section}

\begin{section}{The Radon Transform} \label{S:Radon}

	X-ray tomography seeks to reconstruct a function $f: \mR^m \to \mC$ (e.g.\ describing an electron density) from its line integrals which give the transmitted radiation at different incident angles. Mathematically, this amounts to inverting a two-dimensional \emph{Radon transform} \cite{Radon1917} in the plane of rotation. We restrict the theoretical treatment of this operator to its version in $\mR^2$, being the relevant one to the tomographic applications considered herein, in order to keep the notation simple. For a more general discussion of the Radon transform in $\mR^m$, see for instance \cite{Natterer}.
	\vspace{1em}
	\begin{df1}[2D Radon Transform \text{\cite[p. 9]{Natterer}}]\label{def:Radon}
		Set $Z^m:= [0;2\pi) \times \mR^{m-1}$. For a function $f \in \sS(\mR^2)$, we define its \emph{Radon transform} $\cR  f:  Z^2 \to \mC$ as
		\begin{equation}
		 \cR f(\theta, x) = \int_{\mR}  f(x \bn_{\theta} + y \bn_{\theta}^\perp) \; \D y  \label{eq:Radon}
		\end{equation}
		with $\bn_{\theta} := (\cos\theta, \sin\theta)$ and $\bn_{\theta}^\perp := (\sin\theta, -\cos\theta)$.
		The graph of $\cR(f)$ in the $x$-$\theta$-plane is denoted as the \emph{sinogram} of $f$.
	\end{df1}
	\vspace{1em}
	Note that $\cR$ is linear in $f$ by \eqref{eq:Radon}.
	The principal tool for the further analysis is given by the \emph{Fourier Slice Theorem}, relating the Radon transform to the \emph{polar Fourier transform} $\cF_{\rm p }$. For $f\in \sS(\mR^2)$, the latter is defined by
	\begin{equation}
	 	 \cF_{\rm p }(f) (\theta, \xi_x) := (2 \pi)^{\frac 1 2 } \cF(f) (\xi_x \bn_{\theta}) \MTEXT{for all} (\theta, \xi_x) \in Z^2  \label{eq:PolarFT}.
	\end{equation}
	\begin{th1}[Fourier Slice Theorem \text{\cite[p. 11]{Natterer}}]\label{thm:FourierSlice}
		Let $f \in \sS(\mR^2)$. Then
		\begin{equation}
		   \cF_2(\cR f) = \cF_{\rm p }(f). \label{eq:FourierSlice}
		\end{equation}
		where $\cF_2$ is the 1D Fourier transform with respect to the second variable of $\cR f$.
	\end{th1}
	\vspace{1em}
	By \thmref{thm:FourierSlice} and the isomorphism $\cF: \sS(\mR^2) \to \sS(\mR^2)$, measuring the Radon transform of a function $f$ is equivalent to sampling its Fourier transform $\cF(f)$ \emph{on a polar grid}.
	
	This work is exclusively concerned with the tomography of bounded physical objects, i.e.\ such which are parametrized by functions  supported in a bounded set $\Omega \subset \mR^2$. The path of integration, for which the integrand in \eqref{eq:Radon} is non-zero, is then bounded in length by the finite diameter $\diam(\Omega)$, given by the maximum distance of two points $x,y\in\Omega$. For $f \in \Ckc{\infty}{\Omega} \subset \Lp{p}{\Omega}$, this implies
	\begin{equation*}
	    \norm{\cR f}_{L^\infty(Z^2)} \leq  \sup_{(\theta , x) \in Z^2}\int_{\mR}  \left| f(x \bn_{\theta} + y \bn_{\theta}^\perp) \right| \; \D y\leq  \norm{f}_{L^\infty(\Omega)} \cdot \diam(\Omega),
	\end{equation*}
	i.e.\ $\cR$ is $L^\infty$-bounded. Likewise, bounds in $L^2$-norm may be derived using \thmref{thm:FourierSlice} and \cref{cor:FourierL2}. Since $\Ckc \infty {\Omega} \subset \sS(\mR^2)$ is dense in $\Lp{\infty}{\Omega}$ and $\Lp{2}{\Omega}$ by \thmref{thm:AdjointProps}, these observations permit extensions of $\cR$ to $L^p$-spaces:
		\vspace{1em}
	\begin{th1}[Continuity of the Radon Transform on bounded Domains]\label{thm:RadonContinuity}
		Let $\Omega \subset \mR^2$ bounded. Then \eqref{eq:FourierSlice} holds for all $f \in \Lp 2 \Omega$. Moreover, $\cR$ and $\cF_{\rm p }$, defined  by \eqref{eq:Radon} and \eqref{eq:PolarFT}, have unique  extensions to bounded linear operators
		\begin{equation*}
		   \cR: \Lp \infty {\Omega} \to \Lp \infty {Z^2} \MTEXT{and} \cF_{\rm p }, \cR: \Lp 2 {\Omega} \to \Lp 2 {Z^2}.
		\end{equation*}
	\end{th1}	
		\vspace{1em}
	In fact, an even stronger statement holds true than the $L^2$-continuity stated here \cite[Theorem 1.6]{Natterer}. For the present work, however, it is sufficient that the Radon transform of compactly supported functions is both $L^\infty$- and $L^2$-continuous.
	
	Physical objects in X-ray tomography are always three-dimensional, varying not only in the plane of rotation but also along the axis, by convention taken to be the second variable of a function $f: \mR^3 \to \mC$. It is therefore convenient to define the \emph{cylindrical} Radon- and Fourier transforms of a function $f\in \sS(\mR^m)$ by
	\begin{subequations} \label{eq:CylTransforms}
	\begin{align}
	 \CR(f)(\theta,x,\by) &:= \cR(f(\cdot, \by, \cdot))(\theta, x) \label{eq:CylRadon}  \\
	 \CF(f)(\theta,\xi_x,\bxi_y) &:= (2 \pi)^{\frac 1 2 } \cF(f) (\xi_x \cos(\theta), \bxi_y, \xi_x \sin (\theta) )   \label{eq:CylFourier}
	\end{align}
	\end{subequations}
	for all $x,\xi_x \in\mR, \by,\bxi_y \in \mR^{m-2}, \theta \in [0; 2\pi)$. $\CR$ simply amounts to applying $\cR$ to two-dimensional slices of $f$, parametrized by the first and the $m$-th variable. By generalization of the previous results, this yields the following properties:
	\vspace{1em}
	\begin{cor1}[Properties of the Cylindrical Transforms] \label{cor:CylRadonFourProps}
	 Let $\Omega \subset \mR^m$ be bounded in the first and the last dimension. Then we have for all $f \in \sS(\mR^m) \cup \Lp 2 \Omega$
	 	\begin{equation}
	  \cF_{\overline 2}( \CR f)  = \CF (f).  \label{eq:FourierSliceCyl}
	\end{equation}
	Moreover, $\CR$ and $\CF$ have unique  extensions to bounded linear operators
		\begin{equation}
		   \CR: \Lp \infty {\Omega} \to \Lp \infty {Z^m} \MTEXT{and} \CF, \CR: \Lp 2 {\Omega} \to \Lp 2 {Z^m} \label{eq:CylTransformDoms}.
		\end{equation}
	\end{cor1}
	\vspace{1em}
	
	As the definitions \eqref{eq:CylTransforms} may seem bulky, \figref{fig:RadonVisualization} visualizes the slicewise application of $\cR$ and illustrates the physical interpretation of the Radon transform as shadow images of an object which is illuminated under different incident angles.
	\begin{figure}[hbt!]
	 \centering
	 \includegraphics[width=\textwidth]{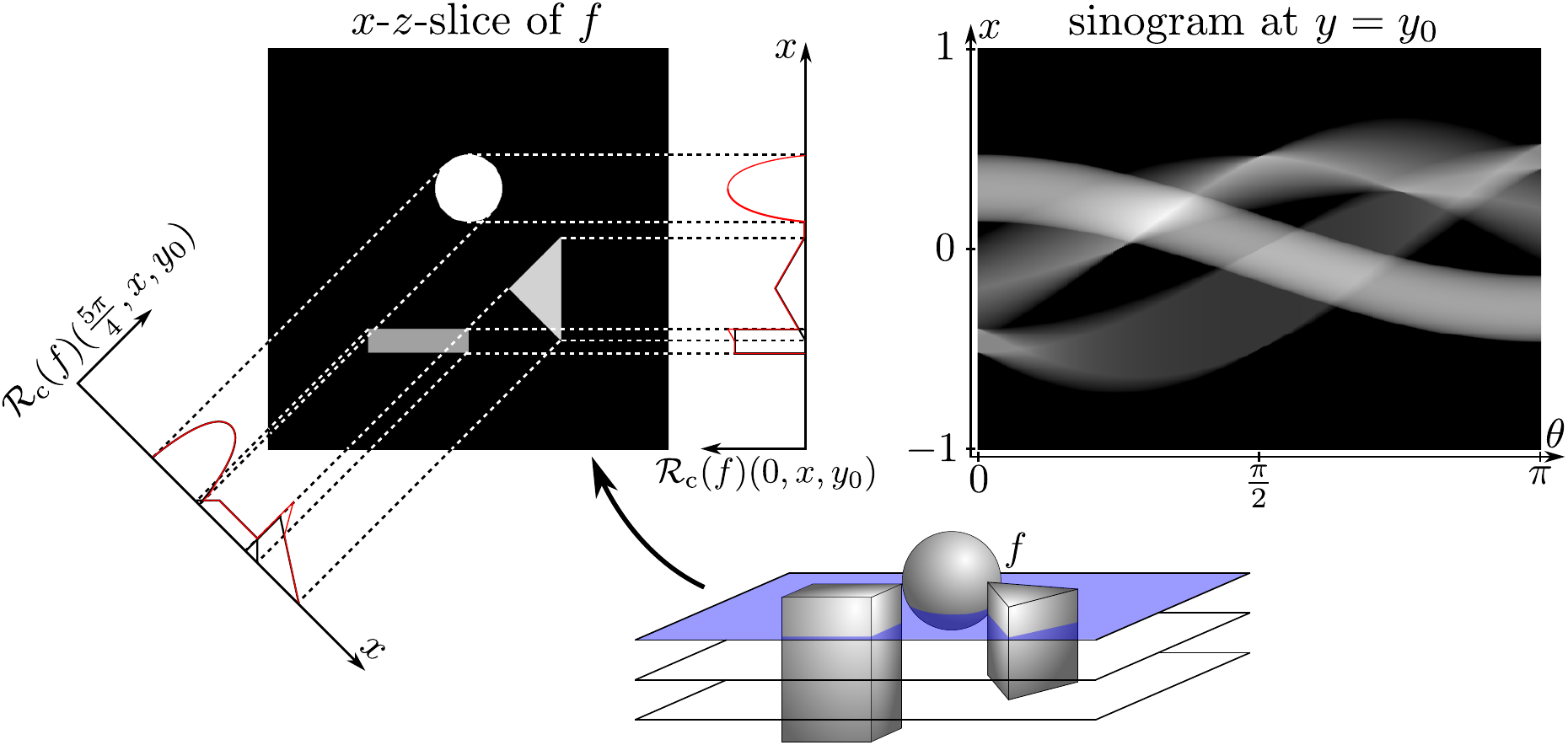}
	 \caption{Illustration of the (cylindrical) Radon transform $\CR$: a 3D-object characterized by a density $f: \mR^3 \to \mC$ (bottom) is projected slice-wise under different rotation angles $\theta$, yielding shadow images. The projections for different $\theta$ form the sinogram. The ensemble of sinograms of the different slices stacked along the $y$-direction defines $\CR(f)$ (Source: \cite{AikeMaster}, modified).\label{fig:RadonVisualization}}
	\end{figure}

\end{section}

\begin{section}{Sobolev Spaces} \label{S:Sobolev}

	In this section, we introduce Sobolev spaces, providing a weak notion of differentiability and smoothness related to the asymptotic behavior of the Fourier transform. This relation is used in the uniqueness analysis of \sref{SS:PhaseRetrFF} and exploited to impose regularity constraints in the Newton methods introduced in \chapref{C:NumMeth}.
	\vspace{1em}
	\begin{df1}[Weak Derivatives \text{\cite[p. 256]{Evans}}]\label{def:WeakD}
		Let $u,v \in L^1_{\rm loc} (\mR^m)$ - locally integrable and $\alpha \in \mN_0^m$ a multi-index. Then, $v$ is the $\alpha$-th \emph{weak derivative} of $u$ and we write $\partial^\alpha u := v $, if
		\begin{equation*}
		 	 \int_{\mR^m} u \cdot \partial^\alpha \phi \; \D x = (-1)^{|\alpha|} \int_{\mR^m} v \cdot   \phi \; \D x  \MTEXT{for all} \phi \in \Ck \infty {\mR^m}.
		\end{equation*}
	\end{df1}
	\vspace{1em}
	Weak derivatives are unique up to a set of measure zero and coincide in this sense with classical derivatives if the latter exist. Moreover, sum- and product rule generalize to weak derivatives (for details, see  \cite[pp. 257, 261]{Evans}). The existence of weak derivatives gives rise to function spaces, called \emph{Sobolev spaces}:
	\vspace{1em}
	\begin{df1}[Sobolev Spaces \text{\cite[p. 258]{Evans}}]\label{def:SobSpace}
		For $k \in \mN_0$, $ 1 \leq p \leq \infty$ the Sobolev space $\WpR{k}{p}$ is defined as the vector space of all functions $u \in L^1_{\rm loc} (\mR^m)$ such that $\partial^\alpha u \in L^p(\mR^m)$ exists for all $\alpha \in \mN_0^m: |\alpha| \leq k$, equipped with the norm
 		\begin{equation*}
 		\norm{u}_{\WpR{k}{p}} := \begin{cases}
 		                          \left( \sum_{\alpha \in \mN_0^m: |\alpha| \leq k} \norm{\partial^\alpha u }_{\Lp p {\mR^m}}^p \right)^{\frac 1 p } &\text{for }p < \infty \\
 		                          \sum_{\alpha \in \mN_0^m: |\alpha| \leq k} \norm{\partial^\alpha u }_{\Lp \infty {\mR^m}}   &\text{for }p = \infty \\
 		                         \end{cases}.
 		\end{equation*}
	\end{df1}
	\vspace{1em}
	$\WpR{k}{p}$ is a \emph{Banach space} for all $k \in \mN_0$, $ 1 \leq p \leq \infty$ \cite[p. 262]{Evans} and thus a \emph{Hilbert space} in the $L^2$-based case $p = 2$. Moreover, the map
	\begin{equation*}
	 \partial^\alpha : \WpR{k}{p} \to \WpR{k- |\alpha|}{p}
	\end{equation*}
	defines a bounded linear operator for all $\alpha \in \mN_0^m: |\alpha| \leq k$.
	
	Owing to the definition of weak derivatives via $L^2$-inner products and \cref{cor:FourierL2}, weak differentiability is related to the decay behavior of the Fourier transform. In fact, using \eqref{eq:FTDeriv}, the following implications can be shown:
	\begin{align*}
	 u \in \WpR{k}{p} \;\;&\Leftrightarrow \;\; \partial^\alpha u \in \Lp 2 {\mR^m} \;\;\;\;\;\;\;\;\MTEXT{for all} |\alpha| \leq k \\
	 &\Leftrightarrow \;\; \xi^\alpha \cdot \cF( u )  \in \Lp 2 {\mR^m} \MTEXT{for all} |\alpha| \leq k \\
	 &\Leftrightarrow \;\; (1+ \norm{\xi}_2^2)^{\frac k 2 } \cdot \cF( u ) \in \Lp 2 {\mR^m}.
	\end{align*}
	Thus, defining $\HsR{s}:= \{ u \in \Lp 2 {\mR^m} : \norm{u}_{\HsR{s}} < \infty \}$ for $s \geq 0$ with
	\begin{align}
	 \norm{u}_{\HsR{s}}  &:= \norm{(1+ \norm{\xi}_2^2)^{\frac s 2 } \cdot \cF( u )}_{\Lp 2 {\mR^m}} ,  \label{eq:DefHs}
	\end{align}
	one arrives at an alternative characterization of Sobolev spaces for $p= 2$:
	\vspace{.5em}
	\begin{th1}[Sobolev Spaces by Fourier Transforms \text{\cite[p. 258]{Evans}}]\label{thm:SobFourier}
		For $k \in \mN_0$, we have
		\begin{equation*}\HsR{k} = \WtwoR{k}\end{equation*}
	      and the norms $\norm{\cdot}_{\HsR{k}}, \norm{\cdot}_{\WtwoR{k}}$ are equivalent. Furthermore,
	      \begin{equation*}
	      \partial^\alpha : \HsR{s} \to \HsR{k- |\alpha|}
	      \end{equation*}
	      is continuous for all $s\geq 0$, $\alpha \in \mN_0^m: |\alpha| \leq s$.
	\end{th1}
	\vspace{.5em}
	Note that the definition of $\HsR{s}$ remains reasonable for non-integer  $s$, other than $\WtwoR{k}$, being based on the existence of an integer number of weak derivatives.
	The following theorem establishes a link between pointwise notions of continuity and differentiability and the introduced Sobolev space characterization:
	\vspace{.5em}
	\begin{th1}[Sobolev Embedding Theorem \text{\cite[p. 160]{Strichartz}}]\label{thm:SobEmbed}
		For $m \in \mN$, $s > \frac m 2$ and $k \in \mN_0$ $1 \leq p < \infty$ such that $k > \frac m p$, we have
		\begin{equation*}
		 \HsR{s} \subset \Lp \infty {\mR^m} \cap \Ck{0}{\mR^m} \MTEXT{and} \WpR{k}{p} \subset \Lp \infty {\mR^m} \cap \Ck{0}{\mR^m}
		\end{equation*}
		with continuous embeddings $\HsR{s} \hookrightarrow \Lp \infty {\mR^m}$, $\WpR{k}{p} \hookrightarrow \Lp \infty {\mR^m}$.
	\end{th1}
	\vspace{.5em}
	
	In the proof of \thmref{thm:UniquePhaseFFMD}, we are concerned with estimating the asymptotic decay of Fourier transforms of \emph{compactly supported} functions. The following, final result of this section shows a certain uniformity of the asymptotic behavior with respect to different dimensions for this class of functions:
	\vspace{.5em}
	\begin{lem1}\label{lem:HsDecay}
		For $s \geq 0$, let $f \in H^s(\mR^{n+m})$ with compact support $\Omega$. Define $\tilde f(\bx,\bxi_y):= \cF(f(\bx, \cdot))(\bxi_y)$ for all $\bx \in \mR^n, \bxi_y \in \mC^m$. Then 
		\begin{equation*} \tilde f(\cdot, \bxi_y) \in H^s(\mR^{n}) \MTEXT{for all} \bxi_y \in \mC^m \end{equation*}
		and the map $\bxi_y \mapsto \norm{\tilde f(\cdot , \bxi_y)}_{H^s(\mR^{n})}$ is continuous in $\mC^m$.
	\end{lem1}
	\begin{pf}
	 Let $\bxi_y \in \mC^m$ be arbitrary. Since $f$ has finite support, so has $f(\bx, \cdot)$ for all $\bx \in \mR^n$, so that $\tilde f(\bx,\cdot) := \cF(f(\bx, \cdot))$ has a unique extension to an entire function in $\mC^m$ by \thmref{thm:PaleyW}. Thus $\tilde f(\cdot, \bxi_y)$ is well defined. Define
	 \begin{align*}
	  \Omega_{x} &:= \{\bx \in \mR^{n}: \exists {\by}\in \mR^m: (\bx,{\by}) \in \Omega \} \\
	  \Omega_{y} &:=  \{{\by} \in \mR^{m}: \exists \bx\in \mR^n: (\bx,{\by}) \in \Omega \}
	 \end{align*}
	  and $y_0 =  \sup \{\norm{{\by}}_2 : {\by} \in \Omega_{y} \}$. By compactness of the support, integrations of $f$ over $\mR^n$ or $\mR^m$ can always be restricted to $\Omega_{x}$ and $\Omega_{y}$, respectively, which are of finite Lebesgue measure $\mu(\Omega_{x}), \mu(\Omega_{y}) < \infty$. Setting $C(\bxi_y)^2:= (2\pi)^{-m } \mu(\Omega_{y})^2 \cdot \exp(2y_0\norm{\Im(\bxi_y)}_2)$ and using the Cauchy--Schwarz inequality, Fubini's theorem and the estimate $\E^{\Im(\bxi_y) \cdot {\by}} \leq \E^{y_0 \norm{\Im(\bxi_y)}_2}$ for $\by \in \Omega_y$ we thus obtain for all $\bxi_x \in\mR^n$
	 \begin{align*}
	   |\cF(f)(\bxi_x, \bxi_y)|^2 &=   (2\pi)^{-(n+m)} \left| \int_{\Omega_x} \cc{\E^{\I \Re(\bxi_y) \cdot {\by}} } \cdot \left( \int_{\mR^n} \E^{\Im(\bxi_y) \cdot {\by}} \E^{-\I \bxi_x \cdot \bx} f(\bx,{\by}) \; \D x \right) \; \D y \right|^2 \\
	   &\leq (2\pi)^{-m } \mu(\Omega_{y})^2  \exp(2y_0\norm{\Im(\bxi_y)}_2) \int_{\mR^m} \left| \cF(f(\cdot, {\by}))(\bxi_x)   \right|^2 \; \D y \\  
	   &\stackrel{\eqref{eq:Parseval}}= C(\bxi_y)^2   \norm{ \cF (f)(\bxi_x, \cdot) }_{L^2(\mR^m)}^2.
	  \end{align*}
	 Accordingly, we have
	 \begin{align*}
	   \infty > \norm{f}^2_{H^s(\mR^{n+m})} &= \int_{\mR^{n+m}} (1+ \norm{\bxi_x}_2^2 + \norm{\bxi_y}_2^2)^s   |\cF(f)(\bxi_x, \bxi_y)|^2 \D \xi_x \; \D \xi_y \\
	   &\geq \int_{\mR^{n}} (1+ \norm{\bxi_x}_2^2 )^s   \norm{ \cF(f)(\bxi_x, \cdot)}_{L^2(\mR^m)}^2  \; \D \xi_x  \\
	   &\geq C(\bxi_y)^{-2} \int_{\mR^{n}} (1+ \norm{\bxi_x}_2^2 )^s   |\cF(f)(\bxi_x, \bxi_y)|^2 \; \D \xi_x \\
	   &= C(\bxi_y)^{-2} \int_{\mR^{n}} (1+ \norm{\bxi_x}_2^2 )^s   |\cF (\tilde f(\cdot, \bxi_y))(\bxi_x)|^2 \; \D \xi_x \\
	   &= C(\bxi_y)^{-2}\norm{\tilde f(\cdot , \bxi_y)}^2_{H^s(\mR^{n})},
	 \end{align*}
	 which proves that $\tilde f(\cdot , \bxi_y) \in H^s(\mR^{n})$.
	 
	 Concerning continuity of $\bxi_y \mapsto \norm{\tilde f(\cdot , \bxi_y)}_{H^s(\mR^{n})}$ note that, by analyticity of $\cF(f)$, $\bxi_y \mapsto|\cF(f)(\bxi_x, \bxi_y)|^2$ is continuous in $\mC^m$ for all $\bxi_x \in \mC^n$ and that for any sequence $\seqn{\bxi}{k} \subset \mC^m$ converging to $\bxi_y$, $\tilde C^2 := \sup_{k \in \mN } C^2(\bxi_k)$ is finite. Hence, 
	 \begin{equation*}\bxi_x \to (1+ \norm{\bxi_x}_2^2 )^s  |\cF(f)(\bxi_x, \bxi_k)|^2\end{equation*}
	 is dominated for all $k \in \mN$ by the integrable function 
	 \begin{equation*}\bxi_x \mapsto \tilde C^2 (1+ \norm{\bxi_x}_2^2 )^s    \norm{ \cF (f)(\bxi_x, \cdot) }_{L^2(\mR^m)}^2 \end{equation*}
	 according to the above estimates. By application of Lebesgue's dominated convergence theorem (see for instance \cite[p. 516]{Werner2007FA}), this yields 
	 \begin{align*}
	    \lim_{k \to \infty} \norm{\tilde f(\cdot , \bxi_k)}_{H^s(\mR^{n})}^2 &= \lim_{k \to  \infty} \int_{\mR^{n}} (1+ \norm{\bxi_x}_2^2 )^s   |\cF_{n}(\tilde f(\cdot, \bxi_k))(\bxi_x)|^2 \D \bxi_x \\
	    &=  \int_{\mR^{n}} (1+ \norm{\bxi_x}_2^2 )^s  \lim_{k \to  \infty} |\cF(f)(\bxi_x, \bxi_k)|^2 \D \bxi_x \\
	    &= \int_{\mR^{n}} (1+ \norm{\bxi_x}_2^2 )^s  |\cF(f)(\bxi_x, \bxi_y)|^2 \D \bxi_x = \norm{\tilde f(\cdot , \bxi_y)}_{H^s(\mR^{n})}^2.
	 \end{align*}
	 By generality, this proves that $\bxi_y \mapsto \norm{\tilde f(\cdot , \bxi_y)}_{H^s(\mR^{n})}^2$ and thus $\bxi_y \mapsto \norm{\tilde f(\cdot , \bxi_y)}_{H^s(\mR^{n})}$ are continuous maps in $\mC^m$.
	\end{pf}
 
\end{section}

\end{chapter}

\begin{chapter}{Appendix: Supplementary Proofs}\label{B}

  \begin{section}{Generalization of \thmref{thm:UniquePhaseFFMD} to $m \geq 2$} \label{A:MDFourPhaseRetr}
      
      In this appendix, we prove \thmref{thm:UniquePhaseFFMD} for in arbitrary dimensions $m \geq 2$ by reduction to the case $m=2$ for which statement is shown in \sref{SS:PhaseRetrFF}.
	\begin{pf}[Proof of \thmref{thm:UniquePhaseFFMD} for $m \geq 2$:] Let the assumptions of \thmref{thm:UniquePhaseFFMD} hold for  $m > 2$.
	  For $(\xi_3, \ldots, \xi_m) \in V:= \left( - \frac \pi 2; \frac \pi 2 \right)^{m-2}$, set $C:=\prod_{j=3}^m  \cF( b_{\rect} )(\xi_j)$. Then
	  \begin{equation*}
	   |C|=\prod_{j=3}^m \left|\cF( b_{\rect} )(\xi_j)\right| = 2^{m-2} \prod_{j=3}^m \underbrace{\frac{|\sin(\xi_j)|}{|\xi_j|}}_{\geq \frac 1 2} \geq 1.
	  \end{equation*}
	  Now define $\psi_{\xi_3, \ldots, \xi_m} := b_{\xi_3, \ldots, \xi_m} + h_{\xi_3, \ldots, \xi_m}$ such that for all $x_1,x_2 \in \mR$, $\xi_3, \ldots, \xi_m \in \mC$
	  \begin{align*}
	   b_{\xi_3, \ldots, \xi_m}(x_1, x_2) &:=b_{\exp}(x_1) b_{\rect}(x_2)  \prod_{j=3}^m \cF( b_{\rect} )(\xi_j)   \\ 
	   h_{\xi_3, \ldots, \xi_m}(x_1, x_2) &:=  \cF (h(x_1,x_2, \cdot))(\xi_3, \ldots, \xi_m).
	  \end{align*}
	  Note that $\supp(h_{\xi_3, \ldots, \xi_m}) \subset [0;1] \times [-1;1]$ and that $h_{\xi_3, \ldots, \xi_m} \in H^{\frac 3 2}(\mR^2)$ according to \lemref{lem:HsDecay}. Moreover, $b_{\xi_3, \ldots, \xi_m} = C \cdot b_{\exp, 2}$ by construction (compare \eqref{eq:DefBaseMD}).
	  
	  Thus, up to the nonzero scaling constant $C$, the setting given by $b_{\xi_3, \ldots, \xi_m}$, $h_{\xi_3, \ldots, \xi_m}$ and $g_{\xi_3, \ldots, \xi_m}$ exactly matches that of \thmref{thm:UniquePhaseFFMD} for $m =2$, where the Fourier intensity data is given for all $ \xi_1,\xi_2, \xi_3, \ldots, \xi_m  \in \mR $ by
	  \begin{equation*}
	   |\cF (\psi_{\xi_3, \ldots, \xi_m})|^2(\xi_1,\xi_2) = | \cF( \psi ) |^2 (\xi_1,\xi_2, \xi_3, \ldots, \xi_m).
	  \end{equation*}
	  Hence, $h_{\xi_3, \ldots, \xi_m}$  can be reconstructed uniquely by the statement for $m=2$. As this holds for all $(\xi_3, \ldots, \xi_m) \in V$, where $V$ is open and $(\xi_3, \ldots, \xi_m) \mapsto h_{\xi_3, \ldots, \xi_m}$ is entire by \thmref{thm:PaleyW}, $(x_1,x_2, \xi_3, \ldots, \xi_m) \mapsto h_{\xi_3, \ldots, \xi_m}(x_1,x_2)$ and thus $h$ are uniquely determined as functions in $\mR^m$.
	\end{pf}
	
   \end{section}
  
\end{chapter}


\cleardoublepage
\bibliography{literature.bib}


\newpage
\thispagestyle{empty}
\section*{Danksagung}
 
An dieser Stelle m\"ochte ich mich bei all denjenigen bedanken, die mich durch zahl-reiche Ideen, motivierende Fragestellungen, weiterf\"uhrende Diskussionen, n\"achtliche Druckdienste oder einfach nur durch willkommene Ablenkung beim Anfertigen dieser Arbeit unterst\"utzt haben. Mein Dank gilt dabei insbesondere den Mitgliedern des IRP f\"ur die geniale Arbeits- und Nicht-Arbeitsatmosph\"are im Institut, die einen wunderbaren Rahmen f\"ur dieses Masterprojekt bildete. Speziell seien dabei erleuchtende Diskussionen mit Aike Ruhlandt zu Eindeutigkeit und Konsistenz erw\"ahnt sowie Martin Krenkel, der mir die Vor- und Nachteile von CTF-Rekonstruktionen gegen\"uber Newton-Verfahren n\"aher brachte, und Johannes Hagemann f\"ur die spontane Beantwortung diverser Spontanfragen. Weiterhin danke ich Matthias Bartels f\"ur die Aufbereitung und Bereitstellung des Nanokolloid-Datensatzes.

Dar\"uber hinaus m\"ochte ich mich herzlich bei Thorsten Hohage bedanken, der mir einerseits das faszinierende interdisziplin\"are Thema dieser Arbeit vorschlug und mich durch immer neue Aspekte inspirierte, andererseits aber auch stets offen war f\"ur meine Vorschl\"age zur Vertiefung oder Kurs\"anderung. Ebenso dankbar bin ich Tim Salditt f\"ur die immer wieder aufgezeigte experimentelle Perspektive, die mich nicht zu weit in den Elfenbeinturm der Idealisierung hat abdriften lassen, und seine wundervoll enthusiastische Einstellung gegen\"uber selbst mathematisch abstraktesten Resultaten. Mein Dank gilt auch Professor Plonka-Hoch f\"ur die Zweitbegutachtung dieser vielleicht nicht immer ganz kurz und b\"undig geratenen Masterarbeit.

Dieses Projekte wurde im Rahmen des SFB 755 Nanoscale Photonic Imaging durch die Deutsche Forschungsgemeinschaft unterst\"utzt. Des Weiteren m\"ochte ich an dieser Stelle der Studienstiftung des Deutschen Volkes meinen Dank aussprechen, deren Stipendiatenf\"orderung mein nun endendes Studium entscheidend begleitet und bereichert hat.

Au{\ss}erdem m\"ochte ich mich noch bei meinen Eltern bedanken, die mir letzteres erm\"oglicht haben - mit der Freiheit selbst zu ergr\"unden, ob und wo die Reise nach f\"unf Jahren endet. Zu guter Letzt gilt mein besonderer Dank noch Jenni, ohne deren Unterst\"utzung diese Masterarbeit wohl eher mich geschafft h\"atte als umgekehrt.

\end{document}